\documentclass[reqno]{amsart}
\usepackage[margin = 1in]{geometry}
\usepackage{amsmath, amssymb, amsthm, bm, fancyhdr, verbatim, graphicx}
\usepackage{enumerate}
\usepackage[all]{xy}
\usepackage[usenames,dvipsnames]{xcolor}
\usepackage{mathrsfs}
\usepackage{tikz-cd}
\usepackage{framed, hyperref}
\usepackage[titletoc]{appendix}
\usepackage{bbm}
\usepackage{lipsum}
\usepackage{adjustbox}

\usepackage{stmaryrd}

\numberwithin{equation}{subsection}

\newcommand{\U}{\mathrm{U}}

\DeclareFontFamily{U}{matha}{\hyphenchar\font45}
\DeclareFontShape{U}{matha}{m}{n}{
      <5> <6> <7> <8> <9> <10> gen * matha
      <10.95> matha10 <12> <14.4> <17.28> <20.74> <24.88> matha12
      }{}
\DeclareSymbolFont{matha}{U}{matha}{m}{n}
\DeclareFontFamily{U}{mathx}{\hyphenchar\font45}
\DeclareFontShape{U}{mathx}{m}{n}{
      <5> <6> <7> <8> <9> <10>
      <10.95> <12> <14.4> <17.28> <20.74> <24.88>
      mathx10
      }{}
\DeclareSymbolFont{mathx}{U}{mathx}{m}{n}

\DeclareMathSymbol{\obot}         {2}{matha}{"6B}

\RequirePackage{xspace}

\newcommand{\BQ}{\ensuremath{\mathbb{Q}}\xspace}

\newcommand{\sK}{\ensuremath{\mathscr{K}}\xspace}
\newcommand{\sL}{\ensuremath{\mathscr{L}}\xspace}

\newcommand{\sV}{\ensuremath{\mathscr{V}}\xspace}
\newcommand{\sW}{\ensuremath{\mathscr{W}}\xspace}

\usepackage{color}

\newcommand{\F}{\mathbf{F}}

\newcommand{\wt}[1]{\widetilde{#1}}

\newcommand{\Q}{\mathbf{Q}}
\newcommand{\Z}{\mathbf{Z}}

\newcommand{\mf}[1]{\mathfrak{#1}}

\newcommand{\R}{\mathbf{R}}

\newcommand{\ul}[1]{\underline{#1}}
\newcommand{\ol}[1]{\overline{#1}}

\newcommand{\Cal}[1]{\mathcal{#1}}
\newcommand{\A}{\mathbf{A}}

\newcommand{\ft}{{}^{\tau}} 
\newcommand{\co}{\colon}
\newcommand{\mrm}[1]{\mathrm{#1}}

\newcommand{\bbm}[1]{\mathbbm{#1}}
\newcommand{\PP}{\mathbf{P}}

\newcommand{\inj}{\hookrightarrow}
\newcommand{\surj}{\twoheadrightarrow}

\newcommand{\tw}[1]{\langle #1 \rangle}

\DeclareMathOperator{\GL}{GL}

\DeclareMathOperator{\Frob}{Frob}

\DeclareMathOperator{\coker}{coker}

\DeclareMathOperator{\Tr}{Tr}

\DeclareMathOperator{\Hom}{Hom}
\newcommand{\cHom}{\Cal{H}om}

\DeclareMathOperator{\Ind}{Ind}

\DeclareMathOperator{\rank}{rank}

\DeclareMathOperator{\Aut}{Aut}

\DeclareMathOperator{\Nm}{Nm}
\DeclareMathOperator{\Spec}{Spec}

\DeclareMathOperator{\End}{End}

\DeclareMathOperator{\Stab}{Stab}
\DeclareMathOperator{\Bun}{Bun}

\DeclareMathOperator{\Pic}{Pic}

\DeclareMathOperator{\Id}{Id}

\DeclareMathOperator{\Sht}{Sht}

\DeclareMathOperator{\Fl}{Fl}

\DeclareMathOperator{\pr}{pr}

\DeclareMathOperator{\Hk}{Hk}

\DeclareMathOperator{\Ch}{CH}

\DeclareMathOperator{\IC}{IC}

\DeclareMathOperator{\Herm}{Herm}
\DeclareMathOperator{\Lagr}{Lagr}

\DeclareMathOperator{\Db}{Db}
\DeclareMathOperator{\Coh}{Coh}
\DeclareMathOperator{\supp}{supp}
\DeclareMathOperator{\Den}{Den}

\DeclareMathOperator{\Prym}{Prym}

\DeclareMathOperator{\ns}{ns}

\DeclareMathOperator{\Spr}{Spr}
\DeclareMathOperator{\HSpr}{HSpr}

\DeclareMathOperator{\AJ}{AJ}

\DeclareMathOperator{\Corr}{Corr}
\DeclareMathOperator{\FYZ}{FYZ}

\newcommand{\Ql}{\Q_{\ell}}

\newcommand{\bu}{\bullet}

\newcommand{\map}[1]{\xrightarrow{#1}}

\newcommand{\iso}{\cong}
\newcommand{\C}{\mathbf C}
\newcommand{\p}{{}^{\mf{p}}}

\newcommand{\nc}{\newcommand}
\nc{\renc}{\renewcommand}
\nc{\ssec}{\subsection}
\nc{\sssec}{\subsubsection}
\nc{\on}{\operatorname}

\nc{\us}[2]{\underset{#1}{#2}}
\nc{\os}[2]{\overset{#1}{#2}}
\nc{\fr}[2]{\frac{#1}{#2}}
\nc{\pfr}[2]{\frac{\partial #1}{\partial #2}}
\nc{\la}{\langle}
\nc{\ra}{\rangle}
\nc{\hra}{\hookrightarrow}

\nc{\Loc}{\on{Loc}}
\nc{\rk}{\on{rk}}
\nc{\Mod}{\on{Mod}}
\nc{\Hhom}{\underline{\on{Hom}}}
\nc{\QCoh}{\on{QCoh}}
\nc{\gr}{\on{gr}}
\nc{\Grpd}{\on{Grpd}}
\nc{\id}{\on{id}}
\nc{\dr}{{\on{dR}}}
\nc{\Maps}{\on{Maps}}
\nc{\oblv}{{\on{oblv}}}
\nc{\alg}{\on{-alg}}
\nc{\comod}{{\on{-comod}}}
\nc{\perf}{\on{perf}}

\nc{\vphi}{\varphi}

\nc{\CA}{{\mathcal{A}}}
\nc{\CB}{{\mathcal{B}}}

\nc{\CD}{{\mathcal{D}}}
\nc{\CE}{{\mathcal{E}}}
\nc{\CF}{{\mathcal{F}}}
\nc{\CG}{{\mathcal{G}}}
\nc{\CI}{{\mathcal{I}}}
\nc{\CJ}{{\mathcal{J}}}
\nc{\CK}{{\mathcal{K}}}
\nc{\CL}{{\mathcal{L}}}
\nc{\CM}{{\mathcal{M}}}
\nc{\CN}{{\mathcal{N}}}
\nc{\CO}{{\mathcal{O}}}
\nc{\CP}{{\mathcal{P}}}
\nc{\CQ}{{\mathcal{Q}}}
\nc{\CR}{{\mathcal{R}}}
\nc{\CS}{{\mathcal{S}}}
\nc{\CT}{{\mathcal{T}}}
\nc{\CU}{{\mathcal{U}}}
\nc{\CV}{{\mathcal{V}}}
\nc{\CW}{{\mathcal{W}}}
\nc{\CX}{{\mathcal{X}}}
\nc{\CY}{{\mathcal{Y}}}
\nc{\CZ}{{\mathcal{Z}}}
\nc{\blambda}{\boldsymbol{\lambda}}
\nc{\bmu}{\boldsymbol{\mu}}

\nc{\fa}{{\mathfrak{a}}} \nc{\fA}{{\mathfrak{A}}}
\nc{\fb}{{\mathfrak{b}}} \nc{\fB}{{\mathfrak{B}}}

\nc{\fd}{{\mathfrak{d}}} \nc{\fD}{{\mathfrak{D}}}
\nc{\fe}{{\mathfrak{e}}} \nc{\fE}{{\mathfrak{E}}}
\nc{\ff}{{\mathfrak{f}}} \nc{\fF}{{\mathfrak{F}}}
\nc{\fg}{{\mathfrak{g}}} \nc{\fG}{{\mathfrak{G}}}
\nc{\fh}{{\mathfrak{h}}} \nc{\fH}{{\mathfrak{H}}}
\nc{\fj}{{\mathfrak{j}}} \nc{\fJ}{{\mathfrak{J}}}
\nc{\fk}{{\mathfrak{k}}} \nc{\fK}{{\mathfrak{K}}}
\nc{\fl}{{\mathfrak{l}}} \nc{\fL}{{\mathfrak{L}}}
\nc{\fm}{{\mathfrak{m}}} \nc{\fM}{{\mathfrak{M}}}
\nc{\fn}{{\mathfrak{n}}} \nc{\fN}{{\mathfrak{N}}}
\nc{\fp}{{\mathfrak{o}}} \nc{\fO}{{\mathfrak{O}}}
\nc{\fq}{{\mathfrak{q}}} \nc{\fQ}{{\mathfrak{Q}}}
\nc{\fR}{{\mathfrak{R}}}
\nc{\fs}{{\mf{s}}} \nc{\fS}{{\mf{s}}}

\nc{\fu}{{\mathfrak{u}}} \nc{\fU}{{\mathfrak{U}}}
\nc{\fv}{{\mathfrak{v}}} \nc{\fV}{{\mathfrak{V}}}
\nc{\fw}{{\mathfrak{w}}} \nc{\fW}{{\mathfrak{W}}}
\nc{\fx}{{\mathfrak{x}}} \nc{\fX}{{\mathfrak{X}}}
\nc{\fy}{{\mathfrak{y}}} \nc{\fY}{{\mathfrak{Y}}}
\nc{\fz}{{\mathfrak{z}}} \nc{\fZ}{{\mathfrak{Z}}}
\nc{\fsl}{{\mathfrak{sl}}}
\nc{\fgl}{{\mathfrak{gl}}}
\nc{\fso}{{\mathfrak{so}}}
\nc{\fsu}{{\mathfrak{su}}}
\nc{\fsp}{{\mathfrak{sp}}}

\nc{\full}{\on{full}}

\newcommand{\csd}[1]{\sigma^* {#1}^\vee}
\nc{\Qd}{\on{Qd}}
\nc{\Lg}{\on{Lagr}}
\nc{\sym}{\mathfrak{sym}}
\nc{\rlr}{\overset{\longrightarrow}{\underset{\longrightarrow}\longleftarrow}}
\nc{\enh}{\on{enh}}

\newcommand{\qbin}[2]{\begin{bmatrix}{#1}\\ {#2}\end{bmatrix}_q}
\nc{\triv}{\on{triv}}
\nc{\TSp}{\on{TSp}}
\nc{\ex}{\on{ex}}
\nc{\Fr}{\on{Fr}}
\nc{\rss}{\on{rss}}
\nc{\kbar}{\ol{k}}
\nc{\ulQl}{\ul{{\BQ_{\ell}}}}

\nc{\type}{\on{type}}

\renewcommand\a\alpha
\renewcommand\b\beta
\newcommand\g\gamma
\renewcommand\d\delta
\newcommand\D\Delta

\newcommand{\s}{\sigma}

\newcommand\dm{\diamondsuit}

\newcommand\sh{\sharp}

\newcommand\cA{\mathcal{A}}
\newcommand\cB{\mathcal{B}}
\newcommand\cC{\mathcal{C}}
\newcommand\cD{\mathcal{D}}
\newcommand\cE{\mathcal{E}}
\newcommand\cF{\mathcal{F}}

\newcommand\cL{\mathcal{L}}
\newcommand\cM{\mathcal{M}}

\newcommand\cO{\mathcal{O}}

\newcommand\cZ{\mathcal{Z}}

\newcommand{\rH}{\ensuremath{\mathrm{H}}\xspace}

\newcommand{\rR}{\ensuremath{\mathrm{R}}\xspace}

\newcommand{\Qll}[1]{\Q_{\ell, #1}}

\newcommand{\BunU}{\Bun_{\U^\dagger(1)}}

\newcommand{\BunUN}{\Bun_{\U^\star(1)}}
\newcommand{\HkU}{\Hk_{\U^\dagger(1)}}
\newcommand{\HkUN}{\Hk_{\U^\star(1)}}

\DeclareMathOperator{\univ}{univ}
\newcommand{\Ud}{\U^\dagger(1)}
\newcommand{\UdN}{\U^\star(1)}
\newcommand{\detd}{{\det}^{\dagger}}
\newcommand{\PrymN}{\Prym_{\mf{N}}}
\DeclareMathOperator{\vir}{vir}
\newcommand{\Sub}{\operatorname{Sub}}
\newcommand{\act}{\operatorname{act}}

\newtheorem{thm}{Theorem}[subsection]
\newtheorem{lemma}[thm]{Lemma}
\newtheorem{prop}[thm]{Proposition}
\newtheorem{cor}[thm]{Corollary}

\newtheorem{conj}[thm]{Conjecture}
\newtheorem{defn-prop}[thm]{Definition-Proposition}

\theoremstyle{remark}
\newtheorem{remark}[thm]{Remark} 
\newtheorem{defn}[thm]{Definition}

\newtheorem{example}[thm]{Example}

\newtheorem{definition}[thm]{Definition}

\makeatletter
\def\th@remark{%
  \thm@headfont{\bfseries}%
  \normalfont 
  \thm@preskip \thm@preskip 
  \thm@postskip\thm@preskip
}
\def\imod#1{\allowbreak\mkern5mu({\operator@font mod}\,\,#1)}
\makeatother

\title[Higher Siegel--Weil formula for unitary groups II]{Higher Siegel--Weil formula for unitary groups II: \\ corank one terms}

\author{Tony Feng}
\address{University of California Berkeley, Department of Mathematics, Berkeley, CA 94720, USA}
\email{fengt@berkeley.edu}

\author{Benjamin Howard}
\address{Boston College, Department of Mathematics, 21 Campanella Way Fifth floor, Chestnut Hill, MA 02467}
\email{howardbe@bc.edu}

\author{Mikayel Mkrtchyan}
\address{Massachusetts Institute of Technology, Department of Mathematics, 77 Massachusetts Avenue, Cambridge, MA 02139, USA}
\email{mikayelm@mit.edu}

\begin{document}

\begin{abstract}We prove the higher Siegel--Weil formula for \emph{corank one} terms, relating (1) the $r^{\rm th}$ central derivatives of corank one Fourier coefficients of Siegel--Eisenstein series, and (2) the degrees of special cycles of virtual dimension 0 on the moduli stack of Hermitian shtukas with $r$ legs. Notably, the formula holds for all $r$, regardless of the order of vanishing of the Eisenstein series. This extends earlier work of Feng--Yun--Zhang, who proved the higher Siegel--Weil formula for the \emph{non-singular} (corank zero) terms. 
\end{abstract}

\maketitle

\tableofcontents

\section{Introduction}

The  \emph{arithmetic Siegel--Weil formula}, conjectured by Kudla and Kudla--Rapoport \cite{Kud04, KRII}, and proved by Li--Zhang \cite{LZ1}, relates degrees of $0$-cycles on integral models of unitary Shimura varieties to central (first)  derivatives of Eisenstein series on  unitary groups over number fields.

A  \emph{higher  Siegel--Weil formula} proposed in \cite{FYZ} relates degrees of $0$-cycles on  moduli spaces of Hermitian shtukas to central derivatives (of all orders) of  coefficients of  Eisenstein series on  unitary groups over  function fields.
The main result of \cite{FYZ} is a proof of this  formula for  non-singular coefficients.

The purpose of this paper is to extend the results of \cite{FYZ} to include some singular coefficients; namely the coefficients of \emph{corank one}, which are the least singular among all singular coefficients.
As an application, we prove a formula relating intersection multiplicities  of special cycles to higher  derivatives of automorphic $L$-functions.

\subsection{Towards a higher Siegel--Weil formula}
\label{ss: intro 1}

Let $\nu \co X' \to X$ be a finite \'etale double cover of smooth, projective, geometrically connected curves over a finite field $k$ of characteristic $p \neq 2$,\footnote{The assumption $p\neq 2$ is made throughout the predecessor paper \cite{FYZ}, so we follow it here.} and denote by $\sigma  \in \Aut(X'/X)$ the nontrivial automorphism.

For every integer $r \geq 0$, there is an associated moduli space $\Sht_{\U(n)}^r$ of rank $n$ Hermitian shtukas on $X$, defined in \cite[\S 6]{FYZ}. 
It is a smooth Deligne-Mumford stack over $k$ of dimension $rn$.

Let $\cE$ be a vector bundle on $X'$ of rank $m>0$,  and let $a: \cE \to \sigma^* \cE^\vee$  be a vector bundle map satisfying the Hermitian condition 
$\sigma^*a^\vee = a$.   Here $\cE^\vee$ is the Serre dual of $\cE$, in the sense of  \S \ref{ss:notation}.
 In \cite{FYZ2} one finds the construction of a  \emph{special cycle} 
 \[
 \cZ_{\cE}^r(a)\to \Sht_{\U(n)}^r,
 \]
 finite and unramified over $ \Sht_{\U(n)}^r$,   and an associated  virtual fundamental class 
\begin{equation}\label{intro virtual}
[\cZ_{\cE}^r(a)]^{\vir} \in \Ch_{ r(n-m)} (  \cZ_{\cE}^r(a)  ) 
\end{equation}
in the Chow group of cycles of dimension $r(n-m)$.  
Using the pushforward of cycle classes,  this virtual fundamental class can be viewed as a cycle on $\Sht_{\U(n)}^r$ of codimension $rm$.

When $\cE$ has rank $n$ and the Hermitian morphism $a:\cE \to \sigma^*\cE^\vee$ is injective (so $\rank(a)= n$, or in other words $a$ is \emph{non-singular}), 
the virtual $0$-cycle class \eqref{intro virtual} is more elementary to construct and was defined earlier in \cite{FYZ}. In \cite[Theorem 1.1]{FYZ}, its degree was calculated and expressed in terms of the $r^\mathrm{th}$ derivative of the $(\cE,a)$-Fourier coefficient of an unramified Siegel Eisenstein series  on the rank $2n$ quasi-split unitary group $H_n=\U(n,n)$  over the field of rational functions on $X$. This constitutes the non-singular part of a \emph{higher Siegel--Weil formula}.

In this paper, we study some of the \emph{singular terms}, where the vector bundle $\cE$ still has rank $n$, but the Hermitian morphism  $a: \cE\to \sigma^*\cE^\vee$ now has rank $n-1$ (i.e., ``corank one'').  In other words, the kernel of $a$ is a line bundle on $X'$. We prove a higher Siegel--Weil formula, again relating the degree of the $0$-cycle \eqref{intro virtual} to the $r^\mathrm{th}$-derivative of the $(\cE,a)$-Fourier coefficient of the same Eisenstein series. To state this corank one higher Siegel--Weil formula, we first introduce some notation. 
\begin{itemize}
\item
For any vector bundle $\cE$ on $X'$, we abbreviate
\begin{equation}\label{eq: d}
d( \cE ) := \frac{   \deg( \cE^\vee )  - \deg(\cE)  }{2}   
=  \mathrm{rank}(\cE)  \deg_X(\omega_X) - \deg(\cE).
\end{equation}
\item 
Let $q=\#k$ be the cardinality of the field of definition of $X$,
 and write $F=k(X)$ and $F'=k(X')$ for the fields of rational functions on $X$ and $X'$.
 \item
Denote by  $\eta : \A_F^\times \to \{ \pm 1\}$  the quadratic character determined by $F'/F$, and let 
 \begin{equation}\label{intro hecke}
 \chi : \A_{F'}^\times \to \C^\times
 \end{equation}
  be an unramified Hecke character whose restriction $ \chi_0=\chi|_{ \A_F^\times}$ is a power of $\eta$.
  We may regard $\chi : \Pic(X') \to \C^\times$ as a character of the group of line bundles on $X'$,   using the conventions of \cite[\S 2.6]{FYZ}.
 \item
 For any integer $m>0$, let $E(g,s,\chi)_m$ be the unramified Siegel Eisenstein series  on the quasi-split unitary group $H_m=\U(m,m)$ over $F$. 
 For a vector bundle $\cE$ of rank $m$ on $X'$ and a Hermitian morphism $a : \cE \to \sigma^* \cE^\vee$, denote by $E_{ ( \cE,a) }(s,\chi)_m$ its $( \cE,a)$-Fourier coefficient.  The precise definitions can be found  in \S \ref{ss:eisenstein}.
 \item
 For any $m>0$ set
\begin{equation}\label{eq: intro L-factor}
\sL_m(s,\chi_0) = \prod_{i=1}^m L(2s+i , \eta^{i-m}\chi_0 ).
\end{equation}
\end{itemize}

Our main result, stated in the text as Theorem \ref{thm:main corank one},  is as follows. 

\begin{thm}[Higher Siegel--Weil formula for corank one terms]\label{thm: intro main} 
Suppose   $\cE$ is a vector bundle on $X'$ of rank $n$,  and $a: \cE \to \sigma^* \cE^\vee$ is a Hermitian morphism of rank $n-1$.  
The stack $\cZ_{\cE}^r(a)$ is proper over $k$, and  the $0$-cycle class \eqref{intro virtual} has degree
\[
\deg\, [\cZ_{\cE}^r(a)]^{\vir}  
=    \frac{1}{ (\log q)^r } \cdot \frac{ q^{  \frac{n}{2} d( \cE)   } }{ \chi(  \det( \cE )) }
\cdot \frac{d^r}{ds^r} \Big|_{s=0} \left(   q^{ns \deg_X(\omega_X) }  \sL_n(s,\chi_0)  E_{ ( \cE,a) } (s,\chi)_n  \right)  
\]
for any unramified Hecke character \eqref{intro hecke} satisfying $\chi_0 = \eta^n$.
\end{thm}

We emphasize that Theorem \ref{thm: intro main} holds for \emph{every} $r$, not just for the leading term in the Taylor expansion of the right hand side of the formula.

Theorem \ref{thm: intro main} is a higher derivative version of the classical Siegel--Weil formula.
The historical context for the Siegel--Weil formula and its variants and generalizations is discussed in \cite[\S 1]{FYZ}. 
As described in \cite[\S 1.1]{FK}, one motivation for a higher Siegel--Weil formula is to relate special cycles on moduli spaces of shtukas to higher derivatives of automorphic $L$-functions, following ideas in the number field setting laid out in \cite{Kud04} and \cite{LL1}, for example.  
See \S \ref{ss:intro applications} below for results in this direction.

\begin{remark}\label{rem: general formulation}
If one replaces the assumption that $a$ has rank $n-1$ with the assumption that $a$ has rank $n$, Theorem \ref{thm: intro main} is precisely the main result of \cite{FYZ}.  One would like to conjecture that this formula holds without any restriction on the rank of $a$, but properness 
 of $ \cZ_{\cE}^r(a)$ is needed to  make sense of the degree  map
\[
\deg: \Ch_0  (  \cZ_{\cE}^r(a)  )  \to \Q.
\]
When $a : \cE \to \sigma^*\cE^\vee$ has corank greater than $1$, no such properness result is known or expected.  
This presents a serious obstruction to  formulating, let alone proving, the analogue of Theorem \ref{thm: intro main} in corank $>1$.
\end{remark}

\begin{remark}
Because of a functional equation on the analytic side, and the emptiness of $\Sht_{\U(n)}^r$ on the geometric side, both sides of Theorem \ref{thm: intro main} vanish when $r$ is odd. 
By incorporating similitude twists as in \cite[\S 2.3]{FYZ2}, one can obtain variants of the theorem when $r$ is odd in which the equality is more interesting than $0=0$: see Section \ref{sec: simil}.
\end{remark}

Now let us suppose that $\cE$ is a vector bundle on $X'$ of rank $n-1$, and that $a: \cE \to \sigma^* \cE^\vee$   is an injective Hermitian morphism.  
In particular, \eqref{intro virtual} is now a cycle class of dimension $r$.  There is a collection of \emph{tautological bundles} $\ell_1 , \ldots, \ell_r$ on   $\Sht^r_{\U(n)}$, defined in \S \ref{ss:GLdegree}.   Using the Chern class map
\[
c_1 : \Pic( \Sht^r_{\U(n)} ) \to \mathrm{CH}^1 (  \Sht^r_{\U(n)}  )
\]
and the intersection pairing
\begin{equation}\label{intro pairing}
\mathrm{CH}^r (  \Sht^r_{\U(n)}  )  \times \mathrm{CH}_r ( \cZ_{\cE}^r(a) ) 
\to \mathrm{CH}_0 ( \cZ_{\cE}^r(a) )
\end{equation}
 of \cite[\S 7.7]{FYZ}, we  form the $0$-cycle class
 \begin{equation}\label{intro tautological intersection}
 \left( \prod_{i=1}^r c_1( \ell_i) \right) \cdot  [\cZ_{\cE}^r(a)]^{\vir}
 \in \mathrm{CH}_0 ( \cZ_{\cE}^r(a) ).
 \end{equation}

Our second main result  expresses the degree of the $0$-cycle \eqref{intro tautological intersection} in terms of an \emph{off-center} $r^\mathrm{th}$  derivative of an Eisenstein series on the lower rank unitary group $H_{n-1}=\U(n-1,n-1)$.

\begin{thm}\label{thm:intro off-center}
If   $\cE$ is a vector bundle on $X'$ of rank $n-1$,  and $a: \cE \to \sigma^* \cE^\vee$   is an injective Hermitian morphism, then 
the stack $\cZ_{\cE}^r(a)$ is proper over $k$. 
If $r$ is odd, then the degree of  \eqref{intro tautological intersection} is $0$. If $r$ is even, then the degree of  \eqref{intro tautological intersection} is
  \begin{equation}\label{intro shifted eisenstein}
 \frac{ 2  \cdot  q^{  \frac{n}{2}   [ d(\cE) +  \deg_X(\omega_X)  ]  }   }{  \chi(\det(\cE)) \cdot  (\log q)^r  }
  \frac{d^r}{ds^r}\Big|_{s=0}
  \Big(      q^{     ns   \deg_X(\omega_X)   } 
   \mathscr{L}_n (s , \chi_0 )    E_{(\cE  , a  ) } (  s + 1/2    ,\chi )_{n-1}  \Big)
\end{equation}
for any unramified Hecke character \eqref{intro hecke} satisfying $\chi_0 = \eta^n$.
\end{thm}

We will deduce Theorem \ref{thm: intro main} as a consequence of Theorem \ref{thm:intro off-center}, using the relation between
 corank one singular  coefficients of $E(g,s,\chi)_n$ and  non-singular coefficients of $E(g,s,\chi)_{n-1}$ found in Proposition \ref{prop:genus drop}.
Theorem \ref{thm:intro off-center} is actually a special case of the slightly more general Theorem \ref{thm:off-center}, which incorporates an auxiliary line bundle $\cE_0$ on $X'$.  This extra generality is needed in the proof of Theorem \ref{thm: intro main}.

\subsection{Applications}
\label{ss:intro applications}

In \S \ref{ss:applications} we give some applications of Theorem \ref{thm: intro main}.
These include an intersection formula for cycles on $\Sht^r_{\U(2)}$,  loosely in the spirit of the Gross--Zagier formula on heights of Heegner points, but heavily influenced by ideas of Kudla.
We explain this intersection formula here.

A special case of the Modularity Conjecture of \cite[Conjecture 4.15]{FYZ2} predicts that the middle codimension cycle classes
\[
[\cZ^r_{\cE}(a) ]^{\vir} \in \mathrm{CH}^r (  \Sht^r_{\U(2)}  ),
\]
indexed by pairs $(\cE,a)$ consisting of a line bundle $\cE$ on $X'$ and a Hermitian morphism $a:\cE \to \sigma^* \cE^\vee$, are the Fourier coefficients of an unramified automorphic form on the quasi-split unitary group $H_1=\U(1,1)$ over $F$,  valued in $\mathrm{CH}^r (  \Sht^r_{\U(2)}  )$.  
The precise statement is Conjecture \ref{conj:modularity}, whose proof will appear in the forthcoming work \cite{FYZ5}.  

Suppose $\pi \subset \mathcal{A}^\mathrm{cusp}(H_1)$ is an unramified irreducible cuspidal automorphic representation.
 Assuming Conjecture \ref{conj:modularity}, one can associate to any unramified  $f\in \pi$ an \emph{arithmetic theta lift}  
 \[
\vartheta^{r,\chi}(f)  \in  \Ch^r( \Sht^r_{\mathrm{U}(2)} ) _\C
\]
in the sense of Kudla \cite{Kud04}; see \eqref{arith theta} for the precise definition.  Here $\chi$ is any unramified Hecke character of $\A_{F'}^\times$ whose restriction to $\A_F^\times$ is trivial.

One might like to express the degree of the self-intersection
$\vartheta^{r,\chi}(f) \cdot \vartheta^{r,\chi}(f)$ in terms of the $r^\mathrm{th}$ central derivative of the twisted base-change $L$-function $ L( s   ,  \mathrm{BC}(\pi) \otimes \chi )$.  
Such a formula would be the higher derivative, function field analogue of the \emph{arithmetic Rallis inner product formula} 
for unitary Shimura varieties, 
 proposed by Kudla \cite{Kud04} and
studied in detail in \cite{Liu11,LL1,LL2}.
Unfortunately, the lack of properness of $\Sht^r_{\mathrm{U}(2)}$ 
means that the degrees of  $0$-cycles are not constant on rational equivalence classes, and so, just as in Remark \ref{rem: general formulation},  the degree of the self-intersection is not well-defined.

However,  one can make sense of  the degree of the intersection of $\vartheta^{r,\chi}(f)$ against any  class
in the Chow group $\Ch^r_c( \Sht^r_{\mathrm{U}(2)} )$ with proper support.
By the properness claim of  Theorem \ref{thm: intro main}, this includes all middle codimension special cycles $[\cZ^r_{\cE}(a)]^{\vir}$ for which $a:\cE \to \sigma^* \cE^\vee$ is injective.

The following consequence of Theorem \ref{thm: intro main} is stated in the text as Corollary \ref{cor:GZformula} (but the provided proof of Corollary \ref{cor:GZformula} is essentially a reference to the relevant parts of \cite{CH}, which contains the details of the argument).  
Again, this Corollary is currently conditional on a special case of the Modularity Conjecture posed in \cite{FYZ2}, whose proof  is expected  to appear in the forthcoming work \cite{FYZ5}. 

 \begin{cor}
Assume Conjecture \ref{conj:modularity}. Fix a pair $(\cE,a)$ consisting of a line bundle $\cE$ on $X'$ and an injective Hermitian morphism $a : \cE \to \sigma^* \cE^\vee$.    For any unramified $f \in \pi$, we have the equality
\begin{align*}   
 \deg \big(   \vartheta^{r,\chi} (f) \cdot [ \cZ^r_{\mathcal{E}}(a)]^{\vir} \big)   =
  f_{( \cE , -a)  }   \cdot   \frac{    q^{ d(\cE) }    }{ (\log q)^r} \cdot 
    \frac{d^r}{ds^r}\Big|_{s=0}\left( q^{  2s \deg_X ( \omega_X)  }    L( s +1/2  ,  \mathrm{BC}(\pi) \otimes \chi )\right),
\end{align*}
where $f_{(\cE,-a)}$ is the $(\cE,-a)$-Fourier coefficient of $f$.
\end{cor}


\subsection{Related results}

The idea of the higher Siegel--Weil formula was inspired by the \emph{arithmetic Siegel--Weil formula} conjectured by Kudla and Kudla--Rapoport \cite{Kud04, KRII}, relating degrees of $0$-cycles on integral models of unitary Shimura varieties to (first) derivatives of Eisenstein series.

For non-singular terms, the arithmetic Siegel--Weil formula was proven by Li--Zhang in \cite{LZ1}. Recently, Ryan Chen \cite{ChenI, ChenII, ChenIII, ChenIV} has proved the arithmetic Siegel--Weil formula for corank one terms, by studying the limiting behavior of the local corank zero results of Li--Zhang towards the singular terms. 

The higher Siegel--Weil formula for non-singular terms was proved in \cite{FYZ}. The proof here of Theorem \ref{thm: intro main},  covering the corank one terms, builds on the constructions and ideas of \emph{op. cit.}, although it does not directly use the higher Siegel--Weil formula of \emph{op. cit.}.  
In particular, \emph{our strategy is completely distinct from Chen's}, who deduces corank one results as consequences of known corank zero results. We will describe the strategy of proof below.


\subsection{Outline of the strategy}

As noted above, Theorem \ref{thm: intro main} as a consequence of Theorem \ref{thm:intro off-center}, and  we sketch some of the ideas that go into the proof of the latter.

The basic format of our strategy resembles that of \cite{FYZ}, which is in turn inspired by the earlier work of Yun--Zhang \cite{YZI}. 
The idea is to realize the degree of the $0$-cycle \eqref{intro tautological intersection} as the trace of a certain cohomological correspondence acting on the cohomology of a space $\cM_{(\cE,a)}$ using the Grothendieck--Lefschetz trace formula, and then to package the spaces $\cM_{(\cE,a)}$ for varying $(\cE,a)$ as the fibers of a map 
\begin{equation}\label{intro hitchin}
 \pi :  \cM \rightarrow \cA.
\end{equation}
We refer to this map as a \emph{Hitchin fibration}, although it is not the fibration that was studied by Hitchin; the language (which follows \cite{YZI}) is chosen to evoke an analogy to Ng\^{o}'s proof \cite{Ngo10} of the Fundamental Lemma.

In this way, the function sending $(\cE,a) \in \cA(k)$ to the degree of  \eqref{intro tautological intersection} is realized as the trace function associated to a certain endomorphism of the complex $\rR \pi_! \Ql$ on $\cA$. To analyze this trace, we decompose 
$\rR \pi_! \Ql$ into (shifted) perverse sheaves. The simple summands, called intersection complexes, can be understood in terms of representation-theoretic data, and then their trace functions calculated and related to the analytic side.

\subsubsection{Comparison to \cite{FYZ}} So far, the strategy has been described in sufficiently vague terms that it applies equally well to the one used in \cite{FYZ} to compute the degree of the special cycle class $[\cZ_\cE^r(a)]^{\vir}$ in the \emph{non-singular case}, where $\rank(\cE)=n$ and $a:\cE\to \sigma^*\cE^\vee$ is injective. Let us now highlight some major \emph{differences} to \cite{FYZ} which arise in the execution of this strategy, thanks to new obstacles such as the following. 
\begin{enumerate}
\item The calculation of the \emph{Hitchin pushforward} $\rR \pi_! \Ql$ in \cite{FYZ} was easy, because $\pi$ was a \emph{small morphism} in that situation. In our case, this step becomes substantially more difficult.
Our solution involves factoring the Hitchin fibration as a composition of two more manageable  (but still not small) maps, and using this factorization to prove an appropriate \emph{Support Theorem}. This will be discussed further below. 

\item The relevant \emph{Hecke operators} on $\rR \pi_!\Ql $ are more complicated in our case, because they must  incorporate the tautological bundles $\ell_1,\ldots, \ell_r$ appearing in \eqref{intro tautological intersection}, and their calculation is correspondingly more involved.

\item The \emph{geometrization of the singular Fourier coefficients}, as trace functions associated to perverse sheaves, involves a more general class of Springer sheaves than appeared in \cite{FYZ}. 
This is because of the restriction $\chi_0 = \eta^n$  in Theorem \ref{thm:intro off-center} on the character  used in the construction of the Eisenstein series $E(s,\chi)_{n-1}$.
The mismatch in parity between $n$ and $n-1$ means that the nonsingular Fourier coefficients appearing in Theorem \ref{thm:intro off-center} are different from those calculated in \cite{FYZ}, and appear as the trace functions of  different perverse sheaves.  
\end{enumerate}

Each of these issues will now be discussed in greater detail. 

\subsubsection{The Hitchin pushforward}
As mentioned above, the Hitchin fibration 
\[
\pi_{\FYZ} : \cM_{\mrm{FYZ}} \rightarrow \cA_{\mrm{FYZ}}
\]
appearing in \cite{FYZ} was small, so one knew immediately that $\rR \pi_{\mrm{FYZ}!} \Ql$ was an intersection complex with full support on $\cA_{\mrm{FYZ}}$. Thanks to the full support property, one can  describe $\rR \pi_{\FYZ!} \Ql$ in terms of the generic geometry of $\pi_{\FYZ}$. Using this, \cite{FYZ} expressed $\rR \pi_{\FYZ!} \Ql$ in terms of Springer theory. 

In the setting of  Theorem \ref{thm:intro off-center}, where $\rank(\cE)=n-1$, the corresponding Hitchin fibration $\pi:  \cM \rightarrow \cA$ is not small, and the determination of $\rR \pi_! \Ql$ is less straightforward.  
 An essential  observation is that $\pi$ can be enhanced in a natural way to a map
\begin{equation}\label{intro enhanced hitchin}
f :  \cM \to  \cA \times \Bun_{\Ud} 
 \end{equation}
 whose geometry is more tractable,
  where $\Bun_{\Ud}$ is the twisted form of the moduli stack of rank $1$ Hermitian bundles on $X'$  defined in \S \ref{sec: Hitchin space}. The morphism $f$ is still not small, but it is more controllable than the original $\pi$, thanks to a crucial idea (suggested to us by Zhiwei Yun) called the \emph{complementary line trick}, explained in \S \ref{ss:cline}. (We remark that the geometry of our situation resembles that of the non-singular terms in the symplectic/orthogonal case, so our arguments may have broader scope for applications.) For example, we show that there are  natural decompositions $\cA = \bigsqcup_{d \ge 0} \cA_d$ and $\cM = \bigsqcup_{d \ge 0} \cM_d$ as disjoint unions of open and closed substacks  in such a way that $f$ restricts to a map
\[
f_d : \cM_d \to \cA_d \times  \Bun_{\Ud}
\] 
that is generically a $( \PP^1)^d$-bundle. We remark that the appearance of each $\PP^1$ has an analogue in the number field setting: see \cite[\S 10.2, \S 10.4, Remark 10.2.1]{LZ1}.

The appearance of the enhanced Hitchin fibration \eqref{intro enhanced hitchin} has an interpretation on the analytic side of Theorem \ref{thm:intro off-center}.  
Proposition \ref{prop:key degree} (and its proof) tells us that when $r$ is even, the expression \eqref{intro shifted eisenstein} can be rewritten as 
\begin{equation}\label{LxDen}
 \frac{ 2 }{  (\log q)^r  }
 \frac{d^r}{ds^r}\Big|_{s=0}
  \Big(  q^{ s \deg_X(\omega_X)  + s d(\cE) }  
 L(2s,\eta)  \Den_\eta( q^{ 1- 2 s}  ,\cE) \Big) ,
\end{equation}
where $\Den_\eta(T,\cE)$ is the twisted density polynomial defined in \S \ref{ss:density polynomials}.
In some sense, it is the geometry of the factor $\cA$ in the enhanced Hitchin base \eqref{intro enhanced hitchin} that controls these twisted density polynomials, while the geometry of the factor $\Bun_{\Ud}$  controls the Dirichlet $L$-function $L(s,\eta)$.

Although the generic geometry of the enhanced Hitchin fibration \eqref{intro enhanced hitchin} is explicit, 
the lack of smallness makes the analysis of 
$\rR f_! \Ql$ quite subtle. The Decomposition Theorem implies that it must be a direct sum of shifted intersection complexes supported on closed substacks, but a priori some of these supports could have positive codimension.  
A key step in our arguments is a ``support theorem'' (Theorem \ref{thm: geometric K}) showing that this does not happen, and so $\rR f_! \Ql$ is determined by the generic geometry of the map $f$.
This allows us to ultimately  express $\rR f_! \Ql$, and hence also $\rR \pi_! \Ql$, in terms of Springer theory.

 \begin{remark}
In Ng\^{o}'s proof of the Fundamental Lemma \cite{Ngo10}, the proof of the appropriate support theorem 
is the most involved part. Our situation is not as subtle, but certainly goes beyond the simple situation of a small morphism. The proof relies on a comparison to the split case, and in that case a new geometric criterion (developed in \S \ref{sec: support}) for a derived direct image to consist of only full support perverse sheaves. 

\end{remark}

\subsubsection{Hecke action} After $\rR \pi_! \Ql$ is adequately characterized, we need describe the endomorphism of it whose trace  on the stalk at $(\cE,a) \in \cA(k)$ computes the degree of the $0$-cycle \eqref{intro tautological intersection}.  This endomorphism of $\rR \pi_! \Ql$  comes from Hecke correspondences on $\cM$ defined in \S \ref{sec: Hecke}. The complicated geometry of those Hecke correspondences makes them difficult to analyze directly. Our solution is to decompose the Hecke operator into pieces, each of which descends to some simpler geometric situation. 

More specifically, we are able to decompose the Hecke operator as a sum of three pieces:
one pulled back from the $\Bun_{\Ud}$ factor in the enhanced Hitchin base \eqref{intro enhanced hitchin}, another pulled back from the $\cA$ factor, and a third pulled back from their product. 
This decomposition roughly mirrors the process of computing the derivative in \eqref{LxDen} using the product rule from freshman calculus. 

The calculation of the Hecke operator is carried out in \S \ref{sec: Hecke} -- \S \ref{sec: geometric side}, and is quite intricate. A delicate analysis of the geometry allows us to descend the calculation of each of the three pieces to the base from which it is pulled back. That breaks the problem down into simpler components, which we study in turn. 

For the piece which is pulled back from $\cA$, the Hecke action is calculated by using the full support property to reduce to the regular semisimple locus, and explicitly analyzing the geometry there. Such a step also occurs \cite{FYZ}, although the geometry itself is quite different. This analysis is carried out in \S \ref{sec: comparison correspondence}. 

The remaining two pieces, which are treated in \S \ref{sec:univ-action} and \S \ref{sec:taut-action}, do not have counterparts in \cite{FYZ}; they are a new feature of this paper. The idea is to explicitly calculate the Hecke action in terms of the Atiyah--Bott description of the cohomology of $\BunU$. The method is inspired by forthcoming work \cite{volume} on arithmetic volumes for general moduli stacks of shtukas.

\subsubsection{Geometrization of twisted density polynomials}
 The calculation of $\rR \pi_! \Ql$ and its Hecke endomorphisms completes the analysis of the geometric side. We then need to relate the answers to the analytic side.

In the analogous step of \cite{FYZ}, the major work was to express local density polynomials, which are used to describe the non-singular Fourier coefficients of the Eisenstein series $E(g,s,\chi)_n$ on $H_n$, in terms of Springer theory  via the sheaf-function correspondence. 
For this the constraint $\chi_0=\eta^n$ on the Hecke character $\chi$, assumed throughout \cite{FYZ}, is important.

If we maintain the constraint $\chi_0=\eta^n$ but work with the Eisenstein series $E(g,s,\chi)_{n-1}$ on the lower rank group $H_{n-1}$, as we do in Theorem \ref{thm:intro off-center}, then the non-singular Fourier coefficients take a different shape, and are no longer expressed in terms of the usual density polynomials.
Instead, they are expressed in terms of  the twisted density polynomials already mentioned above in \eqref{LxDen}, and defined in \S \ref{ss:density polynomials}.  We show that these  twisted density polynomials are again related, via the sheaf-function correspondence, to Springer theory, but are expressed in terms of different Springer sheaves than the ones appearing in \cite{FYZ}.
This analysis is carried out in \S \ref{s:density geometrization}.

\subsubsection{Conclusion of the proofs} 
In \S \ref{ss:first main theorem}  we complete the proof of Theorem \ref{thm:intro off-center} by simply comparing 
Theorem \ref{thm: geometric side} with Theorem \ref{thm: maindensityresult}.  The first of these expresses the degree of the $0$-cycle \eqref{intro tautological intersection} in terms of Frobenius traces on stalks of certain Springer sheaves on $\cA$, while the second expresses the higher derivatives of the twisted density polynomials in the same form.

Switching to Theorem \ref{thm: intro main}, so that  $\cE$ has rank $n$ and $a : \cE \to \sigma^*\cE^\vee$ has corank one, the deduction of that theorem from 
Theorem \ref{thm:intro off-center} exploits the  structure of $[\cZ_{\cE}^r(a)]^{\vir}$ as the product of Chern classes of tautological bundles with a non-singular factor  $[\cZ_{\cE^\flat}^r(a^\flat)]^{\vir}$ defined by an injective Hermitian morphism $a^\flat:\cE^\flat \to \sigma^* (\cE^\flat)^\vee$ on a vector bundle $\cE^\flat$ of rank $n-1$.  More precisely (see Lemma \ref{lem: special cycle digest})
\begin{equation}\label{intro singular cycle decomp} 
 [\cZ_{\cE}^r(a)]^{\vir} =
\left( \prod_{i=1}^r c_1(p_i^* \sigma^* \cE_0^{-1} \otimes \ell_{\U(n),i}) \right) \cdot  [\cZ_{\cE^\flat}^r(a^\flat)]^{\vir},
 \end{equation}
 where $\cE_0 := \ker(a)$ and $\cE^\flat := \cE/\cE_0$.  We can apply Theorem \ref{thm:intro off-center} to the pair $(\cE^\flat,a^\flat)$, and so obtain an expression for the degree of \eqref{intro singular cycle decomp} in terms of non-singular Fourier coefficients of the Eisenstein series $E(g,s,\chi)_{n-1}$.  These are related to corank one coefficients of $E(g,s,\chi)_n$ by Proposition \ref{prop:genus drop}, and Theorem \ref{thm: intro main} follows easily.

\subsubsection{Final remarks} The decomposition \eqref{intro singular cycle decomp} is not a special feature of our corank one pair $(\cE,a)$.
According to \cite[\S 4]{FYZ}, for general $(\cE,a)$ -- meaning for $a:\cE \to \sigma^* \cE^\vee$ of arbitrary corank --
the class $[\cZ_{\cE}^r(a)]^{\vir}$  is a sum of terms which are similarly products of tautological Chern classes with non-singular factors (i.e., defined by injective Hermitian morphisms). 
At the two extremes, we have 
\begin{itemize}
\item virtual classes which are purely non-singular (no Chern class factor), and 
\item virtual classes which are purely products of Chern classes (no non-singular factor), which are analyzed in \cite{volume}. 
\end{itemize}
The present work shows how to combine the techniques used in these two extremes, in order to calculate virtual classes which have a more general form. We may thus hope that it sets the template for a proof of the general higher Siegel--Weil formula, once the issues confounding the formulation of such a formula are resolved (Remark \ref{rem: general formulation}).


\subsection{Acknowledgments} 

Parts of this research were carried out during 
\begin{itemize}
\item the SLMath program on \emph{Algebraic cycles, Euler systems, and L-functions}, and 
\item the  American Institute for Mathematics workshop  \emph{Arithmetic Intersection Theory on Shimura Varieties} held in January of 2024.  
\end{itemize}
The authors thank SLMath and AIM for hosting stimulating environments for collaborative work. We also thank Ryan Chen, Zhiwei Yun, and Wei Zhang for helpful discussions, as well as comments on a draft of the paper. We are especially grateful to Yun for the suggestion of the complementary line trick.

T.F.~was supported by NSF grants DMS-2302520, DMS-2441922, the Simons Foundation, and the Alfred P. Sloan Foundation. B.H.~was supported in part by NSF grant DMS-2101636. M.M. was supported in part by Prof. Zhiwei Yun's Simons Investigator grant.


\subsection{Notation and conventions} 
\label{ss:notation}


Our notation is as in \cite[\S 1.3]{FYZ}, some of which will be reiterated below. 
In particular, we fix an odd prime $p$ and work over a finite field $k = \F_q$ of characteristic $p$.  
A $k$-algebra $R$ is always assumed to be commutative.

\subsubsection{}
 As above, $X$ is a smooth, proper, geometrically connected curve over $k$, and  $\nu : X' \to X$ is  a finite \'etale double cover.  
 Denote by   $\sigma\in \Aut(X'/X)$  the nontrivial automorphism. 

We assume that $X'$ is also geometrically connected, except in \S \ref{sec: Hitchin space} and \S  \ref{sec: support}, where we also allow $X'=X\sqcup X$.  
This is for a particular technical reason: we need to allow the possibility that $X'$ is split in the statement of Theorem \ref{thm: geometric K}, because reduction to that case is an important part of its proof.

We denote by $\omega_X$ and $\omega_{X'} \iso \nu^* \omega_X$ the canonical bundles on $X$ and $X'$. 
The \emph{Serre dual} of  a vector bundle $\cF$ on $X'$ is
\[
\cF^\vee := \cHom(\cF, \omega_{X'}).
\]

For any $i\in \Z$, we denote by $\omega_{X'}^i := \omega_{X'}^{\otimes i}$  the $i^\mathrm{th}$ tensor power of $\omega_{X'}$. A pairing $h:  \cF \otimes \sigma^*\cF  \to  \omega_{X'}^{i} $ is 
\emph{Hermitian} if it is equal to  the composition
\[
\cF \otimes \sigma^*\cF 
\iso
\sigma^*\cF \otimes  \cF \map{ \sigma^* h} \sigma^* \omega_{X'}^{i}  \iso \omega_{X'}^{i}
\]
(the first isomorphism swaps the factors in the tensor product, and the last is induced by the canonical descent $\omega_{X'} \iso \nu^* \omega_X$).  
A morphism $a : \cF  \to \sigma^*\cF^\vee$ is \emph{Hermitian} if $\sigma^* a^\vee =a$.
Thus a Hermitian morphism $\cF \to \sigma^* \cF^\vee$ is the same as a Hermitian pairing $\cF \otimes \sigma^* \cF \to \omega_{X'}$.


\subsubsection{}
The word  \emph{stack}  always means Artin stack.
For a $k$-stack $Y$,  we denote by $D^b_c(Y)$ the bounded derived category of constructible \'etale $\Q_\ell$-sheaves on $Y$, where $\ell \neq p$ is some auxiliary prime.  For an object $E \in D^b_c(Y)$ and an integer $i\in \Z$, we define the \emph{shear}
\begin{equation}\label{shear def}
E\tw{i} : = E[2i](i)
\end{equation}
as the composition of a shift and a Tate twist.

\subsubsection{} For a constructible \'etale $\Q_\ell$-sheaf $E$ on a $k$-stack $Y$,  we denote by $\rH^*(Y, E)$ the \emph{geometric} cohomology $H^*(Y_{\ol{k} }; E_{\ol{k} })$. (Note the different font.)
When $E=\Ql$ we usually omit it from the notation, so that $\rH^*(Y)$ is a shorthand for $H^*(Y_{\ol{k}}; \Ql)$.
We will remind the reader of this convention whenever confusion is possible.

\subsubsection{} We use the language of \emph{cohomological correspondences} from \cite{FYZ3}. In particular, for a correspondence
\[
\begin{tikzcd}
A_0 & \ar[l, "a_0"'] C \ar[r, "a_1"] & A_1
\end{tikzcd}
\]
and $\sK_0 \in D^b_c(A_0)$, $\sK_1 \in D^b_c(A_1)$ we write 
\[
\Corr_C(\sK_0, \sK_1) := \Hom_{ D^b_c(C) } (a_0^* \sK_0, a_1^! \sK_1).
\]
The notions of \emph{pushforward}, \emph{pullback}, and \emph{trace} of cohomological correspondences are defined in \cite[\S 4]{FYZ3}.


\section{Geometry of Hitchin spaces}\label{sec: Hitchin space}


In this section we define our Hitchin fibration \eqref{intro hitchin} and its enhanced version \eqref{intro enhanced hitchin}, and study their geometry.


\subsection{The Hitchin fibration}
\label{ssec: global Hitchin space} 


Fix an integer $n \geq 1$.

As in \cite[ \S 6.1]{FYZ}, denote by $\Bun_{\U(n)}$ the $k$-stack whose $R$-points  (for a $k$-algebra $R$) are rank $n$ vector bundles $\cF$ on $X'_R$ endowed with a Hermitian isomorphism 
\[
h \co \cF \xrightarrow{\sim} \sigma^* \cF^\vee.
\]

We also need a twisted version of $\Bun_{\U(1)}$.  Following the notation of \cite[\S 3.1.1]{FYZ2}, let 
\begin{equation}\label{BunU def}
\Bun_{ \U^\dagger(1)} := \Bun_{ \U(1) , \omega_X^{1-n} }  
\end{equation}
be the moduli stack whose $R$-points are line bundles $\cL$ on $X'_R$ equipped with a Hermitian isomorphism
$\cL \otimes   \sigma^*\cL \iso   \omega_{X'}^{ 2-n }$.  
The Hermitian structure on any  $\cF\in \Bun_{\U(n)}(R)$ determines  a Hermitian isomorphism
\begin{equation}\label{hermdet}
\det(\cF)  \otimes \sigma^* \det(\cF)  \iso  \omega_{X'}^{n},
\end{equation} 
and formation of the  \emph{twisted determinant}
\begin{equation}\label{eq: twisted det}
{\det}^\dagger( \cF) : = \det (\cF) \otimes \omega_{X'}^{1-n}  
\end{equation} 
defines a morphism
$
{\det}^\dagger : \Bun_{\U(n)} \to \BunU.
$

\begin{definition}
Using the notation of \cite[\S 3.2]{FYZ2}, the \emph{Hitchin space} is 
\begin{equation}\label{hitchin space}
\cM := \cM_{\GL(n-1)', \U(n)}^{\ns} .
\end{equation}
Its functor of points assigns to a   $k$-algebra $R$ the groupoid $\cM(R)$  of triples $(\cE, \cF, t)$ in which
\begin{itemize}
\item $\cE$ is a vector bundle of rank $n-1$ on $X'_R$,
\item $\cF \in \Bun_{\U(n)} (R),$ 
\item $t \co \cE \rightarrow \cF$ is a morphism of coherent sheaves on $X'_R$, such that the composition 
\begin{equation}\label{eq: a(t)}
a \co \cE \xrightarrow{t} \cF \xrightarrow{h} \sigma^* \cF^\vee \xrightarrow{\sigma^* t^\vee} \sigma^* \cE^\vee
\end{equation}
is injective fiberwise over $\Spec (R)$. 
\end{itemize}
\end{definition}

\begin{remark}\label{rem: Mdifferent}
The superscript ``ns'' in \eqref{hitchin space} stands for the non-singular substack, and refers to the injectivity condition in the third bullet point above.
Note that this is stronger than the fiberwise injectivity of $t$ imposed in the definition of $\cM(n-1,n)$ found in \cite[\S 8.1]{FYZ}, and so our stack \eqref{hitchin space} is an open substack of the latter. 
\end{remark}

\begin{defn}
The \emph{Hitchin base} $\cA$ is the stack whose functor of points assigns to a   $k$-algebra $R$ the groupoid $\cA(R)$ of pairs $(\cE , a )$ in which
\begin{itemize}
\item $\cE$ is rank $n-1$ vector bundle on $X_R'$, 
\item $a : \cE \to \sigma^* \cE^\vee$ is a Hermitian morphism  that is injective fiberwise over $\Spec(R)$.
\end{itemize}
\end{defn}

\begin{definition}
The \emph{Hitchin fibration} is the morphism
\begin{equation}\label{eq: basic Hitchin}
\pi \co \cM \rightarrow \cA 
\end{equation}
sending $(\cE,\cF,t) \mapsto (\cE,a)$, where $a: \cE \to \sigma^*\cE^\vee$ is the Hermitian morphism \eqref{eq: a(t)}.
\end{definition}

\begin{definition}
The \emph{enhanced Hitchin fibration} is the  morphism
\begin{equation}\label{eq: Hitchin map}
f \co \cM  \rightarrow \cA \times \BunU
\end{equation}
 sending $( \cE, \cF , t ) \in \cM(R)$ to the triple   $( \cE , a , \cL ) \in  (\cA \times \BunU )(R)$  defined by  \eqref{eq: a(t)} and the twisted determinant 
$\cL :={\det}^\dagger (\cF)$.
\end{definition}

The motivation for forming the enhanced Hitchin fibration  will become clearer when we discuss the \emph{complementary line trick} in \S \ref{ss:cline}.


\subsection{Hermitian torsion sheaves}
\label{ss: hermitian torsion}


The geometry of the enhanced Hitchin fibration \eqref{eq: Hitchin map} is complicated. We will study it by restricting to certain dense open substacks
of the source and target,   defined in this subsection  in terms of moduli spaces of Hermitian torsion sheaves on $X'$.

Recall from \cite[\S 4.1]{FYZ}  the moduli stack $\Herm_d$ of length $d$ Hermitian torsion  sheaves on $X'$. 
For a $k$-algebra $R$, a point of $\Herm_d(R)$ is a coherent sheaf $Q$ on $X'_R$ whose pushforward to $\Spec(R)$ is locally free of rank $d$, together with a nondegenerate  $\cO_{X'}$-sesquilinear pairing 
\[
\langle - , - \rangle : Q \times Q \to \omega_{F'}/\omega_{X'}
\]
satisfying $\langle v_1,v_2\rangle = \sigma^*\langle v_2, v_1 \rangle$.  Here $\omega_{F'}$ is the $F'$-module of rational $1$-forms on $X'_R$, viewed as a  quasi-coherent sheaf.  Other formulations of what it means to give such a pairing on $Q$ can be found in \cite[\S 4.1]{FYZ}.

According to \cite[Lemma 4.3]{FYZ}, $\Herm_d$  is a smooth $k$-stack of dimension $0$.
The unramifiedness of  $X' \to X$  implies that   $\Herm_d = \emptyset$ unless $d$ is even.

\subsubsection{The support morphisms}
We regard the scheme-theoretic quotient 
\[
X_d := S_d \backslash X^d
\]
as the moduli space of degree $d$ divisors on $X$, and similarly with $X$ replaced by $X'$.
As in \cite[\S 4.1]{FYZ}, denote by 
\[
\mf{s}_{2d}' \co \Herm_{2d} \rightarrow X_{2d}'
\]
be the map sending a Hermitian  torsion sheaf to its support divisor. 
Using the Hermitian structure, the support divisor descends to $X_d$, and we denote by  \emph{support morphism} 
\begin{equation}\label{eq: support map}
\mf{s}_{2d} \co \Herm_{2d} \rightarrow X_d
\end{equation}
the composition of $\mf{s}_{2d}'$ with this descent.

\subsubsection{Open-closed substacks of the Hitchin space}

As in \cite[\S 8.4]{FYZ}, there is a morphism 
\[
g \co \cA \rightarrow \Herm:= \bigsqcup_{d \geq 0} \Herm_{2d} .
\]
For a $k$-algebra $R$, this  sends   $(\cE, a) \in \cA (R)$ to the torsion coherent sheaf 
\[
Q:= \coker(a : \cE \to \s^* \cE^\vee)
\]
on $X'_R$,  endowed with its natural Hermitian structure.
For $d \geq 0$, denote by 
\begin{equation}\label{eq: g_d}
g_d \co \cA_d \rightarrow \Herm_{2d} 
\end{equation}
 the base change of $g$ along $\Herm_{2d} \inj \Herm$. Equivalently, using the notation of \eqref{eq: d}, $\CA_d \subset \CA$ is the open and closed substack characterized by the relation $d(\cE) = d$.  Similarly, denote by 
$
 \pi_d \co \cM_d \rightarrow \cA_d
$
  the base change of the Hitchin fibration \eqref{eq: basic Hitchin} along  $\cA_d \inj \cA$, and by 
\begin{equation}\label{eq: Hitchin map d}
 f_d \co \cM_d \rightarrow \cA_d \times \BunU
 \end{equation}
  the base change of the enhanced Hitchin fibration \eqref{eq: Hitchin map}.
 We now have morphisms
 \[
\cM_d \map{f_d} \cA_d  \times \BunU \to \cA_d \map{g_d} \Herm_{2d} \map{ \mf{s}_{2d} }  X_d.
 \]

\subsubsection{Regular semisimple loci}
\label{sss:rss}
Denote by 
$
X_d^\circ \subset X_d
$
the open subscheme of multiplicity-free divisors, and by
 $X^{d ,\circ} \subset X^d$   the multiplicity-free locus; i.e., the pre-image of $X_d^\circ$ under the quotient map $X^d \rightarrow X_d$. The map $X^{d ,\circ } \rightarrow X_d^\circ$ is an $S_d$-torsor.  
Taking the preimages of $X_d^\circ \subset X_d$ under the arrows above defines stacks
\begin{equation}\label{all the rss}
\cM^\circ_d \map{f_d} \cA_d^\circ  \times \BunU \to \cA_d^\circ \map{g_d} \Herm_{2d}^\circ \map{ \mf{s}_{2d} }  X_d^\circ.
 \end{equation}

To spell this out a bit more explicitly, the \emph{regular semisimple locus}
 \[
\Herm_{2d}^\circ \subset \Herm_{2d}
\]
is the space of Hermitian torsion sheaves with multiplicity-free support.
An $R$-point $(\cE,a)$ of the  regular semisimple locus
\[
\cA_d^\circ \subset \cA_d
\]
consists of a rank $n-1$ vector bundle $\cE$ on $X'_R$, and a Hermitian morphism $a:\cE \to \sigma^*\cE^\vee$  whose cokernel is a length $2d$ Hermitian torsion sheaf on $X'_R$ with multiplicity-free support.
The regular semisimple locus   
\[
\cM_d^\circ \subset \cM_d
\]
is the preimage  of $\cA_d^\circ \subset \cA_d$ under the Hitchin fibration \eqref{eq: basic Hitchin}. 
All of these are inclusions of open dense substacks.

\subsection{The complementary line trick}
\label{ss:cline}


Given a $k$-algebra $R$ and an $R$-point $(\cE,a,\cL)$ of $\cA \times \BunU$,  we  define the \emph{complementary line bundle} on $X'_R$ by
\begin{equation}\label{eq: cline def}
\cD := \det(   \sigma^*  \cE)   \otimes   \cL    \in \Bun_{\GL(1)'}(R) . 
\end{equation}
Here,  as in \cite[\S 2]{FYZ},  $\Bun_{\GL(1)'}(R)$ is the groupoid of line bundles on $X'_R$.

\begin{lemma}
The complementary line bundle carries  a canonical Hermitian morphism
\begin{equation}\label{compline herm}
\cD \to \sigma^*\cD^\vee,
\end{equation}
which is injective fiberwise on $\Spec(R)$. 
\end{lemma}

\begin{proof}
Taking the determinant of  $a : \cE \to \sigma^* \cE^\vee$ provides a morphism
\[
\det(\cE) \to \sigma^* \det(\cE)^{-1} \otimes \omega_{X'}^{n-1},
\]
which in turn induces  
\begin{equation}\label{compline a}
\cD = 
\det(   \sigma^*  \cE)  \otimes  \cL 
\to  
\det(   \cE)^{-1}    \otimes \omega_{X'}^{n-1} \otimes   \cL  .
\end{equation}
As $\cL$ comes equipped with an isomorphism 
$
\cL \iso ( \sigma^* \cL^\vee) \otimes \omega_{X'}^{1-n} ,
$
we have
\begin{equation}\label{compline b}
\det(   \cE)^{-1}    \otimes \omega_{X'_R}^{n-1} \otimes   \cL
\iso \det(\cE)^{-1} \otimes \sigma^* \cL^\vee = \sigma^*\cD^\vee.
\end{equation}
The morphism \eqref{compline herm} is defined as the composition of \eqref{compline a} and \eqref{compline b}.
\end{proof}

Now  start with a point  $( \cE ,\cF , t) \in \cM(R)$, and form the complementary line bundle
\[
\cD := \det(   \sigma^*  \cE)   \otimes    {\det}^\dagger (\cF) 
\]
associated to its image under the enhanced Hitchin fibration \eqref{eq: Hitchin map}.

\begin{prop}\label{prop:compline}
From the data  $t \co \cE \rightarrow \cF$,  there is a functorially induced map $\cD \to \cF$.
 This map is an isometric embedding, in the sense that the composition
\[
\cD \to \cF \iso \sigma^*\cF^\vee \to \sigma^* \cD^\vee
\]
agrees with \eqref{compline herm}, and its image is orthogonal  to the image of $\cE \to \cF$.\end{prop}

\begin{proof}
The morphism $t : \cE \to \cF$ induces a morphism 
\[
\cF \otimes \det( \cE) \to   \det(\cF)
\]
 sending $f \otimes e_1\wedge \cdots \wedge e_{n-1} \mapsto f \wedge t(e_1) \wedge \cdots \wedge t(e_{n-1})$. 
We may regard this as a morphism
\[
\det( \cF )^\vee  \otimes \det(\cE) \to \cF^\vee,
\]
and then apply  $\sigma^*$ to obtain the first arrow in 
\begin{equation}\label{eq:compline1}
\det(  \sigma^* \cF) ^\vee   \otimes  \det(  \sigma^* \cE) \to  \sigma^* \cF^\vee \map{ h^{-1} } \cF.
\end{equation}
The desired morphism is the composition
\[
\cD  = \det( \cF )   \otimes  \det(   \sigma^*  \cE)   \otimes \omega_{X'}^{1-n } 
\map{\eqref{hermdet}} 
 \det( \sigma^* \cF) ^\vee   \otimes \det(   \sigma^*  \cE)
 \map{ \eqref{eq:compline1} }  \cF.
\]
The proof that this morphism has the stated properties is a tedious but straightforward  exercise in linear algebra (that can be performed over the generic point of $X'$).
\end{proof}

Proposition \ref{prop:compline} associates to any $(t : \cE \to \cF) \in \cM(R)$  a morphism
\begin{equation}\label{eq: augment cline}
\wt{t} \co     \cE  \oplus  \cD \to \cF 
\end{equation}
between vector bundles on $X'_R$ of the same rank $n$.  This morphism is  injective fiberwise over $\Spec(R)$, and its formation is functorial in $R$.  It is an isometric embedding, in the sense of Proposition \ref{prop:compline}, where  $ \cE  \oplus  \cD$ is  endowed with the orthogonal direct sum of the Hermitian morphisms   \eqref{eq: a(t)} and \eqref{compline herm}.
This canonical upgrading of $t : \cE \to \cF$ to \eqref{eq: augment cline} is the \emph{complementary line trick}, suggested to us by Zhiwei Yun.

The complementary line trick allows us to understand the fibers of the enhanced Hitchin fibration \eqref{eq: Hitchin map},   which we recall sends $(t : \cE \to \cF) \in \cM(R)$ to the triple 
\[
(\cE,a,\cL) \in ( \cA \times  \Bun_{\U^\dagger(1)} ) (R),
\] 
where $\cL :={\det}^\dagger (\cF)$. 
 Given an $R$-valued point $(\cE,a,\cL)$  of the enhanced Hitchin base $ \cA \times  \Bun_{\U^\dagger(1)}$, we may form the complementary line bundle \eqref{eq: cline def}.  
 Giving an $R$-valued point of the Hitchin space $\cM$ above 
$(\cE,a,\cL)$ is then the same as giving a  Hermitian vector bundle $\cF \iso \sigma^* \cF^\vee$ together with an isometric embedding \eqref{eq: augment cline}.
 In other words, the points of $\cM$ above $(\cE,a,\cL)$ are the self-dual   vector bundles lying between the source and target of the Hermitian morphism 
\[
\cE \oplus \cD \to \sigma^*\cE^\vee \oplus \sigma^*\cD^\vee.
\]

\begin{remark}
The vector cross product of freshman physics takes two vectors in  $\R^3$ and returns a third vector orthogonal to them. More generally, the Hodge star takes $n-1$ vectors in an $n$-dimensional inner product space and returns an $n^\mathrm{th}$ vector orthogonal to all of them. The complementary line trick simply imitates this construction in the context of  Hermitian vector bundles.
\end{remark}


\subsection{The local Hitchin space}
\label{ssec: local Hitchin space}


Fix an integer $d \ge 0$.

Our next goal is relate the enhanced Hitchin fibration 
\eqref{eq: Hitchin map d} to a simpler fibration, 
expressed entirely in terms of moduli spaces of  
Hermitian torsion sheaves. The following Lemma is 
essential to this.

\begin{lemma}\label{lem: supp coker det}
Let $R$ be a $k$-algebra. For any $R$-valued point 
$(\cE,a,\cL)$ of the enhanced Hitchin base $\cA_d 
\times \BunU$, the Hermitian torsion sheaves
\begin{equation}\label{eq: PQpair}
Q := \coker( \cE \to \sigma^* \cE^\vee) \quad 
\mbox{and} \quad P := \coker( \cD \to \sigma^* 
\cD^\vee)
\end{equation}
on $X'_R$ both have length $2d$, and have the same image under the support map \eqref{eq: support map}.
Here $\cD$ is the complementary line \eqref{eq: cline def}, endowed with its natural Hermitian structure \eqref{compline herm}. 
\end{lemma}

\begin{proof}
As $Q$ has length $2d$ by definition of $\cA_d$,  it 
suffices to show that $P$ and $Q$ have the same image 
under the support map. 
From the construction of \eqref{eq: cline def}, we see that this is a special case of the more general assertion that given an injective map $\varphi \co \cE \rightarrow \cE'$ between vector bundles on $X'$ of the same rank, $\coker(\varphi)$ and $\coker(\det \varphi \co \det \cE \rightarrow \det \cE')$ have the same support. This assertion can be checked over the completion at each point $v \in X'$ where $\varphi$ fails to be an isomorphism, where we see by elementary linear algebra that the multiplicity of $v$ in the support divisor of $\coker(\varphi)$ is the $v$-adic valuation of $\det \varphi$, as desired. 
\end{proof} 

The Lemma allows us to make the following definition.

\begin{defn}\label{def: double}
The \emph{local enhanced Hitchin base} is 
\[
\cB_d^{\Herm}  := \Herm_{2d} \times_{X_d} \Herm_{2d},
\]
where the fiber product  is taken with respect to the support map \eqref{eq: support map}.
The  \emph{doubling map}
\[
\Db \co \cA_d \times \BunU \rightarrow \cB_d^{\Herm}
\]
sends an $R$-point  $(\cE, a, \cL)$ of the enhanced Hitchin base   to the pair 
\[
\Db(\cE, a, \cL) := ( Q , P ) 
\]
of Hermitian torsion sheaves on $X'_R$ defined by \eqref{eq: PQpair}.  
\end{defn}

\begin{remark}\label{rem: double is g}
The terminology \emph{doubling} refers to the fact that  $\Db$ is roughly a doubled version  of \eqref{eq: g_d}. 
A bit more precisely, we have a commutative diagram
\[
\begin{tikzcd}
{ \cA_d \times \BunU } \ar[r, "\Db" ]   \ar[d]    & {  \cB_d^{\Herm} } \ar[d] \ar[dr] \\ 
{  \cA_d } \ar[r, "g_d"']    & { \Herm_{2d}  } \ar[r, "\mf{s}_{2d}"']  & {  X_d }
\end{tikzcd}
\]
in which the vertical arrow on the right is projection to the first factor in the product defining $ \cB_d^{\Herm}$.
\end{remark}

As in \cite[\S 4.5]{FYZ}, there is a stack $\Lagr_{4d} \rightarrow \Herm_{4d}$  whose fiber over an $R$-point $N \in \Herm_{4d}(R)$ is the set of Lagrangian subsheaves $L \subset N$. Pulling back this construction via the map
\[ 
 \cB_d^{\Herm} =  \Herm_{2d} \times_{X_d} \Herm_{2d} \map{ (Q,P) \mapsto Q\oplus P} \Herm_{4d},
\]
 we obtain a Cartesian diagram (this is the definition of the upper left corner)
\[
\xymatrix{
{  \Lagr^{\Db}_{4d}  } \ar[d] \ar[r]  & { \Lagr_{4d}   }  \ar[d]   \\
{ \cB_d^{\Herm} } \ar[r]  & {   \Herm_{4d} . } 
}
\]
By construction, an $R$-point of  $\Lagr^{\Db}_{4d}$ is a  point $(Q,P)\in\cB_d^{\Herm}(R) $, together with a  Lagrangian subsheaf $L \subset Q \oplus P$.

\begin{defn}[Balanced Lagrangian subsheaves]\label{def: balanced}
Suppose $R$ is a $k$-algebra, and $N \in \Herm(R)$ is a Hermitian torsion sheaf on $X'_R$. 
If  $R$ is a field, we say that a Lagrangian subsheaf $L \subset N$ is  \emph{balanced} if for every closed point  $x' \in X'_R$ we have the equality 
\[
\dim_R  ( L_{x'} ) = \dim_R  ( L_{\sigma (x')} )
\]
where $L_{x'}$ and $L_{\sigma (x')}$ are the stalks of $L$ at $x'$ and $\sigma(x')$, respectively. For general $R$ we say that a Lagrangian $L\subset N$ is \emph{balanced} if it is so fiberwise on $\Spec(R)$.
\end{defn}

\begin{remark}\label{rem: balanced support}
If $N \in \Herm_{2d}(R)$, a Lagrangian subsheaf $L \subset N$ is balanced if and only if the support divisor of $L$ is $\sigma$-fixed. When this holds, it is easy to see that the support divisor of $N$ is twice the support divisor of $L$.  
In particular, the support divisor of $L$ is a $\sigma$-fixed divisor of degree $d$ on $X_R'$, which implies that $d$ must be even.
\end{remark}

\begin{defn}
The \emph{local Hitchin space} $\cM_d^{\Herm}$ is the closed substack
 \begin{equation}\label{eq: BandD}
 \cM_d^{\Herm} \subset \Lagr^{\Db}_{4d}
 \end{equation}
parametrizing triples $(Q,P,L)$  consisting of a point $(Q,P)$ of $\cB_d^{\Herm}$  and a balanced Lagrangian $ L \subset Q\oplus P$.  
The \emph{local enhanced Hitchin fibration} 
\[
f_d^{\Herm}  : \cM_d^{\Herm}  \to \cB_d^{\Herm}
\]
is the morphism that forgets the balanced Lagrangian.
\end{defn}

We now come to the main result of this subsection, which tells us that the local enhanced Hitchin fibration is a good approximation to the enhanced Hitchin fibration $f_d$ from \eqref{eq: Hitchin map d}.

\begin{prop}\label{prop: hitchin global to local}
There is a Cartesian diagram
\begin{equation}\label{diag: hitchin global to local}
\begin{tikzcd}
\cM_d \ar[r] \ar[d, "f_d"'] & { \cM_d^{\Herm}}  \ar[d, "f_d^{\Herm}"]  \\
\cA_d  \times \BunU \ar[r, "\Db"] & \cB_d^{\Herm} .
\end{tikzcd}
\end{equation}
\end{prop}

\begin{proof}
We first need to define the top horizontal arrow in the diagram \eqref{diag: hitchin global to local}.
Given a $k$-algebra $R$ and a point  $(t \co \cE \rightarrow \cF)\in \cM_d(R)$,
we use the complementary line trick \eqref{eq: augment cline} to form the morphisms 
\[
\cE \oplus \cD \to \cF \iso \sigma^* \cF^\vee \to \sigma^* \cE^\vee \oplus \sigma^* \cD^\vee
\]
of rank $n$ vector bundles on $X'_R$.
The self-duality of $\cF$ implies that the subsheaf 
\begin{equation}\label{eq: lagrangian subsheaf}
\frac{\cF}{\cE \oplus \cD}  \subset 
\frac{\sigma^* \cE^\vee }{\cE } 
\oplus
\frac{ \sigma^* \cD^\vee}{ \cD} 
\end{equation}
is Lagrangian (see \cite[\S 8.4]{FYZ}), and so defines a point of  $\Lagr^{\Db}_{4d}(R)$.
This defines the top horizontal arrow in a diagram
\begin{equation}\label{Db-sl diagram}
\begin{tikzcd}
\cM_d \ar[d, "f_d"'] \ar[r] & \Lagr^{\Db}_{4d} \ar[d, "f_d^{\Herm}"] \\
\cA_d \times \BunU \ar[r, "\Db"] & {  \cB_d^{\Herm} },
\end{tikzcd}
\end{equation}
and it is clear from the definitions that this diagram commutes.

To complete the construction of the diagram  \eqref{diag: hitchin global to local}, we need the following Lemma.

\begin{lemma}
Lagrangians arising from the construction \eqref{eq: lagrangian subsheaf} are balanced.
Equivalently,  the top horizontal arrow in \eqref{Db-sl diagram}  factors through the closed substack  \eqref{eq: BandD}.  
\end{lemma}

\begin{proof}
Suppose we have $(\cE, \cF,t) \in \cM_d(R)$, and let 
\begin{equation}\label{balanced augment cline}
\cE \oplus \cD \to \cF
\end{equation}
be the isometric embedding from \eqref{eq: augment cline}.  As $\cD$ is a line bundle, taking determinants yields a morphism 
\[
\det(\cE) \otimes \cD \to \det(\cF). 
\]
Recalling that $\cD = \det(   \sigma^*  \cE)   \otimes   \det (\cF) \otimes \omega_{X'}^{1-n} $, we can rewrite this as a  Hermitian morphism
\begin{equation}\label{Edet Herm}
\det(\cE) \otimes \sigma^* \det(\cE) \to \omega_{X'}^{n-1}.
\end{equation}
Tracing through the constructions, one can see that this is none other the Hermitian structure on $\det(\cE)$ induced by the Hermitian structure  $ \cE \otimes \sigma^*\cE \to \omega_{X'}$ from \eqref{eq: a(t)}.

Now fix a point $x' \in X'(R)$.
Because \eqref{Edet Herm} is Hermitian, its cokernel has the same length along $x'$ as it does along $\sigma(x')$.  
Taking determinants and twisting by line bundles does not change the lengths of cokernels (as in the proof of Lemma \ref{lem: supp coker det}), so we deduce that \eqref{balanced augment cline} also has this property.
This is precisely what it means for the Lagrangian \eqref{eq: lagrangian subsheaf} to be balanced.
\end{proof}

It remains to show that the diagram  \eqref{diag: hitchin global to local} is Cartesian.
 Let $\cM_d'$ be the fiber product 
\[
\begin{tikzcd}
\cM'_d \ar[r] \ar[d, "f_d"] & { \cM_d^{\Herm}} \ar[d, "f_d^{\Herm}"]  \\
\cA_d  \times \BunU \ar[r, "\Db"] &  \Herm_{2d} \times_{X_d} \Herm_{2d} .
\end{tikzcd}
\]
Thus we have an induced map $\cM_d \rightarrow \cM_d'$, which we will show to be an isomorphism by constructing an inverse map. 

An $R$-point of $\cM_d'$ consists of an $R$-point $(\cE, a, \cL)$ of  $\cA_d  \times \BunU$ together with a balanced Lagrangian 
\[
L \inj    \frac{\sigma^* \cE^\vee}{ \cE}  \oplus \frac{ \sigma^* \cD^\vee}{ \cD}.
\]
Given such data,  we define $\cF$ to be the pullback 
\begin{equation}\label{eq: reconstruct F}
\begin{tikzcd} 
\cF \ar[r, hook] \ar[d, twoheadrightarrow] & \sigma^*\cE^\vee \oplus \sigma^* \cD^\vee \ar[d, twoheadrightarrow]   \\
L \ar[r, hook] & \frac{\sigma^* \cE^\vee}{ \cE}  \oplus \frac{ \sigma^* \cD^\vee}{ \cD}.
\end{tikzcd}
\end{equation}
By construction, there is an injection $\cE \oplus \cD \inj \cF$, and in particular a map $t \co \cE \inj \cF$.  
This defines a point of $\cM_d(R)$, and so defines a morphism $\cM'_d \to \cM_d$.

Checking that this is  inverse to $\cM_d \to \cM_d'$  amounts to verifying that if we start with a balanced Lagrangian $L$ as above, then the  $\cF$ defined by \eqref{eq: reconstruct F} satisfies ${\det}^\dagger (\cF) \cong \cL$.

To this end, first note that the line bundle 
\[
\cC := \det(\cE) \otimes \cD \iso \det(\cE) \otimes \det(\sigma^* \cE) \otimes \cL
\]
on  $X'_R$ carries a natural Hermitian pairing
\begin{equation}\label{eq: C Herm}
\cC \otimes \sigma^*\cC \to \omega_{X'}^n 
\end{equation}
that can be characterized in two equivalent ways.  One can either take the determinant of the Hermitian map 
$
\cE \oplus \cD \to \sigma^*\cE^\vee \oplus \sigma^* \cD^\vee
$
to obtain 
\[
\cC \iso \det(\cE) \otimes \cD   \to \sigma^* \det(\cE)^{-1} \otimes \sigma^*\cD^{-1}    \otimes \omega_{X'}^n  \iso 
\sigma^*\cC^{-1} \otimes \omega_{X'}^n,
\]
or combine the map  
\begin{equation}\label{eq: CtoL}
\cC \iso  \det(\cE) \otimes \det(\sigma^* \cE) \otimes \cL \map{ \eqref{Edet Herm} } \omega_{X'}^{n-1} \otimes  \cL 
\end{equation}
with the Hermitian  isomorphism   $\cL \otimes \sigma^* \cL  \iso  \omega_{X'}^{2-n}$ to obtain 
\[
\cC \otimes \sigma^* \cC \map{ \eqref{eq: CtoL}}   \omega_{X'}^{2n-2}  \otimes \cL \otimes    \sigma^*\cL \iso \omega_{X'}^n .
\]

Taking  determinants in the composition 
\[
\cE \oplus \cD \to  \cF \to \sigma^*\cE^\vee \oplus  \sigma^*\cD^\vee 
\]
provides us with maps
$
\cC \to \det(\cF) \to \sigma^* \cC^{-1} \otimes \omega_{X'}^n,
$
and the assumption that the $L$ in \eqref{eq: reconstruct F} is a balanced Lagrangian implies that 
\[
\frac{ \det(\cF) }{ \cC  } \subset \frac{\sigma^* \cC^{-1} \otimes \omega_{X'}^n  }{ \cC} 
\]
is a balanced Lagrangian with respect to the Hermitian pairing \eqref{eq: C Herm}.
On the other hand, the map \eqref{eq: CtoL} also realizes 
\[
\frac{  \cL \otimes \omega_{X'}^{n-1} }{ \cC  } \subset \frac{\sigma^* \cC^{-1} \otimes \omega_{X'}^n  }{ \cC} 
\]
as a balanced Lagrangian. 
The cokernel of a Hermitian structure on a line bundle  admits at most one balanced Lagrangian
(any subsheaf of such a cokernel is determined by its support divisor, and all balanced Lagrangians have the same support divisor by Remark \ref{rem: balanced support}), and it follows that  $\det (\cF)  \cong \cL \otimes \omega_{X'}^{n-1}$, as desired. \end{proof}

\begin{cor}\label{cor: f_d proper}
The enhanced Hitchin fibration $f_d \co \cM_d \rightarrow \cA_d \times \BunU$ is a proper morphism between smooth $k$-stacks.
\end{cor}

\begin{proof}
The smoothness of $\cM_d$ is  a special case of \cite[Proposition 3.11(2)]{FYZ2}.
The stack $\cA_d$ is smooth because $\Herm_{2d}$ is smooth by  \cite[Lemma 4.3]{FYZ}, and the morphism $g_d : \cA_d \to \Herm_{2d}$ is smooth by \cite[Proposition 8.12]{FYZ}.  The stack $\BunU$ is smooth by \cite[Lemma 6.8]{FYZ}.

It remains to prove the properness of $f_d$.
The vertical arrow on the right in \eqref{Db-sl diagram} is evidently proper, by the properness of partial flag varieties.
As \eqref{eq: BandD} is a closed substack, the vertical arrow on the right in 
\eqref{diag: hitchin global to local} is also proper, and hence so is the vertical arrow on the left. 
\end{proof}


\subsection{The regular semisimple part of the enhanced Hitchin fibration}


In this subsection we study the regular semisimple part (\S \ref{sss:rss})
\[
f_d : \cM_d^\circ \to  \cA_d^\circ  \times \BunU 
\]
of the enhanced Hitchin fibration, and  give a concrete description of $\rR f_{d *} (\Ql)$.

From Definition \ref{def: double} we have a natural map $\cB_d^{\Herm} \to X_d$.
By pulling back the diagram \eqref{diag: hitchin global to local} along the inclusion $X_d^\circ  \hookrightarrow X_d$ of the open subscheme of multiplicity-free divisors, we obtain a Cartesian diagram
\begin{equation}\label{diag: rss hitchin global to local}
\begin{tikzcd}
{ \cM_d^\circ}  \ar[r] \ar[d, "f_d"'] & { \cM_d^{\Herm, \circ}}  \ar[d, "f_d^{\Herm}"]  \\
{  \cA_d^\circ }   \times \BunU \ar[r, "\Db"] & {  \cB_d^{\Herm,\circ}}  .
\end{tikzcd}
\end{equation}
This is the definition of the two stacks on the right; the two on the left were  defined in \S \ref{ss: hermitian torsion}.

Define an open subscheme $(X')^{d,\circ} \subset (X')^d$  as the fiber product
\[
(X')^{d,\circ} = (X')^d \times_{X_d} X_d^\circ .
\]
For a $k$-algebra $R$, a point of $(X')^{d,\circ} (R)$ consists of an ordered list of points  $x_1',\ldots, x_d' \in X'(R)$ with the property that the $2d$ points 
$
x_1',\ldots, x_d', \sigma(x_1'),\ldots, \sigma(x_d') 
$
are pairwise distinct.
The natural map
\begin{equation}\label{Xtilde torsor}
 (X')^{d,\circ} \to X^\circ_d
\end{equation}
 has the structure of a torsor under the finite group
\[
W_d := (\Z/2\Z)^d \rtimes S_d,
\]
and taking the quotient by  the normal subgroup $(\Z/2\Z)^d$ yields the $S_d$-torsor $X^{d,\circ} \to X_d^\circ$ of  \S \ref{ss: hermitian torsion}.
We now form the Cartesian diagram 
\begin{equation}\label{diag: local rss torsor hitchin}
\begin{tikzcd}
 { \wt \cM_d^{\Herm, \circ}}  \ar[d]   \ar[r, " \wt f_d^{\Herm}" ]   &    {  \wt \cB_d^{\Herm,\circ} }   \ar[r]  \ar[d] &   {  (X')^{d,\circ}  }  \ar[d, "\eqref{Xtilde torsor}" ]   \\
   { \cM_d^{\Herm, \circ}} \ar[r, "f_d^{\Herm}"]    &   {  \cB_d^{\Herm,\circ} }  \ar[r]   & { X_d^\circ} 
 \end{tikzcd}
\end{equation}
(this is the definition of the left and center stacks in the top row),  in which all vertical morphisms are $W_d$-torsors.

\begin{lemma}\label{lem:projective bundle} 
The morphism  $\wt f_d^{\Herm}$  in \eqref{diag: local rss torsor hitchin}  has the structure of a Zariski $(\PP^1)^d$-bundle.
More precisely,  there are line bundles $\mathscr{W}_1,\ldots, \mathscr{W}_d$ and $\mathscr{V}_1,\ldots, \mathscr{V}_d$ on $\wt  \cB_d^{\Herm,\circ}$ such that
\[
\wt \cM_d^{\Herm,\circ} \iso   \PP( \mathscr{W}_1 \oplus \mathscr{V}_1) \times \cdots \times \PP( \mathscr{W}_d \oplus \mathscr{V}_d) ,
\]
where all fiber products are taken over $\wt \cB_d^{\Herm,\circ}$. 
\end{lemma} 

\begin{proof}
Let $R$ be a $k$-algebra.  A point
\[
(Q,P,x_1',\ldots, x_d')  \in  \wt \cB_d^{\Herm,\circ}(R)
\]
consists of a pair $(Q,P) \in \cB_d^{\Herm}(R)$ of Hermitian torsion sheaves on $X'_R$ and points  $x'_1,\ldots, x'_d \in  X'(R)$
such that the common (and multiplicity-free) support divisor of $Q$ and $P$ is 
$
x'_1 + \cdots + x'_d + \sigma(x_1') + \cdots + \sigma(x_d') .
$

The pullbacks of $Q$ and $P$ along $x_i' : \Spec(R) \to X_R'$ are line bundles on $\Spec(R)$, denoted
\[
\mathscr{W}_{i,R} :=  Q|_{ x_i' }  \quad \mbox{and} \quad \mathscr{V}_{i,R} :=  P|_{ x_i' },
\]
and the Hermitian structures on $Q$ and $P$ identify 
\[
\mathscr{W}^{-1}_{i,R} \iso  Q|_{ \sigma(x_i') }  \quad \mbox{and} \quad \mathscr{V}^{-1}_{i,R} \iso  P|_{ \sigma(x_i') }.
\]

By definition, lifts of our point to $\wt \cM_d^{ \Herm, \circ}(R)$ are in bijection with  balanced Lagrangians $L \subset Q \oplus P$. 
To give such an $L$ is to give, for every $1\le i \le d$, local direct summands
\begin{align*}
L|_{x'_i}  & \subset Q|_{x'_i} \oplus P|_{x'_i}    =\mathscr{W}_{i,R} \oplus \mathscr{V}_{i,R}  \\
L|_{ \sigma( x'_i) }  & \subset Q|_{ \sigma( x'_i) } \oplus P|_{ \sigma(x'_i)}=  \mathscr{W}^{-1}_{i,R} \oplus \mathscr{V}^{-1}_{i,R} 
\end{align*}
that are exact annihilators of one another under the Hermitian form (the Lagrangian condition)  and have the same rank (the balanced condition).
Of course this is the same as giving only the rank one local direct summands $L|_{x'_i}$ for $1\le i \le d$.
With these observations, the Lemma follows immediately. 
\end{proof}

Because $X_d^\circ$ carries the $W_d$-torsor  \eqref{Xtilde torsor}, any representation of $W_d$ on a finite dimensional $\Q_\ell$-vector space determines an $\ell$-adic local system on $X_d^\circ$, constant after pullback via \eqref{Xtilde torsor}.
Of special interest to us are the local systems $K_d^{ i , \circ} $ on $X_d^\circ$ defined (for $0 \le i \le d$)  by  the representations $\Ind_{W_i \times W_{d-i}}^{W_d} (\bbm{1})$, and their shears 
(in the sense of \eqref{shear def})
\begin{equation}\label{eq:Kshear}
K_d^{ i , \circ} \tw{-i} \in D_c^b(  X_d^\circ ) .
\end{equation}

\begin{prop}\label{prop: local Hitchin local system}
The  morphism
$
f_d^{\Herm}  :  \cM_d^{\Herm, \circ} \to \cB_d^{\Herm,\circ} 
$
 is smooth and proper of relative dimension $d$, and satisfies 
 \begin{equation}\label{eq: local Hitchin local system}
\rR  f^{\Herm}_{d*} (\Ql) \cong   \bigoplus_{i=0}^d K_d^{i , \circ}   \tw{-i}    |_{\cB_d^{\Herm,\circ}   }   .
\end{equation}
Here  the restriction on the right is pullback via the natural map 
$
\cB_d^{\Herm,\circ}  \to X_d^\circ.
$
\end{prop}

\begin{proof} The morphism $\widetilde{f}_d^{\Herm}$ in the Cartesian diagram  \eqref{diag: local rss torsor hitchin} is smooth and proper by Lemma \ref{lem:projective bundle}. Since the vertical arrows are $W_d$-torsors, $f_d^{\Herm}$ is also smooth and proper.

For the claim about the direct image, if we denote by 
\begin{equation}\label{p_i}
p_i : \PP( \mathscr{W}_i \oplus \mathscr{V}_i )  \to \widetilde{\cB}_d^{\Herm,\circ}
\end{equation}
the $\PP^1$-bundles (for $1\le i \le d$) of  Lemma \ref{lem:projective bundle}, 
then proper base change and the K\"unneth formula give isomorphisms
\begin{equation}\label{direct image bundle factor}
\rR f^{\Herm}_{d*} (\Ql)|_{ \widetilde{\cB}_d^{\Herm,\circ}   } 
 \iso \rR \widetilde{f}^{\Herm}_{d*} (\Ql) 
\iso  \bigotimes_{i=1}^d \rR  p _{i*} ( \Ql ).
\end{equation}

The canonical morphism $p_i^*\Ql \to \Ql$ of sheaves on  $\PP( \mathscr{W}_i \oplus \mathscr{V}_i )$ induces a morphism
\begin{equation*}
\Ql \to \rR p _{i*} ( \Ql ) .
\end{equation*}
On the other hand, the subsheaf $\mathscr{W}_i \oplus 0 \subset  \mathscr{W}_i \oplus \mathscr{V}_i$ defines a section $q_i$ to the bundle map $p_i$, and applying $ \rR p _{i*}$ to the morphism $\Ql \to \rR q _{i*} ( \Ql )$  determines a map
$
 \rR p _{i*} ( \Ql ) \to   \Ql .
$
Applying relative Poincar\'e duality to this last map defines a morphism 
\begin{equation*}
\Ql\tw{-1}  \to \rR p _{i*} ( \Ql ).
\end{equation*}
Combining these,  we obtain a morphism $\Ql \oplus \Ql\tw{-1}  \to \rR p _{i*} ( \Ql )$.  One can verify that this is an isomorphism by checking on the level of stalks of cohomology sheaves, where it becomes the usual description of the cohomology of $\PP^1$.  In summary, each  $\PP^1$-bundle \eqref{p_i}  satisfies
\begin{equation}\label{P1 cohomology}
\rR p _{i*} ( \Ql )  \iso
\rR^0 p _{i*} ( \Ql ) \oplus \rR^2 p _{i*} ( \Ql )[-2]
\iso
\Ql \oplus \Ql\tw{-1} .
\end{equation}

The left hand side of \eqref{direct image bundle factor} is equipped with descent data relative to the $W_d$-torsor 
\begin{equation}\label{B torsor}
\widetilde{\cB}_d^{\Herm,\circ}  \to   \cB_d^{\Herm,\circ} , 
\end{equation}
 encoded by isomorphisms
\[
s^* (  \rR f^{\Herm}_{d*} (\Ql)|_{ \widetilde{\cB}_d^{\Herm,\circ}   } ) \iso \rR f^{\Herm}_{d*} (\Ql)|_{ \widetilde{\cB}_d^{\Herm,\circ}   }
\]
for every $s\in W_d$.  By examining the proof of Lemma \ref{lem:projective bundle}, one can work out what the corresponding isomorphisms are on the right hand side of \eqref{direct image bundle factor}.

If $s \in S_d \subset W_d$ is a permutation, there  are distinguished isomorphisms
\[
s^*\mathscr{W}_{s(i) }  \iso \mathscr{W}_{ i } \quad \mbox{and}\quad s^*\mathscr{V}_{ s(i) }  \iso \mathscr{V}_{ i } 
\]
of line bundles, which induce Cartesian diagrams
\begin{equation*}
\begin{tikzcd}
 {    \PP( \mathscr{W}_i \oplus \mathscr{V}_i )    }  \ar[r]   \ar[d, "p_i"']   &  {  \PP( \mathscr{W}_{s(i)} \oplus \mathscr{V}_{s(i)}  )   } \ar[d, "p_{s(i)}"]    \\
 {  \widetilde{\cB}_d^{\Herm,\circ} }   \ar[r, "s"]   &   {  \cB_d^{\Herm,\circ} }  
\end{tikzcd}
\end{equation*}
for all $1\le i \le d$.  By proper base change we obtain  isomorphisms
\[
s^* \rR  p _{s(i)*} ( \Ql ) \iso \rR  p _{i*} ( \Ql ),
\]
and hence isomorphisms (the first permutes the tensor factors)
\[
s^* \left(  \bigotimes_{i=1}^d \rR  p _{i*} ( \Ql )  \right)  
\iso s^* \left(  \bigotimes_{i=1}^d \rR  p _{s(i)*} ( \Ql )  \right)
\iso  \bigotimes_{i=1}^d \rR  p _{i*} ( \Ql ) .
\]

If $\epsilon_j \in (\Z/2\Z)^d \subset W_d$ is the tuple whose  $j^\mathrm{th}$ coordinate is $1$ and has $0$ in all other coordinates,  there are  distinguished isomorphisms
\[
\epsilon_j^* \mathscr{W}_i \iso 
\begin{cases}
\mathscr{W}_i^{-1} & \mbox{if } i=j \\
\mathscr{W}_i & \mbox{if }i \neq j 
\end{cases}
\quad \mbox{and} \quad 
\epsilon_j^* \mathscr{V}_i \iso 
\begin{cases}
\mathscr{V}_i^{-1} & \mbox{if } i=j \\
\mathscr{V}_i & \mbox{if }i \neq j .
\end{cases}
\]
These isomorphisms determine a Cartesian diagram
\begin{equation*}
\begin{tikzcd}
 {    \PP( \mathscr{W}_i \oplus \mathscr{V}_i )    }  \ar[r]   \ar[d, "p_i"']   &  {  \PP( \mathscr{W}_{i} \oplus \mathscr{V}_{i}  )   } \ar[d, "p_{i}"]    \\
 {  \widetilde{\cB}_d^{\Herm,\circ} }   \ar[r, "\epsilon_j"]   &   {  \cB_d^{\Herm,\circ} }  
\end{tikzcd}
\end{equation*}
for every $1\le i \le d$, where in the case $i=j$ we make use of the  isomorphism
\[
\PP( \mathscr{W}_j \oplus \mathscr{V}_j) \iso \PP( \mathscr{W}^{-1}_j \oplus \mathscr{V}^{-1}_j) 
\]
sending a line in $\mathscr{W}_j \oplus \mathscr{V}_j$ to its annihilator in $\mathscr{W}^{-1}_j \oplus \mathscr{V}^{-1}_j$. 
By proper base change we obtain an isomorphism
$
\epsilon_j^* \rR  p _{i*} ( \Ql )  \iso \rR  p _{i*} ( \Ql ) 
$
for every $i$, and hence an isomorphism
\[
\epsilon_j^*\left(  \bigotimes_{i=1}^d \rR  p _{i*} ( \Ql )  \right)  
\iso  \bigotimes_{i=1}^d \rR  p _{i*} ( \Ql ) .
\]

Combining \eqref{P1 cohomology} with \eqref{direct image bundle factor} defines the first isomorphism in 
\[
\rR f^{\Herm}_{d*} (\Ql)|_{ \widetilde{\cB}_d^{\Herm,\circ}   }  \iso  \bigotimes_{i=1}^d (\Ql \oplus \Ql\tw{-1})  \iso 
\bigoplus_{ i=0}^d \Big( \bigoplus_{  \substack{  S \subset \{ 1, \ldots, d\} \\  \#S=i }  } \Ql \Big) \tw{-i}  .
\]
By our analysis above,  the descent data on the left hand side corresponds to the natural action of $W_d$ (factoring through the quotient $W_d \to S_d$) in the middle by permuting the tensor factors,   and on the right by permuting the subsets $S$.  
Elementary considerations show that 
\[
\bigoplus_{  \substack{  S \subset \{ 1, \ldots, d\} \\  \#S=i }  } \Ql 
\iso
\Ind_{W_i \times W_{d-i}}^{W_d} (\bbm{1}).
\]   

In other words, the two sides of \eqref{eq: local Hitchin local system} become $W_d$-equivariantly isomorphic after pullback along the $W_d$-torsor \eqref{B torsor}. Since $W_d$ is finite, hence has no higher cohomology with rational coefficients, the descent of the isomorphism along this $W_d$-torsor is a condition and not a structure. It follows that the isomorphism descends to $\cB_d^{\Herm,\circ}$. 
\end{proof}

The main result of this subsection is the following consequence of Proposition \ref{prop: local Hitchin local system}.
This Corollary will be strengthened later, in Theorem \ref{thm: geometric K}.

\begin{cor}\label{cor: Hitchin local system}
The regular semisimple part
$
f_d : \cM_d^\circ \to  \cA_d^\circ  \times \BunU 
$
of the enhanced Hitchin fibration is smooth and proper of relative dimension $d$, and  satisfies
\[
\rR f_{d *} (\Ql) \cong  \left( \bigoplus_{i=0}^d  K_d^{i, \circ}  \tw{-i}|_{ \cA_d^\circ  }  \right)  \boxtimes \Qll{\BunU}.
\] 
Here $K_d^{i, \circ}\tw{-i}|_{ \cA_d^\circ  } $ is the sheared local system from  \eqref{eq:Kshear},  pulled back along the composition
\[
 \cA_d^\circ \map{ \eqref{eq: g_d} } \Herm_{2d}^\circ \map{ \eqref{eq: support map} }  X_d^\circ.
 \]
\end{cor}

\begin{proof} 
This follows immediately  from Proposition \ref{prop: local Hitchin local system}  by applying proper base change to the Cartesian diagram \eqref{diag: rss hitchin global to local}, and noting that
\[
 \Db^* \left(  K_d^{ i , \circ}   \tw{- i}    |_{\cB_d^{\Herm,\circ}   } \right)
 \iso
K_d^{ i , \circ}  \tw{-i}|_{ \cA_d^\circ  }  \boxtimes \Ql
 \]
by the  commutativity of the  diagram in Remark \ref{rem: double is g}.
\end{proof}

\section{The support theorem}
\label{sec: support}

Fix a $d \ge 0$, and recall the enhanced Hitchin fibration $f_d$ from \eqref{eq: Hitchin map d}.
In this section we calculate the perverse cohomology sheaves of the direct image  $\rR f_{d*} \Ql$.


\subsection{Statement of the support theorem} 


Recall from the discussion surrounding \eqref{eq:Kshear} that for  $0\le i \le d$ we have constructed a local system $K_d^{i,\circ}$ on $X_d^\circ$  from   the representation $\Ind_{W_i \times W_{d-i}}^{W_d} (\bbm{1})$.

We pull  back $K_d^{i,\circ}$  along the support map $ \Herm_{2d}^\circ \to X_d^\circ$ from \eqref{all the rss}, and denote by 
\begin{equation}\label{ICK}
\sK_d^i :=  \IC_{\Herm_{2d} } ( K_d^{i,\circ} |_{\Herm_{2d}^\circ }   )  \in D_c^b( \Herm_{2d} )
\end{equation}
the $\IC$-extension to $\Herm_{2d}$. 
The  map $g_d : \cA_d \to \Herm_{2d}$ from \eqref{eq: g_d} is smooth by \cite[Proposition 8.12]{FYZ}, and so the pullback 
 \begin{equation}\label{ICK2}
 \sK_d^i|_{\cA_d}   = \IC_{\cA_d} ( K_d^{i,\circ} |_{\cA^\circ_d}  )  \in D_c^b( \cA_d)
\end{equation}
of \eqref{ICK}  agrees with the $\IC$-extension  of the local system on $\cA_d^\circ$ from Corollary \ref{cor: Hitchin local system}. 
Going forward, we abbreviate \eqref{ICK2} simply to $\sK_d^i$, as the stack on which this sheaf lives will always be clear from context.

To lighten notation, we understand that  both \eqref{ICK} and \eqref{ICK2} are $0$ for $i \not\in \{ 0,\ldots, d\}$.

\begin{remark}
Our convention for $\IC$-extension is that  \eqref{ICK2}  is a \emph{shifted} perverse sheaf on $\cA_d$ whose restriction to $\cA^\circ_d$ is identified with the local system  $K_d^{i,\circ} |_{\cA_d^\circ }$.  More precisely, the sheaf \eqref{ICK2} satisfies
\[
\sK_d^i  [\dim \cA_d]  \in \mathrm{Perv}( \cA_d ) .
\]
Note that since $\dim \Herm_{2d} = 0$,  no shift is needed to make \eqref{ICK} perverse.
\end{remark}

\begin{thm}\label{thm: geometric K}
The enhanced Hitchin fibration  $f_d \co \cM_{d} \rightarrow \cA_d \times \BunU$ satisfies the following property:  after base change to $\ol{k}$,  there exists an isomorphism
 \begin{equation}\label{eq: geometric K}
\rR f_{d*} (\Ql) \cong   \left( \bigoplus_{i=0}^{d}   \sK_d^{i}\tw{-i} \right) \boxtimes \Qll{\BunU}.
\end{equation}
\end{thm}

We will not address the subtle question  of whether or not there is an isomorphism \eqref{eq: geometric K} defined over $k$.  
The  Frobenius-equivariance of the isomorphism of the following Corollary,  a key ingredient in the proof of Theorem \ref{thm: geometric side}, is all that we will need.
Recalling the Hitchin fibration 
\[
\pi_d : \cM_d \to \cA_d
\]
 from \eqref{eq: basic Hitchin}, 
 denote by $\p \rR^m \pi_{d*}(\Ql) $  the $m^\mathrm{th}$ perverse cohomology sheaf of $\rR \pi_{d*}(\Ql)$.

\begin{cor}\label{cor:Hitchin-perverse-cohomology}
For every $m\in \Z$, there is a    Frobenius-equivariant isomorphism
\begin{equation}\label{eq:Hitchin-perverse-cohomology}
\p \rR^{m}   \pi_{d!}(\Ql)  \cong \bigoplus_{\dim \cA_d + 2i+j= m} \sK_d^i(-i) \otimes \rH^j(\BunU)[\dim \cA_d]
\end{equation}
of perverse sheaves on $ \cA_{d,\ol{k}}$. 
Here  $\rH^j(\BunU)$  is  geometric cohomology (as in \S \ref{ss:notation}), regarded as a constant sheaf on $\cA_{d,\ol{k}}$,   with its Frobenius  induced by the $k$-structure on $\BunU$.
\end{cor}

\begin{proof}
Denote by $a : \BunU \to \Spec(k)$  the structure map.
By pullback, we regard  the sheaf $\rR^ja_!(\Ql)$  on $\Spec(k)$ also as a sheaf on  $\cA_d$.
Since $\BunU$ is proper\footnote{The definition of properness for stacks does not automatically imply cohomological properness, i.e., that cohomology is isomorphic to compactly supported cohomology. However, it is true in this case, since we work with rational coefficients.}, 
the base change of  $\rR^ja_!(\Ql)$ to $\ol{k}$   is just the constant sheaf on $\cA_{d,\ol{k}}$ associated to 
\[
\rH^j_c(\BunU)  \iso \rH^j (\BunU) .
\]

Theorem \ref{thm: geometric K}  implies the existence of an isomorphism 
\begin{equation}\label{eq:Hitchin-full-cohomology}
 \rR   \pi_{d!}(\Ql)  \iso  \left( \bigoplus_{i=0}^{d}   \sK_d^{i}\tw{-i} \right) \otimes 
  \left( \bigoplus_{j \ge 0}  \rR^ja_!(\Ql) [ -j]  \right)
\end{equation}
in  $D_c^b( \cA_{ d , \ol{k} } ) $.   
As the right hand side is a direct sum of shifted perverse sheaves, we deduce that for every $m\in \Z$ there is an isomorphism 
\begin{equation}\label{eq:Hitchin-perverse-2}
\p \rR^{m}   \pi_{d!}(\Ql)  \cong \bigoplus_{\dim \cA_d + 2i+j= m} \sK_d^i(-i) \otimes \rR^ja_!(\Ql)  [\dim \cA_d]
\end{equation}
of perverse sheaves on the base change $\cA_{d,\ol{k}}$.  
On the other hand,  Corollary \ref{cor: Hitchin local system} provides us with a $k$-rational isomorphism
\eqref{eq:Hitchin-full-cohomology},  and hence also \eqref{eq:Hitchin-perverse-2},  but only over the regular semisimple locus $\cA_d^\circ$.

There is no a priori relation between our two  isomorphisms \eqref{eq:Hitchin-perverse-2}, one defined only after base change, and the other $k$-rational but defined only over the regular semisimple locus.
However, we can modify the first in order to ensure that they agree over the base change to $\ol{k}$ of the regular semisimple locus.  
Indeed,  by composing one with the inverse of the other  we obtain an automorphism of the shifted perverse sheaf
\begin{equation}\label{local system sum}
 \bigoplus_{\dim \cA_d + 2i+j= m} \sK_d^i(-i) \otimes \rR^ja_!(\Ql) ,
\end{equation}
defined over  $\cA^\circ_{d , \ol{k}}$.
Because \eqref{local system sum} is the  $\IC$-extension of a local system on $\cA^\circ_{d , \ol{k}}$, this automorphism extends to an automorphism defined over all of $\cA_{d , \ol{k}}$.
 We use this extension to modify the first isomorphism \eqref{eq:Hitchin-perverse-2}, obtaining the desired compatibility with the second.

What we now have  is an  isomorphism \eqref{eq:Hitchin-perverse-2} defined over $\cA_d^\circ$, which, after base change to $\ol{k}$,  extends to all of $\cA_{d,\ol{k}}$.  We must show that this extension is  Frobenius-equivariant. 
This is a question of whether two endomorphisms of the base change of  \eqref{local system sum}, one induced by its $k$-structure and the other induced by the $k$-structure of $\p \rR^m \pi_{d !}(\Ql)$, coincide. 
As the $\IC$-extension functor  is fully faithful,  this can be checked over the regular semisimple locus, where it is already known.
\end{proof}

The proof of Theorem \ref{thm: geometric K} will occupy the remainder of \S \ref{sec: support}. For the remainder of \S \ref{sec: support} we base change all $k$-stacks to $\ol{k}$, but typically omit the $\ol{k}$ from the notation.

\begin{lemma}\label{lem: springer containment}
The LHS of \eqref{eq: geometric K} non-canonically contains the RHS as a summand. 
In particular, to prove \eqref{eq: geometric K}   it suffices to show that the stalks of both sides have the same Poincar\'e polynomial (whose coefficients are the dimensions of cohomology) at each $\ol{k}$-point of $\cA_d \times \BunU$.
\end{lemma}

\begin{proof}
As $f_d$ is a  proper morphism with smooth domain (Corollary \ref{cor: f_d proper}),  the Decomposition Theorem shows that $\rR f_{d*} (\Ql)$ decomposes as a direct sum of shifts of intersection complexes supported on closed substacks of  $\cA_d \times \BunU$.   Given this, the claim follows from Corollary \ref{cor: Hitchin local system}, which shows that the two sides of 
\eqref{eq: geometric K} become isomorphic after restriction to the dense open substack $\cA^\circ_d \times \BunU$.
\end{proof}

\begin{remark}
As the proof of Lemma \ref{lem: springer containment} demonstrates, the essential new content of Theorem \ref{thm: geometric K} is that the intersection complexes appearing in the decomposition of  $\rR f_{d*} (\Ql)$ have full support, hence  we refer to it as a ``Support Theorem''.
\end{remark}

A key ingredient in the proof of Theorem \ref{thm: geometric K} is the following general observation.

\begin{lemma}\label{lem: stratification support theorem}
Suppose  we are given a proper morphism $f : Y \rightarrow Z$ of smooth stacks over an algebraically closed field, and a decomposition $Y = \bigsqcup_{i=0}^d Y_i$ into locally closed substacks.   Assume that these satisfy:
    \begin{enumerate}
        \item There is an open dense substack $Z_0 \subset Z$ over which the morphism $f$ is smooth, and local systems $L_1 , \ldots, L_d$  on $Z_0$ satisfying
 \[
 \rR f_{*} (\Ql) |_{Z_0} \cong \bigoplus_{i=0}^d L_i \langle -i \rangle.
\]
        \item Each restriction  $f_{i} : Y_i \to Z$ satisfies $\rR f_{i !} (\Ql) \cong \IC_Z( L_i) \langle -i \rangle$.
    \end{enumerate}
Then there is an isomorphism
\begin{equation}\label{eq:full sppt lemma}
\rR f_* (\Ql) \cong \bigoplus_{i=0}^d \IC_Z( L_i) \langle -i \rangle.
\end{equation}
\end{lemma}

\begin{proof}
Using the first hypothesis,  the Decomposition Theorem implies that the right hand side of \eqref{eq:full sppt lemma} is non-canonically a direct summand of the left hand side.  Therefore, to check that equality holds it suffices to prove that the Poincar\'e polynomials of the stalks of both sides are equal at all $z \in Z(\ol{k})$.

    Let $Y_z := f^{-1}(z)$ be the fiber over $z$, with its induced stratification by $Y_{i,z} = f_i^{-1}(z)$. By proper base change, $(\rR f_* \Ql)_z \cong \rH^*(Y_z)$. On the other hand, $(\rR f_* \Ql)_z$ contains as a direct summand 
    \[ \bigoplus_{i=0}^d \IC_Z(L_i)_z \langle -i \rangle,\]
    which by  proper base change and the second hypothesis is isomorphic to
    \[ \bigoplus_{i=0}^d \rH^*_c(Y_{i,z}). \]
    Thus $Y_z$ is a proper scheme with  a locally closed stratification $ Y_z = \bigsqcup_{i=0}^d Y_{i,z} $ such that $\rH^*_c(Y_z)$ contains as a direct summand $\bigoplus_{i=0}^d \rH^*_c(Y_{i,z})$. In particular, for each $n$ we have an inequality
    \begin{equation}\label{eq:dim-inequality}
\dim  \rH^n_c(Y_z)   \geq  \sum_{i=0}^d \dim \rH^n_c(Y_{i,z})
    \end{equation}
 and we want to show that it is an equality. 
 
    By induction, we may reduce to the case where there are only two strata, say with $Y_{0,z}$ closed and $Y_{1,z}$ open. Then we have a long exact sequence
    \[
\ldots \rightarrow    \rH^n_c(Y_{1,z}) \rightarrow \rH^n_c(Y_z) \rightarrow \rH^n_c(Y_{0,z}) \rightarrow \ldots 
    \] 
which supplies the reverse inequality to  \eqref{eq:dim-inequality}, forcing it to be an equality. 
\end{proof}


\subsection{Springer sheaves}\label{ssec: springer sheaves}


We recall some notation from \cite[\S 4]{FYZ}. Let 
\begin{equation}\label{eq: hermpi}
\pi_{2d}^{\Herm} \co \wt{\Herm}_{2d} \rightarrow \Herm_{2d}
\end{equation}
 be the stack whose fiber over $Q \in \Herm_{2d}(R)$ is the set of full flags of the form 
\[
0 \subset Q_1 \subset \cdots \subset Q_i \subset \cdots \subset Q_{2d-1}  \subset Q_{2d}= Q, 
\]
with  $Q_{2d-i} = Q_{i}^\perp$ under the Hermitian form on $Q$.
Thus $\pi_{2d}^{\Herm}$ is the analogue of the Grothendieck--Springer resolution for $\Herm_{2d}$. 
Let 
\begin{equation}\label{SprHerm}
\Spr_{2d}^{\Herm} := \rR  \pi_{2d*}^{\Herm} (\Ql) \in D_c^b( \Herm_{2d}) .
\end{equation}
The map $\pi_{2d}^{\Herm}$ is small, and its restriction to the regular semisimple locus $\Herm_{2d}^\circ$ is a torsor for the finite group $W_d = (\Z/2\Z)^d \rtimes S_d$.
It follows that   $\Spr_d^{\Herm}$ is a perverse sheaf\footnote{No degree shift is needed for this statement since $\Herm_{2d}$ is 0-dimensional.} carrying a Springer action of $W_d$.

 For any $0 \le i \le d$,  the sheaf \eqref{ICK} satisfies
 \begin{equation}\label{ICK3}
\sK_d^i \iso (  \Spr_{2d}^{\Herm })^{W_i \times W_{d-i} } \in D_c^b(  \Herm_{2d} ),
 \end{equation}
 where the right hand side is the largest subsheaf of $\Spr_{2d}^{\Herm}$ on which the subgroup $W_i \times W_{d-i} \subset W_d$ acts trivially.    Indeed, both of these sheaves are $\IC$-extensions from the regular semisimple locus, over which both restrict to the local system determined by the representation $\Ind^{W_d}_{W_i \times W_{d-i}} ( \bbm{1} )$.


\subsection{The split case}
\label{ssec: split support theorem}


 We will first prove Theorem \ref{thm: geometric K} in the case where the double cover 
 \[
 X' = X \sqcup X
 \]
 is split.  Recall from  \S \ref{ss:notation} that this is allowed only in \S \ref{sec: Hitchin space} and \S  \ref{sec: support}.

Denote by $\Coh_d$ the moduli stack of torsion coherent sheaves on $X$ of length $d$, as in  \cite[\S 3.1]{FYZ}.
Restricting a Hermitian torsion sheaf on $X'$ to the first copy of $X$ defines an isomorphism 
$\Herm_{2d} \iso \Coh_d$.  
There is a natural map $\Coh_d \to X_d$ sending a torsion sheaf to its support divisor, and the enhanced local Hitchin base of Definition \ref{def: double} can be identified with 
\begin{equation}\label{coherent B}
\cB_d^{\Coh} \iso \Coh_d \times_{X_d}\Coh_d.
\end{equation}

Similarly,  the local Hitchin space $\cM_d^{\Herm}$ from \eqref{eq: BandD}  can be identified with the moduli stack $\cM_d^{\Coh}$ parametrizing triples  $(Q,P,L)$,  in which $Q$ and $P$ are length $d$ torsion coherent sheaves on $X$ with the same support divisor, and $L \subset Q \oplus P$ is a subsheaf that is \emph{balanced}, in the sense that its support divisor is equal to the common degree $d$ support divisor of $Q$ and $P$.   The map that forgets $L$ is denoted
\begin{equation}\label{fCoh}
f_d^{\Coh} : \cM_d^{\Coh} \to \cB_d^{\Coh}.
\end{equation}

Denote by  $\Coh_d^\mathrm{cyc}$   the $k$-stack of cyclic torsion coherent  sheaves of length $d$ on $X$.
More precisely, for a $k$-algebra $R$, an $R$-valued point of  $\Coh_d^\mathrm{cyc}$  consists of a torsion coherent sheaf $P$ on $X_R$ \emph{together with} a presentation
$
\mathcal{L}_P \twoheadrightarrow P
$
of $P$ as the quotient of a line bundle on $X_R$.
Although the line bundle and the presentation are part of the data of a cyclic torsion coherent  sheaf,  we usually omit explicit mention of them;   the existence of such a presentation of $P$ is more important than the particular choice.  
There is a canonical morphism
$
\Coh_d^\mathrm{cyc}   \to \Coh_d
$
that forgets the presentation, and hence a canonical morphism from
\begin{equation}\label{cyclic B}
\cB_d^\mathrm{cyc}  : = \Coh_d \times_{X_d} \Coh_d^\mathrm{cyc}  
\end{equation}
to \eqref{coherent B}.

The doubling map
 of Definition \ref{def: double}, now viewed as a map taking values in $\cB_d^{\Coh}$, factors naturally  as 
\begin{equation}\label{cyclic factorization}
\cA_d \times \BunU \to \cB_d^\mathrm{cyc} \to \cB_d^{\Coh},
\end{equation}
because the Hermitian torsion sheaf $P$  on $X'$ appearing in Definition \ref{def: double} is cyclic by construction: it is defined in \eqref{eq: PQpair} as the cokernel of a morphism of line bundles.  This allows us to augment the   diagram of Proposition \ref{prop: hitchin global to local}  to a Cartesian diagram
\begin{equation}\label{diag: cyclic hitchin global to local}
\begin{tikzcd}
\cM_d \ar[r] \ar[d, "f_d"'] &  { \cM_d^{\mathrm{cyc}} }    \ar[r]  \ar[d ,  " f_d^{\mathrm{cyc}}" ] &  { \cM_d^{\Coh} }   \ar[d, "f_d^{\Coh}"]  \\
\cA_d  \times \BunU \ar[r] & \cB_d^{\mathrm{cyc}}   \ar[r]  & \cB_d^{\Coh} ,
\end{tikzcd}
\end{equation}
where $\cM_d^{\mathrm{cyc}}$ is defined to make the square on the right Cartesian.  In other words, a point of  $\cM_d^{\mathrm{cyc}}$ is a point  $(Q,P)$ of  \eqref{cyclic B} together with a balanced subsheaf $L \subset Q \oplus P$.

For an $R$-valued point  $L \subset Q\oplus P$  of $\cM^\mathrm{cyc}_d$,   consider the subsheaves 
\begin{equation}\label{JandI}
J := L \cap Q  \subset Q \qquad \mbox{and} \qquad I := \mathrm{Image}( L \map{\pi_P} P)  \subset P.
\end{equation}
Here   $\pi_P: Q\oplus P \to P$ is the projection.
If $R$ is a field (so that $J$ and $I$ are $R$-flat, hence define torsion coherent sheaves on $X'_R$), the exactness of 
\[
0 \to J \to L \map{\pi_P} I \to 0
\]
implies the equality of support divisors  $\supp(J)+\supp(I)= \supp(L)$, which is equal to the common support divisor of $Q$ and $P$ by the balanced condition on $L$.
In particular,  if we denote by $j$ and $i$ the lengths of $J$ and $I$, respectively, then $i+j =d$.
This allows us to define a stratification  
\begin{equation}\label{cyclic stratification}
\cM^\mathrm{cyc}_{d} = \bigsqcup_{i=0}^d \cM^\mathrm{cyc}_{d,i}
\end{equation}
into locally closed substacks, where $\cM^\mathrm{cyc}_{d,i} \subset \cM^\mathrm{cyc}_{d}$ is the locus of points for which the subsheaves $J$ and $I$ in \eqref{JandI} are torsion coherent sheaves of lengths $d-i$ and $i$, respectively.

Denote by  $\Coh_{i \subset d}$  the stack over $\Coh_d$ whose fiber over $P \in \Coh_d(R)$ is the set of length $i$ subsheaves $I   \subset P$.

\begin{lemma}\label{lem:strata bundle}
The  morphism
\[
\pi_{d,i}  : \cM^\mathrm{cyc}_{d,i}  \to \Coh_{d-i\subset d}  \times_{X_d} \Coh_d^\mathrm{cyc}
\]
sending $L \subset Q \oplus P$ to   the pair $(L \cap Q \subset  Q , P) $   has the structure of a  rank $i$ vector bundle.
\end{lemma}

\begin{proof} 
Abbreviate $j = d-i$. Suppose we are given a $k$-algebra $R$ and a point 
\begin{equation}\label{JQP}
 (  J  \subset Q , P)   \in  \left(  \Coh_{j  \subset d}  \times_{X_d} \Coh_d^\mathrm{cyc} \right)(R).
\end{equation}
As $P$ is cyclic with the same support divisor as $Q$, it has a unique subsheaf 
 whose support divisor satisfies    $\supp (Q) = \supp (J) + \supp(I)$.  
 Namely, recalling that the cyclic torsion sheaf $P$ comes with a presentation  $\cL_P \to P$ as a quotient of a line bundle,  $I$ is the image of the composition
 \[
\cL_I :=\cL_P ( -  \supp(J)   )  \to    \cL_P \to P.
 \]

 The preimage of \eqref{JQP} under $\pi_{d,i}$  is the set of all balanced subsheaves $L \subset Q \oplus P$ such that $L \cap Q = J$.   By the discussion following \eqref{JandI}, this is the same as the set of all subsheaves $L \subset Q \oplus P$ satisfying the two conditions 
 $L \cap Q =J $ and  $\mathrm{Image}( L \to P)=I$.
  This is a standard chart in a Grassmannian, and it is well-known that such $L$ are naturally  in bijection with the set $ \Hom_{\cO_X } ( I , Q /J)$.
Explicitly,  to a morphism $f : I \to Q/J$  we  associate the  subsheaf
\[
L = \{  ( y  +   f(x)    , x ) :  y \in J,\, x \in I \}    \subset Q \oplus P . 
\]
 
  We are left to prove that $ \Hom_{\cO_X } ( I , Q /J)$ is a locally free $R$-module of rank $i$.
  By construction, $I$ comes with a presentation $\cL_I \to I$ as a quotient of a line bundle.  
It follows that the $R$-module 
  \[
  \Hom_{\cO_X } ( I , Q /J) 
  \iso  \Hom_{\cO_X } \big( \cL_I   , Q /J \big) 
  \iso  ( Q /J ) \otimes_{\cO_X} \cL_I^{-1} 
  \]
is    locally free of rank $i$,  because $Q/J$ is.
   \end{proof}

Let $\Spr_d$ be the perverse sheaf on $\Coh_d$ from \cite[Corollary 3.4]{FYZ}; this can be thought of as the Springer sheaf associated to the regular representation of $S_d$. The construction is entirely analogous to \eqref{SprHerm}, and there is a canonical isomorphism 
\begin{equation}\label{ICK4}
 (  \Spr_{2d}^{\Herm })^{ (\Z/2\Z)^d } \iso \Spr_d
\end{equation}
of perverse sheaves on $\Herm_{2d} \iso \Coh_d$, equivarlant for the  action of $W_d / (\Z/2\Z)^d \iso S_d$.

\begin{lemma}\label{lem:strata cohomology}
The restriction of  $f_d^{\mathrm{cyc}}$  to a map
$
f_{d , i} ^{\mathrm{cyc}}  : \cM^\mathrm{cyc}_{d , i }   \to \cB^\mathrm{cyc}_d 
$
 satisfies
  \[
    \rR ( f_{d,i}^{\mathrm{cyc}}  )_! \Ql \cong \Spr_d^{S_i \times S_{d-i}} \langle -i \rangle   \boxtimes \Ql .
    \]
    \end{lemma}
    
 \begin{proof}
 Abbreviating $j=d-i$,  the morphism  $f_{d , i} ^{\mathrm{cyc}}$ factors as the composition
\[
\begin{tikzcd}
	{\cM^\mathrm{cyc}_{d, i}}  \ar[  r , "{\pi_{d, i }}" ]
	 & {\Coh_{ j  \subset d}  \times_{X_d} \Coh_d^\mathrm{cyc} } \ar[rr , "{ h_{ j  \subset d }  \times \Id  }"]
	& &   {\Coh_d  \times_{X_d} \Coh_d^\mathrm{cyc}  = \cB_d^\mathrm{cyc} } ,
\end{tikzcd}
\]
where $h_{ j \subset d}$ is the morphism that forgets the subsheaf $J \subset Q$.   
As  $\pi_{d,i}$ is a vector bundle of rank $i$ by Lemma \ref{lem:strata bundle}, it satisfies  
$
\rR (\pi_{ d ,i} )_ ! (\Ql) \cong \Ql \tw{-i}.
$
On the other hand,  the morphism $h_{ j  \subset d}$ is proper and small and satisfies
\begin{equation}\label{eq: partial small}
\rR ( h _{j \subset d})_* \Ql \cong \Spr_d^{S_i  \times S_j}  .
\end{equation}
This last claim  can be checked by noting that both sides are IC-extended from the regular semisimple locus $\Coh_d^\circ$, over which the isomorphism is evident. 
The lemma follows by combining these observations.
\end{proof}

\begin{prop}\label{prop:CohHitchin}
Theorem \ref{thm: geometric K} holds under our assumption that   $X'=X\sqcup X$.
\end{prop}

\begin{proof}
Pulling back the stratification  \eqref{cyclic stratification} along the upper left horizontal map in \eqref{diag: cyclic hitchin global to local} yields a stratification
$
\cM_{d} = \bigsqcup_{i=0}^d \cM_{d,i}.
$
Applying proper base change to the diagram \eqref{diag: cyclic hitchin global to local},  Lemma \ref{lem:strata cohomology} and  the isomorphisms
\begin{equation}\label{ICK5}
 \sK_d^i  \stackrel{\eqref{ICK3}}{\iso} (  \Spr_{2d}^{\Herm })^{W_i \times W_{d-i} } \stackrel{\eqref{ICK4}}{\iso}  \Spr_d^{S_i \times S_{d-i} } 
\end{equation}
 imply that each  restriction $f_{d,i} : \cM_{d,i} \to \cA_d \times \BunU$ satisfies
 \[
    \rR ( f_{d,i}  )_! \Ql \cong \sK_d^i \langle -i \rangle   \boxtimes \Ql .
    \]
    
    The morphism $f_d$ in \eqref{diag: cyclic hitchin global to local} is  proper with smooth source and target by Corollary \ref{cor: f_d proper}, and so Lemma  \ref{lem: stratification support theorem}  shows  that the isomorphism of sheaves on $\cM_d^\circ$ from Corollary \ref{cor: Hitchin local system} extends to an isomorphism \eqref{eq: geometric K} on $\cM_d$.  
    \end{proof}

Now we deduce from the global Proposition \ref{prop:CohHitchin} a purely local statement about the map \eqref{fCoh}.
 Note that the statement only concerns torsion coherent sheaves on $X$, and makes no mention of the  double cover $X' \to X$  (although in the proof we continue to assume that $X'=X\sqcup X$).  This will be essential in the next subsection.

 \begin{cor}\label{cor:Coh support}
Fix a point   $y=(Q,P) \in  \cB_d^{\Coh}(\ol{k})$ corresponding to a pair of  torsion coherent  sheaves on $X_{\ol{k}}$ with the same (degree $d$) support divisors.
 If  $P$ can be realized as a quotient of a line bundle,  there is an isomorphism of stalks
 \[
 \rR f^{\Coh}_{d*} (\Ql)_y \cong   \bigoplus_{i=0}^{d}   \left( \Spr_d^{S_i \times S_{d-i} }  \tw{-i}  \boxtimes \Ql \right)_y.
 \]
 \end{cor}

 \begin{proof}
If  $y$ is in the image of the composition \eqref{cyclic factorization}, the claim is obtained from Proposition \ref{prop:CohHitchin} by applying proper base change to the  outermost square of the diagram \eqref{diag: cyclic hitchin global to local}, and using the isomorphisms \eqref{ICK5}.
 
For  general $y$, we exploit the fact that the statement of the Corollary is purely local, in that it only involves torsion coherent  sheaves on $X$.
In particular, it does not involve the integer $n\ge 1$ fixed in \S \ref{ssec: global Hitchin space}, and used in the definitions of $\cA_d$ and $\BunU$.
Thus we are free to assume that we have chosen $n \ge d$, and we  claim   that under this assumption $y$ lies in the image of the composition \eqref{cyclic factorization}.

To see this, first note that our  assumption on $P$ implies that it is determined uniquely  by $Q$, because a coherent torsion sheaf that can be realized as a quotient of a line bundle is uniquely determined by its support divisor.     It therefore suffices to show that  $Q$ lies in the image of the map
\[
\cA_d  \map{\eqref{eq: g_d}} \Herm_{2d} \iso \Coh_d
\]
sending $a: \cE \to \sigma^* \cE^\vee$ to the restriction of $\coker(a)$ to the first component of  $X'=X\sqcup X$.
As $Q$ has length $d \le n$ it can be expressed as the cokernel of an injective morphism  of rank $n$ vector bundles on $X_{\ol{k}}$, and from this it is easy to construct a point of  $\cA_d(\ol{k})$ above it.
 \end{proof}


\subsection{Reduction to the split case}
\label{ssec: reduction to split case}


We now prove Theorem \ref{thm: geometric K} in full generality.

 Recall from \cite[\S 4.5]{FYZ} the moduli stack $\Lg_{2d}$ of  pairs $(Q,L)$ consisting of a Hermitian torsion sheaf $Q$ on $X'$ and   a Lagrangian subsheaf $L \subset Q$.  The open substack 
 \[
 \Lg_{2d}^{\dm} \subset \Lg_{2d}
 \]
 is defined by the condition that  
\[
 \supp(L) \cap \sigma(\supp(L)) = \emptyset.
\]
 Exactly as in  \cite[\S 4.5]{FYZ}, there are surjective morphisms
 \[
\begin{tikzcd}
{  \Herm_{2d}  } & \ar[l,  "b_0"'] \Lg_{2d}^{\dm} \ar[r,  "b_1"] & \Coh_d 
\end{tikzcd}
\]
defined by $b_0(L \subset Q) = Q$ and  $b_1(L\subset Q) = \nu_* L$.
 We remind the reader that $\nu : X' \rightarrow X$ is our fixed double cover.  
 We will use this correspondence to convert Corollary \ref{cor:Coh support} into an analogous statement about the local enhanced Hitchin fibration 
 \[
 f_d^{\Herm}  : \cM_d^{\Herm} \to \cB_d^{\Herm} = \Herm_{2d} \times_{X_d} \Herm_{2d}
 \]
appearing in the Cartesian diagram of Proposition \ref{diag: hitchin global to local}.

 \begin{prop}\label{prop: Herm support}
Fix a point  $x=(Q,P) \in \cB_d^{\Herm}(\ol{k})$.   If the torsion sheaf $P$ can be expressed as the quotient of a line bundle on $X'$, there is an isomorphism of stalks
  \[
 \rR f^{\Herm}_{d*} (\Ql)_x \cong   \bigoplus_{i=0}^{d}   \left( (\Spr^{\Herm}_d)^{W_i \times W_{d-i} }  \tw{-i}  \boxtimes \Ql \right)_x.
 \]
 \end{prop}

 \begin{proof}
 Define a morphism $\Lagr^{\dm}_{2d} \to X'_d$ by sending $(Q,L)$ to the support divisor of $L$, and  
 consider the diagram 
 \[
\begin{tikzcd}
{   \cM^{\Herm}_d   }  \ar[d , " { f_d^{\Herm} }  "  ']  &  & {  \cM^{\Coh}_d  } \ar[d , " { f_d^{\Coh} }  " ]  \\
{  \cB_d^{\Herm}   } &  {  \Lagr^{\dm}_{2d} \times_{X'_d} \Lagr^{\dm}_{2d}   }    \ar[l,  "{b_0 \times b_0} " ']    \ar[r,  "{b_1 \times b_1}"] & {  \cB_d^{\Coh}   } .
\end{tikzcd}
\]

Fix a point 
\[
z\in   ( \Lagr^{\dm}_{2d} \times_{X'_d} \Lagr^{\dm}_{2d} )(\ol{k})
\]
 above $x\in \cB_d^{\Herm}(\ol{k})$.
 The choice of $z$ consists of choices of 
 Lagrangian subsheaves  $Q^\flat \subset Q$ and $P^\flat \subset P$,
along with a divisor $D^\flat \in X'_d(\ol{k})$ satisfying both
\[
\supp(Q^\flat) = D^\flat = \supp(P^\flat)
\]
and $D^\flat \cap D^\sharp=\emptyset$, where we abbreviate $D^\sharp := \sigma(D^\flat)$.
These conditions imply the existence of unique decompositions
\[
Q = Q^\flat \oplus  Q^\sharp 
\quad \mbox{and} \quad
P = P^\flat \oplus  P^\sharp 
\]
satisfying 
\[
\supp(Q^\sharp) = D^\sharp= \supp(P^\sharp).
\]
Similarly,   any $\cO_{X'}$-submodule  $L \subset Q \oplus P$ decomposes uniquely as 
\[
L = L^\flat \oplus L^\sharp
\]
with $L^\flat \subset Q^\flat \oplus P^\flat$ and $L^\sharp \subset Q^\sharp \oplus P^\sharp$.

The pair $(\nu_*Q^\flat , \nu_* P^\flat)$ of coherent torsion sheaves on $X$ corresponds to the point 
\[
y:=(b_1\times b_1)(z)  \in \cB_d^{\Coh}(\ol{k}),
\]
and the assumption that $P$ is a quotient of a line bundle on $X'$ implies that $\nu_*P^\flat$ is a quotient of a line bundle on $X$ (being a quotient of a line bundle is equivalent to every stalk being generated by a single element as a module over the local ring, and if $P$ has this property then so do $P^\flat$ and $\nu_*P^\flat$).
We claim that there is an isomorphism  
\begin{equation}\label{fiber match}
 \cM^{\Herm}_{d,x}  \iso  \cM^{\Coh}_{d,y} 
\end{equation}
 (as $\ol{k}$-schemes) between the fibers above $x$ and $y$.

A point of the left hand side of \eqref{fiber match} is a balanced Lagrangian $L \subset Q \oplus P$.
The Lagrangian condition implies that 
\begin{equation}\label{lagrangian support}
\supp(L) + \sigma(\supp(L)) = \supp(Q\oplus P).
\end{equation}
  On the other hand,
by Remark \ref{rem: balanced support} the balanced condition is equivalent to the support divisor $\supp(L)$ being $\sigma$-fixed. It follows that the balanced condition on $L$ is equivalent to the  support divisor of $L$ agreeing with the common support divisor of $Q$ and $P$, which is in turn equivalent to the two conditions 
\begin{equation}\label{flat sharp balance}
\supp(L^\flat) = D^\flat \quad \mbox{and} \quad \supp(L^\sharp) = D^\sharp.
\end{equation}
In particular, the first condition implies that $\supp( \nu_*L^\flat) = \nu_*D^\flat$ agrees with the common support divisor of  $\nu_* Q^\flat$ and $\nu_*P^\flat$, which  is precisely the balanced condition needed for the subsheaf 
\[
\nu_* L^\flat \subset \nu_* Q^\flat \oplus \nu_* P^\flat
\]
to define a point on the  right hand side of \eqref{fiber match}.

For the inverse construction, first note that  $L^\flat \mapsto \nu_* L^\flat$ establishes a bijection between subsheaves of $Q^\flat \oplus P^\flat$ and subsheaves of $\nu_* Q^\flat \oplus \nu_* P^\flat$.
A point   on the right hand side of \eqref{fiber match} consists of a subsheaf of $\nu_* Q^\flat \oplus \nu_* P^\flat$ with support divisor $\nu_*D^\flat$, which therefore has the form $\nu_* L^\flat$ for some subsheaf $L^\flat \subset Q^\flat\oplus P^\flat$ with support divisor $D^\flat$.
There is a unique subsheaf $L^\sharp \subset Q^\sharp \oplus P^\sharp$ making 
\[
L:=L^\flat \oplus L^\sharp \subset Q\oplus P
\]
 a Lagrangian, and from the relation \eqref{lagrangian support}  one deduces that $L^\sharp$ has support divisor $D^\sharp$.  
As $L$ satisfies the two conditions \eqref{flat sharp balance}, it is a balanced Lagrangian, and so defines a point of the left hand side of \eqref{fiber match}.

By  \eqref{fiber match} and proper base change, there is  an isomorphism of stalks 
\[
\rR f^{\Herm}_{d*} (\Ql)_x \iso \rR f^{\Coh}_{d*} (\Ql)_y.
\]
On the other hand, by   \cite[Proposition 4.11]{FYZ}  there is a canonical  isomorphism
\[
b_0^* (\Spr^{\Herm}_{2d})^{W_i \times W_{d-i} }    \iso b_1^* \Spr_d^{S_i \times S_{d-i} }  
\]
of sheaves on $\Lagr_{2d}^{\dm}$, and hence a canonical isomorphism of stalks
\[
  \bigoplus_{i=0}^d  \left(  (\Spr^{\Herm}_{2d})^{W_i \times W_{d-i}}    \langle - i \rangle  \boxtimes \Ql \right)_x
 \iso
  \bigoplus_{i=0}^d  \left( \Spr_{d}^{S_i \times S_{d-i}}  \langle - i \rangle  \boxtimes \Ql \right)_y.
\]
The Proposition now follows immediately from Corollary \ref{cor:Coh support}.
 \end{proof}

\begin{proof}[Proof of Theorem \ref{thm: geometric K}]
Fix  a geometric point 
\[
y = (\cE ,a ,\cL) \in (\cA_d \times \BunU)(\ol{k}).
\]
If we denote by  $x\in \cB_d^{\Herm}( \ol{k})$  its image under the doubling map of Definition \ref{def: double},
there are isomorphisms of stalks 
\[
\rR f_{d*} (\Ql)_y  \iso \rR f^{\Herm}_{d*} (\Ql)_x
 \iso 
  \bigoplus_{i=0}^d  \left(  (\Spr^{\Herm}_{2d})^{W_i \times W_{d-i}}    \langle - i \rangle  \boxtimes \Ql \right)_x
\iso 
  \bigoplus_{i=0}^d  \left(  \sK_d^i   \langle - i \rangle  \boxtimes \Ql \right)_y.
\]
The first is obtained by  applying proper base change to the diagram of Proposition \ref{prop: hitchin global to local}.
The second is by Proposition \ref{prop: Herm support}, whose hypotheses are satisfied by the construction of the doubling map.
  The third is by \eqref{ICK3}, where in the final expression $\sK_d^i$ is the sheaf \eqref{ICK2} obtained from   \eqref{ICK} by pullback. 
  
   As observed in Lemma \ref{lem: springer containment}, such an isomorphism on the level of stalks is enough to complete the proof of  Theorem \ref{thm: geometric K}.
\end{proof}

\section{Hecke correspondence for the Hitchin space}
\label{sec: Hecke} 


We continue to fix an integer $n\ge 1$. 
In this section we construct  a Hecke stack for the Hitchin space $\cM$ from \S \ref{ssec: global Hitchin space},
 define several cohomological correspondences on it, and formulate in Theorems \ref{thm: global C action}, \ref{thm: universal bundle action}, and \ref{thm: tautological action} their induced actions on cohomology.  The proofs of these three theorems will then be given in 
\S \ref{sec: comparison correspondence}, \S \ref{sec:univ-action}, and \S \ref{sec:taut-action}.

The central result of \S \ref{sec: Hecke} is  Corollary \ref{cor: tautological correspondence}.  This will be a crucial ingredient in the proof of Theorem \ref{thm: geometric side}, which computes the  geometric side of the higher Siegel--Weil formula


\subsection{The Hecke stack for $\BunU$}\label{ssec: Hecke BunU}


Recall the $k$-stack $\BunU$ from \eqref{BunU def}.
As in  \cite[\S 6.2]{FYZ}, there is a  Hecke correspondence
\begin{equation}\label{eq: HkU correspondence}
\BunU \xleftarrow{h_0^{  \U^\dagger(1) }} \HkU^1 \xrightarrow{h_1^{  \U^\dagger(1) }} \BunU .
\end{equation}
To spell this out, for a $k$-algebra $R$ we define $\HkU^1(R)$ to be the groupoid of  pairs $(x', \cL_{\bu})$ in which  $x' \in X'(R)$,  and $\cL_{\bu}$ is a diagram 
\begin{equation}\label{eq: HkU mod}
\xymatrix{
{  \cL_0 }  & &  {  \cL^\flat_{1/2} } \ar[ll]_{x'} \ar[rr]^{\sigma (x') }   & & { \cL_1 } 
}
\end{equation}
of morphisms  of line bundles on $X'_R$   with   $\cL_0, \cL_1 \in \Bun_{\U^\dagger(1)}(R)$.
In the terminology of \cite[Definition 6.5]{FYZ}, the leftward map is required to be an upper modification of length one along $x'$, the rightward map is required to be an upper modification of length one along $\sigma (x')$, and the Hermitian pairings on  $\cL_0$ and $\cL_1$ are required to  pull back to the same pairing 
\[
\cL^\flat_{1/2} \otimes   \sigma^*\cL^\flat_{1/2}    \to    \nu^* \omega_X^{1-n}  .
\]
The  arrows $h_0^{  \U^\dagger(1) }$ and $h_1^{  \U^\dagger(1) }$ in \eqref{eq: HkU correspondence} are defined by
\[
  \cL_0   \mapsfrom (x', \cL_{\bu}) \mapsto \cL_1.
\]

\begin{remark}\label{rem: HkU}
Any modification as in \eqref{eq: HkU mod} satisfies 
\[
\cL_0( -x')= \cL^\flat_{1/2} =  \cL_1( - \sigma(x') ),
\]
so is determined  by the data  $( x' , \cL_0 )$. This defines an isomorphism 
\[
 \HkU^1 \iso  X' \times  \Bun_{ \U^\dagger(1) } .
 \]
\end{remark}

\subsection{The Hecke stack for $\Bun_{\U(n)}$}

We continue to fix an integer $n\ge 1$, and   define the Hecke correspondence 
\[
\Bun_{\U(n)}\xleftarrow{h_0^{\U(n)}} \Hk_{\U(n) }^1 \xrightarrow{h_1^{ \U(n) }} \Bun_{\U(n)}
\]
for $\Bun_{\U(n)}$ in the same way as  \eqref{eq: HkU correspondence}.  
To spell this out,  for a $k$-algebra $R$, a point 
\[
(x', \cF_{\bu}) \in \Hk_{\U(n) }^1 (R)
\]
consists of an $x' \in X'(R)$,  and a diagram of rank $n$ vector bundles
\begin{equation}\label{Hk U(n) mod}
\begin{tikzcd}
\cF_0 & \cF^\flat_{1/2} \ar[l, "x' "'] \ar[r, "\sigma (x')"] & \cF_1
\end{tikzcd}
\end{equation}
on $X'_R$.    The left and right maps in \eqref{Hk U(n) mod} are modifications of length one along $x'$ and $\sigma (x')$, respectively.
The vector bundles  $\cF_0$ and $\cF_1$ are endowed with Hermitian structures making them points of  $\Bun_{\U(n)}(R)$, and the Hermitian pairings on them pull back to the same Hermitian pairing on $\cF^\flat_{1/2}$.
The maps $h_0^{\U(n)}$ and $h_1^{\U(n)}$  are
\[
  \cF_0 \mapsfrom (x', \cF_{\bu})  \mapsto \cF_1 . 
\]


\subsection{The Hecke stack for $\cM$}


Hecke stacks are defined in \cite[\S 3]{FYZ2} for Hitchin spaces in significant generality. 
Here we spell out the definition of the Hecke stack 
\[
 \cM  \xleftarrow{h_0^{\cM  }} \Hk_{\cM }^1 \xrightarrow{h_1^{ \cM  }} \cM 
\]
for the Hitchin space $\cM = \cM_{\GL(n-1)', \U(n)}^{\ns}$ from \S \ref{ssec: global Hitchin space}.

For a  $k$-algebra $R$, denote by $\Hk_{\cM}^1(R)$  the groupoid of tuples  $(x', \cE, \cF_{\bu}, t_{\bu})$ in which  $( x'  , \cF_{\bu} ) \in \Hk_{\U(n)}^1(R)$,  $\cE$ is a rank $n-1$ vector bundle on $X'_R$, and  $t_{\bu}=(t_0, t_{1/2}, t_1)$ is a tuple of morphisms  fitting into a commutative diagram 
\begin{equation}\label{global hecke point}
\xymatrix{
& &  { \cE } \ar[dll]_{t_0} \ar[drr]^{t_1}  \ar[d]_{t_{1/2}} \\
 { \cF_0 }  & & { \cF^\flat_{1/2} } \ar[ll]^{x'}  \ar[rr]_{\sigma ( x' ) }  &  &  { \cF_1, }  
}
\end{equation}
and making   $( \cE, \cF_0 , t_0)$ and $(\cE, \cF_1, t_1)$  points of $\cM(R)$.
The morphisms $h_0^{\cM}$ and $h_1^{\cM}$ are 
\[
( \cE, \cF_0 , t_0)  \mapsfrom (x', \cE, \cF_{\bu}, t_{\bu})   \mapsto ( \cE, \cF_1 , t_1) .
\]

\begin{remark}\label{rem: MtoBun hecke}
The three Hecke stacks just defined are related by  a  commutative diagram
\begin{equation}\label{MtoBun hecke}
\begin{tikzcd}
{ \cM }   \ar[d]  &  {  \Hk_{\cM}^1   } \ar[l , "h_0^{\cM}" ' ]      \ar[ r , "h_1^{\cM}"  ]   \ar[d] &  { \cM  } \ar[d]  \\
{  \Bun_{\U(n)}   }   \ar[d , " {\det}^\dagger " ']   & {  \Hk_{ \U(n) }^1   }   \ar[l , "h_0^{\U(n)} "  ']      \ar[ r , "h_1^{\U(n)}"  ]    \ar[d ] &  {  \Bun_{\U(n)}  }  \ar[d , " {\det}^\dagger " ]   \\
{  \BunU   }    & {  \Hk_{ \U^\dagger(1) }^1   }   \ar[l , "h_0^{  \U^\dagger(1) } "  ']      \ar[ r , "h_1^{ \U^\dagger(1) }"  ]   &  {  \BunU } ,
\end{tikzcd}
\end{equation} 
where ${\det}^\dagger$ is the twisted determinant \eqref{eq: twisted det}.
In the middle  column, the first arrow is the obvious  $(x', \cE, \cF_{\bu}, t_{\bu}) \mapsto (x', \cF_{\bu} )$.
The second sends an $R$-point $(x', \cF_{\bu})$ to the diagram 
\begin{equation}\label{twisted determinant modification}
\begin{tikzcd}
\cL_0 & \cL^\flat_{1/2} \ar[l, "x' "'] \ar[r, "\sigma (x')"] & \cL_1 ,
\end{tikzcd}
\end{equation}
where $\cL_i = {\det}^\dagger (\cF_i)$, and the line bundle $\cL^\flat_{1/2}$ is constructed from $\cF^\flat_{1/2}$ in the same way.
\end{remark}

\begin{remark}\label{rem: Hk cline}
By the complementary line trick described in \S \ref{ss:cline}, an $R$-point \eqref{global hecke point} of $\mathrm{Hk}_{\mathcal{M}}^1$ can be functorially upgraded to a diagram
 \begin{equation}\label{upgraded hecke}
\begin{tikzcd}
   {  \cE \oplus  \cD_0  } \ar[d , " \wt{t}_0 " ' ]    &  {  \cE \oplus   \cD_{1/2}^\flat  }   \ar[l  ,  " x' "  ' ] \ar[r , "  \sigma ( x' ) "  ]  \ar[d , "  \wt{t}_{1/2} " ]  &  {  \cE \oplus   \cD_1 }  \ar[d ,  " \wt{t}_1 "  ]  \ar[d] \\
 { \cF_0 }   & { \cF^\flat_{1/2} } \ar[l  ,  " x' " ' ]   \ar[r , " \sigma (x' ) "  ]   &    { \cF_1 }  
\end{tikzcd}
\end{equation}
of injections of rank $n$ vector bundles on $X'_R$,  where the modification of line bundles
\begin{equation}\label{Dmod}
\begin{tikzcd}
{  \cD_0  }   &  {    \cD_{1/2}^\flat  }   \ar[l , " x' " ' ]   \ar[ r , " \sigma(x') " ]    &   {  \cD_1    } 
\end{tikzcd}
\end{equation}
is obtained by twisting \eqref{twisted determinant modification} by $\det(   \sigma^*  \cE)$, exactly as in   \eqref{eq: cline def}.
\end{remark}

Taking the product of \eqref{eq: HkU correspondence} with the Hitchin base $\cA$ defines a correspondence 
\[
\begin{tikzcd}
{ \cA  \times \BunU } &  {  \cA \times \HkU^1 } \ar[l] \ar[r] &  {  \cA \times \BunU } .
\end{tikzcd}
\]
Using the morphisms of \eqref{MtoBun hecke},  the    Hitchin fibration  \eqref{eq: basic Hitchin} and its enhancement \eqref{eq: Hitchin map} now fit into a commutative diagram 
 \begin{equation} \label{global hecke}
   \adjustbox{scale = 0.75}{
\begin{tikzcd}
{  \cM   }     \ar[dr , " f  " '  ]   \ar[dd,  bend right,  "  \pi "' ] &&   {   \Hk^1_{  \cM }   }   \ar[ll ] \ar[rr ]  \ar[d]   & &  {  \cM   }    \ar[dl, " f " ]   \ar[dd, bend left,  " \pi " ] \\
& {  \cA  \times \BunU}   \ar[dl  ]   & {  \cA \times  \HkU^1 }     \ar[l] \ar[r]    \ar[d  ] &   {  \cA  \times \BunU}   \ar[ dr ] \\
  {   \cA  }  &    & {  \cA }  \ar[ll, equal ] \ar[rr , equal ] & &  {  \cA.}
\end{tikzcd}
}
\end{equation}

 \begin{remark}\label{rem: d hecke}
Recall from \S \ref{ss: hermitian torsion} the morphism  $\pi_d: \cM_d \to \cA_d$ obtained from $\pi: \cM \to \cA$ by pullback along the open and closed substack  $\Herm_{2d} \subset \Herm$.
There is an obvious version of \eqref{global hecke} for these stacks as well, obtained by the same pullback.
\end{remark}


\subsection{Tautological bundles and the comparison divisor}
\label{ssec: tautological bundles}


We next define a tautological bundle on $\Hk_{\cM}^1$.  This is an object of central importance to us, but it is difficult to analyze directly. We will use the diagram \eqref{MtoBun hecke} to break it down into several pieces, which will be easier to understand.

\begin{defn}\label{def: U(1) tautological bundle}
The  \emph{tautological bundle} $\ell_{\U^\dagger(1)}$ on  $\Hk_{\U^\dagger(1)}^1$ is the line bundle  whose fiber at an $R$-valued  point  $(x', \cL_{\bu}) \in \Hk_{\U^\dagger(1)}^1(R)$ is the pullback of the torsion coherent sheaf
\[
  \mathrm{coker}(   \cL^\flat_{1/2} \to   \cL_1 )  
\]
along  its support divisor $\sigma(x') : \Spec(R) \to X'_R$.
\end{defn}

\begin{defn}\label{def:U(n) taut}
The \emph{tautological bundle} $\ell_{\U(n)}$ on $\Hk_{\U(n)}^1$  is the line bundle  whose fiber over an $R$-point $(x', \cF_{\bu})$ is the pullback of the torsion coherent sheaf
\[
\mathrm{coker}(   \cF^\flat_{1/2} \to   \cF_1 ) 
\]
along its support $\sigma(x') : \Spec(R) \to X_R'$.  
\end{defn}

\begin{defn}\label{def:M taut}
The \emph{tautological bundle} $\ell_{\cM}$  on $\Hk_{\cM}^1$ is the pullback of $\ell_{\U(n)}$ along 
\[
\Hk_{\cM}^1 \map{ \eqref{MtoBun hecke}  }  \Hk^1_{ \U(n)  } .
\] 
 Abusing notation, we denote again by $\ell_{\Ud}$  the pullback of the tautological  bundle from Definition \ref{def: U(1) tautological bundle} along
\[
\Hk^1_{\cM}  \map{ \eqref{MtoBun hecke} } \Hk^1_{\U^\dagger(1) },
\]
 relying on context to make the base clear.   
\end{defn}

\begin{defn} \label{def:comparison divisor}
The  \emph{comparison divisor} is the closed substack
\[
C  \hookrightarrow \Hk_{\cM}^1
\]
whose $R$-points are those tuples $(x',\cE, \cF_{\bu}, t_{\bu}) \in \Hk_{\cM}^1(R)$ for which the morphism of vector bundles
$
\cD_0 \rightarrow \cF_0
$
 from Remark \ref{rem: Hk cline} vanishes along  $\sigma(x') \in  X'(R)$. 
\end{defn}
  
\begin{defn}\label{def:Euniv}
The \emph{universal vector bundle}  $\cE^{\univ}$  on $\Hk_{\cM}^1$, is the vector bundle whose fiber over an $R$-valued point  $(x',\cE, \cF_{\bu}, t_{\bu})  \in \Hk_{\cM}^1(R)$ is the restriction of $\cE$ to  $x' \in  X'(R)$.  In other words, it is the pullback of the universal rank $n-1$ vector bundle on $\cA \times X'$, along the map 
  $\Hk_{\cM}^1  \to \cA \times X'$.
\end{defn}

 The following Proposition gives us a factorization of $\ell_{\cM}$ in terms of the other, simpler,  objects just defined.

  \begin{prop}
  \label{prop:bundle decomp}
 There is a canonical morphism
\begin{equation}\label{taut compare}
\ell_{ \U^\dagger(1) } \otimes  \det( \cE^\mathrm{univ} )   \to \ell_{\cM}
\end{equation}
of line bundles on $\Hk_{\cM}^1$, whose vanishing locus is equal to the comparison divisor $C$.  
In particular, $C$ is a Cartier divisor, and there is an isomorphism 
\[
\ell_{ \U^\dagger(1) } \otimes  \det (\cE^{\univ} ) \otimes  \cO(  C )
 \iso \ell_{ \mathcal{M}}  .
\]
\end{prop}

\begin{proof}
First we define the morphism \eqref{taut compare}.
Given $(x', \cE, \cF_{\bu}, t_{\bu}) \in \Hk_{\cM}^1(R)$, we form the diagram of Remark \ref{rem: Hk cline}, and from it the morphism of   torsion coherent sheaves
\[
\coker(\cD_{1/2}^\flat \rightarrow \cD_1) \rightarrow \coker(\cF_{1/2}^{\flat} \rightarrow  \cF_1)
\]
on $X_R'$.
Unraveling the definition of the $\cD_i$, we see that 
\[
\coker(\cD_{1/2}^\flat \rightarrow \cD_1) \cong 
\mathrm{coker} \big(    \cL^\flat_{1/2}   \to  \cL_1   \big)    \otimes    \det( \sigma^* \cE ) 
\]
where $\cL_i= {\det}^\dagger( \cF_i)$  as in Remark \ref{rem: MtoBun hecke}.   Pulling back the resulting morphism 
\[
\mathrm{coker} \big(    \cL^\flat_{1/2}   \to  \cL_1   \big)    \otimes    \det( \sigma^* \cE ) \to 
 \coker(\cF_{1/2}^{\flat} \rightarrow  \cF_1)
\]
along  $\sigma(x')  : \Spec(R) \to  X'_R $ defines the morphism \eqref{taut compare}; more precisely, defines its restriction to the $R$-valued point of $\Hk_{\cM}^1$ determined by the tuple $(x', \cE, \cF_{\bu}, t_{\bu})$.

Denote by $C' \hookrightarrow  \Hk_{\cM}^1$ the vanishing divisor of the morphism \eqref{taut compare}. 
Directly from the construction above,  one sees that $(x', \cE, \cF_{\bu}, t_{\bu}) \in  C'(R)$  precisely when 
 there exists a (necessarily unique) dashed arrow 
\[
\begin{tikzcd}
   {  \cE \oplus  \cD_0  } \ar[d , " t_0"  '  ]     & {  \cE \oplus  \cD_{1/2}^\flat  }   \ar[l , " x'  "  ' ]   \ar[r , "  \sigma (x') " ]   \ar[d]&    { \cE \oplus \cD_1 }  \ar[d , " t_1"  ]   \ar[d]   \ar[dl , dashrightarrow] \\
 { \cF_0 }   & { \cF^\flat_{1/2} } \ar[l  ,  " x'  "  ' ]    \ar[r , "  \sigma ( x' )  " ]   &    { \cF_1 }  
\end{tikzcd}
\]
(see Remark \ref{rem: Hk cline}) making the diagram commute. 
In other words, $C'$ is the locus of points \eqref{upgraded hecke} for which the map 
\[
 \cD_0(  - x'  ) \iso  \cD_{1/2}^\flat \to \cF^\flat_{1/2} 
\]
 extends to 
\[
 \cD_0( \sigma (x')  - x' ) \iso  \cD_1 \to \cF^\flat_{1/2} .
 \]
 This is equivalent to the map
 \[
 \cD_0( - x' ) \to \cF^\flat_{1/2} 
 \]
 vanishing along $\sigma(x')$, which is equivalent to    $\cD_0 \to \cF_0$ vanishing along $\sigma (x')$.
 This is precisely  the condition that cuts out $C$,  which is therefore equal to $C'$.
\end{proof}


\subsection{Hecke action on cohomology}
\label{ssec:hecke-action-on-cohomology}


Fix an integer $d \ge 0$, and  return to the diagram
\begin{equation} \label{global hecke d}
   \adjustbox{scale = 0.75}{
\begin{tikzcd}
{  \cM_d   }     \ar[dr , " f_d  " '  ]   \ar[dd,  bend right,  "  \pi_d "' ] &&   {   \Hk^1_{  \cM_d }   }   \ar[ll , "h_0^{\cM}" ' ] \ar[rr, "h_1^{\cM}" ]  \ar[d]   & &  {  \cM_d   }    \ar[dl, " f_d " ]   \ar[dd, bend left,  " \pi_d " ] \\
& {  \cA_d  \times \BunU}   \ar[dl  ]   & {  \cA_d \times  \HkU^1 }     \ar[l] \ar[r]    \ar[d  ] &   {  \cA_d  \times \BunU}   \ar[ dr ] \\
  {   \cA_d  }  &    & {  \cA_d }  \ar[ll, equal ] \ar[rr , equal ] & &  {  \cA_d ,}
\end{tikzcd}
}
\end{equation}
from  \eqref{global hecke} and  Remark \ref{rem: d hecke}.

\subsubsection{Cohomological correspondences}\label{sssec:cohomological-correspondences} The morphism  labeled  $h_1^{\cM}$   is LCI  (the source and target are smooth, as established in \cite[\S 3.4, \S 3.5]{FYZ2})  and proper, of relative dimension $1$. Hence the Gysin map $\Ql\tw{1} = h_1^{\cM *} (\Ql) \tw{1} \rightarrow h_1^{\cM!} \Ql$ is an isomorphism. Hence, referring to the notions of cohomological correspondences introduced in \cite[\S 3]{FYZ}, we have canonical isomorphisms 
\begin{align*}
\Corr_{\Hk_{\cM_d}^1}(\Ql, \Ql ) & :=  \Hom_{  D_c^b(\Hk_{\cM_d}^1)  } (  h_0^{\cM *} \Q_\ell, h_1^{\cM !} \Q_\ell  ) \\ 
& \iso  \Hom_{  D_c^b(\Hk_{\cM_d}^1 )  }    (   \Q_\ell,  \Q_\ell  \tw{ 1}  )    \nonumber \\
& \iso H^2(\Hk_{\cM_d}^1) (1)  .   \nonumber
\end{align*}
Here $H^2(\Hk_{\cM_d}^1) (1) $ refers to the absolute cohomology over $k$, which admits a pullback map to geometric cohomology (recall the notational conventions from \S \ref{ss:notation}) landing in the subspace of Frobenius invariants $\rH^2(\Hk_{\cM_d}^1)(1)^{\Frob}$. Any cohomology class $\mf{c} \in H^{2}( \Hk_{\cM_d}^1 )  (1) $ therefore induces a Frobenius-equivariant cohomological correspondence in $\Corr_{\Hk_{\cM_d}^1, \ol k}(\Ql, \Ql)$.

\subsubsection{Definition of endomorphisms}

 The Hitchin map 
 \[
 \pi_d : \cM_d \to \cA_d
 \]
 is proper, by Corollary \ref{cor: f_d proper} and the properness of $\BunU$ as a $k$-stack. 
 Therefore, for any 
 \[
 \mf c \in H^2(\Hk_{\cM_d}^1) (1) \iso \Corr_{\Hk_{\cM_d}^1}(\Ql, \Ql ), 
 \]
  there is a pushforward cohomological correspondence
\begin{equation}\label{hitchin correspondence endo}
\Gamma_{\mf c} := \pi_{d!} (\mf c) \in \Corr_{\cA_d}( \rR \pi_{d!} \Ql,  \rR \pi_{d!} \Ql ) \iso \End_{D_c^b(\cA_d)}( \rR \pi_{d!}  \Ql).
\end{equation}
We  denote again  by $\Gamma_{\mf c}$ the  endomorphisms of the perverse cohomology sheaves $\p \rR^m \pi_{d!} ( \Ql ) $ induced by \eqref{hitchin correspondence endo}.

\begin{remark}\label{rem:stalk hecke}
For any $m\in \Z$, the action of the endomorphism \eqref{hitchin correspondence endo} on the stalk 
\[
( \rR^m \pi_{d!}  \Ql)_\mathfrak{a} \iso \rH_c^m( \cM_{d,\mf{a}} ) 
\]
at a point $\mathfrak{a} \in \cA_d(\ol{k})$ is given by the explicit and familiar formula
$\alpha \mapsto h^{\cM} _{1 !} ( \mf{c} \cdot h^{\cM *}_0 \alpha)$.
\end{remark}

 Applying the construction \eqref{hitchin correspondence endo} to the cycle class 
 \[
 [C] \in H^2(\Hk_{\cM_d}^1)( 1 )
 \]
 of the comparison divisor (Definition \ref{def:comparison divisor}), we obtain endomorphisms
 \begin{equation}\label{eq:Gamma_C}
 \Gamma_{[C]} \in \End_{ D_c^b(\cA_d)}(\p\rR^m \pi_{d!}(\Ql))
 \end{equation}
 for all $m$.  Similarly, the Chern classes 
 \[
 c_1(\ell_{\U^\dagger(1) }),\ c_1(    \det \cE^{\univ}   ), \ c_1(\ell_{\cM}) \in H^2(\Hk_{\cM_d}^1)( 1 )
 \]
 of the various line bundles introduced in \S \ref{ssec: tautological bundles}
 induce endomorphisms
 \begin{equation*}
 \Gamma_{c_1(\ell_{\U^\dagger(1) })}, \
 \Gamma_{c_1(    \det \cE^{\univ}   )}, \ \Gamma_{c_1(\ell_{\cM})} \in \End_{ D_c^b(\cA_d)}(\p\rR^m \pi_{d!}(\Ql)).
 \end{equation*}

Recall  from Corollary \ref{cor:Hitchin-perverse-cohomology} the  Frobenius-equivariant isomorphisms 
\begin{equation}\label{eq: hitchin direct image}
\p \rR^{m}  \pi_{d!}(\Ql)  \cong \bigoplus_{\dim \cA_d + 2i+j= m} \sK_d^i(-i) \otimes \rH^j(\BunU)[\dim \cA_d]
\end{equation}
of perverse sheaves on $\cA_{d,\ol{k}}$.
Our goal is to explain how the above endomorphisms interact with this direct sum decomposition.  
They do not quite respect the decomposition, but we  will see that they do respect the   Frobenius-stable filtration defined by 
\begin{equation}\label{eq: filtration}
F_a \left(\p \rR^{m} \pi_{d!}(\Ql)\right) := \bigoplus_{\substack{\dim \cA_d + 2i+j = m \\ i \leq a}} \sK_d^i(-i) \otimes \rH^j(\BunU) [\dim \cA_d]. 
\end{equation}
Note that the associated graded sheaf is again just \eqref{eq: hitchin direct image}, and so 
 any endomorphism of $\p \rR^{m}  \pi_{d!}(\Ql)$ that respects the filtration has a \emph{semisimplification}, which is an endomorphism of  \eqref{eq: hitchin direct image} respecting the direct sum decomposition.

 \subsubsection{The automorphism $w$} Let $\Prym$ be the Prym stack of $X'/X$, parametrizing line bundles $\cL$ on $X'$ plus an isomorphism $\mathrm{Nm}_{X'/X}( \cL )  \cong \cO_X$.  
 Giving such an isomorphism is equivalent to giving a Hermitian isomorphism $\cL \otimes \sigma^* \cL \cong \cO_{X'}$, and using this one sees that  $\Prym$ acts on $\Bun_{ \U^\dagger(1)}$ by taking the tensor product of Hermitian line bundles.

Consider the Abel-Jacobi map $\AJ \co X' \rightarrow \Prym$ sending $x' \mapsto \cO(\sigma (x') - x')$. Composing $\AJ$ with the tensoring action of $\Prym$ on $\BunU$, each $x'\in X'( \ol{k})$ determines an automorphism of $\Bun_{ \U^\dagger(1), \ol k}$ given explicitly by $
\cL \mapsto \cL( \sigma(x')-x')
$. Since $X'$ is connected, $\AJ$ lands in a connected component of $\Prym$, hence the induced action on $\rH^*(\BunU)$ is in fact independent of $x'$; we denote it 
\begin{equation}\label{eq:w-defn}
w : \rH^*( \Bun_{\U^\dagger(1)}) \to \rH^*( \Bun_{\U^\dagger(1)}).
\end{equation}
This is an order 2 automorphism, which we will describe in a different way in \S\ref{ssec:cohomology-prym}.

\subsubsection{Action of the comparison divisor} 
The following Theorem will be proved in \S \ref{sec: comparison correspondence} below.  
 
 \begin{thm}\label{thm: global C action}
The endomorphism $\Gamma_{[C]}$ from \eqref{eq:Gamma_C} preserves the direct sum decomposition 
\eqref{eq: hitchin direct image}, and acts by
\[
\text{$(2d-2i)  \otimes w$ on the summand $\sK_d^i(-i) \otimes \rH^j(\BunU)[\dim \cA_d]$.}
\]    
  \end{thm}

\subsubsection{Action of the universal bundle} 
Next we explicate the action of $\Gamma_{c_1(    \det \cE^{\univ}   )}$. 
   The following theorem will be proved in \S \ref{sec:univ-action} below.

 \begin{thm}\label{thm: universal bundle action} The endomorphism $\Gamma_{c_1(    \det \cE^{\univ}   )} \in \End_{D^b_c(\cA_d)}(\p\rR^m \pi_{d!} (\Ql))$ preserves the filtration \eqref{eq: filtration} on $\p \rR^m \pi_{d!} (\Ql)$, and its semisimplification acts by 
\[
\text{$(-d + (n-1)\deg_X(\omega_X)) \otimes w $ on the summand $\sK_d^i(-i) \otimes \rH^j(\BunU)[\dim \cA_d]$.}
\]
 \end{thm}

 \subsubsection{Action of the $\U^\dagger(1)$ tautological bundle} 
	Next we explicate the action of $\Gamma_{c_1(\ell_{\U^{\dagger}(1)})}$. 
    The following Theorem will be proved in \S \ref{sec:taut-action} below. 
	
 \begin{thm}\label{thm: tautological action}
The endomorphism $\Gamma_{c_1(   \ell_{\U^\dagger(1) } ) } \in \End_{D^b_c(\cA_d)}(\p\rR^m \pi_{d!}(\Ql))$ preserves the filtration \eqref{eq: filtration} on $\p \rR^m \pi_{d!} (\Ql)$, and its semisimplification acts by 
\[
\text{$(-2j + (2-n) \deg_X(  \omega_X ))  \otimes  w$ on the summand $\sK_d^i(-i) \otimes \rH^j(\BunU)[\dim \cA_d]$.}
\] 
 \end{thm}

 \subsubsection{Action of the  tautological bundle} 
The following Corollary (of the three preceding Theorems) will be a crucial ingredient of the proof of Theorem \ref{thm: geometric side} below.

\begin{cor}\label{cor: tautological correspondence}The endomorphism $\Gamma_{c_1(\ell_{\cM})} \in \End_{D^b_c(\cA_d)}(\p\rR^m \pi_{d!}(\Ql))$ preserves the filtration \eqref{eq: filtration}, and its semisimplification acts by 
\[
\text{$(d-2i -2j + \deg_X(\omega_X) )  \otimes w$  on the summand $\sK_d^i(-i) \otimes \rH^j(\BunU)[\dim \cA_d]$.}
\]
\end{cor}

\begin{proof}
By Proposition \ref{prop:bundle decomp},  there is a  decomposition of cohomology classes
\[
 c_1(    \ell_{\cM}  )    =    [ C ] + c_1(    \det (\cE^{\univ} )  )    +  c_1(   \ell_{\U^\dagger(1) } )         \in  H^2(\Hk^1_{\cM_d})(1).
\]
This gives a corresponding decomposition of the endomorphism $\Gamma_{  c_1(    \ell_{\cM}  )}$ as a sum of   three terms.  The first term is computed by Theorem \ref{thm: global C action}, the second term is computed by Theorem \ref{thm: universal bundle action}, and the third term is computed by Theorem \ref{thm: tautological action}. 
\end{proof}


\section{Hecke action of the comparison divisor}
\label{sec: comparison correspondence}


This section is devoted to the proof of Theorem \ref{thm: global C action}.
We work with a fixed $d \ge 0$.


\subsection{Hecke stacks for the local Hitchin space}


Recall from \S \ref{ssec: local Hitchin space}  the  local enhanced Hitchin fibration 
\[
 f_d^{\Herm}   : 
\cM_d^{\Herm} \to
\cB_d^{\Herm} = \Herm_{2d} \times_{X_d} \Herm_{2d} .
\]
Our first task is to upgrade this to a  commutative diagram 
\begin{equation}\label{local hecke}
\begin{tikzcd}
{ \cM_d^{\Herm } }  \ar[d , " f_d^{\Herm} " ' ]    & { \Hk_{\cM_d^{\Herm}}^1 }  \ar[l]   \ar[r]      \ar[d] & { \cM_d^{\Herm} }  \ar[d ,  " f_d^{\Herm} " ]   \\
 {  \cB_d^{\Herm  }     }   & {   \Hk^1_{ \cB_d^{\Herm   } }  }   \ar[l]  \ar[r]   &    {   \cB_d^{\Herm }   }
\end{tikzcd}
\end{equation}

The Hecke stack in the bottom row   is  defined by 
\begin{equation}\label{eq: B Hecke}
  \Hk^1_{  \cB_d^{\Herm}  } : = X' \times \cB_d^{\Herm}  .
\end{equation}

Next we will describe the functor of points for $\Hk_{\cM_d^{\Herm}}^1$. First note that given a torsion coherent sheaf $P$ on $X'$ and a divisor $D$ on $X'$, we can define the twist $P(D) := P \otimes_{\cO_{X'}} \cO(D)$. Moreover, using the canonical Hermitian structure on $\cO(\sigma (x')-x')$, such twisting by $D := \sigma (x')-x'$ will carry torsion Hermitian sheaves to torsion Hermitian sheaves. Now suppose we are given a $k$-algebra $R$ and a point  $(x' , Q , P ) \in \Hk^1_{  \cB_d^{\Herm}  }  (R)$.
 From the Hermitian torsion sheaf $P$ we construct a triple $P_{\bu}=(P_0,P_{1/2}, P_1)$ of torsion coherent sheaves 
\begin{equation}\label{Ptwists}
P_0 := P ,  \quad P_{1/2} := P(-x') , \quad P_1 := P( \sigma(x') - x' ) .
\end{equation}
on $X'_R$, and use this construction to regard points of  \eqref{eq: B Hecke} as tuples 
\begin{equation}\label{hecke B point}
( x' , Q,P_\bu  ) \in \Hk^1_{  \cB_d^{\Herm}  }  (R).
\end{equation}
The leftward and rightward arrows in the bottom row of \eqref{local hecke}  are
\[
\begin{tikzcd}
{ (Q,P_0) }   & { (x' , Q,P_\bu )   }    \ar[l, mapsto ]   \ar[r , mapsto ]  & {   (Q,P_1)  } .
\end{tikzcd}
\]

Turning to the top row of \eqref{local hecke}, the Hecke stack in the middle  is defined by its $R$-points
\[
(  x' , Q, P_{\bu}, L_{\bu}  ) \in \Hk_{\cM_d^{\Herm}}^1(R) ,
\]
 which consist of a point  \eqref{hecke B point},  together with a tuple  $L_{\bu} = (L_0, L_{1/2}, L_1)$  of torsion coherent sheaves on $X'_R$ related by a diagram
\begin{equation}\label{eq: L half diagram}
\begin{tikzcd}
L_0 \ar[d, hook] & L_{1/2} \ar[l] \ar[d, hook] \ar[r] & L_1 \ar[d, hook] \\
Q \oplus P_0  & Q \oplus P_{1/2} \ar[r] \ar[l] & Q \oplus P_1 .
\end{tikzcd}
\end{equation}
We demand that the vertical arrows are inclusions of  $R$-module local direct summands of rank  $2d$, 
and that $L_i \subset Q \oplus P_i$ is a balanced Lagrangian (Definition \ref{def: balanced})  for $i\in \{ 0,1\}$.
The leftward and rightward arrows in the top row of  \eqref{local hecke} are
\[
\begin{tikzcd}
{ (Q, P_0, L_0) }   & { (  x' , Q, P_{\bu}, L_{\bu}  ) }    \ar[l, mapsto ]   \ar[r , mapsto ]  & {  (Q, P_1, L_1)  } .
\end{tikzcd}
\]

\begin{lemma}\label{lem: cartesian cube}
There is a commutative diagram 
\begin{equation}\label{eq:  global-semilocal Hecke}
\adjustbox{scale = 0.75}{
\begin{tikzcd}[column sep = tiny]
{  \cM_d }  \ar[dr, dashed] \ar[ddd, dashed] & & {   \Hk_{\cM_d}^1  }  \ar[rr]  \ar[ll] \ar[dr, dashed] \ar[ddd, dashed] & &  { \cM_d }  \ar[dr, dashed]  \ar[ddd, dashed] \\
&  \cA_d  \times \BunU  \ar[ddd, dashed] & & { \cA_d \times \HkU^1 }   \ar[rr] \ar[ll] \ar[ddd, dashed] &  &  {  \cA_d \times \BunU } \ar[ddd, dashed] \\  
\\
\cM_d^{\Herm}  \ar[dr, dashed] & & \Hk_{\cM_d^{\Herm}}^1 \ar[rr] \ar[ll] \ar[dr, dashed] & & \cM_d^{\Herm}   \ar[dr, dashed]  \\
& \cB_d^{\Herm}   & &  \Hk^1_{  \cB_d^{\Herm}  }  \ar[ll] \ar[rr] &  & \cB_d^{\Herm} 
\end{tikzcd}
}
\end{equation}
in which the three squares with all dashed sides are Cartesian, and also all the squares in the front and back sides are Cartesian. 
\end{lemma}

\begin{proof}
The leftmost and rightmost squares with all dashed arrows are the Cartesian diagram of Proposition \ref{prop: hitchin global to local}.  The top face of the diagram is the top half of  \eqref{global hecke}, and the bottom face is \eqref{local hecke}.
It remains to define  the vertical dashed arrows in the middle of the diagram. 

To define the dashed vertical arrow in the middle of the back face, start with a point 
$(x' , \cE , \cF_\bu , t_\bu) \in   \Hk_{\cM_d}^1 (R)$ 
as in \eqref{global hecke point}, and use the complementary line trick to upgrade it to a diagram \eqref{upgraded hecke}.
   From that diagram we construct morphisms of torsion coherent sheaves 
\begin{equation}\label{stray subsheaves}
\begin{tikzcd}
 { \frac{ \cF_0 }{  \cE \oplus \cD_0 }  }   \ar[d, hookrightarrow]  & {   \frac{ \cF^\flat_{1/2} }{  \cE \oplus \cD_{1/2}^\flat  }  }   \ar[d, hookrightarrow]  \ar[l] \ar[r]  &    { \frac{ \cF_1 }{  \cE \oplus \cD_1 }  } \ar[d, hookrightarrow]  \\
    {  Q \oplus  P_0 }    &  {   Q \oplus   P _{1/2}  }    \ar[l] \ar[r] &  {   Q  \oplus     P_1 }     
\end{tikzcd}
\end{equation}
where  $Q :=   \sigma^* \cE^\vee / \cE$,   and the $P$'s are defined by the equalities in 
\[
 \begin{tikzcd}
 { \frac{ \sigma^* \cD_0^\vee } {  \cD_0  }  = P_0 }   &   {  P_0 ( -x' )  =  P_{1/2} = P_1( - \sigma(x') )   } \ar[r]  \ar[l]  & P_1 = \frac{  \sigma^* \cD_1^\vee }{  \cD_1 } .
   \end{tikzcd}
\]
This diagram has the form  \eqref{eq: L half diagram}, so defines a point of $\Hk^1_{ \cM_d^{\Herm}}(R)$.

To define the dashed vertical arrow in the middle of the front face, start with a point 
$
( \cE , x'  ,  \cL_\bu) \in   ( \cA_d \times \Hk_{\U^\dagger(1) }^1)  (R).
$
As in Remark \ref{rem: Hk cline},  twist the morphisms in  \eqref{eq: HkU mod} by $\det(\sigma^*\cE)$  to obtain  a modification of line bundles $\cD_\bu$ as in  \eqref{Dmod}, and then define  $Q$ and $P_{\bu}=(P_0,P_{1/2} , P_1)$ exactly as above.  The data $(x', Q , P_{\bu})$ defines a point \eqref{hecke B point}.

The claimed Cartesian properties not already proved in Proposition \ref{prop: hitchin global to local} follow directly from the constructions.
\end{proof}


 \subsection{Local comparison divisors} 
 \label{ss:local comparison divisors}
 

 As all of the stacks in \eqref{eq:  global-semilocal Hecke}  admit canonical morphisms to $X_d$,   we can pull back the entire diagram along the composition
\[
(X')^{d, \circ} \map{  \eqref{Xtilde torsor} } X_d^\circ \hookrightarrow X_d
\]
to obtain a commutative diagram 
\begin{equation}\label{eq:  global-semilocal Hecke torsor}
\adjustbox{scale = 0.75}{
\begin{tikzcd}[column sep = tiny]
{  \wt \cM^\circ_d }  \ar[dr, dashed] \ar[ddd, dashed] & &   {   \Hk^1_{ \wt \cM^\circ_d }   }   \ar[rr]  \ar[ll] \ar[dr, dashed] \ar[ddd, dashed] & &  {  \wt \cM^\circ_d }  \ar[dr, dashed]  \ar[ddd, dashed] \\
&  { \wt  \cA_d^\circ  \times \BunU}   \ar[ddd, dashed] & & { \wt  \cA_d^\circ  \times  \HkU^1 }    \ar[rr] \ar[ll] \ar[ddd, dashed] &  & { \wt  \cA_d^\circ  \times \BunU}  \ar[ddd, dashed] \\  
\\
{   \wt \cM_d^{\Herm, \circ}  }  \ar[dr, dashed, "\wt f_d^{\Herm}  " ' ] & & \Hk_{   \wt \cM_d^{\Herm, \circ}  }^1 \ar[rr] \ar[ll] \ar[dr, dashed] & & {   \wt \cM_d^{\Herm, \circ}  }   \ar[dr, dashed  , "\wt f_d^{\Herm} " ]  \\
&  { \wt \cB_d^{\Herm,\circ} }   & &  \Hk^1_ { \wt \cB_d^{\Herm,\circ} }  \ar[ll] \ar[rr] &  &  { \wt \cB_d^{\Herm,\circ} } ,
\end{tikzcd}
}
\end{equation}
in which the three squares made of dashed arrows, and all squares in the front and back faces, are Cartesian.
By construction, the two dashed arrows labeled $\wt f_d^{\Herm}$ agree with the arrow of the same name in 
 \eqref{diag: local rss torsor hitchin}.

To spell out the moduli interpretations of the two Hecke stacks in the bottom face,  an $R$-valued  point 
\begin{equation}\label{hecke B rss point}
( x' , Q, P_{\bu}  , x'_1,\ldots, x_d' ) \in \Hk_{ \wt \cB_d^{\Herm,\circ}}^1(R) 
\end{equation}
is a point \eqref{hecke B point}, together with a presentation  of the (common) support of the torsion sheaves $Q$, $P_0$, $P_{1/2}$, and $P_1$ as a multiplicity-free Cartier divisor
\begin{equation}\label{comp-comp sppt}
x'_1+ \cdots + x'_d + \sigma(x'_1) + \cdots +  \sigma(x'_d) 
\end{equation}
for points  $x_1',\ldots, x_d' \in X'(R)$.

To give a lift of \eqref{hecke B rss point} to  $\Hk_{\wt \cM_d^{\Herm,\circ}}^1(R)$ is to further give a diagram 
 \eqref{eq: L half diagram} of torsion coherent sheaves on $X_R'$.  
 This diagram is uniquely  determined by the collection (indexed by $1\le i \le d$) of diagrams
\begin{equation}\label{heckeP diagram}
\begin{tikzcd}
L_0  |_{x_i'}     \ar[d, hook] & L_{1/2}   |_{x_i'}    \ar[l] \ar[d, hook] \ar[r] & L_1 |_{x_i'}    \ar[d, hook]  \\
( Q \oplus P_0)|_{x_i'}   &  (  Q \oplus P_{1/2}  ) |_{x_i'}   \ar[r] \ar[l] &  (   Q \oplus P_1 )  |_{x_i'}  
\end{tikzcd}
\end{equation}
of projective $R$-modules obtained from \eqref{eq: L half diagram} by restricting to  each  $x'_i \in  X'(R)$.  
Each  vertical arrow here is  the  inclusion of a rank one local direct summand into a projective $R$-module of rank two.

\begin{definition}\label{def: incidence}
For $1 \le r \le d$, the \emph{incidence divisors} 
\[
I_r , I_r^\star \hookrightarrow \Hk_{\wt \cB_d^{\Herm, \circ}}^1
\]
 are the closed substacks cut out by the conditions $x'= x'_r$ and $x'=\sigma(x'_r)$, respectively, on points \eqref{hecke B rss point}.  The \emph{local comparison divisors}
\[
C_r , C_r^\star \hookrightarrow \Hk_{\wt \cM_d^{\Herm,\circ}}^1
\]
are defined as follows: 
\begin{itemize}
\item 
$C_r$ is the closed substack defined by the condition  $x'=x_r'$ on the underlying point \eqref{hecke B rss point}, and the condition $L_0|_{x'_r} = Q|_{x'_r}$ in the diagram \eqref{heckeP diagram} indexed by $i=r$.  
There is no condition imposed on the diagrams \eqref{heckeP diagram} with $i \neq r$.
\item
$C_r^\star$ is the closed substack defined by the condition  $x'=\sigma(x_r')$ on the underlying point \eqref{hecke B rss point}, and the condition $L_0|_{x'_r} = P_0|_{x'_r}$ in the diagram \eqref{heckeP diagram} indexed by $i=r$.  
There is no condition imposed on the diagrams \eqref{heckeP diagram} with $i \neq r$.
\end{itemize}
  \end{definition}
The map $\Hk_{\wt \cM_d^{\Herm,\circ}}^1 \rightarrow \Hk_{\wt \cB_d^{\Herm, \circ}}^1$ carries $C_r^{(\star)}$ to $I_r^{(\star)}$ for each $1 \leq r \leq d$. Note that the incidence divisors $I_1,\ldots, I_d, I_1^\star, \ldots, I_d^\star$ are pairwise disjoint, and hence the same is true of the local comparison divisors $C_1,\ldots, C_d, C_1^\star, \ldots, C_d^\star$.

\begin{lemma}\label{lem:de-comp divisor}
Consider the morphisms
\[
\begin{tikzcd}
 {  \Hk^1_{ \cM_d }   }  &  {  \Hk^1_{ \wt \cM_d^\circ}  }  \ar[l] \ar[r] & { \Hk_{   \wt \cM_d^{\Herm, \circ}  }^1  }. 
\end{tikzcd}
\]
The pullback of the comparison divisor $C$ from Definition \ref{def:comparison divisor} along the leftward arrow is equal to the pullback of 
 $C_1 \sqcup \cdots \sqcup C_d \sqcup C_1^\star \sqcup \cdots \sqcup C^\star_d $
along the rightward arrow.
\end{lemma}

\begin{proof} 
Recall from \eqref{upgraded hecke} that for a $k$-algebra $R$, a point
$
( x' , \cE , \cF_{\bu} , t_\bu )  \in  \Hk^1_{  \cM_d}   (R)
$
determines injective morphisms of rank $n$ vector bundles
\begin{equation}\label{comp-comp sandwich}
\cE \oplus \cD_0 \to \cF_0 \to \sigma^* \cE^\vee \oplus \sigma^*\cD_0^\vee
\end{equation}
on $X'_R$.
A lift of  the above point to  
\begin{equation}\label{hecke torsor point}
( x' , \cE , \cF_{\bu} , t_\bu, x_1,\ldots, x_d)  \in  \Hk^1_{ \wt \cM_d^\circ}   (R)
\end{equation}
consists of a presentation  of the (common) support of  $Q=\sigma^*\cE^\vee / \cE$ and  $P_0=\sigma^*\cD_0^\vee / \cD_0$ as a multiplicity-free Cartier divisor \eqref{comp-comp sppt}.

Denote by 
\[
\wt C \subset \Hk^1_{\wt \cM_d^\circ}
\]
 the pullback of $C\subset   \Hk^1_{ \cM_d }$.
 By definition, a point \eqref{hecke torsor point} lies on  $\wt C$ when the morphism
$\cD_0 \to \cF_0$ in \eqref{comp-comp sandwich} vanishes along $\sigma(x')$.
This can only happen if $\sigma(x')$ is one of the points  in  \eqref{comp-comp sppt}, for otherwise both morphisms in \eqref{comp-comp sandwich} would restrict to isomorphisms along $\sigma(x')$.  
This allows us to  decompose 
\[
\wt C   = \wt C_1 \sqcup \cdots \sqcup \wt C_d \sqcup \wt C_1^\star \sqcup \cdots \sqcup \wt C_d^\star,
\]
where $\wt C_r \subset \wt C$ is the locus where $x' = x_r'$, and $\wt C^\star_r \subset \wt C$ is the locus where $x' = \sigma(x_r')$.

Equivalently,  $\wt C_r$ is the locus  in $\Hk^1_{\wt \cM_d^\circ}$  where $x'=x_r'$, and the Lagrangian subsheaf 
\[
\frac{ \cF_0 }{  \cE \oplus \cD_0 } \subset    \frac{ \sigma^* \cE^\vee } {  \cE  } \oplus  \frac{ \sigma^* \cD_0^\vee } {  \cD_0  } 
\]
in \eqref{stray subsheaves}  satisfies 
\[
\frac{ \cF_0 }{  \cE \oplus \cD_0 }  \Big|_{ \sigma(x_r')}  =    \left(  0 \oplus  \frac{ \sigma^* \cD_0^\vee } {  \cD_0  }  \right)  \Big|_{ \sigma(x_r')}   .
\]
By the Lagrangian condition, this is equivalent to 
\[
\frac{ \cF_0 }{  \cE \oplus \cD_0 }  \Big|_{x_r'}  =    \left(  \frac{ \sigma^* \cE^\vee } {  \cE  }  \oplus 0 \right)  \Big|_{x_r'}   ,
\]
which is exactly the condition defining the pullback of $C_r$ to $\Hk^1_{\wt \cM_d^\circ}$.

Similar reasoning shows that $\wt C_r^\star$ agrees with the pullback of $C_r^\star$ to $\Hk^1_{\wt \cM_d^\circ}$.
\end{proof}


\subsection{Hecke stacks for projective bundles}
\label{ss:Hecke proj}


Let us now reconsider Lemma \ref{lem:projective bundle} and its proof.  
There we saw that the map $\wt f_d^{\Herm}$ in \eqref{eq: global-semilocal Hecke torsor} could be identified with 
\begin{equation}\label{recall isabundle}
 \wt \cM_d^{\Herm, \circ}  \iso \prod_{ 1 \le i \le d}  \PP_i \to  \wt \cB_d^{\Herm , \circ },
\end{equation}
where each 
\[
\PP_i  : =  \PP( \mathscr{W}_i \oplus \mathscr{V}_i )  \map {p_i }  \wt \cB_d^{\Herm , \circ }  
\]
is  the $\PP^1$-bundle defined by distinguished line bundles $\mathscr{W}_i$ and $\mathscr{V}_i$ on $ \wt \cB_d^{\Herm , \circ }$,
and the fiber products in \eqref{recall isabundle} are taken over the base $\wt \cB_d^{\Herm , \circ }$.
Each bundle map $p_i$  has two distinguished sections, defined by the two summands  in  $\mathscr{W}_i \oplus \mathscr{V}_i$. 
To spell this out,  recall that a point of $\PP_i(R)$ consists of a 
\[
( Q, P   , x'_1,\ldots, x_d' ) \in \wt \cB_d^{\Herm,\circ}(R) 
\]
together with a rank one $R$-module local direct summand of $(Q \oplus P)|_{x'_i}$, which we denote by $L|_{x_i'} \subset (Q \oplus P)|_{x'_i}$.  
We define 
\begin{itemize}
\item
 $\bm{0}_i \hookrightarrow \PP_i$  to be the closed substack cut out by the condition $L |_{x'_i} = Q|_{x'_i} $,
\item
$\bm{\infty}_i \hookrightarrow \PP_i$  to be the closed substack cut out by the condition $L |_{ x'_i} =  P|_{x'_i}$.
\end{itemize}

From the above discussion, we see that the bottom face of \eqref{eq:  global-semilocal Hecke torsor} can be identified with 
\begin{equation}\label{eq: local proj hecke}
\begin{tikzcd}
{   \prod_{i=1}^d \PP_i  } \ar[d, " \prod p_i " '] & { \prod_{i=1}^d  \Hk_{ \PP_i } ^1 }   \ar[l] \ar[r] \ar[d] & {   \prod_{i=1}^d \PP_i  }     \ar[d, "  \prod p_i " ] \\
 { \wt \cB_d^{\Herm , \circ }     }  & {  \Hk_{  \wt \cB_d^{\Herm , \circ } }   }  \ar[l ] \ar[r]  &    { \wt  \cB_d^{\Herm , \circ }  ,  } 
\end{tikzcd}
\end{equation}
in which  $\Hk_{ \PP_i } ^1$ is the moduli stack of points \eqref{hecke B rss point}, together with inclusions of rank one local direct summands as in \eqref{heckeP diagram}. Here and throughout the subsection, the products are formed over the appropriate base (which in \eqref{eq: local proj hecke} means the term vertically below). For any $1 \le r \le d$ we denote by 
\[
\big( \prod_{i=1}^d \Hk^1_{ \PP_i}  \big) |_{I_r}
  =  \prod_{i=1}^d  \big( \Hk^1_{ \PP_i} |_{I_r} \big)   
  \quad \mbox{and}\quad 
 \big( \prod_{i=1}^d \Hk^1_{ \PP_i}  \big) |_{I^\star_r}
  =  \prod_{i=1}^d  \big( \Hk^1_{ \PP_i} |_{I^\star_r} \big)
\]
  the fibers of $\prod_{i=1}^d \Hk^1_{ \PP_i} $ over the incidence divisors $I_r$ and $I^\star_r$, respectively (Definition \ref{def: incidence}).

\begin{lemma}\label{lem:incidence fiber}
Fix $1 \le r \le d$.
If $i \neq r$ then all horizontal arrows in the diagram
\[
\begin{tikzcd}
{  \PP_i  } \ar[d, " p_i " ']    &  {    \Hk_{ \PP_i }^1 |_{I_r}    }   \ar[l, " h_0" '] \ar[r, " h_1" ]   \ar[d]   &   { \PP_i  }   \ar[d, "  p_i "]     \\
{  \wt \cB_d^{\Herm , \circ }}   & { I_r }   \ar[l , " \iso" ' ] \ar[r, "\iso"] &   {   \wt \cB_d^{\Herm , \circ }}
\end{tikzcd}
\] 
are isomorphisms.  If $i=r$ then only the lower horizontal arrows are isomorphisms,  and there is a decomposition 
\[
 \Hk_{ \PP_i }^1 |_{I_r}   = h_0^{-1}( \bm{0}_r) \cup h_1^{-1}( \bm{\infty}_r)
\]
as a union of closed substacks. Moreover,  the morphisms $h_0$ and $h_1$ restrict to isomorphisms
\[
h_0 : h_1^{-1}( \bm{\infty}_r) \iso \PP_i  \quad \mbox{and} \quad h_1 : h_0^{-1}( \bm{0}_r) \iso \PP_i .
\]
  A similar result holds with $I_r$ replaced by $I_r^\star$, except now 
    \[
 \Hk_{ \PP_i }^1 |_{I^\star_r}   =h_0^{-1}( \bm{\infty}_r)   \cup h_1^{-1}( \bm{0}_r).
\]
\end{lemma}

\begin{proof}
We show first that the bottom horizontal arrows in the diagram are isomorphisms.  This is done most easily by exhibiting their inverses.  Suppose we have a point 
\[
(   Q, P   , x'_1,\ldots, x_d' ) \in  \wt \cB_d^{\Herm,\circ}  (R) .
\]
The inverse to the leftward arrow sends  this to the point \eqref{hecke B rss point} defined by $x' = x'_r$ and 
\[
P_0 = P ,  \quad P_{1/2} = P(-x') , \quad P_1 = P( \sigma(x') - x' ),
\]
while the inverse to the rightward  sends it to the point \eqref{hecke B rss point} defined by $x' = x'_r$ and 
\[
P_0 = P(  x' - \sigma(x')  )  ,  \quad P_{1/2} = P( \sigma(x') ) , \quad P_1 = P.
\]

Now we turn to the upper horizontal arrows in the diagram.  An $R$-valued point of the divisor  $\Hk_{ \PP_i }^1 |_{I_r}$ is given by a point \eqref{hecke B rss point} with $x'=x'_r$ and inclusions of rank one local direct summands 
\begin{equation}\label{heckeP diagram 2}
\begin{tikzcd}
L_0  |_{x_i'}     \ar[d, hook] & L_{1/2}   |_{x_i'}    \ar[l] \ar[d, hook] \ar[r] & L_1 |_{x_i'}    \ar[d, hook]  \\
( Q \oplus P_0)|_{x_i'}   &  (  Q \oplus P_{1/2}  ) |_{x_i'}   \ar[r] \ar[l] &  (   Q \oplus P_1 )  |_{x_i'}  
\end{tikzcd}
\end{equation}
as in \eqref{heckeP diagram}.  In this diagram  the bottom horizontal arrows are induced by the natural maps
\[
 \begin{tikzcd}
 { P_0 }   &   {  P_0 ( -x_r ' )  =  P_{1/2} = P_1( - \sigma(x_r ') )   } \ar[r]  \ar[l]  & P_1  .
   \end{tikzcd}
\]

If $i \neq r$ then the divisors $x'_i$ and $x'=x_r'$ are disjoint, the bottom horizontal arrows in \eqref{heckeP diagram 2} are isomorphisms, and  each of the local direct summands in the top row is uniquely determined by the other two.  From this it is clear that the arrows labeled $h_0$ and $h_1$ in the statement of the Lemma are isomorphisms.

Now suppose $i=r$, so that $x'=x'_i=x'_r$.   In this case the bottom row of  \eqref{heckeP diagram 2} becomes
\[
\begin{tikzcd}
( Q \oplus P_0)|_{x_i'}   &  (  Q \oplus P_{1/2}  ) |_{x_i'}   \ar[r , "\iso" ] \ar[l , "\Id \oplus 0" '] &  (   Q \oplus P_1 )  |_{x_i'}  .
\end{tikzcd}
\]
The line  $L_{1/2}|_{x'_i}$ is uniquely determined by $L_1|_{x'_i}$, and its image under the leftward arrow must be contained in $L_0|_{x'_i}$. 

The divisor $h_1^{-1}( \bm{\infty}_i)$ consists of those diagrams \eqref{heckeP diagram 2} for which $L_1|_{x'_i} = P_1|_{x'_i}$.  When this happens, $L_{1/2}|_{x'_i} = P_{1/2}|_{x'_i}$ is in the kernel  of the leftward arrow, and  $L_0|_{x'_i}$ can be chosen arbitrarily. It follows that $h_0$ restricts to an isomorphism $h_1^{-1}( \bm{\infty}_i) \iso \PP_i$.

The divisor $h_0^{-1}( \bm{0}_i)$ consists of those diagrams \eqref{heckeP diagram 2} for which $L_0|_{x'_i} = Q|_{x'_i}$.  When this holds,  the line $L_{1/2}|_{x'_i}$ can be chosen arbitrarily, and then $L_1|_{x'_i}$ is determined by that choice.  It follows that $h_1$ restricts to an isomorphism $h_0^{-1}( \bm{0}_i ) \iso \PP_i$.

No matter what the direct summands are in the top row of  \eqref{heckeP diagram 2}, they must satisfy either 
$L_0|_{x'_i} = Q|_{x'_i}$ or $L_1|_{x'_i} = P_1|_{x'_i}$.  Indeed, suppose the second equality fails.
 The image of  $L_{1/2}|_{x'_i} \neq P_{1/2}|_{x'_i}$ under the leftward arrow is then $Q|_{x'_i}$,  which must therefore be contained in  $L_0|_{x'_i}$.  It follows that $L_0|_{x'_i} = Q|_{x'_i}$.  This shows that the union of 
 $h_0^{-1}( \bm{0}_i)$ and $h_1^{-1}( \bm{\infty}_i)$ is all of  $\Hk_{ \PP_i }^1 |_{I_i}$.
 
 The analysis with $I_r$ replaced by $I_r^\star$ is entirely similar.
 \end{proof}

\begin{lemma}\label{lem:untangle}
Fix $1 \le r \le d$.   The local comparison divisors of Definition \ref{def: incidence}, when regarded as closed substacks of 
\[
\prod_{i=1}^d\Hk^1_{ \PP_i} \iso  \Hk_{   \wt \cM_d^{\Herm, \circ}  }^1 ,
\]
  admit factorizations 
$C_r =  \prod_{i=1}^d  C_{ r ,i }$ and $C^\star_r =   \prod_{i=1}^d  C^\star_{ r , i }$,  in which 
\[
C_{ r  , i } = \begin{cases}
 \Hk_{ \PP_i }^1|_{I_r}   &  \mbox{if  } i \neq r \\
 h_0^{-1} ( \bm{0}_r ) & \mbox{if } i=r
\end{cases} 
\quad \mbox{and} \quad
C^\star_{ r , i } = \begin{cases}
 \Hk_{ \PP_i }^1|_{I^\star_r }   &  \mbox{if  } i \neq r \\
 h_0^{-1} ( \bm{\infty}_r ) & \mbox{if } i=r .
\end{cases} 
\]
\end{lemma}

\begin{proof}
This is  a direct translation of Definition \ref{def: incidence} into the notation of Lemma \ref{lem:incidence fiber}
\end{proof}


\subsection{Proof of Theorem \ref{thm: global C action}}\label{ssec:ComparisonAction}


At last we combine the  results proved above to prove Theorem \ref{thm: global C action}. 
 Recall that we are trying to prove that
\[
\Gamma_{[C]} \in \End_{D_c^b(\cA_{d, \ol{k} })}( \rR \pi_{d!}  \Ql)
\]
acts as $(2d-2i)  \otimes w$ on the summand
\[
\sK_d^i(-i) \otimes \rH^j(\BunU) [ \dim \cA_d] \subset  \p \rR^{m} \pi_{d!}(\Ql).
\]
The proof takes place after base changing all $k$-stacks to $\ol{k}$, so we henceforth omit the base change from the notation.

\subsubsection{Reduction to the regular semisimple locus}
As a first reduction, recall from \eqref{ICK2} that 
\[
\sK_d^i  =  \IC_{\cA_d } ( K_d^{i,\circ}  |_{ \cA^\circ_d} ) .
\]
As an endomorphism of an IC sheaf is determined by its restriction to the generic locus, it suffices to prove the analogous result over the regular semisimple locus $\cA_d^\circ \subset \cA_d$.

\subsubsection{Reduction to a covering} 
Define $\wt \cA_d^\circ$ as the fiber product
\begin{equation}\label{eq:tildeA}
\begin{tikzcd}
\wt \cA_d^\circ \ar[r] \ar[d]  & (X')^{d, \circ} \ar[d , "\eqref{Xtilde torsor}"] \\
\cA_d^\circ \ar[r, " \eqref{all the rss} " '] & X_d^\circ. 
\end{tikzcd}
\end{equation}

If we  pull back the  entire diagram  \eqref{global hecke d} along $\wt \cA_d^\circ \to \cA_d$, the result is  a diagram of correspondences
  \begin{equation}\label{eq:torsor full hecke}
  \adjustbox{scale = 0.75}{
\begin{tikzcd}
{  \wt \cM^\circ_d  }     \ar[dr , " \wt f_d " '  ]   \ar[dd,  bend right,  " \wt \pi_d "' ] &&   {   \Hk^1_{ \wt \cM^\circ_d }   }   \ar[ll ] \ar[rr ]  \ar[d]   & &  {  \wt \cM^\circ_d  }    \ar[dl, " \wt f_d" ]   \ar[dd, bend left,  " \wt \pi_d " ] \\
& { \wt  \cA_d^\circ  \times \BunU}   \ar[dl  ]   & { \wt  \cA_d^\circ  \times  \HkU^1 }     \ar[l] \ar[r]    \ar[d  ] &   { \wt  \cA_d^\circ  \times \BunU}   \ar[ dr ] \\
  {   \wt \cA_d^\circ }  &    & {   \wt \cA_d^\circ }  \ar[ll, equal ] \ar[rr , equal ] & &  {   \wt \cA_d^\circ .}
\end{tikzcd}
}
\end{equation}
The comparison divisor $C \hookrightarrow \Hk^1_{  \cM } $ pulls back to a divisor 
$
\wt C  \hookrightarrow  \Hk^1_{ \wt \cM^\circ_d } 
$
 whose fundamental class determines, exactly as in \eqref{hitchin correspondence endo},  a cohomological correspondence
\[
\wt \pi_{d !} [ \wt C  ] \in  \Corr_{ \wt \cA^\circ_d }(   \rR  \wt \pi_{d !} \Ql ,  \rR  \wt \pi_{d !} \Ql )  =  \End_{  D_c^b(  \wt \cA^\circ_d  ) } (  \rR  \wt \pi_{d !} \Ql) .
\]

By pulling back the isomorphism  \eqref{eq:Hitchin-perverse-cohomology} along $\wt \cA_d^\circ \to \cA_d$, we obtain 
\[
  \left( \p\rR^m  \wt \pi_{d !} \Ql  \right) [-m] 
  \iso     
  \bigoplus_{  \dim \cA_d + 2i + j =m   }    K_d^{ i ,\circ}  \tw{-i}  |_{  \wt \cA^\circ_d }   \otimes     \rH^j( \BunU )  [ -j] 
 \in D_c^b(  \wt \cA^\circ_d  ).
\]
We note that each $K_d^{ i ,\circ}   |_{  \wt \cA^\circ_d }$ is a constant local system, by construction.

Thanks to Lemma \ref{lem: cartesian cube}, the Base Change Theorem for cohomological correspondences 
 \cite[Theorem 5.1.3]{FYZ3}\footnote{We are actually invoking a relatively easy special case of the Theorem, since all the relevant commutative squares are in fact Cartesian.} applies, and says that the pullback of $\Gamma_{ [C] } =\pi_{d!} [C]$ along $\wt \cA_d^\circ \rightarrow \cA_d^\circ$ agrees with $\wt \pi_{d!} [\wt C]$. The assertion that an endomorphism of a local system is a particular scalar can be checked after pullback along a cover (or indeed, along any map whose image intersects every connected component non-emptily). Hence, in order to prove Theorem \ref{thm: global C action}, it suffices to prove that  
 $\wt \pi_{d !} [ \wt C  ] $ acts as $(2d-2i) \otimes w$ on the summand
\begin{equation}\label{untangled summand}
 K_d^{ i ,\circ} \tw{-i}  |_{  \wt \cA^\circ_d }   \otimes \rH^j( \BunU )[-j] \subset ( \p \rR^m \wt \pi_{d!}\Ql)[-m] .
 \end{equation}

  \subsubsection{Completion of the proof} 
Return to the diagram \eqref{eq: global-semilocal Hecke torsor}.  In \S \ref{ss:local comparison divisors} and \S \ref{ss:Hecke proj} we gave a detailed analysis of the lower face of that diagram, and that analysis allows us to understand the structure of the top face.

For example, the identification of the bottom face with \eqref{eq: local proj hecke} shows that the top face has a similar structure, allowing us to identify \eqref{eq:torsor full hecke} with 
 \begin{equation}\label{untangled global hecke}
  \adjustbox{scale = 0.75}{
\begin{tikzcd}
{  \prod_{i=1}^d \PP_i  }     \ar[dr , " p=\prod p_i " '  ]   \ar[dd,  bend right,  " \wt \pi_d "' ] &&  {  \prod_{i=1}^d  \Hk^1_{ \PP_i}   }   \ar[ll ] \ar[rr ]  \ar[d]   & &   {  \prod_{i=1}^d \PP_i   }   \ar[dl, "  p= \prod p_i " ]   \ar[dd, bend left,  " \wt \pi_d " ] \\
& { \wt  \cA_d^\circ  \times \BunU}   \ar[dl  , "a"']   & { \wt  \cA_d^\circ  \times  \HkU^1 }      \ar[l] \ar[r]    \ar[d  ] &   { \wt  \cA_d^\circ  \times \BunU}   \ar[ dr, "a" ] \\
  {   \wt \cA_d^\circ }  &    & {   \wt \cA_d^\circ }  \ar[ll, equal ] \ar[rr , equal ] & &  {   \wt \cA_d^\circ  .}
\end{tikzcd}
}
\end{equation}
Again, products are formed over their respective bases, found by following the arrows in the diagram. Here each $\PP^1$-bundle  $p_i : \PP_i \to \wt  \cA_d^\circ  \times \BunU$ is defined as the pullback along the morphism 
\[
\wt  \cA_d^\circ  \times \BunU \to  \wt \cB_d^{\Herm , \circ }
\]
in \eqref{eq: global-semilocal Hecke torsor} of the bundle of the same name  in \eqref{eq: local proj hecke}, and similarly for  $\Hk^1_{ \PP_i}$. In particular, each $\PP_i$ comes equipped  with two distinguished sections to the bundle map, whose images we again denote by $\bm{0}_i , \bm{\infty}_i \hookrightarrow \PP_i$.

For any $1\le r \le d$ there are incidence divisors $I_r, I_r^\star \hookrightarrow \wt  \cA_d^\circ  \times \BunU$, defined as the pullbacks of the incidence divisors of Definition \ref{def: incidence}.
 Equivalently, $I_r$ is the pullback of the diagonal along 
\[
\wt  \cA_d^\circ  \times  \HkU^1       \to (X')^{d,\circ} \times X' \map{  \mathrm{pr}_r \times \Id } X' \times X',
\]
and similarly for $I_r^\star$ but with $ \mathrm{pr}_r \times \Id$ replaced by $ \mathrm{pr}_r \times \sigma$.

Similarly, we may pull back the local comparison divisors of Definition \ref{def: incidence} to obtain divisors (for all $1\le r \le d$)
\[
C_r, C_r^\star \hookrightarrow   \Hk^1_{ \wt \cM^\circ_d }  \iso \prod_{i=1}^d  \Hk^1_{ \PP_i} ,
\]
which factor  as  $C_r \iso \prod_{i=1}^d C_{r,i}$ and $C^\star_r \iso \prod_{i=1}^d C^\star_{r,i}$ by 
 Lemma \ref{lem:untangle}.
 The  diagram \eqref{untangled global hecke} restricts to 
 \begin{equation}\label{eq:hecke punchline}
  \adjustbox{scale = 0.75}{
\begin{tikzcd}
{  \prod_{i=1}^d \PP_i  }     \ar[dr , " p " '  ]   \ar[dd,  bend right,  " \wt \pi_d "' ] &&  {  C_r    }   \ar[ll ] \ar[rr , " \iso" ]  \ar[d]   & &   {  \prod_{i=1}^d \PP_i   }   \ar[dl, " p " ]   \ar[dd, bend left,  " \wt \pi_d " ] \\
& { \wt  \cA_d^\circ  \times \BunU}   \ar[dl  , "a"']   & { I_r  }    \ar[l  , " \iso" ' ] \ar[r, " \iso" ]    \ar[d  ] &   { \wt  \cA_d^\circ  \times \BunU}   \ar[ dr, "a" ] \\
  {   \wt \cA_d^\circ }  &    & {   \wt \cA_d^\circ }  \ar[ll, equal ] \ar[rr , equal ] & &  {   \wt \cA_d^\circ .}
\end{tikzcd}
}
\end{equation}
The isomorphisms in the middle row are clear from the definition of the incidence divisor and the morphisms  in \eqref{eq: HkU correspondence}.   
The top row is the product of morphisms
\[
\begin{tikzcd}
{  \PP_i }  & {   C_{r,i}  }   \ar[l]   \ar[r] &  {  \PP_i  } .
\end{tikzcd}
\]
These are isomorphisms if $i\neq r$.  When $i=r$ the rightward arrow is an isomorphism, but the leftward arrow contracts $C_{r,r}$ to the divisor $\bm{0}_r \subset \PP_r$.  These claims are clear from Lemmas \ref{lem:incidence fiber} and \ref{lem:untangle}.
There is a similar picture with  $C_r$,  $I_r$, and $\bm{0}_r$ replaced by $C_r^\star$,  $I_r^\star$, and $\bm{\infty}_r$.
 Lemma \ref{lem:de-comp divisor} tells us that 
\[
\wt C  = C_1 \sqcup \cdots \sqcup C_d \sqcup C_1^\star \sqcup \cdots \sqcup C_d^\star,
\]
 and  we are reduced to computing the endomorphism  of  $\rR  \wt \pi_{d !} \Ql$ associated to each divisor on the right hand side.

The endomorphism $\wt \pi_{d !} [ \wt C_r  ] $ can be computed as the image of the fundamental class $[C_r]$ under the composition from upper left to lower right in the commutative diagram
 \[
 \begin{tikzcd}
{   \Corr_{  C_r    }(   \Ql,   \Ql )   }  \ar[r]  \ar[d, "\rR p_!" ' ] &  {    \Corr_{ \Hk^1_{ \wt \cM^\circ_d } }(   \Ql,   \Ql )   }  \ar[d , "   \rR \wt f_{d !}   " ] \\
 { \Corr_{ I_r }  (   \rR p_! \Ql  ,  \rR p_! \Ql   )   }  \ar[r] \ar[d]   &  {    \Corr_{   \wt  \cA_d^\circ  \times \HkU^1   }(    \rR \wt f_{d !}  \Ql    ,   \rR \wt f_{d !}  \Ql  )   }   \ar[d]    \\
 {    \Corr_{ \wt \cA^\circ_d }(   \rR  \wt \pi_{d !} \Ql ,  \rR  \wt \pi_{d !} \Ql )   }    \ar[ r, equal ]  & {  \End_{  D_c^b(  \wt \cA^\circ_d  ) } (  \rR  \wt \pi_{d !} \Ql) } .
 \end{tikzcd}
 \]
For this computation, regard the three horizontal rightward isomorphisms in \eqref{eq:hecke punchline} as identifications.  The diagram then takes the form 
\begin{equation}\label{eq:hecke punchline simplified}
  \adjustbox{scale = 0.75}{
\begin{tikzcd}
{  \prod_{i=1}^d \PP_i  }     \ar[dr , " p  " '  ]   \ar[dd,  bend right,  " \wt \pi_d "' ] &&  { \prod_{i=1}^d \PP_i  }   \ar[ll , " \beta_r " '  ] \ar[rr, equal ]  \ar[d]   & &   {  \prod_{i=1}^d \PP_i   }   \ar[dl, " p " ]   \ar[dd, bend left,  " \wt \pi_d " ] \\
& { \wt  \cA_d^\circ  \times \BunU}   \ar[dl , "a"  ' ]   & { \wt  \cA_d^\circ  \times \BunU}    \ar[l ,  " w_r" ' ] \ar[r, equal ]    \ar[d, "a"  ] &   { \wt  \cA_d^\circ  \times \BunU}   \ar[ dr  , "a" ] \\
  {   \wt \cA_d^\circ }  &    & {   \wt \cA_d^\circ }  \ar[ll, equal ] \ar[rr , equal ] & &  {   \wt \cA_d^\circ ,}
\end{tikzcd}
}
\end{equation}
where  $a$ is projection to the first factor,  and $w_r$ is the automorphism  
\[
 w_r ( \cE , a ,  x'_1, \ldots, x'_d, \mathscr{L} ) =  ( \cE , a ,  x'_1, \ldots, x'_d,\mathscr{L}(  x'_r - \sigma( x'_r  )  )  .
\]
There is a factorization  $\beta_r = \prod_i \beta_{r,i}$ in which  $\beta_{r,i}  : \PP_i \to \PP_i$ is an isomorphism if $i\neq r$, and has image $\bm{0}_r$ if $i=r$.

These identifications allow us to regard $ \rR p_! [C_r]$ as an element of
\begin{align*}
 \Corr_{ I_j }  (   \rR p_! \Ql  ,  \rR p_! \Ql   ) 
 & \iso  \Corr_{  \wt  \cA_d^\circ  \times \BunU     }  (   \rR p_! \Ql  ,  \rR p_! \Ql   )  \\
& = {  \Hom_{  \wt  \cA_d^\circ  \times \BunU     }  (  w_r^* \rR p_!   \Ql ,  \rR p_! \Ql  )  } ,
\end{align*}
where it is easy to make explicit.
Let us abbreviate 
\[
K := \bigotimes_{i=1}^d   (   \Ql  \oplus  \Ql \tw{-1} ) \in D_c^b(  \wt  \cA_d^\circ  ) ,
\]
and denote by $b_r : K \to K$ the endomorphism that is the  identity on the tensor factors indexed by $i \neq r$, but is 
\[
\Ql \oplus  \Ql \tw{-1}  \map{ \mathrm{id} \oplus 0 } \Ql \oplus  \Ql \tw{-1} 
\]
on the $r^\mathrm{th}$ tensor factor.
 If we use  the canonical $\rR p_{i!} \Ql  \iso \Ql \oplus \Ql\tw{-1}$ to identify
\[
\rR  p_! \Ql  \iso   \bigotimes_{i=1}^d  \rR  p_{i ! }  \Ql   \iso  K \boxtimes \Ql \in D_c^b(  \wt  \cA_d^\circ  \times \BunU ) 
\] 
(as in the proof of Proposition \ref{prop: local Hitchin local system}), then $ \rR p_! [C_r]$ is identified with the morphism
\begin{equation}\label{explicit c_r}
w_r^* \rR p_!   \Ql  \iso   K \boxtimes w_r^* \Ql   \map{  b_r \boxtimes i_r } K \boxtimes \Ql \iso  \rR p_!   \Ql .
\end{equation}
Here $i_r: w_r^* \Ql \iso \Ql$  is the canonical isomorphism.

Applying proper base change  to the bottom half of \eqref{eq:hecke punchline simplified}, 
we have a canonical  isomorphism of functors $\rR a_! \iso \rR a_! \circ w_r^*$, and  the composition 
\[
\bigoplus_{j\ge 0} \rH^j_c( \BunU)  [ -j ] = \rR  a_! \Ql \iso  \rR a_!  ( w_r^* \Ql ) \map{ (\rR a_!) ( i_r)  } \rR  a_! \Ql = \bigoplus_{j\ge 0} \rH^j_c( \BunU)  [ -j ]
\]
is the involution $w$ from \eqref{eq:w-defn}.   Applying $\rR a_!$ throughout \eqref{explicit c_r} therefore shows that the endomorphism 
$\wt \pi_{d!} [ C_r]$ has the explicit form
\[
 \rR  \wt \pi_{d !} \Ql  \iso K \otimes \bigoplus_{j\ge 0} \rH^j_c( \BunU)  [ -j ]
  \map{ b_r \otimes w} K \otimes \bigoplus_{j\ge 0} \rH^j_c( \BunU)  [ -j ]   \iso  \rR  \wt \pi_{d !} \Ql .
\]
The same holds with $C_r$ replaced by $C_r^\star$.

As in the proof of Proposition \ref{prop: local Hitchin local system}, there are canonical isomorphisms
\begin{equation}\label{Kombinatorial}
 \bigoplus_{i=1}^d   K_d^{ i ,\circ}  \tw{-i}  |_{  \wt \cA^\circ_d }    \iso K \iso   \bigoplus_{  S \subset \{1,\ldots, d\}   }     \bigotimes_{ i \in S} \Ql \tw{-1} .
\end{equation}
The endomorphism $b_r \in \End(K)$  kills those summands  $\bigotimes_{ i \in S} \Ql \tw{-1}$ for which $r \in S$, and acts as the identity on those for which $r \not\in S$.    Elementary combinatorics then shows that  $b_1 + \cdots + b_d \in \End(K)$ acts as multiplication by $d-i$ on the summand of \eqref{Kombinatorial} supported in cohomological degree $2i$, namely $ K_d^{ i ,\circ}  \tw{-i}  |_{  \wt \cA^\circ_d }$.

The conclusion is that  the endomorphism  $\wt \pi_{d!} [C_1] + \cdots +\wt \pi_{d!} [C_d]$ of
\[
\rR  \wt \pi_{d !} \Ql = \bigoplus_{m\in \Z} ( \p \rR^m  \wt \pi_{d !} \Ql ) [ -m ] 
\]
acts as $(d-i) \otimes w$ on  the summand \eqref{untangled summand}.  The same is true  if every  $C_r$ is replaced by $C_r^\star$, completing the proof.
\qed


\section{Hecke action of the universal bundle}\label{sec:univ-action}


In this section, we will prove Theorem \ref{thm: universal bundle action}. 
The key step is to understand how cohomological correspondences behave under pullback along the vertical arrows in the diagram
\[
\begin{tikzcd}
{ \cM }   \ar[d]  &  {  \Hk_{\cM}^1   } \ar[l]      \ar[ r   ]   \ar[d] &  { \cM  } \ar[d]  \\
{  \BunU   }    & {  \Hk_{ \U^\dagger(1) }^1   }   \ar[l  ]      \ar[ r ]   &  {  \BunU } 
\end{tikzcd}
\]
extracted from \eqref{MtoBun hecke}.  Although the cohomological correspondence $c_1( \cE^{\univ})$ appearing in Theorem \ref{thm: universal bundle action} is not such a pullback, it can be related to one that is. 

The cohomological correspondence $ c_1(\ell_{\U^\dagger(1) })$ appearing in Theorem \ref{thm: tautological action} does arise as such  a pullback, and the results proved here, especially Proposition \ref{prop: hecke compatible}, will be essential for proving Theorem \ref{thm: tautological action} in \S \ref{sec:taut-action}.

 
 \subsection{Hecke action on $\BunU$}
 \label{ssec:univ-action-Bun}


Let us consider in more detail the Hecke correspondence 
\[
\BunU \xleftarrow{h_0} \HkU^1 \xrightarrow{h_1} \BunU 
\]
from \eqref{eq: HkU correspondence}. 
Both morphisms are smooth and proper  of relative dimension $1$, and so, exactly as 
 in \eqref{sssec:cohomological-correspondences}, there is an induced isomorphism 
\begin{align*}
\Corr_{\HkU^1}(\Ql, \Ql ) &  \cong H^2(\HkU^1) (1)  .   \nonumber
\end{align*}

Since $\BunU$ is proper, the pushforward of cohomological correspondences
\[
\Corr_{\HkU^1}(\Ql , \Ql ) 
\to  \Corr_{\Spec (k)}(a_! \Qll{\BunU}, a_! \Qll{\BunU}) \rightarrow \End(\rH^*(\BunU))^{\Frob}
\]
is defined, where $a : \BunU \to \Spec(k)$ is the structure morphism. 
In particular, any class $\mf{c} \in H^2(\HkU^1 ) (1) $ induces a Frobenius-equivariant endomorphism 
\begin{equation}\label{BunU Hecke endomorphism}
\Gamma_{\mf c} \co \rH^*(\BunU) \rightarrow \rH^*(\BunU)
\end{equation}
given explicitly, as in Remark \ref{rem:stalk hecke}, by 
\begin{equation}\label{explicit cohomological endomorphism}
\alpha   \mapsto  h_{1!}(\mf{c} \cdot h_0^*\alpha).
\end{equation}
Here we again remind the reader of the conventions for cohomology from \S \ref{ss:notation}.

Let us compute the endomorphism \eqref{BunU Hecke endomorphism} for a simple choice of $\mf{c}$.
Denote by 
\begin{equation}\label{eq:HkUlegmap}
 p : \HkU^1 \to X'
\end{equation}
  the morphism sending a point $(x' ,\cL_\bu )$ as in \eqref{eq: HkU mod} to the underlying $x'$. 
Equivalently, using the identification of Remark \ref{rem: HkU}, it is the projection map
\[
\HkU^1  \iso X' \times \BunU \to X'. 
\]

\begin{prop}\label{prop:pullback bundle correspondence}
If  $\mf{M}$ is any line bundle on $X'$,  the endomorphism 
\[
\Gamma_{c_1( p^* \mf{M} )} \co \rH^*(\BunU) \rightarrow \rH^*(\BunU),
\]
induced by $c_1( p^* \mf{M} ) \in  H^2( \HkU^1 )  (1)$,  is equal to $ \deg_{X'}(\mf{M} ) \cdot w$, where $w$ is the order two automorphism \eqref{eq:w-defn}.
\end{prop}

\begin{proof}Let us base change all of our $k$-stacks to $\ol{k}$, and omit the base change from the notation. By additivity, we are immediately reduced to the case where $\mf{M}=\cO_{X'}(x')$ for a single closed point $x' \in X'$.  
The first Chern class $c_1(p^* \mf{M})$ is then equal to the cycle class of the divisor 
\[
\{ x' \} \times \BunU \subset  X' \times \BunU \iso \HkU^1.
\]

It is clear from Remark \ref{rem: HkU}  that the two arrows in 
\eqref{eq: HkU correspondence} restrict to \emph{isomorphisms}
\[
\BunU \xleftarrow{h_0} \{ x' \} \times \BunU \xrightarrow{h_1} \BunU,
\]
and so the endomorphism \eqref{BunU Hecke endomorphism} induced by $\mf{c}=c_1( p^* \mf{M} )$ is simply the map on cohomology induced by the automorphism 
\[
 h_1 \circ h_0^{-1} : \BunU \to \BunU. 
\]
This automorphism sends $\cL \mapsto \cL ( \sigma(x')  - x')$, and   so induces \eqref{eq:w-defn} on cohomology.
\end{proof}


\subsection{Compatibility of Hecke actions}
\label{ss:compatibility statement}


Fix $d \ge 0$.

For a  cohomology class  $\mf{c} \in H^2( \HkU^1)(1)$,  let us again consider the endomorphism  $\Gamma_{\mf{c}}$ from 
\eqref{BunU Hecke endomorphism}.  If we denote by 
$
\mf{c}_{\cM} \in H^2(  \Hk_{\cM_d}^1 )(1) 
$
the  pullback of $\mf{c}$ along the morphism  
\[
 \Hk_{\cM_d}^1  \map{\eqref{MtoBun hecke}}   \HkU^1 ,
 \]
the construction \eqref{hitchin correspondence endo} defines an endomorphism
\[
\Gamma_{\mf{c}_{\cM}}  \in \End_{  D_c^b ( \cA_d) } (  \rR \pi_{d!}  \Ql  )  .
\]

Now return to the  Frobenius-equivariant isomorphism
\[
\p \rR^{m}   \pi_{d!}(\Ql)  \cong \bigoplus_{\dim \cA_d + 2i+j= m} \sK_d^i(-i) \otimes \rH^j(\BunU)[\dim \cA_d]
\]
of perverse sheaves on $ \cA_{d,\ol{k}}$ from Corollary \ref{cor:Hitchin-perverse-cohomology}.
It would be natural to expect that the endomorphism $\Gamma_{\mf{c}_{\cM}}$ of the left hand side respects this decomposition,  and acts as $\mathrm{Id} \otimes \Gamma_{\mf{c}}$ on each summand on the right, but this expectation is  too naive.
Instead, we have the following weaker statement.

\begin{prop}\label{prop: hecke compatible}
The endomorphism $\Gamma_{\mf{c}_{\cM}}$ of $\p \rR^{m}  \pi_{d!}(\Ql)$ 
respects the filtration \eqref{eq: filtration}, and its semisimplification acts by
\[
\text{$\mathrm{Id}  \otimes \Gamma_{\mf{c}} $ on the summand $\sK_d^i(-i) \otimes \rH^j(\BunU)[\dim \cA_d]$.}
\]
\end{prop}

 \begin{remark}
The statement and proof of Proposition \ref{prop: hecke compatible} are the same if one starts instead with a \emph{geometric} cohomology class $\mf{c} \in \rH^2( \HkU^1)(1)$, except of course the endomorphism  $\Gamma_{\mf{c}}$ need not be Frobenius-equivariant, and $\Gamma_{\mf{c}_{\cM}}$ is only defined after base change to $ \cA_{d,\ol{k}}$.
 This minor variant of Proposition \ref{prop: hecke compatible} will be used in \S \ref{ss:proofUBA}.
\end{remark}

 We will prove Proposition \ref{prop: hecke compatible}  by reducing it to Proposition \ref{prop:reduced hecke comparison} below.
 That latter Proposition  will then be proved in \S \ref{ssec:proof-prop-RHC}.  
 The proof  is not  formal, and uses the detailed analysis  of  the geometry of the Hecke correspondences  \eqref{MtoBun hecke}  that was carried out in  \S \ref{sec: comparison correspondence}.
 In \S \ref{ss:proofUBA} we will deduce  Theorem \ref{thm: universal bundle action} as a  consequence of Propositions \ref{prop:pullback bundle correspondence} and  \ref{prop: hecke compatible}.


\subsubsection{A corollary}


 Before all of this, we record a Corollary of Proposition \ref{prop: hecke compatible} that will be needed in the proof of Theorem \ref{thm: geometric side}.
Let $\mf{M}$ be a line bundle on $X'$, and denote by $p$ the composition 
\[
\Hk_{\cM_d}^1 \map{\eqref{MtoBun hecke}} \HkU^1 \map{  \eqref{eq:HkUlegmap}  } X' .
\]
By the construction of \S \ref{sssec:cohomological-correspondences}, we obtain an endomorphism $\Gamma_{c_1( p^* \mf{M} )} \in \End_{D_c^b (\cA_d)}( \rR \pi_{d!} (\Ql))$.

\begin{cor}\label{cor:pullback-bundle-action}
The endomorphism 
\[
\Gamma_{c_1(p^*\mf M)} \in \End_{D^b_c(\cA_d)}( \p \rR^m \pi_{d!} (\Ql))
\]
 preserves the filtration \eqref{eq: filtration}, and its semisimplification acts by $\deg_{X'}(\mf{M} ) \cdot w$.
\end{cor}

\begin{proof}
This follows from  Propositions \ref{prop:pullback bundle correspondence} and \ref{prop: hecke compatible}, by taking 
$
\mf{c} \in H^2( \HkU^1)(1)
$
 to be the first Chern class of the pullback of $\mf{M}$ along \eqref{eq:HkUlegmap}.
\end{proof}


\subsection{A series of reductions}


We  reduce  Proposition \ref{prop: hecke compatible}  to the more tractable Proposition \ref{prop:reduced hecke comparison} stated below.

\subsubsection{Reduction to $\ol{k}$}
The statement of  Proposition \ref{prop: hecke compatible} only involves the base change of $\p \rR^{m}  \pi_{d!}(\Ql)$ to $\cA_{d,\ol{k}}$.  Henceforth we base change all $k$ stacks to $\ol{k}$, and omit the base change from the notation.

We  fix an isomorphism $\Ql(1) \iso \Ql$, to avoid cluttering the notation with Tate twists.

\subsubsection{Reduction to the regular semisimple locus}
We are examining an endomorphism of a perverse sheaf on $\cA_d$, which is expressed as  the direct sum of intermediate extensions of (shifted) local systems on the regular semisimple locus $\cA_d^\circ \subset \cA_d$. 

The filtration \eqref{eq: filtration} is the intermediate extension of a filtration on the local system on the regular semisimple locus. Therefore, it suffices to prove Proposition \ref{prop: hecke compatible} after restriction to the  regular semisimple locus $\cA_d^\circ$. 

\subsubsection{Reduction to a covering}
Since the assertion is now about an endomorphism of a local system, it can be checked after pulling back along the cover $\wt{\cA}_d^\circ \rightarrow \cA_d^\circ$ from \eqref{eq:tildeA}.

\subsubsection{Reduction to a stalk}
 Since the assertion is  about an endomorphism of a local system,  it suffices to show that $\Gamma_{  \mf{c}_{\cM} }$ acts in the correct way on the stalk of $\p \rR^{m}  \pi_{d!}(\Ql)$ at every  geometric point
\begin{equation}\label{eq:mfa point}
\mf{a}:= (\cE, a,  \ul{x}')  \in \wt{\cA}_d^\circ  (\ol{k}),
\end{equation}
where $\ul{x}' = (x_1', \ldots, x_d') \in (X')^{d,\circ}(\ol k)$ is a tuple of points whose 
associated (multiplicity-free) Cartier divisor 
\begin{equation}\label{eq:CompareSupp}
x_1' + \ldots + x_d' + \sigma(x_1') + \ldots + \sigma(x_d')
\end{equation}
is the  support divisor of the Hermitian torsion sheaf  $Q := \coker(a : \cE \to \sigma^* \cE^\vee)$ on $X'$.

The endomorphism $ \Gamma_{ \mf{c}_{\cM} } $ we are studying arises from the diagram  \eqref{global hecke d},
and by our reductions so far we may replace this diagram  with the fiber over $\mf{a}$ of the diagram \eqref{eq:torsor full hecke}.
We once again identify this last diagram with \eqref{untangled global hecke}, so that the fiber over $\mf{a}$ takes the form 
 \begin{equation}\label{untangled hecke fiber}
  \adjustbox{scale = 0.75}{
\begin{tikzcd}
{  \prod_{i=1}^d \PP_i  }     \ar[dr , " p=\prod p_i " '  ]   \ar[dd,  bend right,  " \wt \pi_d "' ] &&  {  \prod_{i=1}^d  \Hk^1_{ \PP_i}   }   \ar[ll , "h_0" ' ] \ar[rr , " h_1" ]  \ar[d]   & &   {  \prod_{i=1}^d \PP_i   }   \ar[dl, "  p= \prod p_i " ]   \ar[dd, bend left,  " \wt \pi_d " ] \\
& { \BunU}   \ar[dl  , "a"']   & {  \HkU^1 }      \ar[l] \ar[r]    \ar[d  ] &   { \BunU}   \ar[ dr, "a" ] \\
  {   \Spec(\ol{k}) }  &    & {   \Spec(\ol{k}) }  \ar[ll, equal ] \ar[rr , equal ] & &  {    \Spec(\ol{k}) . }
\end{tikzcd}
}
\end{equation}
Here  each $\PP_i = \PP( \sW_i \oplus \sV_i)$ is the projective bundle associated to a pair of line bundles on $\BunU$.
Tracing back to their origins in Lemma \ref{lem:projective bundle}, the fibers of $\sW_i$ and $\sV_i$ at a point
$\cL \in \BunU(\ol{k})$ are the $\ol{k}$-vector spaces  $Q|_{ x_i '}$ and $P|_{x'_i}$, respectively, where $Q$ and $P$ are the Hermitian torsion sheaves on $X'$ from \eqref{eq: PQpair}.  Note that $P$ depends on both $\cE$ and $\cL$, but $Q$ depends only on $\cE$.  In particular, as $\cE$ is fixed,  $\sW_i$ is (noncanonically)  trivial.

As in \S \ref{ss:Hecke proj}, denote by 
\begin{equation}\label{eq:comp0}
\bm{0}_i : \BunU \to \PP_i
\end{equation}
 the section determined by $\sW_i \subset \sW_i \oplus \sV_i$.
Regarding this also as a Cartier divisor on $\PP_i$, we denote its cycle class by 
$
 [ \bm{0}_i  ]  \in \rH^2 (\PP_i)  ,
$
and denote by 
\[
\beta_i \in \rH^*(  \textstyle\prod \PP_i )
\]
the pullback of $ [ \bm{0}_i  ] $ along the projection to the $i^\mathrm{th}$ factor in $\prod \PP_i$.
For any  $u\in \rH^*( \BunU )$ and any subset 
$
S = \{ i_1,\ldots, i_a \} \subset \{1,\ldots d\},
$
 we abbreviate
\begin{equation}\label{eq:ubeta def}
u\beta_S := p^*(u) \cdot \beta_{i_1}\cdots \beta_{i_a} \in \rH^* (  \textstyle\prod \PP_i ).
\end{equation}

Fix a class $\mf{c} \in \rH^2( \HkU^1 )$, and denote by  
$
\mf{c}_{\PP} \in \rH^2(  \textstyle\prod \Hk_{\PP_i}^1  ) 
$
its pullback. Exactly as in the definition \eqref{explicit cohomological endomorphism} of 
\[
\Gamma_{\mf{c}} : \rH^*( \BunU ) \to \rH^*( \BunU ), 
\]
there is an associated  endomorphism
\[
\Gamma_{ \mf{c}_{\PP} } : \rH^*( \textstyle\prod \PP_i ) \to \rH^*( \textstyle\prod \PP_i ) .
\]
The following Proposition, which we will prove in \S \ref{ssec:proof-prop-RHC}, relates these.

\begin{prop}\label{prop:reduced hecke comparison}
For all $u \in \rH^*( \BunU )$ and all $S \subset \{ 1, \ldots, d\}$, we have
\[
\Gamma_{ \mf{c}_{\PP} }
(  u  \beta_S )  
\in 
   \Gamma_{\mf{c}}(u) \beta_S   + \sum_{ T\subsetneq S} u_T\beta_T
\]
for some classes $u_T \in  \rH^*( \BunU )$.
\end{prop}

\begin{proof}[Proof that Proposition \ref{prop:reduced hecke comparison} implies Proposition \ref{prop: hecke compatible}]
Exactly as in the proof of Proposition \ref{prop: local Hitchin local system}, especially \eqref{P1 cohomology}, we have isomorphisms
\begin{align*}
\rR p_{i!}(\Ql)   &  \cong \rR^0 p_{i!}(\Ql) \oplus \rR^2 p_{i!}(\Ql) [-2]   \\
& \cong  \Q_{\ell,\BunU} \oplus  \Q_{\ell,\BunU} \tw{-1} 
\end{align*}
in $D_c^b( \BunU)$ for all $1\le i \le d$.    Using  the K\"{u}nneth formula   
\[
\rR p_! (\Ql) \cong \bigotimes_{i=1}^d  \rR p_{i!}(\Ql) 
\]
and applying $\rR a_!$, we obtain a canonical isomorphism of graded vector spaces
\begin{equation}\label{eq: full P1 factorization}
\rR_{  \wt \pi_d  !} ( \Ql )  \iso \left( \bigotimes_{i=1}^d   \big( \Ql \oplus \Ql \tw{-1} \big) \right) \otimes
  \bigoplus_{j\ge 0} \rH^j(  \BunU) [ -j]  .
\end{equation}
The expression in parenthesis is canonically identified with the stalk of  $\sK_d^i\tw{-i} $ at our fixed geometric point \eqref{eq:mfa point}, by the proof of Proposition \ref{prop: local Hitchin local system}.

In  concrete terms, what \eqref{eq: full P1 factorization}  says is  that any cohomology class in 
\begin{equation}\label{eq:stalk is projective cohomology}
\rR^m_{  \wt \pi_d  !} ( \Ql )  \iso  \rH^m (  \textstyle\prod \PP_i )
 \end{equation}
can be expressed as a linear combination of  products of pullbacks from $\rH^*( \BunU )$ and the classes $\beta_1,\ldots, \beta_d$.  In other words, $ \rH^m (  \textstyle\prod \PP_i )$ is spanned by classes of the form $u \beta_S$.

The  filtration \eqref{eq: filtration} induces\footnote{Over the regular semisimple locus of $\cA_d$,  $\rR \pi_{d!} \Q_\ell$  is a direct sum of shifted local systems on a smooth space, and so its perverse cohomology sheaves and usual cohomology sheaves agree  up to a shift of indices.} a filtration on   \eqref{eq:stalk is projective cohomology}, which is the one for which 
\[
F_a  \rH^m (  \textstyle\prod \PP_i )  \subset \rH^m ( \textstyle \prod \PP_i )
\]
is the subspace spanned by  all $u \beta_S$ with $\#S \le a$.

By construction, the  endomorphism  $\Gamma_{ \mf{c}_{\cM}}$ on the left hand side of \eqref{eq:stalk is projective cohomology} agrees with the $\Gamma_{\mf{c}_\PP}$ defined on the right hand side, and combining all of this reduces Proposition \ref{prop: hecke compatible} to Proposition \ref{prop:reduced hecke comparison}.
\end{proof}


\subsection{Proof of Proposition \ref{prop:reduced hecke comparison}}
\label{ssec:proof-prop-RHC}


For a smooth proper $\ol{k} $-stack $Y$ of pure  dimension $N$, we have the trace map 
\begin{equation}\label{eq:trace integral}
\int_{Y} \co   \rH^{2N}( Y)    \rightarrow \Ql
\end{equation}
 of Poincar\'e duality (recall we have fixed $\Ql(1) \iso \Ql$, to suppress Tate twists).
 We extend this to a linear functional on the full cohomology ring of $Y$, vanishing outside of top degree.

\subsubsection{Geometry of Hecke stacks}

For each $1\le i \le d $,  we obtain from \eqref{untangled hecke fiber} a diagram 
\begin{equation}\label{untangled hecke fiber factor}
\begin{tikzcd}
{   \PP_i  }     \ar[d , "  p_i " '  ]  &  {  \Hk^1_{ \PP_i}   }   \ar[l , "h_0" ' ] \ar[r , " h_1" ]  \ar[d]   &   {  \PP_i   }   \ar[d, "   p_i " ]   \\
 { \BunU}      & {  \HkU^1 }      \ar[l ] \ar[r  ]   &    { \BunU}   . 
\end{tikzcd}
\end{equation}

To make this a bit more explicit,  fix a $\ol{k}$-valued  point $(x',\cL)$ of 
\[
\HkU^1 \iso X' \times \BunU,
\]
 as in  Remark \ref{rem: HkU}.
Using the $\cE$ fixed in  \eqref{eq:mfa point}, we  first  use $\cL$ to construct the complementary line bundle $\cD$ of  \eqref{eq: cline def}, and then the Hermitian torsion sheaves $Q$ and $P$ on $X'_{ \ol{k}}$ from \eqref{eq: PQpair}.  
Set $P_0=P$, and then use $x'$ to form the diagram of torsion coherent sheaves 
\begin{equation}\label{eq:ComparePs}
 \begin{tikzcd}
 { P_0 }   &   {  P_0 ( -x' )  =  P_{1/2} = P_1( - \sigma(x') )   } \ar[r]  \ar[l]  & P_1  ,
   \end{tikzcd}
\end{equation}
 exactly as in \eqref{Ptwists}, all having the same multiplicity-free support divisor \eqref{eq:CompareSupp}.
To give a  $\ol{k}$-valued point of $\Hk^1_{ \PP_i}$ above $(x',\cL)$ is to give 
inclusions 
\begin{equation}\label{eq:CompareHeckeP}
\begin{tikzcd}
L_0  |_{x_i'}     \ar[d, hook] & L_{1/2}   |_{x_i'}    \ar[l] \ar[d, hook] \ar[r] & L_1 |_{x_i'}    \ar[d, hook]  \\
( Q \oplus P_0)|_{x_i'}   &  (  Q \oplus P_{1/2}  ) |_{x_i'}   \ar[r] \ar[l] &  (   Q \oplus P_1 )  |_{x_i'}  
\end{tikzcd}
\end{equation}
of rank one local direct summands as in  \eqref{heckeP diagram}.

\subsubsection{Pullback of the zero sections}

We want to understand the pullback of the cycle class $[\bm{0}_i]$ to $\Hk_{\PP_i}^1$ along the arrows $h_0$ and $h_1$ in \eqref{untangled hecke fiber factor}.  Regarding points of $\HkU^1 \iso X' \times \BunU$ as pairs 
  $(x',\cL)$, and recalling that we have fixed the points in  \eqref{eq:CompareSupp}, there are  incidence divisors 
\begin{equation}\label{eq:IncidenceComp}
I_i , I_i^\star \hookrightarrow \HkU^1 
\end{equation}
  defined by  $x'=x_i'$ and  $x'=\sigma(x'_i)$, respectively.

\begin{lemma}\label{lem:heartsection}
There is a section 
 \[
\bm{0}_i^\heartsuit :  \HkU^1  \to   \Hk^1_{ \PP_i} 
  \]
 whose associated divisor satisfies
 \[
h_0^{-1}( \bm{0}_i )  =  \bm{0}_i^\heartsuit  + D_i 
\qquad \mbox{and} \qquad 
h_1^{-1}( \bm{0}_i )  =  \bm{0}_i^\heartsuit  + D^\star_i 
 \]
 for smooth divisors  $D_i$ and $D_i^\star$  on $\Hk^1_{ \PP_i}$ supported above $I_i$ and $I_i^\star$, respectively.
 Moreover, the preimages of $D_i$ and $D_i^\star$ under $\prod \Hk^1_{\PP_j} \to \Hk^1_{\PP_i}$ are again smooth divisors.
\end{lemma} 

\begin{proof}
This is similar to the analysis of the diagram  \eqref{eq: local proj hecke} carried out in the proof of Lemma \ref{lem:incidence fiber}.  
Using the  above description of \eqref{untangled hecke fiber factor},  the divisor  $h_0^{-1} (\bm{0}_i)$ is cut out by the condition $L_0|_{x_i'} =Q$ on the diagram \eqref{eq:CompareHeckeP},  while  $h_1^{-1}(\bm{0}_i)$ is cut out by the condition $L_1|_{x_i'} =Q$.   

The  desired section $\bm{0}_i^\heartsuit$ is the one defined by the diagram
\[
L_0|_{x_i'} = L_{1/2}|_{x_i'}=L_1|_{x_i'}  =Q|_{x_i'} .
\]
The divisor $D_i$ is defined by the conditions $x'= x'_i$ and $L_0|_{x_i'} =Q|_{x_i'}$, while   $D_i^\star$ is defined by 
$x'= \sigma(x'_i)$ and $L_1|_{x_i'} =Q|_{x_i'}$.

The incidence divisor $I_i$ is smooth because it is carried isomorphically to the smooth $k$-stack $\BunU$ by   projection to the second factor in $\HkU^1 \iso X'  \times \BunU$.
As $D_i$ has the structure of a $\PP^1$-bundle over $I_i$, it is  smooth.

The preimage  of $D_i$ under $\prod \Hk^1_{\PP_j} \to \PP_i$ is 
\begin{equation}\label{D_i pullback}
D_i \times \prod_{j\neq i}  \left( \Hk_{\PP_j}^1|_{I_i} \right),
\end{equation}
where all fiber products are taken over $I_i$.   From the above description of the diagram \eqref{untangled hecke fiber factor}, it is easy to see that  $\Hk_{\PP_i}^1 \to \HkU$ is a   $\PP^1$-bundle away  from $I_i \cup I_i^\star$, and it follows that for $j\neq i$ the fiber $\Hk_{\PP_j}^1|_{I_i}$ is a $\PP^1$-bundle over $I_i$.  We deduce that \eqref{D_i pullback} is a product of $\PP^1$-bundles over $I_i$, so is smooth.

Similar arguments prove the smoothness claims concerning $D_i^\star$.
\end{proof}

\subsubsection{Smoothness of the Hecke stack}

We will need to integrate, in the sense of \eqref{eq:trace integral}, cohomology classes on the stack $\prod \Hk^1_{\PP_i}$ from \eqref{untangled hecke fiber}.
In order to make sense of this we need to know that this latter stack is smooth.
This  is a consequence of the smoothness of $\BunU$ and the following Lemma.

\begin{lemma}
The composition 
\begin{equation}\label{eq:smooth hecke}
\prod \Hk^1_{\PP_i}  \map{\eqref{untangled hecke fiber} } \HkU^1 \iso X' \times \BunU \map{\mathrm{pr}} \BunU
\end{equation}
is smooth, where the isomorphism is that of Remark \ref{rem: HkU}.
\end{lemma}

\begin{proof}
First we check  that the fibers of \eqref{eq:smooth hecke} are smooth, using the  above description of \eqref{untangled hecke fiber factor}.
A geometric point $\cL \in \BunU( \ol{k} )$ determines torsion coherent sheaves $Q$ and $P$ on $X'_{\ol{k}}$, both with multiplicity-free support divisor \eqref{eq:CompareSupp}. 
 To give a $\ol{k}$-point of the fiber  $(\prod \Hk^1_{\PP_i})_{\cL}$  is to give a point $x' \in X'(\ol{k})$,  together with a $d$-tuple of diagrams \eqref{eq:CompareHeckeP}, indexed by $1\le i \le d$.  In particular,   there is a natural morphism 
\begin{equation}\label{blowup}
(\prod \Hk^1_{\PP_i})_{\cL} \to X' \times \prod   \PP(   (Q \oplus P_0)|_{ x'_i}   ) .
\end{equation}

For each index $1\le r \le d$, let $U_r \subset X'$  be the open subscheme of points $x'$ that are disjoint from all points in the support divisor \eqref{eq:CompareSupp}, except possibly $x_r'$.   
Let $U^\star_r  = \sigma(U_r) \subset X'$ be the open subscheme of points  $x'$ that are disjoint from all points in the support divisor \eqref{eq:CompareSupp}, except possibly $\sigma(x_r')$.    Denote by 
$V_r , V_r^\star \subset (\prod \Hk^1_{\PP_i})_{\cL}$ their preimages under the map to $X'$ in \eqref{blowup}.
As the $U_r$ and $U_r^\star$ cover $X'$, in order to prove the smoothness of $(\prod \Hk^1_{\PP_i})_{\cL}$  it suffices to prove the smoothness of each  $V_r$ and $V_r^\star$.

For ease of notation, let us explain why  $V_1$ is smooth, the other cases being entirely similar.   
At a point of $V_1$ the bottom horizontal arrows in \eqref{eq:CompareHeckeP} are isomorphisms for all $i>1$, because $x'$ is distinct from both $x'_i$ and $\sigma(x_i')$.  
Hence each  diagram indexed by $i>1$ is completely determined by the line $L_0|_{x'_i} \subset  (Q \oplus P_0)|_{ x'_i} $.  
The diagram \eqref{eq:CompareHeckeP}  indexed by $i=1$ has the form
\[
\begin{tikzcd}
L_0  |_{x_1'}     \ar[d, hook] & L_{1/2}   |_{x_1'}    \ar[l] \ar[d, hook] \ar[r,equal] & L_1 |_{x_1'}    \ar[d, hook]  \\
( Q \oplus P_0)|_{x_1'}   &  (  Q \oplus P_0(-x'  ) ) |_{x_1'}   \ar[r, equal ] \ar[l] &  (   Q \oplus P_0(\sigma(x')-x') )  |_{x_1'}  ,
\end{tikzcd}
\]
and such diagrams are parametrized by the blowup of $X' \times  \PP(   (Q \oplus P_0)|_{ x'_1}   )$ at the ordered pair consisting of the point  
$x'_1 \in X'$ and the line $  Q|_{x_1'} \subset  (Q \oplus P_0)|_{ x'_1}$.
The conclusion is that the restriction of \eqref{blowup}  to a morphism 
\[
V_1\to U_1 \times \prod   \PP(   (Q \oplus P_0)|_{ x'_i}   ) ,
\]
 identifies $V_1$ with the blowup of the codomain along the smooth codimension two subscheme
\[
\{ x_1' \} \times \big\{   Q|_{x_1'} \subset  (Q \oplus P_0)|_{ x'_1}   \big\} \times \prod_{i>1}  \PP(   (Q \oplus P_0)|_{ x'_i}   ),
\]
and therefore $V_1$ is smooth.

Having proved that each fiber of \eqref{eq:smooth hecke} is smooth, it suffices to prove that \eqref{eq:smooth hecke}  has the structure of a fiber bundle.  This amounts to understanding how the Hermitian torsion sheaves $Q$ and $P$ vary with  $\cL \in \BunU( \ol{k} )$.  Certainly $Q$ does not vary at all, as it depends only on the fixed  \eqref{eq:mfa point}.  
By construction, $P$ depends only on the fibers  of $\cL$ at the points of  \eqref{eq:CompareSupp}.  
The fiber  $\cL|_{x'_i}$ is the restriction of the universal line bundle  to 
\[
\BunU \iso \BunU \times \{ x'_i\} \subset \BunU \times X',
\]
so is locally constant as $\cL$ varies.  The fiber  $\cL|_{\sigma(x_i')}$ is locally constant for the same reason.
Because $Q$ and $P$ are locally constant as $\cL$ varies, so are the fibers of  \eqref{eq:smooth hecke}.
\end{proof}

\subsubsection{Comparison of integrals}
 
For a subset $T \subset \{ 1,\ldots, d\}$, let $T^c$ denote its complement.

\begin{lemma} \label{lem:first hecke integral}
Fix  $u , v  \in \rH^*( \BunU )$ and subsets $S , T \subset \{ 1,\ldots, d\}$. If $S  = T $ then 
\[
\int_{ \prod \PP_i } u  \beta_{S } \cdot v \beta_{T^c} = \int_{\BunU} u  \cdot v .
\]
If $S \cup T^c \subsetneq \{ 1,\ldots, d\}$,  the integral on the left is $0$.
\end{lemma}

\begin{proof}
Case (1):  If   $S = T$ then 
\[
\int_{ \prod \PP_i } u  \beta_{S } \cdot v \beta_{T^c} =
\int_{ \prod \PP_i } p^*(u)  \cdot p^*(v)  \cdot \beta_1 \cdots \beta_d
= \int_{\BunU} u  \cdot v.
\]
The first equality is immediate from the definition \eqref{eq:ubeta def}.
The second  follows from the observation that  $\beta_1 \cdots \beta_d \in \rH^{2d}( \prod \PP_i )$ is the cycle class of a section to $p : \prod \PP_i \to \BunU$, namely the product    the zero sections \eqref{eq:comp0}.

Case (2): If $S \cup T^c \subsetneq \{1,\ldots, d\}$,  then fix an index $j \not\in S\cup T^c$.  
The  cohomology class  $u  \beta_{S } \cdot v \beta_{T^c}$  arises as a pullback along the morphism 
\[
\prod \PP_i \to \prod_{ i \neq j} \PP_i .
\]
The target of this morphism has lower dimension than the source, and so if $u  \beta_{S } \cdot v \beta_{T^c}$ is in the top degree cohomology of $\prod \PP_i$, it must vanish.  Hence the integral on the left vanishes.
\end{proof}

\begin{lemma} \label{lem:second hecke integral}
Fix  $u,v \in \rH^*( \BunU )$ and subsets $S,T \subset \{ 1,\ldots, d\}$.  If $S  = T$ then 
\[
\int_{ \prod \PP_i }  \Gamma_{ \mf{c}_{\PP} } ( u  \beta_{S} )  \cdot v  \beta_{T^c}
 = \int_{\BunU} \Gamma_{\mf{c}} ( u  )  \cdot  v .
\]
If $S \cup T^c \subsetneq \{ 1,\ldots, d\}$,  the integral on the left is $0$.
\end{lemma}

\begin{proof}
Recall that 
$
\Gamma_{ \mf{c}_{\PP} } ( u  \beta_{S} )  = h_{1!} ( \mf{c}_{\PP}  \cdot  h_0^*(u  \beta_{S} ) ),
$
and so 
\begin{equation} \label{eq:GammaCPIntegral} 
\int_{ \prod \PP_i }  \Gamma_{ \mf{c}_{\PP} } ( u  \beta_{S} )  \cdot v  \beta_{T^c} = 
\int_{ \prod \Hk_{\PP_i}^1 } \mf{c}_{\PP}  \cdot  h_0^*( u  \beta_{S} ) \cdot h_1^* (v \beta_{T^c} ) .
\end{equation}

If  $S \cup T^c \subsetneq \{1,\ldots, d\}$ then the  integral on the right vanishes by the same reasoning as in Case (2) of the proof of Lemma \ref{lem:first hecke integral}.  
We may therefore  assume that  $S=T$, and that the class $u\cdot v$ has top degree on $\BunU$.

Denote by $[  \bm{0}_i^\heartsuit  ] \in \rH^2(  \Hk^1_{ \PP_i}  )$ the cycle class of the divisor   from Lemma \ref{lem:heartsection}, and by 
\[
\beta_i^\heartsuit \in \rH^2( \textstyle{\prod} \Hk^1_{\PP_i} )
\]
its   pullback.  By the same reasoning as Case (1) of the proof of Lemma \ref{lem:first hecke integral}, we have 
\begin{align}
\int_{ \prod \Hk_{\PP_i}^1 } \mf{c}_{\PP}  \cdot  h_0^*( p^* u)  \cdot h_1^*( p^* v)  \cdot   \beta_1^\heartsuit \cdots \beta_d^\heartsuit
 &  =\int_{\HkU^1} \mf{c} \cdot h_0^*( u) \cdot h_1^*(v) \nonumber \\
 & = \int_{\BunU} \Gamma_{\mf{c}} ( u  )  \cdot  v   .  \label{GammaHeartIntegral}
\end{align}
(By abuse of notation, the maps $h_0$ and $h_1$ in the integral over $\HkU^1$  are the unlabeled arrows in the middle row of \eqref{untangled hecke fiber}.)

Denoting  by $q : \prod \Hk_{\PP_i}^1 \to \rH^*(\HkU^1)$ the morphism in \eqref{untangled hecke fiber}, we may write
\[
\mf{c}_{\PP}  \cdot  h_0^*( p^* u)  \cdot h_1^*( p^* v) = q^* (\alpha)
\]
for some top degree class $\alpha \in \rH^*( \HkU^1 )$.
Comparing  \eqref{eq:GammaCPIntegral} with  \eqref{GammaHeartIntegral}, we see that the proof is complete once we show that 
\begin{equation}\label{heart to noheart}
\int_{ \prod \Hk_{\PP_i}^1 } q^*(\alpha) \cdot  h_0^*(   \beta_{S} ) \cdot h_1^* ( \beta_{S^c} )
=
\int_{ \prod \Hk_{\PP_i} } q^*(\alpha)  \cdot   \beta_1^\heartsuit \cdots \beta_d^\heartsuit .
\end{equation}

 Denote by 
 \[
 E = E_1 \cup \cdots \cup E_d \cup  E_1^\star \cup \cdots \cup E_d^\star \subset \prod \Hk_{\PP_i}^1
 \]
 the union of the preimages of the divisors    $D_1,\ldots, D_d$ and $D_1^\star,\ldots, D_d^\star$ from Lemma \ref{lem:heartsection} 
(the union is disjoint, because these divisors lie above the incidence divisors \eqref{eq:IncidenceComp}, which are pairwise disjoint),  and denote by 
\[
 f : E  \inj  \prod \Hk_{\PP_i}^1
\]
the inclusion.  It follows from Lemma \ref{lem:heartsection} that $E$ is smooth, and that
\[
h_0^*( \beta_i )  =   \beta_i^\heartsuit +   f_! [E_i] 
\qquad \mbox{and} \qquad
h_1^*( \beta_i ) =  \beta_i^\heartsuit +  f_! [E^\star_i]  
\]
where $[E_i]  ,[E^\star_i]   \in \rH^0(E)$ are the fundamental classes of $E_i$ and $E_i^\star$.
Repeatedly applying the relation $(f_! x) \cdot y= f_!(x \cdot f^* y)$, 
we find that we may write 
\[
h_0^*(   \beta_{S} ) \cdot h_1^* ( \beta_{S^c} )  = \beta_1^\heartsuit \cdots \beta_d^\heartsuit   +  f_! \epsilon
\]
for some $\epsilon \in \rH^*(E)$.

To complete the proof of \eqref{heart to noheart}, note that 
  \[
 \int_{ \prod \Hk_{\PP_i}^1 } q^*(\alpha) \cdot  f_!\epsilon =  \int_{E}  f^* ( q^* (\alpha) ) \cdot \epsilon =0 ,
 \]
 because the composition  $q \circ f : E \to \HkU^1$ factors through the union of all incidence divisors \eqref{eq:IncidenceComp}, and the pullback of $\alpha$ to this union vanishes for degree reasons.
\end{proof}

\subsubsection{Completion of the proof}

Fix  $u \in \rH^*( \BunU )$ and a subset $S \subset \{ 1,\ldots, d\}$.

\begin{proof}[Proof of Proposition \ref{prop:reduced hecke comparison}]
The discussion surrounding  \eqref{eq:stalk is projective cohomology} allows us to  write 
\[
 \Gamma_{ \mf{c}_{\PP} } ( u \beta_{S} )  = \Gamma_{\mf{c}}(u ) \beta_{S}  + \sum_{  R \subset \{1,\ldots, d\}  } u_R \beta_R 
 \]
 for  some classes $u_R \in \rH^*( \BunU )$.  
 We want to show that $u_R= 0$ for any  $R$ that is not properly contained in $S$, so let us rewrite this as
 \begin{equation}\label{eq:Integral Rsplit}
  \Gamma_{ \mf{c}_{\PP} } ( u \beta_{S} )  = \Gamma_{\mf{c}}(u ) \beta_{S}  + 
   \sum_{  R \subsetneq S } u_R \beta_R
  +
  \sum_{  R \in \mathscr{R} } u_R \beta_R,
 \end{equation}
where
\[
  \mathscr{R}  := \{ R \subset \{1,\ldots, d\}  : u_R\neq 0 \mbox{ and  $R$ is not a proper subset of  $S$} \}.
\]
To show that $\mathscr{R}$ is empty, we argue by contradiction.
Fix a   $T\in \mathscr{R}$ of maximal cardinality.

As $T$ is not a proper subset of $S$,  either $S=T$ or $S \cup T^c \subsetneq \{ 1,\ldots, d\}$.
In either case Lemmas \ref{lem:first hecke integral} and \ref{lem:second hecke integral} imply 
\[
 \int_{ \prod \PP_i }  \Gamma_{ \mf{c}_{\PP} } ( u \beta_{S} )  \cdot v \beta_{T^c}= \int_{ \prod \PP_i }  \Gamma_{\mf{c}}(u ) \beta_{S}  \cdot v \beta_{T^c} 
\] 
for every   $v \in \rH^*( \BunU )$.  
Another application of  Lemma \ref{lem:first hecke integral} shows that 
\[
\sum_{  R \subsetneq S  }   \int_{ \prod \PP_i } u_R \beta_R   \cdot v \beta_{T^c} =0,
\]
because each $R$ in the sum satisfies $R \cup T^c \subsetneq \{ 1,\ldots d\}$.
Combining \eqref{eq:Integral Rsplit} with the last two equalities, we deduce that 
\begin{equation}\label{eq:BadR 1}
\sum_{  R \in \mathscr{R}  }   \int_{ \prod \PP_i } u_R \beta_R   \cdot v \beta_{T^c} =0 
\end{equation}
for all $v \in \rH^*( \BunU )$.

On the other hand,  it follows from Lemma \ref{lem:first hecke integral}  that
\begin{equation}\label{eq:BadR 2}
\sum_{  R \in \mathscr{R}  }   \int_{ \prod \PP_i } u_R \beta_R   \cdot v \beta_{T^c}  
=  \int_{ \prod \PP_i } u_T \beta_T   \cdot v \beta_{T^c} 
=\int_{\BunU} u_T  \cdot v.
\end{equation}
Indeed,  the terms in the sum  vanish for all $R\in \mathscr{R}$ with  $R \cup T^c  \subsetneq \{1,\dots, d\}$, and the only  $R\in \mathscr{R}$ for which $R \cup T^c = \{1,\dots, d\}$ is  $R=T$ (because the maximality of $T$ implies  $\# R + \# T^c \le d$, which forces the union to be disjoint).

Comparing \eqref{eq:BadR 1} with \eqref{eq:BadR 2} shows that 
\[
\int_{\BunU} u_T  \cdot v  = 0 
\]
for all choices of $v \in \rH^*( \BunU )$, which implies $u_T=0$.  Of course this contradicts $T \in \mathscr{R}$, completing the proof of Proposition \ref{prop:reduced hecke comparison}
 \end{proof}


\subsection{Proof of Theorem \ref{thm: universal bundle action} }
\label{ss:proofUBA}


Having completed the proof of Proposition \ref{prop: hecke compatible}, we now use it to prove Theorem \ref{thm: universal bundle action}.
Theorem \ref{thm: universal bundle action} is a statement about perverse sheaves on $\cA_{d , \ol{k}}$, so we  base change all $k$-stacks to $\ol{k}$, and omit the base change from the notation.

The proof  will make use of the following observation.
Consider the morphism 
\begin{equation}\label{eq:proof-universal-hitchin}
f_d : \Hk^1_{\cM_d} \to  \cA_d \times \HkU^1
\end{equation} 
from  \eqref{global hecke d}.  Any line bundle $\mf{M}$ on the target  determines an endomorphism 
\[
\Gamma_{c_1(   f_d^*\mf{M}  )} \in \End_{D^b_c(\cA_d)}( \rR \pi_{d!} \Ql )
\] 
by the construction \eqref{hitchin correspondence endo}.  If we fix a geometric point 
$\mf{a}  = (\cE ,a)  \in \cA_d(\ol{k})$, 
 the action of this endomorphism on the stalk $(\rR \pi_{d!} \Ql)_{\mf{a}}$ only depends on the restriction of $\mf{M}$ to the fiber 
\begin{equation}\label{eq:proof-universal-fiber}
\HkU^1 \iso  \{ \mf{a} \} \times \HkU^1 \subset  \cA_d \times \HkU^1 .
\end{equation}

\begin{proof}[Proof of Theorem \ref{thm: universal bundle action}]
Fix a geometric point $\mf{a}=(\cE,a)$ of  $\cA_d$ as above, and 
recall from \eqref{eq:HkUlegmap} the morphism 
$
p : \HkU^1 \to X' .
$
Applying Propositions \ref{prop:pullback bundle correspondence} and \ref{prop: hecke compatible} to 
\[
\mf{c} = c_1 ( \det p^* \cE ) \in \rH^2( \HkU^1)(1) 
\]
and its pullback  $\mf{c}_{\cM} \in \rH^2(  \Hk_{\cM_d}^1 )(1)$, we find that 
\[
\Gamma_{ \mf{c}_{\cM} }  \in \End_{ D_c^b( \cA_{d,\ol{k}}) }  (  \p \rR^{m}  \pi_{d!}(\Ql) )
\]
 respects the filtration \eqref{eq: filtration}, and its semisimplification acts  as
\[
\mathrm{Id} \otimes \Gamma_{\mf{c}} = 
\deg(\cE) \otimes w 
=(-d + (n-1)\deg_X(\omega_X)) \otimes w
\]
on the summand  $\sK_d^i(-i) \otimes \rH^j(\BunU)[\dim \cA_d]$.

An endomorphism of an $\IC$-extension of a local system  is determined by its action on stalks.  
Therefore,  because $\mathfrak{a}$ is arbitrary,  to prove Theorem \ref{thm: universal bundle action} it suffices to show that 
\[
\Gamma_{c_1(    \det \cE^{\univ}   )} \in \End_{ D_c^b( \cA_d)  } (\p\rR^m \pi_{d!} (\Ql))
\] 
induces the same endomorphism as $\Gamma_{ \mf{c}_{\cM} }$  on the stalk  at $\mathfrak{a}$.

Recalling Definition \ref{def:Euniv},  the universal vector bundle  $\cE^{\univ}$  is the pullback of another (universal)  vector bundle along the composition
\[
\Hk^1_{\cM_d} \map{f_d}  \cA_d \times \HkU^1 \map{ \mathrm{Id} \times p}  \cA_d \times X' .
\]
Pulling back the determinant of this bundle along only the second arrow, we obtain a line bundle $\mf{M}$ on $\cA_d \times \HkU^1$  satisfying  $\det \cE^{\univ} \iso f_d^*\mf{M}$.
For tautological reasons,  the restriction of $\mf{M}$   to the fiber \eqref{eq:proof-universal-fiber}  is    $\det p^*\cE$.   

Pulling back $\det \cE$ along the composition 
\[
\cA_d \times \HkU^1\to   \HkU^1   \map{p}  X' 
\]
therefore yields  another  line bundle $\mf{M'}$ whose restriction to the fiber \eqref{eq:proof-universal-fiber} agrees with the restriction of $\mf{M}$, and whose further pullback along \eqref{eq:proof-universal-hitchin} has first Chern class  $\mf{c}_{\cM}$.

 By the observation made at the beginning of \S \ref{ss:proofUBA}, the two endomorphisms
\[
\Gamma_{c_1(    \det \cE^{\univ}   )}  = \Gamma_{c_1(    f_d^*\mf{M}  )}  
\quad \mbox{and}\quad 
 \Gamma_{ \mf{c}_{\cM}} = \Gamma_{c_1(    f_d^*\mf{M}'  )}
\]
of $\p\rR^m \pi_{d!} (\Ql)$ induce the same endomorphism on the stalk at $\mf{a}$, completing the proof.
\end{proof}


\section{Hecke action of the tautological bundle}
\label{sec:taut-action}


In this section we prove Theorem \ref{thm: tautological action}. 
The proof follows the same strategy as \S \ref{sec:univ-action}, insofar as it uses Proposition \ref{prop: hecke compatible} to relate an endomorphism of $\p \rR^{m}  \pi_{d!}(\Ql)$ to a simpler endomorphism of 
$\rH^*( \BunU)$.  The difference is that in this case the simpler endomorphism requires a significant amount of effort to  compute.


\subsection{A more general setting}


We remind the reader  that $\eta : \A_F^\times \to \{ \pm 1\}$ is the unramified Hecke character determined by the \'etale double cover $\nu :X'\to X$.  Using the convention of \cite[\S 2.6]{FYZ}, we also regard 
\[
\eta : \Pic(X) \to \{\pm 1\}.
\]

In \S \ref{sec:taut-action} only, we consider  the following  generalization of    \eqref{BunU def}. 
Fix a line bundle
\[
\mf{N} \in \Pic(X).
\]
  Using the notation of \cite[\S 3.1.1]{FYZ2},  let 
\begin{equation}\label{new U(1) twist}
\BunUN := \Bun_{ \U(1) ,  \mf{N} \otimes \omega_X^{-1} }  
\end{equation}
be the moduli stack whose $R$-points, for any  $k$-algebra $R$, are line bundles $\cL$ on $X'_R$ equipped with a Hermitian isomorphism
$\cL \otimes   \sigma^*\cL \iso   \nu^*\mf{N}$.

The obvious generalization of \S \ref{ssec: Hecke BunU} gives a Hecke correspondence
\begin{equation}\label{eq: HkU N correspondence}
\BunUN \xleftarrow{h_0} \HkUN^1 \xrightarrow{h_1} \BunUN ,
\end{equation}
and everything in \S \ref{ssec:univ-action-Bun} holds verbatim with $\BunU$ replaced by $\BunUN$.
This includes the construction of the order two automorphism $w$ from \eqref{eq:w-defn}.

\begin{remark}
Taking $\mf{N} = \omega_X^{2-n}$ in \eqref{new U(1) twist}  recovers  the stack $\BunU$ of \eqref{BunU def}.
The extra generality obtained by allowing $\mf{N}$ to vary will never be used in the proofs of our main results.
We allow it here solely because Proposition \ref{prop:better bun} and Corollary \ref{cor: trace on U(1)} are of independent interest, so we formulate and prove them in this greater generality.
\end{remark}


\subsection{Cohomology of Prym varieties}\label{ssec:cohomology-prym}


Let $\Prym_{\mf{N}}$ be the fiber over $\mf{N}\in \Pic_X(k)$ of the norm map $\Nm \co \Pic_{X'} \rightarrow \Pic_X$ between Picard stacks.  
 Generalizing Remark \ref{rem: HkU}, there are canonical identifications
\begin{equation}\label{BunUtoPrym}
\begin{tikzcd}
{  \BunUN }  \ar[d, equal]  &   {   \HkUN^1    }    \ar[l , "h_0"  '] \ar[r , "h_1"]    \ar[d, equal]  &  {    \BunUN  }   \ar[d , equal]  \\ 
{   \PrymN    }   &  {  X' \times \PrymN  }   \ar[l , "h_0"  '] \ar[r , "h_1"]    &  {  \PrymN       } 
\end{tikzcd}
\end{equation}
defined as follows.  The outside equalities identify a point of $ \BunUN$, consisting of a line bundle $\cL \in \Pic(X')$ and a Hermitian isomorphism 
$
\cL \otimes     \sigma^*\cL  \iso  \nu^*\mf{N} , 
$
with the point of $\PrymN$ defined by $\cL$ with its induced  isomorphism $\mathrm{Nm}_{X'/X}( \cL) \iso \mathfrak{N}$.
The middle equality sends a point of $  \HkUN^1  $, consisting of a modification \eqref{eq: HkU mod},  to the point $(x' , \cL)$ of $X' \times \PrymN$.
The bottom horizontal arrows  are
\[
 \cL    \mapsfrom  (x' , \cL )  \mapsto \cL ( \sigma(x') - x' ) .
\]

For any choice of $\fN$, $\PrymN$ has two geometric connected components, which we denote by $\PrymN^+$ and $\PrymN^-$.
The labeling is arbitrary, except in the special case $\mf{N}=\cO_X$, in which case
\[
\Prym:= \Prym_{\cO_X}
\]
has a group structure, and  we always take $\Prym^+$ to be the neutral component.

As explained in \cite[Lemma 2.13]{FYZ}, the $k$-rationality of the components depends on $\eta(\mf{N})$. If $\eta(\mf{N}) = 1$, the two components $\PrymN^{\pm}$ are both defined over $k$; if $\eta(\mf{N}) = -1$, there is only one connected component over $k$.

The tensor product of line bundles makes each $\PrymN^\pm$ into a  $\Prym^+$-torsor.   
As explained in the discussion leading to \cite[(10.16)]{FYZ2}, any trivialization of this torsor over $\ol{k}$ determines a \emph{canonical} isomorphism
\begin{equation}\label{Prym torsors}
\rH^m (\PrymN^\pm) \iso \rH^m ( \Prym^+ )
\end{equation}
 independent of the choice of trivialization.
As a special case, we find a canonical   isomorphism 
\[
\rH^1( \Prym^- ) \iso \rH^1( \Prym^+ ).
\]

Abbreviating  
\[
V := \rH^1(X')^{\sigma=-1},
\]
the Abel-Jacobi map $\mathrm{AJ} : X' \to \Prym^-$,  normalized by 
$
\mathrm{AJ}(x' )=  \cO_{X'}( \sigma ( x' )  - x') ,
$
 determines an isomorphism
\[
 \rH^1( \Prym^-)  \map{\mathrm{AJ}^*} \rH^1(X')^{\sigma=-1} =V.
\]
Combining all of this, we obtain  canonical isomorphisms
\begin{align}
\rH^m (  \BunUN )  & \iso   \rH^m (\Prym^+_\mathfrak{N} ) \oplus \rH^m (\Prym^-_\mathfrak{N} ) \nonumber  \\ 
& \iso 
\rH^m (\Prym^+ ) \oplus \rH^m (\Prym^+)  \nonumber \\
& \iso 
 \left( {\bigwedge}^m V  \right) \oplus  \left( {\bigwedge}^m V  \right) . \label{BunUtoV}
\end{align}

Tracing through the construction of \eqref{BunUtoV}, one finds that  the automorphism 
\begin{equation}\label{w def}
w  : \rH^* (  \BunUN )  \to \rH^* (  \BunUN )  
\end{equation}
from \eqref{eq:w-defn}  acts as $(a,b) \mapsto (b,a)$ on the final expression in \eqref{BunUtoV}.

Let us also explicate the Frobenius action on $\rH^*(\PrymN)$. If $\eta(\mf{N}) = 1$, then the Frobenius preserves the direct sum decomposition in \ref{BunUtoV}, and acts on each via its natural action on $V$. If $\eta(\mf{N}) = -1$, the Frobenius exchanges the two summands. In fact, we have the following proposition:

\begin{prop}\label{prop: inert Prym Frobenius action}
    Let $\fN \in \Pic_X(k)$ satisfy $\eta(\fN) = -1$. Then for any $i \geq 0$, $\Tr(\Fr_{\PrymN}, \rH^i(\PrymN)) = 0$, and
    \[\Tr(\Fr_{\PrymN} \circ w, \rH^i(\PrymN)) = \Tr(\Fr_{\Prym}, \rH^i(\Prym))\]
\end{prop}

\begin{proof}
    It is proved in \cite[Lemma 2.13]{FYZ2} that under the assumption $\eta(\mf{N}) = -1$, $\Fr_{\PrymN}$ exchanges the two geometric connected components of $\PrymN$. This proves the first claim.

    For the second, as in the proof of \cite[Lemma 2.14]{FYZ2} it suffices to consider the case $\fN = \CO(x)$ for an inert closed point $x \in X$. Following the notation from {\it{loc. cit}}, let $x \times_k \kbar = \{x_1,\dots,x_d\}$ be the set of geometric points above $x$, and $x' \times_k \kbar = \{x_1',\dots,x_{2d}'\}$ be the set of geometric points above $x'$, where $x'$ is the unique closed point of $X'$ above $x$; moreover, we enumerate the $\{x_i\}$ and $\{x'_i\}$ in a manner compatible with the Frobenius and Galois actions as in {\it{loc. cit}}. In particular, $\CF = \CO_{X'}(x'_1 + \dots + x'_d)$ defines an element of $\Prym_{\fN}(\kbar)$.
    
    As $\Prym_{\fN}$ is a $\Prym$-torsor, there is an isomorphism of $\kbar$-schemes $\act_{\CF}: \Prym_{\kbar} \simeq \Prym_{\fN,\kbar}$ defined by acting on $\CF$. Both $\Prym_{\kbar}$ and $\Prym_{\fN,\kbar}$ carry Frobenius endomorphisms, however they are not intertwined under $\act_{\CF}$. Instead, for any $\CL \in \Prym_{\kbar}$, we have
    \[\Fr_{\Prym_{\fN}}(\CL \otimes \CF) = \Fr_{\Prym}(\CL) \otimes \Fr_{\Prym_{\fN}}(\CF) = \act_{\CF}(\Fr_{\Prym}(\CL) \otimes \Fr_{\Prym_{\fN}}(\CF) \otimes \CF^{-1}).\]
    Therefore under the isomorphism $\act_{\CF}$, $\Fr_{\Prym_{\fN}}$ corresponds to the endomorphism of $\Prym$ given by
    \[\CL \mapsto \Fr_{\Prym}(\CL) \otimes \Fr_{\Prym_{\fN}}(\CF) \otimes \CF^{-1}.\]
    It remains to prove that the endomorphism of $\Prym_{\kbar}$ given by translation by $\Fr_{\Prym_{\fN}}(\CF) \otimes \CF^{-1} \in \Prym_{\kbar}$ acts on $\rH^*(\Prym)$ as the automorphism $w$. It is proved in {\it{loc. cit}} that $\Fr_{\Prym_{\fN}}(\CF) \otimes \CF^{-1} \in \Prym^-_{\kbar}$, from which the claim follows, as translation by elements in $\Prym^-_{\kbar}$ acts on cohomology as $w$.
\end{proof}

 
 \subsection{Hecke action of tautological bundles on $\BunUN$}


Recall from \S \ref{ssec:univ-action-Bun} that for each $\mf c \in H^2(\HkUN^1)(1)$, there is a Frobenius-equivariant endomorphism 
\[
\Gamma_{\mf c} \in \End(\rH^*(\BunUN)).
\]
In Proposition \ref{prop:pullback bundle correspondence}, we calculated this action when $\mf c$ is the first Chern class of a line bundle from $X'$. 
In this subsection, we undertake the more difficult calculation where $\mf c$ is instead the first Chern class of the tautological bundle $\ell_{\UdN}$ from Definition \ref{def: U(1) tautological bundle}. This calculation can, with enough work, be extracted as a special case of the method of \cite{volume}.

\subsubsection{Relation to the Poincar\'e bundle} As a first step, consider the Poincar\'e line bundle $\mathscr{P}$ on $X' \times \PrymN$, whose pullback along an $R$-valued point $(x' , \cL)$ is the rank one projective $R$-module  $\cL|_{x'}$.

\begin{lemma}\label{lem: taut to poincare}
The isomorphism
$
\HkUN^1 \iso X' \times  \PrymN
$
from  \eqref{BunUtoPrym}  identifies 
\[
\ell_{\UdN}   \iso \mathscr{P}^{-1} \otimes   \pr_0^* (\omega_{X'}^{-1} \otimes  \nu^* \mathfrak{N}   ) ,
\]
where $\pr_0 : X' \times  \PrymN \to X'$ is the projection to the first factor. 
\end{lemma}

\begin{proof}
The pullback of $\ell_{\UdN}$ along an $R$-valued point $(x' , \cL) \in X' \times  \PrymN$ is the projective $R$-module
\[
\cL(\sigma(x') - x') |_{ \sigma(x')} \iso   \cL(\sigma(x')  ) |_{ \sigma(x')} 
 \iso ( \cL \otimes \omega_{X'}^{-1} )|_{\sigma(x')}  \iso ( \sigma^* \cL)|_{x'} \otimes   \omega_{X'}^{-1} |_{ x'} .
\]
Using the Hermitian isomorphism $\cL \otimes \sigma^* \cL \iso \nu^*\mf{N}$,  this last line bundle is isomorphic to
\[
 \cL^{-1} |_{x'} \otimes  (   \omega_{X'}^{-1}   \otimes \nu^*\mf{N}  )  |_{ x'},
\]
which is what we needed to prove.
\end{proof}

\subsubsection{The class $\beta$} Recall that we abbreviate 
$
V := \rH^1(X')^{\sigma=-1} .
$
Denote by 
\[
 \psi : V (1)  \otimes V  \to \Q_\ell
 \]
the pairing
\[
 V (1) \otimes V  \subset \rH^1(X' )  (1)   \otimes \rH^1(X')  \to  \rH^2(X ' ) (1)  \to  \Ql
 \]
of Poincar\'e duality. 
Define a cohomology class
 \[
 \beta = \sum_j v_j \otimes v^j \in V(1)  \otimes V ,
  \]
 where $\{ v_j\}$ is any  basis of $V(1)$,  and $\{v^j\}$ is the dual basis of $V$ with respect to $\psi$.  
 By the  discussion  leading to \eqref{BunUtoV}, we may regard $\beta$ as an element of 
  \begin{equation}\label{Vtensor}
 \beta \in V (1) \otimes V 
  \subset
 \rH^1( X')(1)  \otimes \rH^1( \Prym^\pm_\mathfrak{N}  ) 
 \subset \rH^2( X' \times  \Prym^\pm_\mathfrak{N}  )(1) . 
 \end{equation}

Consider the  diagram 
\begin{equation}\label{Prym component Hecke}
\begin{tikzcd}
{   \Prym^\pm_{\mf{N}}     }   &  {  X' \times \Prym^\pm_{\mf{N}}     }  \ar[l , "h_0" ' ] \ar[r , " h_1 "]    &  {  \Prym^\mp_{\mf{N}}       } 
\end{tikzcd}
\end{equation}
obtained by restricting \eqref{BunUtoPrym} to connected components.
 Exactly as in the discussion leading to \eqref{BunU Hecke endomorphism}, the  cohomology class $\beta$ in \eqref{Vtensor}
can   be viewed as a cohomological correspondence, and so defines (for any $m\ge 0$) a homomorphism 
\[
\rH^m( \PrymN^\pm) \map{\Gamma_\beta}    \rH^m(  \PrymN^\mp).
\]
 As we have identified both sides with $\bigwedge^m V$ in \eqref{BunUtoV}, we may view it as an \emph{endomorphism}
\begin{equation}\label{beta endomorphism}
{\bigwedge}^m V  \iso \rH^m( \PrymN^\pm) \map{\Gamma_\beta}    \rH^m(  \PrymN^\mp) \iso {\bigwedge}^m V  .
\end{equation}

 \begin{lemma}\label{lem:beta calculation}
The endomorphism \eqref{beta endomorphism} is multiplication by $m-\deg_X(\omega_X)$.
 \end{lemma}
 
 \begin{proof}
 We calculate the effect of the arrows in \eqref{Prym component Hecke}  on cohomology.
 
For any  $a \in \rH^1(  \PrymN^\pm   ) \iso V$,  it is clear that 
\[
h_0^*(a) = 1 \otimes a \in \rH^0( X') \otimes \rH^1(  \PrymN^\pm   ) \subset \rH^1( X' \times \PrymN^\pm ).
\]
It follows that we have a commutative diagram
\[
\begin{tikzcd}
{ \rH^m( \PrymN^\pm)   }  \ar[r, " h_0^*" ] &  { \rH^0(  X')  \otimes  \rH^m( \PrymN^\pm)    }  \ar[r , " \beta \cdot" ] & {  \rH^1(X')(1) \otimes \rH^{m+1}(   \PrymN^\pm)   }   \\
{  \bigwedge^m V  }  \ar[u ,equal]  \ar[rr] &    & {  V(1) \otimes \bigwedge^{m+1} V ,  }  \ar[ u , hook ]  
\end{tikzcd}
\]
in which the bottom horizontal arrow sends $a= a_1\wedge \cdots \wedge a_m$ to 
\begin{equation}\label{half taut Hecke}
 \beta \cdot h_0^*(a) = 
 \sum_j v_j \otimes (v^j \wedge a_1\wedge \cdots \wedge a_m) .
\end{equation}

 Now we turn to the analysis of $h_1$. 
The map $h_1$ factors  as 
\[
X' \times \PrymN^\pm \xrightarrow{ \AJ  \times \Id} \Prym^{-} \times \PrymN^{\pm} \to  \PrymN^{\mp}
\]
where the final  arrow   is the tensor product of line bundles on $X'$.  
Using the isomorphisms \eqref{Prym torsors}, the induced pullback
\[
\rH^1( \PrymN^{\mp}  )  \to \rH^1 (  \Prym^{-} \times \PrymN^{\pm}  ) 
\]
on cohomology is identified with the coproduct 
\[
\rH^1( \Prym^+ )  \to \rH^1 (  \Prym^{+} \times \Prym^{+}  ) 
\]
determined by the group law on $\Prym^+$.  
 For any abelian variety $B$, the coproduct  
\[
\rH^1(B) \to \rH^1(B\times B) =\left( \rH^0(B) \otimes \rH^1(B)   \right) \oplus \left( \rH^1(B) \otimes \rH^0(B)   \right)
\]
 is given by $b \mapsto 1\otimes b + b\otimes 1$. Tracing through our identifications, it follows that for any class
 $b \in \rH^1(  \PrymN^\mp  ) \iso V$ we have
\begin{align*}
h_1^*(b) = 1 \otimes b + b \otimes 1 & \in (\Q_\ell \otimes V)   \oplus ( V  \otimes \Q_\ell)  \\
& \iso 
\left( \rH^0(X') \otimes \rH^1( \PrymN^\pm)   \right)  
\oplus  
\left( \rH^1(X') \otimes \rH^0( \PrymN^\pm)   \right)  
\\
&  \iso  \rH^1(X' \times \PrymN^\pm ) .
\end{align*}

Abbreviate 
$
d : =  \dim V= \deg_X(\omega_X).
$
From the above description of  $h_1^*$ in cohomological degree $1$, it follows that  we have a commutative diagram
\[
\begin{tikzcd}
{   \rH^{d-m}(  \PrymN^\mp   )  }   \ar[dd ,equal ]
\ar[rr , "h_1^* "]    &   & 
{  \rH^{d-m}(  X' \times \PrymN^\pm   )  }   \ar[d , " \mathrm{pr}" ]    \\
{  }    & & {    
\rH^1(  X'    )\otimes \rH^{d -m-1}(  \PrymN^\pm   )  }   \\
{\bigwedge^{d-m} V  }  \ar[rr]  & & { V \otimes \bigwedge^{d-m-1} V  }  \ar[u, hook] 
\end{tikzcd}
\]
in which $\mathrm{pr}$ is projection to the Kunneth factor, and the bottom horizontal arrow sends 
$b =  b_1\wedge \cdots \wedge b_{d-m}$ to 
\[
h_1^*(b)= 
\sum_{  i=1}^{ d-m} (-1)^{ i -1 }   
b_i \otimes( b_1\wedge \cdots \wedge \widehat{b}_i \wedge \cdots \wedge b_{d-m} ) .
\]

The above calculation shows that \eqref{half taut Hecke} satisfies
\begin{align*}
 &    \beta \cdot h_0^*(a)   \cdot h_1^*(b)  \\
 & =
\sum_{  i=1}^{ d-m}   \sum_j (-1)^{m+ i   }   \psi ( v_j , b_i)   (v^j \wedge a_1\wedge \cdots \wedge a_m  \wedge 
 b_1\wedge \cdots \wedge \widehat{b}_i \wedge \cdots \wedge b_{d-m} )  \\
 & =
\sum_{  i=1}^{ d-m}    (-1)^{m+ i   }   b_i \wedge a_1\wedge \cdots \wedge a_m  \wedge 
 b_1\wedge \cdots \wedge \widehat{b}_i \wedge \cdots \wedge b_{d-m}   \\
  & =
 - \sum_{  i=1}^{ d-m}    a_1\wedge \cdots \wedge a_m  \wedge 
 b_1\wedge \cdots  \wedge b_{d-m}   \\
 & = (m-d)  ( a \cdot b) ,
\end{align*}
where we have used  $ \sum_j   \psi(  v_j ,  b _i)    v^j  = b_i$, and the above equalities are understood in 
\[
\rH^{d+2}( X' \times  \PrymN^\pm)  (1)  \iso 
\rH^2(X')(1) \otimes \rH^d( \PrymN^\pm)   \iso   {\bigwedge}^d V .
\]
By definition of $h_{1*}=h_{1!}$, we must therefore have
\[
h_{1! } (   \beta \cdot h_0^*(a)  ) \cdot b  = (m-d) ( a \cdot b) ,
\]
as elements of $\rH^{d}(   \PrymN^\pm)     \iso   {\bigwedge}^d V$, and this must hold for all choices of $b$.  This can only happen if    
\[
h_{1!} (   \beta \cdot h_0^*(a)  ) = (m-d)  a.
\]
Comparing with \eqref{explicit cohomological endomorphism}, this is what we needed to prove.
 \end{proof}

\subsubsection{The action on $\BunUN$} Recall the tautological bundle $\ell_{\UdN}$ from Definition \ref{def: U(1) tautological bundle}.
By the construction \eqref{BunU Hecke endomorphism}, its  first Chern class 
 \[
  c_1(  \ell_{\UdN} )  \in  H^2(\HkUN^1 )  (1)
  \]
induces an endomorphism 
\begin{equation}\label{taut gamma}
 \Gamma_{ c_1(  \ell_{\UdN} ) }  : \rH^*( \BunUN )  \to \rH^*( \BunUN ) .
\end{equation}

\begin{prop}\label{prop:better bun}
 For any $m\ge 0$, we have the equality 
 \[
 \Gamma_{ c_1(  \ell_{\UdN} ) }  =  (   \deg_X(\mathfrak{N})  -2m  ) \cdot w   
 \]
 of endomorphisms of $ \rH^m( \BunUN )$. 
 Here $w$ is the  automorphism \eqref{w def}, and $\mf{N} \in \Pic(X)$ is the line bundle fixed in \eqref{new U(1) twist}.
\end{prop}

\begin{proof}
According to \cite[Lemma 10.5]{FYZ2},  we have the equality 
\[
  c_1 \left(  \mathscr{P} |_{  X' \times  \PrymN^\pm } \right) =
 2\beta +   \deg_X( \mf{N}) \cdot   \gamma 
\]
 in $H^2( X' \times  \PrymN^\pm )(1)$,   where $\gamma$  is  the cycle class of the divisor 
 \[
 \{ x'\} \times \PrymN^\pm \subset X' \times \Prym^\pm_\mathfrak{N}
 \]
  determined by any  closed point $x' \in X'_{ \ol{k}} $.
Using  Lemma \ref{lem:beta calculation}, and arguing as in the proof of Proposition \ref{prop:pullback bundle correspondence} to compute the action of $\gamma$,  the induced homomorphism 
\[
 {\bigwedge}^m V  \iso \rH^m( \PrymN^\pm) \to   \rH^m(  \PrymN^\mp) \iso {\bigwedge}^m V 
\]
 is multiplication by $2m-2\deg_X( \omega_X )  + \deg_X(\mf{N}) $.
 Using the identification
 \[
 \rH^m ( \BunUN ) \iso \left( {\bigwedge}^m V  \right) \oplus  \left( {\bigwedge}^m V  \right)
 \]
from \eqref{BunUtoV}, we find that  the cohomological correspondence 
\[
c_1( \mathscr{P}) \in     H^2( X' \times  \PrymN  )(1) \iso  H^2(  \HkUN^1    )(1)
\]
 induces the endomorphism  
\begin{equation}\label{eq: BunU end 1}
(2m-2\deg_X( \omega_X )  + \deg_X(\mf{N}) ) \cdot w :   \rH^m ( \BunUN )  \to \rH^m ( \BunUN ).
\end{equation}

By Lemma \ref{lem: taut to poincare} we have the equality of first Chern classes
\[
c_1(\ell_{\UdN})    = - c_1( \mathscr{P})  +     c_1(  p^* (\omega_{X'}^{-1} \otimes  \nu^* \mf{N}  ) )  .
\]
Converting from cohomological correspondences to endomorphisms as above, the first term $- c_1( \mathscr{P})$ induces the negative of the endomorphism \eqref{eq: BunU end 1}, while from Proposition \ref{prop:pullback bundle correspondence} we see that the second term $c_1(  p^* (\omega_{X'}^{-1} \otimes  \nu^* \mf{N}  ) ) $ induces the endomorphism 
\[
\deg_{X'} ( \omega_{X'}^{-1} \otimes  \nu^* \mf{N}  ) \cdot w  
= (-2 \deg_X (\omega_X) + 2 \deg_X (\mf{N} ) ) \cdot w .
\]
of $\rH^m ( \BunUN )$.
Adding these together completes the proof of the Proposition. 
\end{proof}

\begin{cor}\label{cor: trace on U(1)}
The endomorphism $ \Gamma_{ c_1(\ell_{\UdN} )  } $ from  \eqref{taut gamma} satisfies
\begin{align}\label{eq: dirichlet L-function as trace}
\sum_{m \ge 0} (-1)^m  \cdot   \mathrm{Tr} \big(   \Gamma_{c_1(  \ell_{\UdN} ) }  ^r  \circ \mathrm{Fr} , \rH^m ( \BunUN )  \big)
  =   \frac{ 2 }{ (\log q)^r } \cdot \frac{d^r}{ds^r} 
  \left[  
  q^{s \deg_X(\mathfrak{N}) }    L( 2s,\eta )  
   \right]_{s=0} 
\end{align}
if $(-1)^r = \eta(\mf{N})$. Otherwise, the sum on the left hand side is $0$.
\end{cor}

\begin{proof}
If $(-1)^r \neq \eta(\fN)$, it is clear from the characterization of $w$ in   \eqref{w def} and Proposition \ref{prop: inert Prym Frobenius action} that 
\[
 \mathrm{Tr} (   w^r  \circ \mathrm{Fr} , \rH^m ( \BunUN ) ) =0  ,
\]
and so the vanishing of the left hand side  follows from Proposition \ref{prop:better bun}.

Now  assume that $(-1)^r = \eta(\mf{N})$ holds, so that Propositions \ref{prop:better bun} and \ref{prop: inert Prym Frobenius action} imply
\begin{align*}
& \sum_{m \ge 0} (-1)^m  \cdot   \mathrm{Tr} \big(   \Gamma_{ c_1(  \ell_{\UdN} ) }  ^r  \circ \mathrm{Fr} , \rH^m ( \BunUN )  \big)  \\
& =
   \sum_{m\ge 0} (-1)^m   \cdot   (  \deg_X(\mathfrak{N})   -2m    )^r    \cdot      \mathrm{Tr}(  \mathrm{Fr} ,  \rH^m ( \Prym)  )   .
\end{align*}
On the other hand, writing $\alpha_1,\ldots, \alpha_{2g-2}$ for the eigenvalues of $\mathrm{Fr}$ acting on  
$
V = \rH^1(X')^{\sigma=-1}
$, we have  
 \begin{align*}
 L(s,\eta)  & = \prod_{i=1}^{2g-2}(1 - \alpha_i q^{-s} )   \\
&  = \sum_{m \ge 0} (-1)^m \cdot  \mathrm{Tr} \big( \mathrm{Fr}, {\bigwedge}^m V \big) \cdot q^{-ms} \\
& =  2^{-1}  \sum_{m  \ge 0} (-1)^m \cdot  \mathrm{Tr}( \mathrm{Fr},  \rH^m( \Prym )   ) \cdot q^{- ms },
 \end{align*}
where we have used \eqref{BunUtoV} for the final equality. The claim follows. 
\end{proof}

\begin{remark}
Upon applying the Lefschetz trace formula of \cite[Proposition 11.8]{FYZ}, Corollary \ref{cor: trace on U(1)} recovers \cite[Theorem 10.2]{FYZ}. The original proof of \emph{loc. cit.} was more difficult and complicated than the one given here. Moreover, we shall see that the present approach is more flexible, as it can be incorporated into the calculation of the degrees of corank one special cycles.
\end{remark}


\subsection{Proof of Theorem \ref{thm: tautological action}}


We now take  $\mf{N} = \omega_X^{2-n}$, so that $\BunUN=\BunU$.
Abbreviate 
\[
\mf{c} = c_1 ( \ell_{\Ud} ) \in H^2(  \HkU )(1)  
\]
for the first Chern class of the tautological bundle of Definition \ref{def: U(1) tautological bundle}.
As in  \S \ref{ss:compatibility statement},  we pull it back to a class
\[
\mf{c}_{\cM} = c_1 ( \ell_{\Ud} ) \in H^2(  \Hk_{\cM_d}^1  )(1)  .
\]
This latter class is, by abuse of notation,  the first Chern class of the line bundle $\ell_{\Ud}$  on $\Hk_{\cM_d}^1$ from Definition \ref{def:M taut}.

\begin{proof}[Proof of Theorem \ref{thm: tautological action}]
As  in \S \ref{ss:compatibility statement},  the cohomology class $\mf{c}$ determines an endomorphism $\Gamma_{\mf{c}}$ of $\rH^j( \BunU)$.  
This is the same as the endomorphism  in Proposition \ref{prop:better bun}, so is given by the explicit formula
\[
\Gamma_{\mf{c}}=  (  -2j +  (2-n) \deg_X(\omega_X)   ) \cdot w.
\]

Again as  in \S \ref{ss:compatibility statement}, the cohomology class $\mf{c}_{\cM}$ determines an  endomorphism $\Gamma_{\mf{c}_{\cM}}$ of $\p \rR^{m}  \pi_{d!}(\Ql)$.  This is the same as the endomorphism  in the statement of Theorem \ref{thm: tautological action}.

Proposition \ref{prop: hecke compatible} tells us that  $\Gamma_{\mf{c}_{\cM}}$  respects the filtration \eqref{eq: filtration}, and its semisimplification acts by
\[
\text{$\mathrm{Id}  \otimes \Gamma_{\mf{c}} $ on the summand $\sK_d^i(-i) \otimes \rH^j(\BunU)[\dim \cA_d]$,}
\]
so we are done.
\end{proof}

 
\section{Cycles on moduli spaces of shtukas}\label{sec: geometric side}


In this section we recall from \cite{FYZ} the construction of corank one special cycles on the moduli space of Hermitian shtukas.
We then prove the fundamental Theorem \ref{thm: geometric side}, which expresses the intersection multiplicities of these cycles against  tautological bundles in terms of sheaves on the Hitchin base $\cA$.


\subsection{Compositions of  correspondences}

Given any correspondence  of the form 
\[
\begin{tikzcd}
Y & \ar[l, "h_0"']  \Hk_Y^1 \ar[r, "h_1"]  & Y,
\end{tikzcd}
\]
we define the ``$r$-legged version'' $\Hk_Y^r$ as the iterated fibered product
\begin{equation}\label{eq: Hk_Y^r}
\underbrace{\Hk_Y^1 \times_{h_1, Y, h_0}  \Hk_Y^1 \times_{h_1, Y, h_0} \ldots  \times_{h_1, Y, h_0} \Hk_Y^1 }_{r \text{ factors}}.
\end{equation}
We regard this as a correspondence 
\begin{equation}\label{eq:HkY-corr}
Y \xleftarrow{h_0} \Hk_Y^r \xrightarrow{h_r} Y.
\end{equation}
where $h_0$ is the projection from \eqref{eq: Hk_Y^r} to the leftmost factor of $\Hk_Y^1$ followed by $h_0$, and $h_r$ is the projection from \eqref{eq: Hk_Y^r} to the rightmost factor $\Hk_Y^1$ followed by $h_1$.   
When $r=0$,  the correspondence $\Hk_Y^r$ is understood to mean $Y$, and \eqref{eq:HkY-corr} becomes
\[
Y \xleftarrow{\Id} Y \xrightarrow{\Id} Y.
\]

Applying the above construction to the Hecke correspondences \eqref{MtoBun hecke} yields a similar $r$-legged analogue of that diagram, and in particular morphisms
\begin{equation}\label{three Hkr stacks}
\Hk^r_{\cM} \to \Hk^r_{\U(n)} \to \Hk^r_{\U^\dagger(1)} .
\end{equation}
For example,  a point  
$
(x_{\bu}', \cF_{\bu} )  \in \Hk^r_{\U(n)}(R)
$
consists of  a diagram of modifications
\begin{equation}\label{Hkr point}
\begin{tikzcd}
 & {  \cF^\flat_{1/2}  }   \ar[dl , "x'_1" ' ] \ar[dr,  " \sigma(x'_1) " ]  & &   {    \cF^\flat_{3/2}   }    \ar[dl,  " x'_2" ' ]   & {  \cdots  } &   {   \cF^\flat_{ r - 1/2 }   }   \ar[dl ,  "  x'_r  "  '  ] \ar[dr , "  \sigma(x'_r)  " ]     \\ 
{ \cF_0 }  &  & {  \cF_1 }  &{  \cdots}    & { \cF_{r-1} }  &  &     { \cF_r  } 
\end{tikzcd}
\end{equation}
of rank $n$ vector bundles on $X'_R$, each as in \eqref{Hk U(n) mod}, with $x_1',\ldots, x_r' \in X'(R)$ and $\cF_0,\ldots, \cF_r \in \Bun_{\U(n)}(R)$.
To simplify notation, we will shorten \eqref{Hkr point} to 
\[
\begin{tikzcd}
\cF_0 \ar[r, dashed] & \cF_1 \ar[r, dashed] & \ldots \ar[r, dashed] & \cF_r . 
\end{tikzcd}
\]
We will not similarly spell out the points of the other two iterated Hecke stacks in \eqref{three Hkr stacks},   but we do note that the commutativity of \eqref{global hecke} implies that the Hitchin fibration \eqref{eq: basic Hitchin} sits in a  commutative diagram
\begin{equation}\label{global hecke r legs}
\begin{tikzcd}
& { \Hk_{\cM }^r }  \ar[dl, "h_0^{\cM}"'] \ar[dr, "h_r^{\cM}"]    \\
{  \cM  }  \ar[dr, "\pi"']   & &    {  \cM }  \ar[dl, "\pi"]  \\ 
& { \cA  }.
\end{tikzcd}
\end{equation}


\subsection{Shtukas and corank one special cycles}
\label{ssec: corank 1 special cycles} 


Recall from \cite[\S 6.3]{FYZ} the Deligne-Mumford  stack  of rank $n$ Hermitian shtukas with $r \ge 0$ legs.
A point 
\begin{equation}\label{shtuka point}
(x_{\bu}', \cF_{\bu}, \varphi) \in \Sht_{\U(n)}^r(R)
\end{equation}
defined over a $k$-algebra $R$ consists of an underlying point $(x_{\bu}', \cF_{\bu} )  \in \Hk^r_{\U(n)}(R)$ as in \eqref{Hkr point},  together with an isomorphism  $ \varphi : \cF_r \cong \ft \cF_0$ respecting Hermitian structures.  Here  $\ft(-)$ denotes Frobenius twist. 
More generally, associated to each of the iterated Hecke stacks in \eqref{three Hkr stacks} is a similarly defined moduli space of shtukas, related by morphisms
\[
\begin{tikzcd}
{ \Sht_{\cM}^r } \ar[r]   \ar[d]  &  {   \Sht_{\U(n)}^r  }  \ar[r, " \det^{\dagger} "  ]   \ar[d]   &  {  \Sht_{\U^\dagger(1) }^r  }  \ar[d]   \\
{ \Hk^r_{\cM} } \ar[r]   &  {   \Hk^r_{\U(n)} }  \ar[r, " \det^{\dagger} "  ]    &  {  \Hk^r_{\U^\dagger(1) }  .}   
\end{tikzcd}
\]
The arrows labeled $\det^\dagger$  are defined by taking the twisted determinant \eqref{eq: twisted det} of each $\cF_i$.

For any $(\cE, a) \in \cA(k)$, consisting of a rank $n-1$ vector bundle $\cE$ on $X'$ and an injective Hermitian morphism $a : \cE \to \sigma^* \cE^\vee$,  we define  the corank one special cycle $\cZ_{\cE}^r(a)$ as the fiber product
\[
\begin{tikzcd}
{ \cZ_{\cE}^r(a) }  \ar[r] \ar[d] & {  \Sht_{\cM}^r }   \ar[d, " \eqref{global hecke r legs}" ]  \\
{  \Spec(k) } \ar[r, "{(\cE,a)}"] & {  \cA  .}  
\end{tikzcd}
\]
More explicitly, a point $(x_{\bu}', \cF_{\bu}, t_{\bu}, \varphi) \in \cZ_{\cE}^r(a)(R)$ consists of an underlying point \eqref{shtuka point}  and a commutative diagram 
\begin{equation}\label{special cycle diagram}
\begin{tikzcd}
\cE_R \ar[r, equals]\ar[d, "t_0"] & 
\cE_R \ar[r, equals]\ar[d, "t_1"]   & \ldots 
 \ar[r, equals]  & \cE_R \ar[r, equals] \ar[d, "t_r"]  \ar[r, equals] & \cE_R \ar[d, "\ft t_0"] \\
\cF_0 \ar[r, dashed] & \cF_1 \ar[r, dashed] & \ldots \ar[r, dashed] & \cF_r \ar[r, "\varphi"] & \ft \cF_0 ,
\end{tikzcd}
\end{equation}
in which  each $t_i$ respects Hermitian structures, in the sense that $\sigma^* t_i^\vee  \circ t_i = a$.

\begin{remark}
Although our $\cM$ does not agree (Remark \ref{rem: Mdifferent}) with the stack of the same name in \cite{FYZ},
our $\cZ^r_{\cE}(a)$ does agree with the stack  of the same name  in \cite[\S 7.2]{FYZ}.  
This is  because the morphism $a$ in the pair $(\cE, a) \in \cA(k)$ is injective, by definition of $\cA$.
\end{remark}

\begin{prop}\label{prop: corank one Z is proper}
The $k$-stack  $\cZ_{\cE}^r(a)$ is proper.   
\end{prop}

\begin{proof}
The twisted determinant \eqref{eq: twisted det} defines a morphism 
\begin{equation}\label{special detd}
\detd \co \cZ_{\cE}^r(a) \rightarrow \Sht_{\U^\dagger(1)}^r 
\end{equation}
sending a point \eqref{special cycle diagram} to 
\[
\begin{tikzcd}
\detd(\cF_0) \ar[r, dashed] & \detd(\cF_1) \ar[r, dashed] & \ldots \ar[r, dashed] & \detd(\cF_r) \ar[rr , " \detd(\varphi)" ]   & & { {}^\tau \detd(\cF_0). } 
\end{tikzcd}
\]
We will exhibit \eqref{special detd} as a composition of  proper morphisms. 
Because $\Sht^r_{\U^\dagger(1)}$ is a proper $k$-stack (an immediate consequence of the properness of $\BunU$), this will prove the claim.

Recall that for a $k$-algebra $R$, a point  $( \cL_\bullet) \in \Sht^r_{ \U^\dagger(1)}(R)$ consists of a sequence of modifications
\begin{equation}\label{Lshtuka}
\begin{tikzcd}
\cL_0 \ar[r, dashed] & \cL_1 \ar[r, dashed] & \ldots \ar[r, dashed] & \cL_r \ar[r , " \iso" ]   &  { {}^\tau \cL_0 } 
\end{tikzcd}
\end{equation}
of line bundles $\cL_0,\ldots, \cL_r  \in \BunU(R)$, which each modification having the form  
\[
\xymatrix{
{  \cL_i }  & &  {  \cL^\flat_{ i+ 1/2} } \ar[ll]_{x_i'} \ar[rr]^{\sigma (x_i') }   & & { \cL_{i+1} } 
}
\]
as in \eqref{eq: HkU mod}.  The line bundles 
\begin{equation}\label{Dshtuka}
\cD_i : = \det(\sigma^* \cE_R) \otimes \cL_i \quad \mbox{and} \quad 
\cD_{i+1/2}^{\flat} := \det (\sigma^* \cE_R) \otimes  \cL_{i+1/2}^{\flat}
\end{equation}
on $X_R'$ are then similarly related by modifications
\[
\xymatrix{
{  \cD_i }  & &  {  \cD^\flat_{ i+ 1/2} } \ar[ll]_{x_i'} \ar[rr]^{\sigma (x_i') }   & & { \cD_{i+1} } 
}
\]
with the property that the Hermitian pairings  on $\cD_i$ and $\cD_{i+1}$ from \eqref{compline herm} restrict to the same Hermitian pairing on $\cD^\flat_{ i+ 1/2} $.
Recalling the doubling map from Definition \ref{def: double}, we  now have morphisms of torsion coherent sheaves
\begin{equation}\label{eq: seq torsion 1}
  \frac{\sigma^* ( \cE_R \oplus \cD_{i+1/2}^\flat)^\vee  }{\cE_R \oplus \cD_{i+1/2}^\flat} 
  \hookleftarrow 
  \frac{\sigma^*   ( \cE _R \oplus \cD_i) ^\vee }{\cE_R \oplus \cD_{i+1/2}^\flat} 
  \surj   \frac{\sigma^* ( \cE_R \oplus \cD_i)^\vee }{\cE_R \oplus \cD_{i}} = \Db(\cE, a, \cL_i)
\end{equation}
and
\begin{equation}\label{eq: seq torsion 2}
  \frac{\sigma^* (  \cE _R \oplus \cD_{i+1/2}^\flat ) ^\vee }{\cE_R \oplus \cD_{i+1/2}^\flat} 
  \hookleftarrow \frac{\sigma^* ( \cE _R \oplus  \cD_{i+1})^\vee }{\cE_R \oplus \cD_{i+1/2}^\flat} \surj \frac{\sigma^* (  \cE_R \oplus \cD_{i+1} )^\vee }{\cE_R \oplus \cD_{i+1}}
   = \Db(\cE, a, \cL_{i+1}).
\end{equation}

Consider the stack $\mathcal{C}  \rightarrow \Sht_{\U^\dagger(1)}^r$ whose fiber over $(\cL_{\bu}) \in\Sht_{\U^\dagger(1)}^r$ parametrizes coherent subsheaves  
\[
L_i \subset \frac{\sigma^* ( \cE_R \oplus \cD_i)^\vee }{\cE_R \oplus \cD_{i}}
\quad \mbox{and}\quad 
L_{i+1/2}^\flat \subset \frac{\sigma^* ( \cE  \oplus \cD_{i+1/2}^\flat)^\vee}{\cE_R \oplus \cD_{i+1/2}^\flat}
\]
(for all $i\in \{0,\ldots, r\}$ and $i\in \{0,\ldots, r-1\}$, respectively) such that 
\begin{itemize}
\item
Each $L_i$ is a Lagrangian subsheaf.
\item 
Each $L_{i+1/2}^\flat$ is contained in the image of the injection $\hookleftarrow$ from \eqref{eq: seq torsion 1}, and maps to $L_i$ under the surjection $\surj$ from \eqref{eq: seq torsion 1}. 
\item 
Each $L_{i+1/2}^\flat$ is contained in the image of the injection $\hookleftarrow$ from \eqref{eq: seq torsion 2}, and maps to $L_{i+1}$ under the surjection $\surj$ from \eqref{eq: seq torsion 2}. 
\end{itemize}
The map $\mathcal{C} \rightarrow \Sht_{\U^\dagger(1)}^r$ is proper, being cut out by closed conditions in a product of partial flag varieties.

Now we claim that \eqref{special detd} can be factored as
\begin{equation}\label{special proper factorization}
\cZ_{\cE}^r(a) \rightarrow  \mathcal{C}  \rightarrow \Sht_{\U^\dagger(1)}^r,
\end{equation}
in which the first arrow is a closed immersion.  
To define the first arrow, fix an $R$-valued point \eqref{special cycle diagram} of $\cZ_{\cE}^r(a)$.
Recall that the modifications in $\cF_{\bu}$ have the form 
\[
\begin{tikzcd}
\cF_i & \cF^\flat_{i+ 1/2} \ar[l, "x' "'] \ar[r, "\sigma (x')"] & \cF_{i+1}
\end{tikzcd}
\]
described in \eqref{Hk U(n) mod}.   The image of this point under \eqref{special detd} is a $\U^\dagger(1)$-shtuka \eqref{Lshtuka}, from which we construct the line bundles \eqref{Dshtuka}.  
The  complementary line trick from \S \ref{ss:cline}  associates to our point \eqref{special cycle diagram}, functorially in $R$,    a sequence of injections
\[
\cE_R \oplus \cD_{i} \inj \cF_{i} \iso \csd \cF_{i} \inj  \sigma^*(  \cE_R \oplus   \cD_i)^\vee.
\]
The first arrow in the sequence restricts to a map $\cE_R \oplus \cD_{i+1/2}^\flat \inj \cF_{i+1/2}^\flat$, and so we obtain a second sequence of injections
\[
\cE_R \oplus \cD_{i+1/2}^\flat \inj \cF_{i+1/2}^\flat \inj  \sigma^*( \cF_{i+1/2}^\flat)^\vee  \inj \sigma^* ( \cE_R \oplus \cD_{i+1/2}^\flat)^\vee,
\]
related to the first by a commutative diagram
\begin{equation}\label{eq: proper diagram}
\begin{tikzcd}
\cE_R \oplus \cD_i \ar[d, hook] & \cE_R \oplus \cD_{i+1/2}^\flat  \ar[l, hook'] \ar[r, hook] \ar[d, hook] & \cE_R \oplus \cD_{i+1}  \ar[d, hook]  \\
\cF_i \ar[d, hook] & \cF_{i+1/2}^\flat \ar[r, hook] \ar[l, hook']  \ar[d, hook]  & \cF_{i+1}  \ar[d, hook]  \\
\sigma^*( \cE_R  \oplus  \cD_i)^\vee  \ar[r, hook] & \sigma^* ( \cE_R  \oplus \cD_{i+1/2}^\flat)^\vee  & \sigma^*( \cE_R  \oplus  \cD_{i+1})^\vee  \ar[l, hook'] 
\end{tikzcd}
\end{equation}
The first arrow in \eqref{special proper factorization} is now defined by sending the point \eqref{special cycle diagram} of $\cZ_{\cE}^r(a)(R)$ to the   $\U^\dagger(1)$-shtuka \eqref{Lshtuka},  together with the subsheaves
\[
L_i := \frac{\cF_i}{\cE_R \oplus \cD_i} \subset \frac{\sigma^* ( \cE_R \oplus \cD_i)^\vee }{\cE_R \oplus \cD_{i}}
\quad \text{ and } \quad  
L_{i+1/2}^\flat := \frac{\cF_{i+1/2}^\flat}{\cE_R \oplus \cD_{i+1/2}^\flat} \subset \frac{\sigma^* ( \cE_R  \oplus \cD_{i+1/2}^\flat)^\vee}{\cE_R \oplus \cD_{i+1/2}^\flat}.
\]

It remains to show that the first arrow in \eqref{special proper factorization}   is a closed immersion. 
From a point $(\cL_{\bu}, L_{\bu}) \in \mathcal{C}(R)$, we can reconstruct  vector bundles $\cF_i$ and $\cF_{i+1/2}^\flat$ fitting into a diagram of the form \eqref{eq: proper diagram} by pulling back along the quotient maps. 
Then by the same argument as in the proof of Proposition \ref{prop: hitchin global to local}, the point $(\cL_{\bu}, L_{\bu}) $ lies in the image of $\cZ_{\cE}^r(a)(R)$ if and only if the following conditions are satisfied: 
\begin{itemize}
\item Each $L_i$ is balanced in $\Db(\cE, a, \cL_i)$.  
\item The isomorphism $\cD_r \cong \ft \cD_0$ carries $L_r$ to $\ft L_0$.
\item Denoting by $\wt{L}_i$ the preimage of $L_i$ under the surjection in \eqref{eq: seq torsion 1}, each torsion sheaf $ \wt{L}_i/L_{i+1/2}^\flat$ is of length 1 supported at $x_i'$, and each $\wt{L}_{i+1}/ L_{i+1/2}^\flat$ is of length $1$ supported at $\sigma (x_i')$.
\end{itemize}
As these are each evidently closed conditions, we are done. 
\end{proof}


\subsection{Calculation of Grothendieck-Lefschetz trace} 
\label{ss:GLdegree}


Fix a $d\ge 0$, and a point 
\[
\mf{a} :=(\cE, a) \in \cA_d(k)
\]
consisting of a rank $n-1$ vector bundle $\cE$ on $X'$ and a Hermitian morphism $a : \cE \to \sigma^* \cE^\vee$.  The associated special cycle
\[
  \cZ_{\cE}^r(a)  \to   \Sht_{\U(n) }^r 
\]
has expected dimension $r$, but this expectation need not be fulfilled.  To correct for this, in  \cite[\S 4]{FYZ2} one finds the construction of a virtual fundamental class 
\[
[\cZ_{\cE}^r(a)]^{\vir} \in \Ch_r( \cZ_{\cE}^r(a) ) .
\]

Denote by $\ell_i$  the pullback of the line bundle $\ell_{\U(n)}$ from Definition \ref{def:U(n) taut} along the $i^\mathrm{th}$ projection  
\[
\Sht_{\U(n)}^r \to \Hk^1_{\U(n)} .
\]
These are the \emph{tautological bundles} $\ell_1,\ldots, \ell_r$ defined in  \cite[\S4.3]{FYZ2}.

Fix a line bundle   $\cE_0$ on $X'$, and abbreviate
\[
d_0 := d(\cE_0)=\deg_X(\omega_X) - \deg(\cE_0)
\]
as in \eqref{eq: d}.
We  now use the intersection pairing of \eqref{intro pairing} to form the $0$-cycle class
\[
\left( \prod_{i=1}^r c_1(p_i^* \sigma^* \cE_0^{-1} \otimes \ell_i ) \right) \cdot  [\cZ_{\cE}^r(a)]^{\vir} \in \Ch_0 (\mathcal{Z}^r_\cE(a)) ,
\]
where
\[
p_i \co \Sht_{\U(n)}^r \rightarrow X'
\]
is the $i^\mathrm{th}$  leg map.
By the properness of $\mathcal{Z}^r_\cE(a)$ proved in Proposition \ref{prop: corank one Z is proper}, the above $0$-cycle has a well-defined degree.

\begin{thm}\label{thm: geometric side}
 Recall the shifted perverse sheaf $\mathscr{K}_d^i \in D_c^b( \cA_d)$ from  \eqref{ICK2}. The degree
\begin{equation}\label{eq: deg chern cycle}
\deg  \left(\left( \prod_{i=1}^r c_1(p_i^* \sigma^* \cE_0^{-1} \otimes \ell_i ) \right) \cdot  [\cZ_{\cE}^r(a)]^{\vir} \right)
\end{equation}
vanishes for odd $r$, and for even $r$ is given by  
\begin{align*}
\eqref{eq: deg chern cycle} =  \frac{ 2 }{  (\log q)^r  }
 \frac{d^r}{ds^r}\Big|_{s=0}
  \Big[  q^{ s (d_0  +  d)  } 
   \cdot    L(2s,\eta) \cdot 
\sum_{  i = 0 }^d  
 \mathrm{Tr} \big( \mathrm{Fr}_\mf{a} ,   ( \mathscr{K}_d^i) _\mf{a}  \big)  \cdot  q^{ i ( 1 - 2  s )    } 
 \Big] .
\end{align*}
\end{thm}

\begin{proof} 
Recalling the tautological bundle  $ \ell_{\cM}$ on $\cM_d$ from Definition \ref{def:M taut},  abbreviate
\[
\mf{c} := c_1( \ell_{\mathcal{M}   }  )  - c_1( p^*  \sigma^* \cE_0  ) \in   
 H^{2}(\Hk_{\cM_d}^1) (1) \ .
\] 
Here $p :  \Hk_{\cM_d}^1 \to X'$ is the natural map.
As explained in \S \ref{ssec:hecke-action-on-cohomology}, this class may be viewed as a cohomological correspondence, and may be   pushed forward along the Hitchin fibration $\pi_d: \cM_d \to \cA_d$ to obtain an endomorphism
\[
\Gamma_{ \mf{c} }
=  \Gamma_{ c_1( \ell_{\mathcal{M}   }  )   } 
-  \Gamma_{ c_1( p^*  \sigma^* \cE_0  ) } 
\in  \End_{   D_c^b( \cA_d )  }( \rR \pi_{d!} \Ql) .
\]
It follows from the Lefschetz trace formula of \cite[Proposition 11.8]{FYZ} that
\begin{align*}
\eqref{eq: deg chern cycle}  & = 
 \Tr \left( \Gamma_{\mf{c} }^r   \circ \mathrm{Fr}_{\mf{a}} , ( \rR \pi_{d!} \Ql)_{\mf{a}} \right)  \\ 
  & = 
 \sum_{m \in \Z} (-1)^m  \Tr \left( \Gamma_{\mf{c} }^r   \circ \mathrm{Fr}_{\mf{a}} , ( \p \rR^m \pi_{d!} \Q_\ell)_{\mf{a}} \right) 
\end{align*}
where $\Gamma_{\mf{c} }^r  = \Gamma_{\mf{c} } \circ \cdots \circ \Gamma_{\mf{c} }$ is the $r$-fold composition of endomorphisms, and $ \p \rR^m \pi_{d!} \Ql$ are the perverse cohomology sheaves of $\rR \pi_{d!} \Ql$.

On the other hand, we know from Corollary \ref{cor:Hitchin-perverse-cohomology} that for every $m$ there is a Frobenius-equivariant isomorphism
\[
\p \rR^{m}  \pi_{d!}\Ql  \cong \bigoplus_{\dim \cA_d + 2i+j= m} \sK_d^i(-i) \otimes \rH^j(\BunU)[\dim \cA_d] .
\]
 The endomorphisms $ \Gamma_{ c_1( \ell_{\mathcal{M}   }  )   }$ and $\Gamma_{ c_1( p^*  \sigma^* \cE_0  ) }$ were computed, up to semisimplification, in Corollary \ref{cor: tautological correspondence} and Corollary \ref{cor:pullback-bundle-action}, respectively.  
 Combining those results, we find that the endomorphism $\Gamma_{ \mf{c} }$ respects the filtration \eqref{eq: filtration}, and its semisimplification acts as 
\[
\Gamma_{ \mf{c} }    = (d+d_0-2i -2j  )   \otimes w 
\]
on the summand  $\sK_d^i ( -i )_{\mf{a}}  \otimes H^j(\Bun_{\U^\dagger(1)})[\dim \cA_d]$.
Recall that  $w$ is the order two automorphism \eqref{eq:w-defn}.

We are only interested in  the \emph{traces} of these endomorphisms, composed with Frobenius. Passage to the semisimplification does not affect the trace, so upon combining this information, we find that 
\begin{align*}
\eqref{eq: deg chern cycle} 
& =  
 \sum_{  i,j \ge 0  }  (-1)^j  \cdot 
 \Tr \left(  \Gamma_{ \mf{c} }^r   \circ \mathrm{Fr}_{\mf{a}} ,      \sK_d^i ( -i )_{\mf{a}}  \otimes \rH^j(\Bun_{\U^\dagger(1)})   \right)  \\
& =  
\sum_{ i ,j  \ge 0} (-1)^j    \cdot   
( d+ d_0  -2 i    - 2j  )^r   \cdot 
  \mathrm{Tr} ( \mathrm{Fr}_{\mf{a}} ,    \mathscr{K}_d^i ( -i )_{ \mf{a}}  )  
  \cdot  \mathrm{Tr} ( w^r \circ  \mathrm{Fr} ,   \rH^j ( \Bun_{\U^\dagger(1)} )    )   .
  \end{align*} 
 If $r$ is odd then, as in the proof of Corollary \ref{cor: trace on U(1)}, we have
 \[
 \mathrm{Tr} ( w^r \circ  \mathrm{Fr} ,   \rH^j( \Bun_{\U^\dagger(1)} )    )=0
 \]
  for all $j\ge 0$, and we are done.  This leaves the case where $r$ is even, so that $w^r=\Id$. In this case, we have 
 \begin{align*}
\eqref{eq: deg chern cycle} 
& =  
\sum_{ i ,j  \ge 0} (-1)^j    \cdot   
( d+ d_0  -2 i    - 2j  )^r   \cdot  
  \mathrm{Tr} \big( \mathrm{Fr}_{\mf{a}} ,     \mathscr{K}_d^i (-i)_{ \mf{a}}  \big)  
  \cdot  \mathrm{Tr} \big(   \mathrm{Fr} ,   \rH^j( \Bun_{\U^\dagger(1)} )    \big)   ,
  \end{align*} 
and the claim follows from the identity
 \[
 L(s,\eta) 
 =  \frac{1}{2} \sum_{j  \ge 0} (-1)^j \cdot  \mathrm{Tr}( \mathrm{Fr},  \rH^j( \Bun_{\U^\dagger(1)} )   ) \cdot q^{- sj }
 \]
already used in the proof of Corollary \ref{cor: trace on U(1)}.
\end{proof}


\section{Eisenstein series and density polynomials}


We now turn to the analytic side of our calculations.  We recall the construction of Siegel Eisenstein series on quasi-split unitary groups, the expression of their nonsingular Fourier coefficients in terms of representation densities, and the expression of their corank one Fourier coefficients in terms of nonsingular coefficients of Eisenstein series on a lower rank group.  Nothing here is really new, but the precise formulas we need are not easily found in the literature.


\subsection{Eisenstein series on quasi-split unitary groups}
\label{ss:eisenstein}


Recall from \S \ref{ss: intro 1} the quadratic character  $\eta : \A_F^\times \to \{ \pm 1\}$ corresponding to the extension of function fields $F'/F$.
Fix an unramified Hecke character 
\begin{equation}\label{hecke}
\chi :   \A_{F'}^\times \to \C^\times,
\end{equation}
whose  restriction  $ \chi_0=\chi|_{ \A_F^\times}$ is a power of $\eta$.

We recall some notation from \cite[\S 2]{FYZ}.
For an integer $m\ge 1$,  we denote by 
\[
H_m := \U(m,m) \subset \mathrm{Res}_{F'/F} \GL_{2m}
\]
the rank $2m$ quasi-split unitary group over $F$, as in \cite[\S 2.1]{FYZ}, and by $P_m=M_m N_m$  its standard Siegel parabolic.  Thus the Levi factor and unipotent radical have adelic points
\begin{align*}
M_m(\A_F) & =  \left\{
m(\alpha) = \begin{pmatrix}  \alpha & \\ & {}^t \sigma(\alpha)^{-1} \end{pmatrix}  : \alpha \in \GL_m(\A_{F'}) 
\right\} \\
 N_m(\A_F) & = \left\{
n(\beta) = \begin{pmatrix}  1_m &  \beta \\ & 1_m  \end{pmatrix}  : \beta \in \Herm_m(\A_F) 
\right\} .
 \end{align*}

In the usual way, we associate  an Eisenstein series
 \[
 E(g,s,\chi)_m  = \sum_{ \gamma \in P_m(F) \backslash H_m(F) } \Phi (\gamma g ,s, \chi)_m
 \]
 on $H_m(\A_F)$  to the   unramified  section 
 \[
 \Phi(g,s,\chi)_m \in I_m(s,\chi) = \mathrm{Ind}_{P_m(\A_F)}^{H_m(\A_F) } \big( \chi \cdot |\cdot |_F^{ m+ \frac{s}{2} } \big),
 \]
 normalized by $ \Phi( 1_{2m} ,s,\chi)_m =1$.  Although we will never explicitly need it,  this Eisenstein series satisfies a functional equation in $s\mapsto -s$: see \cite[Remark 10.1]{FYZ2}.
 
 After fixing a nontrivial additive character $\psi_0 : k \to \C^\times$,  this Eisenstein series has a Fourier expansion
\[
E(g,s,\chi)_m = \sum_{  T \in \Herm_m(F,\omega_F)   }     E_T(g,s,\chi)_m, 
\]
exactly as in \cite[\S 2.2]{FYZ}, whose coefficients are indexed by $m\times m$ Hermitian matrices with entries in the one-dimensional $F'$-vector space of rational $1$-forms on $X'$.  
 
 As  in \cite[\S 2.6]{FYZ},  given a rank $m$ vector bundle $\cE$  on $X'$ and a (not necessarily injective) Hermitian morphism $a :  \cE  \to \sigma^* \cE^\vee$, we define a Fourier coefficient
 \begin{equation}\label{geometric coefficient}
 E_{ ( \cE ,a) }(s,\chi)_m  := E_T( m(\alpha) ,s,\chi)_m 
 \end{equation}
 by fixing a trivialization $ (F')^m \iso \cE_{F'}$ of the space of rational sections of $\cE$, taking $T$ to be the  matrix of the pairing  $a : \cE \otimes \sigma^*\cE \to \omega_{X'}$ with respect to this trivialization, and choosing $\alpha\in \GL_m(\A_{F'})$ in such a way that for every closed point  $w\in X'$ the above trivialization identifies the lattice $ \alpha_w \cO_{F'_w}^m \subset (F_w')^m$ with the completed stalk of $\cE$ at $w$.


\subsection{Density polynomials and non-singular coefficients}
\label{ss:density polynomials}


 When the character \eqref{hecke} satisfies $\chi_0=\eta^m$, it is explained in  \cite[\S 2]{FYZ} how to express the  Fourier coefficient 
 \eqref{geometric coefficient} associated to a non-singular pair $(\cE,a)$  in terms of certain density polynomials.   
 In this subsection, we  recall these formulas and explain how to extend them to  the case  $\chi_0=\eta^{m+1}$.

Denote by $|X|$  the set of closed points of $X$.  Fix a $v\in |X|$,  a uniformizer $\varpi_v \in F_v$ of the corresponding completion of $F$, and denote by $q_v=\# k_v$  the cardinality of residue field of $v$.
For an integer $a\ge 0$, define a polynomial
\begin{align} \label{m and n}
\mathfrak{m}_v(a , T)  := \prod_{ i=0 }^{ a-1 } \big( 1 - \eta_v(\varpi_v)^i q_v^i T \big) \in \Z[T]
\end{align}
of degree $a$.

Given   a torsion Hermitian sheaf $Q \in \Herm (k)$  on $X'$, the completed stalk $Q_v$ is a torsion Hermitian module for $\cO_{F'_v}/ \cO_{F_v}$.
Define the \emph{local density polynomial}
\begin{equation}\label{local density def}
 \mathrm{Den}( T , Q_v)  
 : = 
  \sum_{   I \subset I^\perp \subset  Q_v   } 
 T^{  \mathrm{length}_{\cO_{F_v}}(I) } \mathfrak{m}_v(   t_v'(  I^\perp / I )  , T)   ,
\end{equation}
where  $I$ runs over all totally isotropic $\cO_{F'_v}$-submodules of $Q_v$, and 
\[
t_v'( I^\perp / I  )  : = \frac{1}{2} \dim_{k_v} \big( (I^\perp / I )\otimes_{\cO_{F_v}} k_v \big) \in \Z .
\]
Here we are effectively taking the Cho-Yamauchi formula of \cite[Theorem 2.3 (3)]{FYZ}, which is the unitary version of the formula pioneered in the quadratic setting in \cite{CY}, as our working definition of the local density polynomial; compare with   \cite[Remark 2.4]{FYZ}. 
We will also need the \emph{twisted local density polynomial} 
\begin{equation}\label{local twisted density}
 \mathrm{Den}_\eta( T , Q_v) : = \mathrm{Den} (  \eta_v(\varpi_v) T , Q_v)
 \end{equation}
which does not appear in \cite{FYZ}.

The \emph{global density polynomial} is defined by
\begin{equation}\label{glob density}
 \mathrm{Den} ( T , Q)  
 := \prod_{v\in |X|}   \mathrm{Den} (   T^{  \deg(v) }   , Q_v),
\end{equation}
and similarly the \emph{twisted global density polynomial} is defined by 
\begin{equation}\label{twisted glob density}
 \mathrm{Den}_\eta( T , Q)  
: = \prod_{v\in |X|}   \mathrm{Den}_\eta(   T^{  \deg(v) }   , Q_v) .
\end{equation}

Given a pair $(\cE,a)$ consisting of a  vector bundle $\cE$ on $X'$ and an injective  Hermitian morphism $a: \cE \to \sigma^*\cE^\vee$, set 
\begin{align}\label{Edensity}
\mathrm{Den} ( T , \cE)  &:= \mathrm{Den} ( T , Q)  \\
 \mathrm{Den}_\eta ( T , \cE)  & := \mathrm{Den}_\eta ( T , Q) \nonumber
\end{align}
where $Q := \coker( a )$ is the  torsion Hermitian sheaf associated to $(\cE,a)$.  
Our interest in these polynomials stems from their close connection with the Eisenstein series of the previous subsection.

\begin{prop}\label{prop:eisenstein density}
For any pair  $(\cE,a)$  consisting of a rank $m$ vector bundle $\cE$ on $X'$ and an injective Hermitian morphism
$a : \cE \to \sigma^*\cE^\vee$, the Fourier coefficient \eqref{geometric coefficient} is given by  
\[
E_{(\cE  , a ) } \left(   s   ,\chi \right)_m
 = \frac{   \chi(\det(\cE) ) q^{   (  s  -   \frac{m}{2}  ) d(\cE)  }        q^{    -ms    \deg_X(\omega_X)   }    }{ \sL_m(s,\chi_0)  }    \times
\begin{cases}
\mathrm{Den}(q^{-2s} , \cE )  & \mbox{if }\chi_0 =\eta^m \\
\mathrm{Den}_\eta(q^{-2s} , \cE )  & \mbox{if }\chi_0 =\eta^{m+1} ,
\end{cases}
\]
where $d(\cE)$ and $ \sL_m(s,\chi_0)$ are defined in the introduction, cf. \eqref{eq: d} and \eqref{eq: intro L-factor}.  
\end{prop}

\begin{proof}
The case $\chi_0=\eta^m$  is \cite[Theorem 2.8]{FYZ}, and we explain the  changes needed for the other case.  
By  \cite[(2.15)]{FYZ}, there is a factorization 
\[
E_{(\cE  , a ) } \left(   s   ,\chi \right)_m
 =  \chi(\det(\cE) ) q^{   (  s  -   \frac{m}{2}  ) d(\cE)  }        q^{    -s  m  \deg_X(\omega_X)   }    
  \prod_{v \in |X|}  W_{ T_v}  ( 1 , s,\Phi_v) 
\]
as a product of the local Whittaker functions of \cite[(2.2)]{FYZ}, where $T_v \in \Herm_m(F_v, \omega_{F_v})$  is  the matrix of the local Hermitian $\cO_{F'_v}$-lattice $\cE_v$ obtained by completing the stalks of $\cE$ at the points of $X'$ above $v \in |X|$.
  This holds for both  $\chi_0=\eta^m$ and $\chi_0=\eta^{m+1}$.
  
If $\chi_0=\eta^m$ then \cite[Lemma 2.7]{FYZ} tells us that
\[
W_{ T_v } ( 1 , s,\Phi_v)
= \mathrm{Den}(q_v^{-2s} , \cE_v) \cdot 
 \prod_{i=1}^m \frac{1}{ L_v(2s+i ,\eta^i)}.
\]
If $\chi_0=\eta^{m+1}$ we can apply \cite[Lemma 2.3.1]{ChenIV} to obtain
\begin{align*}
W_{ T_v} ( 1 , s,\Phi_v)
& = \mathrm{Den}_\eta(  q_v^{-2s} , \cE_v) \cdot 
\prod_{i=1}^m  \frac{1}{ L_v(2s+i ,\eta^{i+1})}  .
\end{align*}
The proposition follows immediately.
\end{proof}

The following functional equation tells us that  both polynomials in \eqref{Edensity} have palindromic coefficients.

\begin{prop}[Functional equation]\label{prop:functional equation}
For any pair $(\cE, a)$ consisting of a vector bundle $\cE$ on $X'$ and an injective Hermitian morphism $a : \cE \to \sigma^* \cE^\vee$, the global density polynomial $\Den(T,\cE)$  satisfies the functional equation
\[
  \Den  (   T  , \cE)  =  T^{d(\cE)}     \Den ( 1 / T   , \cE ).
  \]
  The same functional equation holds for the twisted global density polynomial $\Den_\eta(T,\cE)$.
\end{prop}

\begin{proof}
Fix a $Q \in \Herm_{2d}(k)$.   There are local functional equations
\begin{align*}
\mathrm{Den}(T , Q_v) 
&= ( \eta(\varpi_v)T)^{ \frac{1}{2}  \mathrm{length} (   Q_v ) }  
\mathrm{Den}( 1/T ,  Q_v )  \\
\mathrm{Den}_\eta(T ,Q_v) 
&= T^{ \frac{1}{2}  \mathrm{length} (  Q_v ) }   \mathrm{Den}_\eta( 1/T , Q_v ) ,
\end{align*}
where length is taken as an $\cO_{F_v}$-module, so that $2 d= \sum_{v\in |X|} \mathrm{length} (  Q_v )$.  
The first of these is  \cite[Remark 2.6]{FYZ}, and the second is immediate from the definition \eqref{local twisted density}.

The local functional equations imply the global functional equations
\begin{align*}
  \Den  (   T  , Q) &  =  T^d     \Den  ( 1 / T   , Q ) \cdot   \prod_{ v \in |X| } \eta(\varpi_v)^{  \frac{1}{2}  \mathrm{length} (  Q_v )} \\
    \Den_\eta  (   T  , Q)  &=  T^d     \Den_\eta ( 1 / T   , Q ) ,
 \end{align*}
and \cite[Lemma 11.13]{FYZ} implies that the product over $v\in |X|$ on the right hand side of the first equality is  $1$ whenever $Q= \coker(a : \cE \to \sigma^*\cE^\vee)$ is the cokernel of a Hermitian morphism of vector bundles. 
  \end{proof}

\begin{remark}These local and global functional equations can also be deduced from results stated in \S \ref{ss:density geometrization statements}.
See especially \eqref{eq: inert local density 1} and Theorem \ref{thm: maindensityresult}.
\end{remark}

\begin{example}\label{ex:simple densities}
Suppose that  $Q \in \Herm_2(k)$  is supported above a single closed point $v\in |X|$ of degree one.
Directly from the definitions one finds  $ \mathrm{Den}_\eta ( T , Q)   = 1 +   T$,  while 
\[
\Den(T, Q)= \begin{cases} 1+ T & \mbox{if $v$ is split}   \\ 1-T & \mbox{if $v$ is inert.}
\end{cases} 
\]
In the case where $v$ is inert in $X'$ the coefficients of $\Den(T, Q)$ are not palindromic, but this does not contradict 
Proposition \ref{prop:functional equation}.   
Instead, it demonstrates that it is impossible to express such a $Q$  as the cokernel of a Hermitian morphism of vector bundles on $X'$
(it would be possible if we relaxed our standing assumption that $X' \to X$ is unramified).
\end{example}


\subsection{Singular Fourier coefficients} 


We now explain how to express the corank one Fourier coefficients of the Eisenstein series $E(g,s,\chi)_m$ on $H_m(\A_F)$ in terms of the non-singular  Fourier coefficients of the Eisenstein series $E(g,s,\chi)_{m-1}$ on $H_{m-1}(\A_F)$.
The latter are known by Proposition \ref{prop:eisenstein density}.

Suppose we are given an exact sequence of vector bundles
\[
0 \to \cE_0 \to \cE \to \cE^\flat \to 0
\]
with $\cE$ and $\cE^\flat$ having ranks $m$ and $m-1$, respectively.
Suppose we are also given  an injective Hermitian morphism $a^\flat : \cE^\flat \to \sigma^* (  \cE^\flat )^ \vee $, and  denote by $a$ the composition 
\[
\cE  \to \cE^\flat \map{a^\flat}  \sigma^* (  \cE^\flat )^ \vee   \to  \sigma^*  \cE^\vee .
\]

\begin{prop}\label{prop:genus drop}
We have the equality
\begin{align*}
&
 q^{ ms \deg_X(\omega_X) }  \mathscr{L}_m (s , \chi_0 ) \cdot  E_{ ( \cE , a ) } (s,\chi)_m  \\
& =
  \chi(\cE_0)  q^{(\frac{m}{2} +s) \deg(\cE_0)}   q^{ m s \deg_X(\omega_X) }   \mathscr{L}_m(s,\chi_0)  \cdot  E _{(\cE^\flat  , a^\flat ) } (  s+ 1/2   ,\chi  )_{m-1}   \\
& \quad   +   \chi(\cE_0)  q^{(\frac{m}{2}-s) \deg(\cE_0)}   q^{ -ms \deg_X(\omega_X) }   \mathscr{L}_m(-s,\chi_0)  \cdot  E_{(\cE^\flat  , a^\flat ) } (    - s + 1/2   ,\chi  )_{m-1} ,
\end{align*}
where $\sL_m(s,\chi_0)$ is defined by \eqref{eq: intro L-factor}.
\end{prop}

 \begin{proof}
 There is a general method to express the corank $r$ Fourier coefficients of $E(g,s,\chi)_m$ in terms of non-singular coefficients of $E(g,s,\chi)_{m-r}$,  explained in \cite[\S 5.2.2]{GS19}.  
The particular case of corank one is spelled out in detail in \cite[\S 2.4]{ChenIV}, at least in the  number field setting.  In the function field setting, the case $m=2$ is spelled out in detail in \cite[\S 2.5]{CH}. 
As the proof of the proposition in the generality stated here requires no new ideas, we leave the details of the calculation  to the interested reader.
 \end{proof}


\section{Geometrization of twisted density polynomials}
\label{s:density geometrization}


One of the main technical results of \cite{FYZ} is the geometrization of the density polynomials $\Den(T, Q)$ of \S \ref{ss:density polynomials} in the sense of the sheaf-function correspondence. In other words, \cite{FYZ} constructs a polynomial with coefficients being perverse sheaves on $\Herm_{2d}$, such that the trace of Frobenius on the stalk at $Q \in \Herm_{2d}(k)$ is $\Den(T, Q)$.  In this section we do the same for the twisted density polynomials $\Den_\eta(T, Q)$.


\subsection{Statement of the result}
\label{ss:density geometrization statements}


Recall from \S  \ref{ssec: springer sheaves}  the perverse sheaf
\[
\Spr_{2d}^{\Herm}  \in D_c^b( \Herm_{2d}),
\]
with its  Springer action of the group $W_d = (\Z/2\Z)^d \rtimes S_d$. 

Given a pair of integers $a,b \ge 0$ such that $a+b=d$, we  regard $W_a \times W_b \subset W_d$ in the obvious way.
If we denote by 
\[
\chi_a \co (\Z/2\Z)^a \rightarrow \{\pm 1 \}
\]
  the character obtained by summing the coordinates and then applying the identification $(\Z/2\Z) \cong \mu_2$, then the  $(\chi_a \boxtimes \triv_b)$-isotypic component
  \begin{equation}\label{HSprab}
  \HSpr_{a,b} := (\Spr_{2d}^{\Herm})^{  \chi_a \boxtimes \triv_b  }  \subset \Spr_{2d}^{\Herm} 
  \end{equation}
for the action of   $(\Z/2\Z)^a \times (\Z/2\Z)^b \subset W_d$   is a perverse sheaf on $\Herm_{2d}$, and  carries an action of $W_a \times W_b$.  In particular, it carries an action of the subgroup $S_a \times S_b$.  
According to   \cite[Theorem 5.3 and Proposition 12.3]{FYZ}, for any $Q \in \Herm_{2d}(k)$  we have the equality of polynomials
    \begin{equation}\label{eq: inert local density 1}
    \Den(T,Q) = \sum_{a+b=d}  \Tr( \Fr_Q , ( \HSpr_{a,b}^{S_a \times S_b})_Q ) \cdot   T^a.
    \end{equation}

In the special case where $(a,b)=(0,d)$, the perverse sheaf $\HSpr_{a,b}$ agrees with the sheaf 
 \begin{equation}\label{HSpr}
 \HSpr_d :=  (\Spr_{2d}^{\Herm})^{(\Z/2\Z)^d}
 \end{equation}
  from  \cite[Definition 4.7]{FYZ}.  
The following variant of \eqref{eq: inert local density 1} expresses the  twisted global density polynomials of \S \ref{ss:density polynomials} in terms of the perverse sheaves
    \[
  \sK_d^i = (\Spr_{2d}^{\Herm})^{W_i \times W_{d-i}} = \HSpr_{d}^{S_i \times S_{d-i}}  \in D_c^b( \Herm_{2d}) 
  \] 
from \eqref{ICK} and \eqref{ICK3}.

\begin{thm}\label{thm: maindensityresult}
  For any $Q \in \Herm_{2d}(k)$,   we have the equality of polynomials
    \[
    \Den_{\eta}(T,Q) = \sum_{i=0}^d \Tr(\Fr_Q, (\sK_d^i)_Q)   \cdot T^i.
    \]
\end{thm}

The proof of Theorem \ref{thm: maindensityresult} will be given in \S
\ref{ss:maindensityresultproof}, after we prove some preliminary results.



\subsection{Comparison of Springer fibers}
\label{ss: SpringerComp}

We need to strengthen some results from \cite[\S 4]{FYZ} on the structure of the stalks of $\Spr_{2d}^{\Herm}$.  This requires first recalling a substantial amount of notation from \emph{loc.~cit.}.

\subsubsection{Hermitian Springer fibers} For $Q \in \Herm_{2d}(\ol{k})$, let $\cB_Q^{\Herm}$ be the \emph{Hermitian Springer fiber} defined in \cite[\S 4.4]{FYZ}.  It is the $\ol{k}$-scheme classifying complete flags of $\cO_{X'}$-modules 
\[
0 \subset Q_1 \subset Q_2 \subset \ldots \subset Q_{2d-1} \subset  Q_{2d} = Q 
\]
satisfying $Q_{2d-i} = Q_{i}^\perp$ under the Hermitian form on $Q$.   
This is just the fiber over $Q$ of the morphism \eqref{eq: hermpi}, and so by proper base change there are isomorphisms (as graded $\Ql$-vector spaces)
\begin{equation}\label{eq: Spr fiber stalk}
(\Spr_{2d}^{\Herm})_Q =   ( \rR  \pi_{2d*}^{\Herm} \Ql )_Q \iso  \rH^*( \cB_Q^{\Herm}  )  .
\end{equation}

Denote by $\supp(Q)$ the support divisor of $Q$, a divisor on $X'_{\ol{k}}$ of degree $2d$, and by  
\[
Z' :=  \{ x\in X'(\ol{k}) : Q|_x \neq 0   \}  
\]
the set of points in the support of this divisor (awkwardly, the support of the support divisor).
As explained in  \cite[\S 4.4]{FYZ}, there is a natural decomposition 
\begin{equation}\label{BHermQ decomp}
\mathcal{B}^{\Herm}_Q  =  \bigsqcup_{ y'  \in \Sigma(Z')  } \mathcal{B}^{\Herm}_Q(y')
\end{equation}
as a disjoint union of open and closed subschemes, indexed by the set  $\Sigma(Z')$  of functions 
\[
y' \co \{1,2, \ldots, 2d\} \rightarrow Z'
\]
satisfying $y'(2d+1-i) = \sigma    y'(i)   $  and the equality $\sum_{i=1}^{2d} y'(i) = \supp(Q)$ of divisors on $X'_{\ol{k}}$.

\subsubsection{Correspondence of Springer bases} Recall from \S \ref{ssec: reduction to split case}  and \cite[\S 4.5]{FYZ} the diagram 
\begin{equation}\label{LagrCorr}
\begin{tikzcd}
\Herm_{2d}
  & 
  \Lagr_{2d}^{\dm} 
    \arrow[l, "b_0"', twoheadrightarrow] 
    \arrow[r, "b_1", twoheadrightarrow] 
  & 
  \Coh_d,
\end{tikzcd}
\end{equation}
where $\Coh_d$ is the moduli space of length $d$ torsion coherent sheaves on $X$.
The points of  $\mathrm{Lagr}_{2d}^{\dm}(\bar{k})$ above $Q$ are in bijection with subsets $Z^\sharp \subset Z'$ satisfying 
\begin{equation}\label{Zsplit}
Z' = Z^\sharp \sqcup \sigma( Z^\sharp ) ,
\end{equation}
by sending $Z^\sharp$ to the Lagrangian subsheaf $Q|_{Z^\sharp} \subset Q$.
The map $b_1$ sends this Lagrangian subsheaf to the torsion coherent sheaf
 \begin{equation}\label{Qflatdef}
 Q^\flat := \nu_* (Q|_{Z^\sharp}) \in \mathrm{Coh}_d(\bar{k})
 \end{equation}
on $X_{\ol{k}}$,  supported on the set  
\[
Z : = \nu(Z')  \subset X(\ol{k}).
\]
 (The notation $Q^\flat$ here is not consistent with the  usage in the proof of Proposition \ref{prop: Herm support}, but is consistent with \cite[\S 4.6.1]{FYZ}.)

\subsubsection{General linear Springer fibers} Fixing one $Z^\sharp$ as above, and hence a $Q^\flat \in \mathrm{Coh}_d(\bar{k})$, 
 let $\cB_{Q^\flat}$ be the \emph{Springer fiber} defined in \cite[\S 3.3]{FYZ}.
 This is the $\ol{k}$-scheme parametrizing  complete flags of $\cO_X$-modules 
\[
0 \subset Q_1 \subset Q_2 \subset \ldots \subset Q_{d-1} \subset Q_d =  Q^\flat. 
\]
As explained in   \cite[\S 3.3]{FYZ}, there is a natural decomposition 
 \begin{equation}\label{BQ decomp}
\mathcal{B}_{Q^\flat} 
= 
\bigsqcup_{ y \in \Sigma(Z)  } \mathcal{B}_{Q^\flat} (y)
\end{equation}
as a disjoint union of open and closed substacks, indexed by the set $\Sigma(Z)$ of functions 
\[
y  \co \{ 1,2, \ldots, d \} \rightarrow Z
\]
satisfying the equality  $\sum_{i=1}^d y(i) = \supp(Q^\flat)$ of degree $d$ divisors on $X_{\ol{k}}$.
The support divisor $\supp(Q^\flat)$ is actually independent of the choice of $Z^\sharp$, as it agrees with the image of $Q$ under the support map \eqref{eq: support map}.

\subsubsection{Isomorphisms of Springer fibers} Because of the relation \eqref{Zsplit}, the double cover $\nu : X' \to X$ restricts to a bijection $\nu : Z^\sharp \to Z$.
For each $y \in \Sigma(Z)$, we get a $y^{\sh} \in \Sigma(Z')$ by defining $y^{\sh}(i) = y(i)^{\sharp}$ to be the unique lift of $y(i)$ in $Z^\sharp$ for $i \in \{1, \ldots, d\}$, and then extending $y^{\sh}$ to $\{1, \ldots , 2d\}$ by the equation $y^{\sh}(2d+1-i) = \sigma y^{\sh}(i)$. This defines  an injection 
\[
\Sigma(Z) \map{ y \mapsto y^\sharp}  \Sigma(Z') 
\]
whose image we denote by  $\Sigma(Z^\sharp) \subset \Sigma(Z')$. In other words, $\Sigma(Z^\sharp) \subset \Sigma(Z')$ is the subset of those $y'\in \Sigma(Z')$ for which $y'(i) \in Z^\sharp$ for all $i \in \{ 1,2,\ldots, d\}$.

According to  \cite[(4.1)]{FYZ}, for every $y\in \Sigma(Z)$ there is a canonical isomorphism
\[
\gamma_{Z^\sharp,y^\sharp} : 
\mathcal{B}_{Q^\flat} (y) \iso \mathcal{B}^{\Herm}_Q(y^\sharp).
\]
By varying $y$ and using \eqref{BHermQ decomp} and \eqref{BQ decomp}, we obtain a canonical open and closed immersion
\begin{equation}\label{gamma}
\gamma_{Z^\sharp} : \mathcal{B}_{Q^\flat} \inj \mathcal{B}^{\Herm}_Q
\end{equation}
with image 
$ \bigsqcup_{ y'  \in \Sigma(Z^\sharp)  } \mathcal{B}^{\Herm}_Q(y')$.
In particular, there is an induced injection
\begin{equation}\label{eq: comparison injection} 
i_{Z^\sharp} : \rH^*(  \mathcal{B}_{ Q^\flat }  )  \inj \rH^*( \mathcal{B}^{\Herm}_Q ),
\end{equation}
whose image is (canonically) a direct summand.  The pullback 
\begin{equation}\label{eq: comparison surjection} 
\gamma_{Z^\sharp}^* : \rH^*( \mathcal{B}^{\Herm}_Q )  \twoheadrightarrow
 \rH^*(  \mathcal{B}_{ Q^\flat }  ) ,
\end{equation}
is just the projection to this  summand.

\subsubsection{Equivariance for Springer action} The sheaf $\Spr^{\Herm}_d$ carries a Springer action of  $W_d = (\Z/2\Z)^d \rtimes S_d$.
Hence, by \eqref{eq: Spr fiber stalk}, so does $\rH^* (\mathcal{B}^{\Herm}_Q  )$.  
Analogously,  the discussion of  \cite[\S 3.2]{FYZ} provides us with a perverse sheaf
\[
\Spr_d \in D_c^b( \Coh_d) 
\]
carrying a Springer action of $S_d$, and related to $\rH^*(  \mathcal{B}_{ Q^\flat }  )$ by the obvious analogue of \eqref{eq: Spr fiber stalk}. 
 In particular,  there is a Springer action of $S_d$ on $\rH^*(  \mathcal{B}_{ Q^\flat }  )$. 

We realize $W_d$ as the group of permutations of $\{ 1, \ldots, 2d\}$ commuting with the involution $i \mapsto 2d+1-i$. The subgroup $S_d \subset  W_d$ is realized as the subgroup of permutations that stabilize $\{1, \ldots, d\}$.  We then let $W_d$ act on $\Sigma(Z')$ by precomposition.

We want to show that \eqref{eq: comparison injection} and \eqref{eq: comparison surjection} are $S_d$-equivariant for the Springer actions on source and target, and induce (using Frobenius reciprocity) an isomorphism of $W_d$-representations
\[
\mathrm{Ind}_{S_d}^{W_d} \rH^*(  \mathcal{B}_{ Q^\flat }  ) \iso \rH^*( \mathcal{B}^{\Herm}_Q ).
\]
Note that taking $(\Z/2\Z)^d$-invariants on both sides yields the weaker result 
\[
\rH^*(  \mathcal{B}_{ Q^\flat }  ) \iso \rH^*( \mathcal{B}^{\Herm}_Q )^{ (\Z/2\Z)^d} 
\]
already proved in \cite[\S 4.6.1]{FYZ},  but this weaker result will not suffice for our purposes.
The following Lemma is a first step towards the stronger result. 

\begin{lemma}\label{lem:springer pre-invariance}
The image of \eqref{eq: comparison injection}
is stable under the Springer action of the subgroup $S_d \subset W_d$, and satisfies
\begin{equation}\label{LagrangianSpringer}
\rH^* (\mathcal{B}^{\Herm}_Q  )   =
\bigoplus_{ u \in ( \Z/2\Z)^d } u \cdot i_{Z^\sharp}(  \rH^*(  \mathcal{B}_{ Q^\flat }  )  ) . 
\end{equation}
\end{lemma}

\begin{proof}
Recalling the decomposition \eqref{BHermQ decomp}, the image of \eqref{eq: comparison injection} is 
\[
\bigoplus_{  y' \in  \Sigma(Z^\sharp)  } 
\rH^*( \mathcal{B}^{\Herm}_Q(y') )  \subset \rH^* (\mathcal{B}^{\Herm}_Q  )  .
\] 
Given this, the claim is  a formal consequence of \cite[Lemma 4.9]{FYZ}, the stability of $\Sigma(Z^\sharp) \subset \Sigma(Z')$ under $S_d$, and the decomposition
\[
\Sigma(Z') = \bigsqcup_{  u \in (\Z/2\Z)^d }  u \cdot \Sigma(Z^\sharp). \qedhere
\]
\end{proof}

The next result improves \cite[Proposition 4.11]{FYZ}.

\begin{prop}\label{prop: U-GL comparison}
Recalling \eqref{LagrCorr}, there is a canonical $S_d$-equivariant morphism
\begin{equation}\label{SpringerCompare}
b_0^* \mathrm{Spr}^{\Herm}_{2d} 
 \to  b_1^* \mathrm{Spr}_d 
\end{equation}
of perverse sheaves on $\mathrm{Lagr}_{2d}^{\dm}$.  For every character  $\chi$ of the subgroup $(\Z/2\Z)^d \subset W_d$, this isomorphism  restricts to an isomorphism 
 \begin{equation}\label{chiSpringerCompare}
( b_0^* \mathrm{Spr}^{\Herm}_{2d}  ) ^{ \chi }
\iso 
b_1^* \mathrm{Spr}_d ,
\end{equation}
where the left hand side is the $\chi$-isotypic component.
\end{prop}

\begin{proof}
The desired morphism \eqref{SpringerCompare} is constructed in  the proof of \cite[Proposition 4.11]{FYZ}, but the  $S_d$-equivariance is not established there.

Directly from the construction in \emph{loc.~cit.}, we see that at a geometric point $L \in  \mathrm{Lagr}_{2d}^{\dm}(\bar{k})$ with images 
  \[
  Q \in \mathrm{Herm}_{2d} (\ol{k})
  \quad \mbox{and} \quad 
  Q^\flat \in \mathrm{Coh}_d ( \ol{k}) ,
  \]
   the induced morphism 
  \begin{equation}\label{SpringerStalkMap}
  \rH^* (\mathcal{B}^{\Herm}_Q  )  
  \iso   (  b_0^*\mathrm{Spr}^{\Herm}_{2d})_{ L } 
    \to 
  (b_1^* \mathrm{Spr}_d )_{L}
    \iso 
     \rH^*(  \mathcal{B}_{ Q^\flat }  )
  \end{equation}
  on stalks  agrees with \eqref{eq: comparison surjection} 
  Equivalently, this morphism is projection to the $u=\mathrm{id}$ factor in the decomposition
  \[
  \rH^* (\mathcal{B}^{\Herm}_Q  )   =
\bigoplus_{ u \in ( \Z/2\Z)^d } u \cdot   \rH^*(  \mathcal{B}_{ Q^\flat }    ) 
\]
of  \eqref{LagrangianSpringer}, where we now omit the injection $i_{Z^\sharp}$ from the notation for simplicity.
  
 Let $\chi$ be any character of  $(\Z/2\Z)^d \subset W_d$.  Projection to the $u=\mathrm{id}$ summand above restricts to an isomorphism
 \[
 \rH^* (\mathcal{B}^{\Herm}_Q  )^{  \chi   }   \iso
 \rH^*(  \mathcal{B}_{Q^\flat}  ),
 \]
 and hence (because it is true on stalks) the morphism \eqref{SpringerCompare} restricts to an isomorphism 
 \eqref{chiSpringerCompare}. 
The proof of \cite[Proposition 4.11]{FYZ} shows that the right hand side of \eqref{chiSpringerCompare} is the middle extension of a local system on a dense open substack of $\mathrm{Lagr}_{2d}^{\dm}$, hence the same is true of the left hand side.  As this holds for all characters $\chi$,  the same is also true of  $b_0^* \mathrm{Spr}^{\Herm}_{2d}$.  

To check that   \eqref{SpringerCompare} is $S_d$-equivariant, it now suffices to check this on the level of stalks over the dense open substack of $\mathrm{Lagr}_{2d}^{\dm}$ on which $L$ has multiplicity-free support divisor.   In other words, it suffices to prove the $S_d$-equivariance of  \eqref{SpringerStalkMap} under the assumption that  $Q$ and $Q^\flat$ each have 
multiplicity-free support divisors.  But in this case the source and target of  \eqref{gamma} carry actions of $S_d$ and $W_d$, respectively, and the map \eqref{gamma}  is $S_d$-equivariant by construction.
  \end{proof}

\begin{remark}
Proposition \ref{prop: U-GL comparison} can be restated as saying that there is a canonical $W_d$-equivariant isomorphism of perverse sheaves
\[
b_0^* \mathrm{Spr}^{\Herm}_{2d} 
\iso 
 \mathrm{Ind}_{S_d}^{W_d} ( b_1^* \mathrm{Spr}_d ).
\]
The  right hand side is defined, as a perverse sheaf, by taking the direct sum of copies of $b_1^*\mathrm{Spr}_d$ indexed by $(\Z/2\Z)^d \iso W_d / S_d$.
The action of  $W_d$ is then defined in such a way that Frobenius reciprocity holds. 
\end{remark}

\begin{cor}\label{cor:springer invariance}
The morphisms  
\[
i_{Z^\sharp}:  \rH^*(  \mathcal{B}_{Q^\flat}  )  \hookrightarrow  \rH^* (\mathcal{B}^{\Herm}_Q  )  
\quad \mbox{and} \quad  
\gamma^*_{Z^\sharp}:  \rH^* (\mathcal{B}^{\Herm}_Q  )   \twoheadrightarrow \rH^*(  \mathcal{B}_{ Q^\flat }  )
\]
induced by \eqref{gamma} are equivariant for the Springer $S_d$-actions on source and target.
They induce an isomorphism of $W_d$-representations
\[
\rH^* (\mathcal{B}^{\Herm}_Q  ) 
\iso 
\mathrm{Ind}_{S_d}^{W_d} (   \rH^*(  \mathcal{B}_{Q^\flat}  ) )  
\]
(which depends on the choice of $Z^\sharp$).
\end{cor}

\begin{proof}
The map $\gamma^*_{Z^\sharp}$ agrees with \eqref{SpringerStalkMap}, so its $S_d$-equivariance is a consequence of Proposition \ref{prop: U-GL comparison}.   
The equivariance of $i_{Z^\sharp}$ is  a formal consequence of  the fact that its image is $S_d$-stable (Lemma \ref{lem:springer pre-invariance}), and that it admits the $S_d$-equivariant section  $\gamma^*_{Z^\sharp}$.

The final claim is a consequence of the first and the decomposition \eqref{LagrangianSpringer}.
\end{proof}


\subsection{Comparison of Frobenius actions}
\label{ss:FrobComp}


The results of the previous subsection were over $\bar{k}$. Now we work over $k$ to track Frobenius actions. Fix a Hermitian torsion sheaf  $Q \in \mathrm{Herm}_{2d}(k)$ defined over $k$.

We continue with the notation of the previous subsection, so that $Z' \subset X'(\ol{k})$ is the set of points in the support of the base change $Q_{\ol{k}}$.  As in \eqref{Qflatdef}, each subset $Z^\sharp \subset Z'$ satisfying \eqref{Zsplit} determines a torsion coherent sheaf 
\[
\nu_* (Q_{\ol{k}} |_{Z^\sharp}) \in \Coh_d(\bar{k}).
\]
and  \cite[Lemma 4.12]{FYZ} and its proof show  that, up to isomorphism, there is a unique  $Q^\flat \in \Coh_d(k)$ whose base change to $\ol{k}$ is isomorphic to  $\nu_* (Q_{\ol{k}} |_{Z^\sharp})$.

Recall the $S_d$-equivariant morphisms 
\[
i_{Z^\sharp}:  \rH^*(  \mathcal{B}_{Q^\flat}  )  \hookrightarrow  \rH^* (\mathcal{B}^{\Herm}_Q  )  
\quad \mbox{and} \quad  
\gamma^*_{Z^\sharp}:  \rH^* (\mathcal{B}^{\Herm}_Q  )   \twoheadrightarrow \rH^*(  \mathcal{B}_{ Q^\flat }  )
\]
of  Corollary \ref{cor:springer invariance}.
Because $Q$ and $Q^\flat$ are defined over $k$, the Springer fibers $\mathcal{B}^{\Herm}_Q$ and $\mathcal{B}_{Q^\flat}$ are $k$-varieties, and so the source and target of each map carries an action of  Frobenius.

\subsubsection{Choice of section} To describe how these maps interact with the Frobenius actions, we choose $Z^\sharp$ as in \cite[\S 4.6.1]{FYZ}.
That is, set $Z = \nu(Z') \subset X(\bar{k})$, and let   $|Z| \subset X$ and $|Z'| \subset X'$ be the sets of closed points below the points of $Z$ and $Z'$, respectively.  For each $v \in |Z|$,  fix an $x'_v \in Z'$ above it,  and then set $x_v = \nu(x'_v) \in Z$.  Having made such choices, set 
\begin{equation}\label{good Zsharp}
Z^\sharp = \{ \mathrm{Fr}^i ( x_v') : v\in |Z| \mbox{ and } 0 \le i < \deg(v) \} .
\end{equation}

\subsubsection{Twisting element} For any $y\in \Sigma(Z)$,  let $t_y \in (\Z/2\Z)^d$ be the element determined by 
\begin{equation}\label{t_y def}
(\mbox{$i^\mathrm{th}$-coordinate of $t_y$}) 
= \begin{cases}
1  & \mbox{if $y(i) = x_v$ for some inert $v\in |Z|$} \\
0  & \mbox{otherwise.}
\end{cases}
\end{equation}
For a character $\chi : (\Z/2\Z)^d \to \{ \pm 1\}$, let
\begin{equation}\label{t_chi}
t_\chi :  \rH^* ( \mathcal{B}_{ Q^\flat }   )  \to \rH^* ( \mathcal{B}_{ Q^\flat }   )
\end{equation}
be the automorphism that acts as multiplication by $\chi(t_y)$ on the direct summand indexed by $y$ in the decomposition
\[
\rH^* ( \mathcal{B}_{ Q^\flat }  ) = \bigoplus_{ y \in \Sigma(Z)  }  \rH^* ( \mathcal{B}_{ Q^\flat } (y)  )
\]
induced by \eqref{BQ decomp}.

The group $S_d$ acts as automorphisms of the  character group  of  $(\Z/2\Z)^d$.  
To be precise, if $\epsilon_i$ is the character that is nontrivial on the $i^\mathrm{th}$ component and trivial on all other components, then $s\in S_d$ acts as $s\cdot \epsilon_i = \epsilon_{s(i)}$.  
It is immediate from the definitions that the diagram
\[
\xymatrix{
{  \rH^*( \mathcal{B}_{Q^\flat} )  }  \ar[r]^{t_\chi}  \ar[d]_s  &  {  \rH^*( \mathcal{B}_{Q^\flat} )  }   \ar[d]^s \\
{  \rH^*( \mathcal{B}_{Q^\flat} )  }  \ar[r]_{t_{s\cdot \chi}}     &  {  \rH^*( \mathcal{B}_{Q^\flat} )  } 
}
\]
commutes for any character $\chi : (\Z/2\Z)^d \to \{\pm 1\}$.
In particular, $t_\chi$ commutes with the action of the stabilizer $\mathrm{Stab}(\chi) \subset S_d$.

\subsubsection{Frobenius action} We now prove a generalization of   \cite[Proposition 4.13]{FYZ}.
More precisely \emph{loc.~cit.} is the  case $\chi=1$ of the following   Proposition.

\begin{prop}\label{prop:t-twist}
For any character $\chi :  (\Z/2\Z)^d \to \{ \pm 1\}$ we have a commutative diagram
\[
\begin{tikzcd}
\rH^* (\mathcal{B}^{\Herm}_Q  )^{   \chi    } 
    \arrow[d, "{ \gamma^*_{Z^\sharp} }"'] 
    \arrow[rr, "\mathrm{Fr}"] 
    &  & 
    \rH^* (\mathcal{B}^{\Herm}_Q  )^{  \chi    } 
    \arrow[d, "{ \gamma^*_{Z^\sharp} }"] 
    \\
\rH^*( \mathcal{B}_{ Q^\flat } ) 
    \arrow[r, "{ \theta \circ t_\chi }"'] 
    & 
    \rH^*( \mathcal{B}_{ Q^\flat } ) 
    \arrow[r, "\mathrm{Fr}"'] 
    & 
    \rH^*( \mathcal{B}_{ Q^\flat } ) 
\end{tikzcd}
\]
in which the vertical arrows are isomorphisms of vector spaces,  equivariant with respect to the action of the stabilizer $\mathrm{Stab}(\chi) \subset S_d$.   
The $S_d$-equivariant operator $\theta : \rH^*(  \mathcal{B}_{ Q^\flat }  ) \to  \rH^*(  \mathcal{B}_{ Q^\flat }  )$, which depends on $Z^\sharp$ but not on $\chi$, is that of \cite[Proposition 4.13]{FYZ}.
\end{prop}

\begin{proof}
The claim that the vertical arrows are $\mathrm{Stab}(\chi)$-equivariant isomorphisms follows immediately from Corollary \ref{cor:springer invariance}.   Thus  it suffices to prove the commutativity of the diagram.

 Denote the inverse of the bijection  $\nu : Z^\sharp \xrightarrow{\sim} Z$  by $x\mapsto x^\sharp$, and similarly for the inverse of the induced bijection $\Sigma(Z^\sharp) \to \Sigma(Z)$.
Because of our particular choice of $Z^\sharp$, for any $x\in Z$ we have (cf. the second paragraph of the proof of \cite[Proposition 4.13]{FYZ})
\begin{equation}\label{sharp flip}
 \mathrm{Fr}( x^\sharp ) = 
 \begin{cases}
  \sigma( \mathrm{Fr}( x)^\sharp ) & \mbox{if  $\mathrm{Fr}(x) = x_v$ for some inert $v\in |Z|$}  \\ 
\mathrm{Fr}( x)^\sharp  & \mbox{otherwise}.
\end{cases} 
\end{equation}

As in the proof of \cite[Proposition 4.13]{FYZ}, for each $y\in \Sigma(Z)$ 
 there is a unique $\tau_y \in (\Z/2\Z)^d$ for which 
\[
 \mathrm{Fr} (y)^\sharp = \mathrm{Fr}(y^\sharp) \cdot \tau_y
\]
as elements of $\Sigma(Z')$.    Using  \eqref{sharp flip}, we see that  $\tau_y$ is the tuple whose $j^\mathrm{th}$ component is $1$ if 
$
 \mathrm{Fr} (  y(j) )  =  x_v 
$
for some inert $v\in |Z|$, and is $0$ otherwise. In other words, we have
\begin{equation}\label{eq: tau t}
\tau_y = t_{ \mathrm{Fr}(y) } .
\end{equation}

Using \eqref{eq: tau t} and Corollary \ref{cor:springer invariance}, the desired commutativity follows from the commutativity of the diagram
\[
\xymatrix{
{  \rH^* ( \mathcal{B}_Q^{\Herm}( \mathrm{Fr}(y)^\sharp) )    }  \ar[d]_{ \gamma^*_{ Z^\sharp , \mathrm{Fr}(y)^\sharp} }   \ar[r]^{\tau_y}   & 
{  \rH^* ( \mathcal{B}_Q^{\Herm}( \mathrm{Fr}(y^\sharp) ) )    }   \ar[r]^{\mathrm{Fr} }   \ar[d]_{ \gamma^*_{ \mathrm{Fr}(Z^\sharp)  , \mathrm{Fr}(y^\sharp)} }     & 
{  \rH^* ( \mathcal{B}_Q^{\Herm} (y^\sharp ) )    }  \ar[d]_{ \gamma^*_{ Z^\sharp ,  y^\sharp} }    \\ 
{  \rH^* ( \mathcal{B}_{Q^\flat}( \mathrm{Fr}(y)) )    }   \ar[r]_{\theta_{ \mathrm{Fr}(y) }}   & 
{  \rH^* ( \mathcal{B}_{Q^\flat}( \mathrm{Fr}(y) ) )    }   \ar[r]_{\mathrm{Fr} }    & 
{  \rH^* ( \mathcal{B}_{Q^\flat} (y ) )    } ,
}
\]
found in the  proof of \cite[Proposition 4.13]{FYZ}. 
See especially the diagram (4.7) of \emph{loc.~cit.}.

\end{proof}


\subsection{Proof of Theorem \ref{thm: maindensityresult}}
\label{ss:maindensityresultproof}


For a fixed  $Q \in \Herm_{2d}(k)$,   Theorem  \ref{thm: maindensityresult} is equivalent to the equality
\begin{equation}\label{altmaindensity}
  \Den_{\eta}(T,Q) = \sum_{a+b=d}  \Tr(\Fr_Q, (\HSpr_{d}^{S_a \times S_b}  )_Q)   \cdot T^a.
\end{equation}
The basic idea  for proving \eqref{altmaindensity} is to  deduce it from the known equality  \eqref{eq: inert local density 1}, by using Proposition \ref{prop:t-twist} to directly compare the traces of Frobenius on the stalks  of the perverse sheaves $\HSpr_{a,b}^{S_a \times S_b}  $ and $\HSpr_{d}^{S_a \times S_b} $.

\subsubsection{The case of a split point} The following two Propositions carry out this comparison in the special case where $Q$ is supported above a single closed point of $X$.

\begin{prop}\label{prop:split trace comparison}
Assume that $Q$ is supported above a single closed point $v \in X$ that splits in $X'$.
For all $a,b \ge 0$ satisfying $a+b=d$, we have 
\[
\Tr( \Fr_Q , ( \HSpr_{a,b}^{S_a \times S_b})_Q )=\Tr( \Fr_Q , ( \HSpr_d^{S_a \times S_b})_Q ). 
\]
\end{prop}

\begin{proof}
Fix a geometric point $x'_v \in X'(\ol{k})$ above $v$, and let  $Z^\sharp \subset X'(\ol{k})$ be the subset   \eqref{good Zsharp}.

Now choose a  character   $\chi : (\Z/2\Z)^d \to \{ \pm 1\}$.
 Because $|Z|=\{v\}$ contains only split places, 
the element $t_y \in (\Z/2\Z)^d$ from \eqref{t_y def} is trivial for all $y \in \Sigma(Z)$.
In particular,  the automorphism $t_\chi$ from \eqref{t_chi}  is trivial.  
Applying Proposition \ref{prop:t-twist} to both $\chi$ and the trivial character of $ (\Z/2\Z)^d $, we find an isomorphism
\[
\rH^* (\mathcal{B}^{\Herm}_Q  )^{   \chi    }  \iso \rH^* (\mathcal{B}^{\Herm}_Q  )^{ (\Z/2\Z)^d  }
\]
that respects the Frobenius operators on and source and target, as well as the action of the stabilizer $\Stab(\chi)\subset S_d$.  

Applying this to the character $\chi = \chi_a \boxtimes \triv_b$ of $(\Z/2\Z)^a \times (\Z/2\Z)^b$, and using \eqref{eq: Spr fiber stalk}, 
we obtain an isomorphism
\[
( \HSpr_{a,b})_Q \iso ( \HSpr_d)_Q,
\]
equivariant with respect to both $\Fr_Q$ and  the Springer action of $S_a \times S_b$.  The claim follows immediately.
\end{proof}

\subsubsection{The case of an inert point} Next we turn to the more difficult calculation where $Q$ is supported over a single inert point.

\begin{prop}\label{prop:inert trace comparison}
Assume that $Q$ is supported above a single closed point $v \in X$ that is inert  in $X'$, and  fix  $a,b \ge 0$ satisfying $a+b=d$.
If $\deg(v)$ divides $a$, then 
\[
\Tr( \Fr_Q , ( \HSpr_{a,b}^{S_a \times S_b})_Q )= (-1)^{a / \deg(v)}  \Tr( \Fr_Q , ( \HSpr_d^{S_a \times S_b})_Q ). 
\]
If $\deg(v)$ does not divide $a$, then both traces are $0$.
\end{prop}

\begin{proof}
Fix a geometric point $x'_v \in X'(\ol{k})$ above $v$, and let  $Z^\sharp \subset X'(\ol{k})$ be the subset   \eqref{good Zsharp}.

The decomposition \eqref{BQ decomp} induces a direct sum decomposition
\begin{align*}
\rH^*( \mathcal{B}_{ Q^\flat } ) 
& =
\bigoplus_{ y \in \Sigma(Z)}   \rH^*( \mathcal{B}_{Q^\flat}(y) ) \\
& =
\bigoplus_{ y \in ( S_a \times S_b ) \backslash \Sigma(Z) } 
\left(   \bigoplus_{  u \in (S_a\times S_b) /   (S_a\times S_b )   \cap \mathrm{Stab}(y)    }  u \cdot \rH^*( \mathcal{B}_{Q^\flat}(y) )  \right) .
\end{align*}
Here we are using \cite[Lemma 3.5]{FYZ}, which proves that the Springer action of $S_d$ on $\rH^*( \mathcal{B}_{ Q^\flat } )$ permutes the summands $\rH^*( \mathcal{B}_{Q^\flat}(y) )$ in a way compatible with the natural action of $S_d$ on $\Sigma(Z)$ through precomposition.
By taking fixed points for the Springer action of $S_a \times S_b \subset S_d$ throughout,  we obtain an isomorphism
\begin{equation}\label{Mackey}
\rH^*( \mathcal{B}_{  Q^\flat } ) ^{   S_a \times S_b }  \iso \bigoplus_{ y \in ( S_a \times S_b )  \backslash \Sigma(Z) } V(y)
\end{equation}
in which we have set
\[
V(y) : = \rH^*( \mathcal{B}_{Q^\flat} (y) )^{ (S_a \times S_b)  \cap  \mathrm{Stab}(y)   } .
\]

Each summand $V(y)$ in \eqref{Mackey} is stable under the operator $\theta$ from Proposition \ref{prop:t-twist}.
Pullback along Frobenius  permutes the summands in \eqref{Mackey} according to the rule 
\[
\Fr \co V(y) \xrightarrow{\sim} V(\Fr^{-1}(y)),
\]
where  the action of $\Fr$ on  $\Sigma(Z)$ is  induced by the natural permutation action of $\Fr$ on the set   $Z \subset X( \ol{k})$ of all geometric points above $v$. 
 It follows that 
\begin{equation}\label{first flat trace}
    \mathrm{Tr}\left(  \mathrm{Fr} \circ \theta ,
\rH^*( \mathcal{B}_{  Q^\flat } ) ^{  S_a \times S_b  } 
 \right)  = \sum_{ \substack{  y \in ( S_a \times S_b ) \backslash \Sigma(Z)   \\ \mathrm{Fr}(y) =y   }   } 
   \mathrm{Tr}\left(  \mathrm{Fr} \circ \theta , V(y) 
 \right) .
\end{equation}
Each summand $V(y)$ is also    stable under the automorphism  $t_{ \chi_a \boxtimes \triv_b }$ from \eqref{t_chi}, which  acts   as multiplication by the scalar $(\chi_a \boxtimes  \triv_b)( t_y )$.  Hence the same reasoning also shows 
\begin{align}\label{second flat trace}
   & \mathrm{Tr}\left(  \mathrm{Fr} \circ \theta \circ   t_{ \chi_a \boxtimes \triv_b },
\rH^*( \mathcal{B}_{  Q^\flat } ) ^{  S_a \times S_b  } 
 \right)  \\
 &  = \sum_{ \substack{  y \in ( S_a \times S_b ) \backslash \Sigma(Z)   \\ \mathrm{Fr}(y) =y   }   } 
(\chi_a \boxtimes \triv_b)( t_y) \cdot   \mathrm{Tr}\left(  \mathrm{Fr} \circ \theta , V(y) 
 \right) . \nonumber
\end{align}

The crux of the proof of Proposition \ref{prop:inert trace comparison} is the following combinatorial Lemma.

\begin{lemma}\label{lem:combinatorics}
If  there is a $y\in \Sigma(Z)$ whose $ (S_a \times S_b)$-orbit is stable under the action of $\mathrm{Fr}$,  then $\deg(v)$ divides $a$ and 
\begin{equation}\label{explicit tchi}
(\chi_a \boxtimes \triv_b)(t_y) = (-1)^{ a / \deg(v) }.
\end{equation}
\end{lemma}

\begin{proof}
Let $v' \in X'$ be the unique closed point above $v \in X$.
Let $e$ be the length of $Q$ as an $\cO_{X',v'}$-module, so that  $d=e \deg(v)$.
Partition 
\[
\{ 1, \ldots, d\} = B_1 \sqcup \cdots \sqcup  B_{\deg(v)}
\]
where $B_1=\{ 1,\ldots, e\}$, $B_2=\{ e+1,\ldots, 2e\}$, and so on.
Let $y_0 \in \Sigma(Z)$ be the element defined by
\[
y_0 (i) = \begin{cases}
x_v & \mbox{if }i \in B_1 \\
\mathrm{Fr}(x_v) & \mbox{if } i \in B_2 \\
\quad \vdots \\
\mathrm{Fr}^{\deg(v) -1}(x_v) & \mbox{if } i \in B_{ \deg(v) } .
\end{cases}
\]

The stabilizer of $y_0$ is 
\[
 \mathrm{Stab}(y_0) = \underbrace{S_e \times \cdots  \times S_e}_{ \deg(v) \,  \mathrm{times}} \subset S_d,
\]
and every element of $\Sigma(Z)$\ has the form $y_0 \circ s^{-1}$ for a unique $s \in  S_d / \mathrm{Stab}(y_0)$.
The tuple   $t_{ y_0 \circ s^{-1} }  \in (\Z/2\Z)^d$ from \eqref{t_y def} has the explicit form 
\begin{equation}\label{explicit t}
(\mbox{$i^\mathrm{th}$-coordinate of $t_{ y_0 \circ s^{-1} } $}) 
= \begin{cases}
1  & \mbox{if } i  \in s(B_1)  \\
0  & \mbox{otherwise,}
\end{cases}
\end{equation}
and the action of $\mathrm{Fr}$ on $\Sigma(Z)$ is given by
$
\mathrm{Fr}( y_0 \circ s^{-1}) = y_0 \circ  (s\gamma)^{-1}
$
where  $\gamma \in S_d$ is any permutation restricting to bijections
\begin{align*}
B_{\deg(v)}  \map{\gamma} B_{\deg(v)-1}  \map{\gamma} \cdots  \map{\gamma}   B_2   \map{\gamma} B_1  .
\end{align*}

If  $y= y_0 \circ s^{-1} \in    (S_a  \times S_b ) \backslash  \Sigma(Z)$ is fixed by $\mathrm{Fr}$, then 
$s = s \gamma$ as elements of $(S_a  \times S_b ) \backslash S_d / \mathrm{Stab}(y_0)$, and hence 
 there is a $\mu \in S_a  \times S_b$ for which 
\[
 \mu s  \in s  \gamma  \mathrm{Stab}(y_0)   .
\]
The permutation $\mu$ restricts to bijections
\[
s(B_{k+1})  \cap \{ 1 ,\ldots, a \} \map{\mu}  s(B_k)  \cap \{ 1 ,\ldots, a \}
\]
for all $1 \le k < \deg(v)$, and it follows that 
\[
\{1,\ldots, a \} = \bigsqcup_{k=1}^{\deg(v)} 
\big( s(B_k)  \cap \{ 1 ,\ldots, a \} \big)   
\]
is  a partition of $\{1,\ldots, a \}$ into $\deg(v)$ subsets, each of  the same cardinality.

In particular $\deg(v)$ divides $ a$, and  $s(B_1) \cap  \{ 1 ,\ldots, a \}$ has $a / \deg(v)$ elements.
Recalling \eqref{explicit t}, we see that among the first $a$ coordinates of $t_y=t_{y_0 \circ s^{-1}}$,  exactly $a/ \deg(v)$ of them are equal to $1$, proving \eqref{explicit tchi}.
\end{proof}

We can now complete the proof of Proposition \ref{prop:inert trace comparison}.
Combining \eqref{first flat trace} and \eqref{second flat trace} with  Lemma \ref{lem:combinatorics} shows that 
\begin{equation}\label{flat trace comparison}
 \mathrm{Tr}  (  \mathrm{Fr} \circ \theta ,
\rH^*( \mathcal{B}_{  Q^\flat } ) ^{  S_a \times S_b  }  )   = (-1)^{ a / \deg(v) }  \mathrm{Tr} (   \mathrm{Fr} \circ \theta \circ   t_{ \chi_a \boxtimes \triv_b },
\rH^*( \mathcal{B}_{  Q^\flat } ) ^{  S_a \times S_b  } )
\end{equation}
if $\deg(v)$ divides $a$, and otherwise both traces are $0$.

Applying Proposition \ref{prop:t-twist} twice, once with $\chi$ the trivial character of $(\Z/2\Z)^d$ and once with $\chi=\chi_a \boxtimes \triv_b$, shows that there are $(S_a \times S_b)$-equivariant isomorphisms
\[
\rH^* (\mathcal{B}^{\Herm}_Q  )^{ (\Z/2\Z)^d  }   \iso \rH^*( \mathcal{B}_{ Q^\flat } )  \iso  \rH^* (\mathcal{B}^{\Herm}_Q  )^{   \chi_a \boxtimes \triv_b    },
\]
under which the operator $\Fr$ on the first and last cohomology spaces agree with the operators $\mathrm{Fr} \circ \theta$ and $\mathrm{Fr} \circ \theta \circ   t_{ \chi_a \boxtimes \triv_b }$, respectively,  on $\rH^*( \mathcal{B}_{ Q^\flat } )$.  Proposition \ref{prop:inert trace comparison} follows by combining  these isomorphisms, the equality \eqref{flat trace comparison}, and the isomorphism \eqref{eq: Spr fiber stalk}.
\end{proof}

\subsubsection{The general case} At last, we prove the main result of \S \ref{s:density geometrization}. 

\begin{proof}[Proof of Theorem \ref{thm: maindensityresult}]
Recall that we must prove the equality of polynomials \eqref{altmaindensity}.
First suppose that $Q$ is supported above a single closed point $v \in X$.
In particular, \eqref{glob density} and \eqref{twisted glob density} simplify to 
\[
\mathrm{Den} ( T , Q)  =    \mathrm{Den} (   T^{  \deg(v) }   , Q_v) \quad \mbox{and} \quad  \mathrm{Den}_\eta( T , Q)  
 =    \mathrm{Den}_\eta(   T^{  \deg(v) }   , Q_v).
\]

If $v$ is split in $X'$, it is clear from the definitions of \S \ref{ss:density polynomials} that   $  \Den_{\eta}(T,Q)=  \Den(T,Q)$.
Hence  \eqref{altmaindensity}  follows by combining  \eqref{eq: inert local density 1} and Proposition \ref{prop:split trace comparison}.

If $v$ is inert in $X'$,  the equality  \eqref{eq: inert local density 1} implies that the local density polynomial of $Q_v$ is 
\begin{align*}
 \Den (T,Q_v)   = \sum_{ \substack{ a+b=d \\   \deg(v) \mid a  } }  \Tr( \Fr_Q , ( \HSpr_{a,b}^{S_a \times S_b})_Q ) \cdot   T^{ a / \deg(v) }  .
\end{align*} 
By Proposition \ref{prop:inert trace comparison}, the twisted local density polynomial \eqref{local twisted density} is therefore equal to
\begin{align*}
 \Den_\eta (T,Q_v)    & =  \Den (-T,Q_v)  \\
 & = \sum_{ \substack{ a+b=d \\   \deg(v) \mid a  } }  \Tr( \Fr_Q , ( \HSpr_d^{S_a \times S_b})_Q ) \cdot   T^{ a / \deg(v) }   .
\end{align*} 
The desired equality \eqref{altmaindensity}  follows by replacing $T$ with $T^{\deg(v)}$, and using  the  final claim of Proposition \ref{prop:inert trace comparison}.

We have now proved \eqref{altmaindensity} when $Q$ is supported above a single closed point of $X$.
For the general case, one decomposes $Q$ as a direct sum of Hermitian torsion sheaves on $X'$ with pairwise disjoint supports,  each supported above a single closed point of $X$.  
The left hand side of \eqref{altmaindensity} is multiplicative with respect to this decomposition, as one sees immediately from the definition \eqref{twisted glob density}.  Hence we are reduced to proving that the right hand side has this same  multiplicative property.  

Using \eqref{eq: Spr fiber stalk}, and applying Proposition \ref{prop:t-twist} to the trivial character $\chi$, we have isomorphisms of graded vector spaces
\[
(\HSpr_d^{S_a \times S_b})_Q 
\iso 
\rH^* (\mathcal{B}^{\Herm}_Q  )^{ W_a \times W_b} 
\iso 
\rH^*( \mathcal{B}_{ Q^\flat } )^{S_a \times S_b}.
\]
Recalling the operator $\theta$ from Proposition \ref{prop:t-twist}, it follows that the right hand side of \eqref{altmaindensity} is equal to 
\begin{equation}\label{multiplicative polynomial}
\sum_{a+b=d}  \Tr(\Fr \circ \theta , \rH^*( \mathcal{B}_{ Q^\flat } )^{S_a \times S_b}  )   \cdot T^a.
\end{equation}

The proof that \eqref{multiplicative polynomial} is multiplicative with respect to the support decomposition of $Q^\flat$  is virtually identical to similar multiplicative properties proved in \cite{FYZ}.  More precisely, one first proves this multiplicativity for 
\[
\sum_{a+b=d}  \Tr(\Fr  , \rH^*( \mathcal{B}_{ Q^\flat } )^{S_a \times S_b}  )   \cdot T^a
\]
using the same argument as in \cite[Lemma 5.4]{FYZ}; in fact, the proof in our case can be simplified somewhat,  because we are taking $(S_a \times S_b)$-invariants, rather than taking isotypic components with respect to nontrivial characters.  Then one deduces the multiplicativity of \eqref{multiplicative polynomial} in the same way that   \cite[Proposition 5.5(1)]{FYZ} is deduced from  \cite[Lemma 5.4]{FYZ}.   As these arguments involve no new ideas, we omit the details.
\end{proof}

  
\section{The main results}


Fix an integer $n \ge 1$,  and an unramified Hecke character 
$
\chi :   \A_{F'}^\times \to \C^\times,
$
whose restriction  $ \chi_0=\chi|_{ \A_F^\times}$ satisfies  $\chi_0 = \eta^n$.

  
\subsection{Intersection with the tautological bundles} 
\label{ss:first main theorem}


Fix  a pair $(\cE,a)$  consisting of a rank $n-1$ vector bundle on $X'$ and an injective  Hermitian morphism $a : \cE \to \sigma^* \cE^\vee$.  
As in \S \ref{ssec: corank 1 special cycles}, associated to this data there is a  finite and unramified special cycle 
\[
\cZ_{\cE}^r(a)  \to \Sht^r_{\U(n)}   ,
\]
proper over $k$ by Proposition \ref{prop: corank one Z is proper}, 
and a virtual fundamental  class 
\[
[\cZ_{\cE}^r(a)]^{\vir} \in \mathrm{CH}_r ( \cZ_{\cE}^r(a) )
\]
in the Chow group of $r$-dimensional cycles.
For any line bundle $\cE_0$ on $X'$,  we consider the degree of the $0$-cycle class
 \begin{equation}\label{tautological intersection}
 \left( \prod_{i=1}^r c_1(p_i^* \sigma^* \cE_0^{-1} \otimes \ell_i) \right) \cdot  [\cZ_{\cE}^r(a)]^{\vir}
 \in \mathrm{CH}_0 ( \cZ_{\cE}^r(a) )
 \end{equation}
 of Theorem \ref{thm: geometric side}.

 \begin{prop}\label{prop:key degree}
The degree of  \eqref{tautological intersection} is $0$ if $r$ is odd, and is 
 \[
 \frac{ 2 }{  (\log q)^r  }
 \frac{d^r}{ds^r}\Big|_{s=0}
  \Big[  q^{ s d ( \cE_0 )  + s d(\cE) }  
 L(2s,\eta)  \Den_\eta( q^{ 1- 2 s}  ,\cE) \Big] 
\]
if $r$ is even.  
 \end{prop}
 
 \begin{proof}
Abbreviate  $d:=d(\cE)$,  and let  $Q := \mathrm{coker}(a ) \in \Herm_{2d}(k)$
be the image of $(\cE,a)\in \cA_d(k)$ under  \eqref{eq: g_d}.
Recall also the perverse sheaf 
\[
\sK_d^{i}   = ( \mathrm{Spr}_{2d}^{\mathrm{Herm}})^{W_i \times W_{d-i}  } \in D_c^b(\Herm_{2d}) 
\]
from \eqref{ICK} and \eqref{ICK3}.
By  Theorem \ref{thm: geometric side}, the degree of \eqref{tautological intersection}  vanishes if $r$ is odd, and is equal to 
\begin{align*}
 \frac{ 2 }{  (\log q)^r  }
 \frac{d^r}{ds^r}\Big|_{s=0}
  \Big[  q^{ s d ( \cE_0 )  + s d }  
 L(2s,\eta) 
\sum_{ i = 0 }^d   
 \mathrm{Tr} ( \mathrm{Fr}_{Q} ,    (\mathscr{K}_d^i )_{Q}   )   q^{i (1- 2 s)  } 
 \Big] 
\end{align*}
if $r$ is even.  The claim therefore follows from  Theorem \ref{thm: maindensityresult}.   
\end{proof}

Using Proposition \ref{prop:eisenstein density}, we can rewrite the equality of Proposition \ref{prop:key degree} in terms of the $(\cE,a)$-Fourier coefficient of the  Eisenstein series  $E(g,s,\chi)_{n-1}$  on $H_{n-1}$ from \S \ref{ss:eisenstein}.

\begin{thm}\label{thm:off-center}
The degree of  \eqref{tautological intersection} is $0$ if $r$ is odd, and is 
  \begin{equation}\label{shifted eisenstein}
 \frac{ 2  \cdot  q^{  \frac{n}{2}   [ d(\cE) +  \deg_X(\omega_X)  ]  }   }{  \chi(\det(\cE)) \cdot  (\log q)^r  }
  \frac{d^r}{ds^r}\Big|_{s=0}
  \Big[      q^{   n s   \deg_X(\omega_X)  + s \deg ( \cE_0 )  } 
   \mathscr{L}_n (s , \chi_0 )    E_{(\cE  , a  ) } (  s + 1/2    ,\chi )_{n-1}  \Big] 
\end{equation}
if $r$ is even.
\end{thm}

\begin{proof}
If we abbreviate 
\begin{equation}\label{third lambda}
\Lambda(s,\cE) = 
    q^{   -s   d(\cE)   }         L(-2s ,  \eta)     \mathrm{Den}_\eta(  q^{ 2s + 1 } , \cE  ) ,
  \end{equation}
 then Proposition \ref{prop:key degree} implies that the degree of  \eqref{tautological intersection} vanishes if $r$ is odd, and is 
 \begin{equation}\label{lambda derivative 1}
 \frac{ 2     }{  (\log q)^r  }
  \frac{d^r}{ds^r}\Big|_{s=0}
  \Big[  q^{ - s d ( \cE_0 )   } \Lambda( s,\cE)  \Big] 
\end{equation}
if $r$ is even.
The functional equation  
\[
q^{    (2s + \frac{1}{2} )   \deg_X(\omega_X)  } L(2s+1,  \eta )  =   L(-2s ,  \eta ) 
\]
together with the functional equation for twisted local densities (Proposition \ref{prop:functional equation}) allow us to rewrite \eqref{third lambda} as 
 \[
\Lambda(s,\cE  )    = 
q^{   ( s + 1 )    d(\cE )     }   
 q^{  (   2 s  + \frac{1}{2}   )    \deg_X(\omega_X)     }     L(2s+1 ,  \eta )  
  \mathrm{Den}_\eta(  q^{-2s -1 } , \cE ) ,
  \]
 and taking  $m=n-1$ in Proposition \ref{prop:eisenstein density}  shows that
\begin{equation}\label{Lambda}
    \Lambda(s,\cE )  = \chi(\det(\cE))^{-1}
  q^{  \frac{n}{2}   d(\cE) } q^{   \left[  \frac{n}{2}  + (n+1)s \right]  \deg_X(\omega_X)}  
   \mathscr{L}_n (s , \chi_0 )    E_{(\cE  , a  ) } (  s + 1/2    ,\chi )_{n-1}.
\end{equation}
Using this one sees that  \eqref{shifted eisenstein} is equal to \eqref{lambda derivative 1}.
\end{proof}


\subsection{Corank one higher Siegel--Weil formula}


We will now recast Theorem \ref{thm:off-center} as a relation between corank one special cycles and corank one Fourier coefficients of the Eisenstein series $E(g,s,\chi)_n$ on $H_n$ defined in \S \ref{ss:eisenstein}.

Fix a pair $(\cE,a)$ as before, except now assume that $\cE$ is a rank $n$ vector bundle on $X'$, and that the  Hermitian map $a: \cE \to \sigma^* \cE^\vee$ has rank $n-1$.
Once again, there is an associated finite and unramified morphism
\[
\cZ^r_{\cE}(a) \to \Sht^r_{\U(n)} ,
\]
but in contrast to the situation of \S \ref{ss:first main theorem}, the virtual fundamental class from \cite[\S 4]{FYZ2}
\[
[\cZ_{\cE}^r(a)]^{\vir} \in \mathrm{CH}_0 ( \cZ_{\cE}^r(a) )
\]
now lies in the Chow group of $0$-cycles.

First, we will digest the definition of the virtual fundamental class $[\cZ_{\cE}^r(a)]^{\vir} $. 
It is easy to see that   $\cE_0  := \ker(a)$ and $\cE^\flat := \cE/\cE_0$ are vector bundles on $X'$ of ranks $1$ and $n-1$, respectively, and that there is a unique injective Hermitian morphism $a^\flat: \cE^\flat \to \sigma^*(\cE^\flat)^\vee$ for which the composition
\[
\cE \to \cE^\flat \map{a^\flat} \sigma^*(\cE^\flat)^\vee \to \sigma^* \cE^\vee
\]
is equal to $a$.  In particular, the entire discussion of \S \ref{ss:first main theorem} can be applied to the pair $(\cE^\flat, a^\flat)$.

\begin{lemma}\label{lem: special cycle digest} 
There is a canonical isomorphism of $\Sht^r_{\U(n)}$-stacks
\[
 \cZ_{\cE^\flat}^r(a^\flat)  \iso  \cZ_{\cE}^r(a),
\]
and the induced isomorphism on Chow groups of $0$-cycles identifies   
 \begin{equation}\label{corank one is tautological intersection}
\left( \prod_{i=1}^r c_1(p_i^* \sigma^* \cE_0^{-1} \otimes \ell_i) \right) \cdot  [\cZ_{\cE^\flat}^r(a^\flat)]^{\vir} 
=
 [\cZ_{\cE}^r(a)]^{\vir} .
\end{equation}
Here the left hand side is the intersection  \eqref{tautological intersection}, but with $(\cE,a)$ replaced by $(\cE^\flat, a^\flat)$.
\end{lemma}

\begin{proof}
This is \cite[Lemma 3.2.5]{CH}. 
In that reference it is assumed that $\cE$ splits as the orthogonal direct sum $\cE_0 \oplus \cE^\flat$, with $\cE_0$ endowed with the zero map $\cE_0 \to \sigma^*\cE_0^\vee$, but the same proof works essentially verbatim in the slightly greater generality asserted here.

\end{proof}

\begin{thm}\label{thm:main corank one}
The stack $\cZ_\cE^r(a)$ is proper over $k$, and the $0$-cycle $ [\cZ_{\cE}^r(a)]^{\vir}$ has degree
\begin{equation}\label{final corank one}
\deg\,   [\cZ_{\cE}^r(a)]^{\vir} =
\frac{1}{ (\log q)^r } \cdot \frac{ q^{  \frac{n}{2} d( \cE)   } }{ \chi(  \det( \cE )) }
\cdot \frac{d^r}{ds^r} \Big|_{s=0} \left(   q^{ns \deg_X(\omega_X) }  \sL_n(s,\chi_0)  E_{ ( \cE,a) } (s,\chi)_n  \right) .
\end{equation}
Both sides are $0$ if  $r$ is odd.
\end{thm} 

\begin{proof}
Proposition \ref{prop: corank one Z is proper} tells us that  $\cZ_{\cE^\flat}^r(a^\flat)$ is a proper $k$-stack, and by Lemma \ref{lem: special cycle digest}  the same is true of $\cZ_{\cE}^r(a)$.  
 By Theorem \ref{thm:off-center} (and its proof),  the degree of \eqref{corank one is tautological intersection} vanishes if $r$ is odd, and is 
  \begin{equation}\label{lambda derivative}
 \frac{ 2     }{  (\log q)^r  }
  \frac{d^r}{ds^r}\Big|_{s=0}
  \Big[  q^{ - s d ( \cE_0 )   } \Lambda( s,\cE^\flat)  \Big] 
\end{equation}
if $r$ is even.  Here  $\Lambda(s,\cE^\flat)$ is defined by \eqref{Lambda}, but with $(\cE,a)$ replaced by $(\cE^\flat,a^\flat)$ throughout.
On the other hand, Proposition \ref{prop:genus drop} can be rewritten as 
\begin{align*}
&
 q^{ ns \deg_X(\omega_X) }  \mathscr{L}_n (s , \chi_0 ) \cdot  E_{ ( \cE , a ) } (s,\chi)_n  \\
& =
  \chi(\cE_0)  q^{ - \frac{n}{2}  d(\cE_0)} q^{ - s d(\cE_0)} 
  q^{   \left[  \frac{n}{2}  + (n+1)s \right]  \deg_X(\omega_X)}  
    \mathscr{L}_n(s,\chi_0)   E _{(\cE^\flat  , a^\flat ) } (  s+ 1/2   ,\chi  )_{n-1}   \\
& \quad   +  \chi(\cE_0)  q^{ - \frac{n}{2}  d(\cE_0)} q^{  s d(\cE_0)} 
  q^{   \left[  \frac{n}{2}  - (n+1)s \right]  \deg_X(\omega_X)}  
    \mathscr{L}_n(-s,\chi_0)   E _{(\cE^\flat  , a^\flat ) } (  -s+ 1/2   ,\chi  )_{n-1}   \\
    & = 
   \frac{  \chi( \det(\cE) )  }{  q^{  \frac{n}{2}  d(\cE)}  } 
   \big[  q^{ - s d(\cE_0)}   \Lambda(s,\cE^\flat)  
   + 
 q^{  s d(\cE_0)}   \Lambda(-s,\cE^\flat)   \big].
\end{align*}
The final expression is visibly an even function of $s$,  and it follows that right hand side of \eqref{final corank one} vanishes if $r$ is odd, and is equal to \eqref{lambda derivative} if $r$ is even.
\end{proof}


\subsection{Applications}
\label{ss:applications}


We derive some  corollaries of Theorem \ref{thm:main corank one}, following \cite{CH}.
Because  Theorem \ref{thm:main corank one}  only applies to pairs $(\cE,a)$ of corank one, the corollaries will be limited  to statements about special cycles on the moduli space $\Sht^r_{\U(2)}$, and their connections with automorphic forms on the quasi-split unitary group
$H_1 = \U(1,1)$ over $F$ from \S \ref{ss:eisenstein}.  

Let $K_1 \subset H_1(\A_F)$ be the standard maximal compact open subgroup.
Exactly as in \eqref{geometric coefficient}, any unramified automorphic form 
\[
f :  H_1(F)   \backslash H_1 (\A_F) / K_1 \to \C
\]
is determined by its Fourier coefficients $f_{ (\cE,a)} \in \C$, which are  indexed by pairs  $(\cE,a)$ consisting of a line bundle $\cE$ on $X'$ and a Hermitian morphism $a : \cE \to \sigma^* \cE^\vee$.

Fix an  unramified Hecke character 
$
\chi :   \A_{F'}^\times \to \C^\times
$
whose restriction   $\chi_0 : = \chi|_{\A_F^\times}$ is the trivial character, and let $E (g,s,\chi)_2$ be the  Eisenstein series on $H_2(\A_F)$  from \S \ref{ss:eisenstein}.    The \emph{doubling kernel}
\[
D  (g_1,g_2 , s,\chi) := E  ( i_0( g_1,g_2) , s,\chi)_2
\]
is the unramified automorphic form on $H_1(\A_F) \times H_1(\A_F)$ obtained by pulling back this Eisenstein series 
along the standard embedding $i_0 : H_1 \times H_1 \to H_2$ from \cite[(2.4.1)]{CH}. 
As in \cite[(4.2.5)]{CH}, the doubling kernel has Fourier coefficients 
 \begin{equation}\label{eq:geometric doubling coefficients}
D_{  (\mathcal{E}_1,a_1) , ( \mathcal{E}_2,a_2) }  ( s,\chi)
=
\sum_{  a = \left( \begin{smallmatrix} a_1 & * \\ * & a_2 \end{smallmatrix} \right) }
E_{ (\mathcal{E} , a ) } ( s , \chi)_2 
\end{equation}
indexed by pairs  $( \cE_1,a_1)$  and $(\cE_2,a_2)$,  in which each $\cE_i$ is a line bundle on $X'$, and each $a_i : \cE_i \to \sigma^* \cE_i^\vee$ is a Hermitian morphism.  
On the right hand side $\mathcal{E} : = \mathcal{E}_1\oplus \mathcal{E}_2$, and
 the sum is over all Hermitian morphisms $a\colon\mathcal{E} \to \sigma^*\mathcal{E}^\vee $ for which  the composition
\[
 \mathcal{E}_i  \map{\mathrm{inc.}} \mathcal{E}_1 \oplus \mathcal{E}_2 \map{a}  \sigma^*\mathcal{E}_1^\vee \oplus  \sigma^*\mathcal{E}_2^\vee \map{\mathrm{proj.}}   \sigma^*\mathcal{E}_i^\vee
\]
 agrees with $a_i$ for both $i \in \{1,2\}$.

Our first corollary of Theorem \ref{thm:main corank one} expresses the Fourier coefficients \eqref{eq:geometric doubling coefficients}  in terms of intersection multiplicities of  special cycles  on $\Sht^r_{\U(2)}$.  
Let us fix pairs $(\cE_1,a_1)$ and $(\cE_2,a_2)$ as above, and assume that  $a_2 : \cE_2 \to \sigma^*\cE_2^\vee$ is injective.  
By Proposition \ref{prop: corank one Z is proper},  the special cycle $\cZ_{\cE_2}^r(a_2)$  on $\Sht^r_{\U(2)}$ is a proper $k$-stack, and  we obtain cycle classes in middle codimension
\[
[\cZ_{\cE_1}^r(a_1)]^{\vir} \in \mathrm{CH}^r( \Sht^r_{\U(2)} ) 
\qquad \mbox{and}\qquad 
[\cZ_{\cE_2}^r(a_2)]^{\vir} \in \mathrm{CH}_c^r( \Sht^r_{\U(2)} ) .
\]
The second Chow group here is the Chow group with proper support, and  this allows us to  use the composition 
\[
\mathrm{CH}^r( \mathrm{Sht}^r_{\U(2)}  ) \times \mathrm{CH}_c^r( \mathrm{Sht}^r_{\U(2)}  ) \to \mathrm{CH}_c^{2r}( \mathrm{Sht}^r_{\U(2)}  )  \map{\deg} \Q
\]
to attach an intersection multiplicity to these two cycle classes.

\begin{cor}\label{cor: doubling coefficients}
For any pairs $(\cE_1,a_1)$ and $(\cE_2,a_2)$ as above (in particular, with $a_2$  injective), we have 
\begin{align*}
 & \deg   \big(  [\cZ_{\cE_1}^r(a_1)]^{\vir}  \cdot [\cZ_{\cE_2}^r(a_2)]^{\vir} \big) \\
 &=
\frac{1}{ (\log q)^r } \cdot \frac{ q^{   d( \cE_1)   }q^{  d( \cE_2)  } }{ \chi(   \cE_1 ) \chi( \cE_2  ) }
\cdot \frac{d^r}{ds^r} \Big|_{s=0} \left(   q^{2s \deg_X(\omega_X) }  \sL_2(s,\chi_0)  D_{  (\mathcal{E}_1,a_1) , ( \mathcal{E}_2,a_2) }  ( s,\chi)  \right) .
\end{align*}
Both sides are $0$ if  $r$ is odd.
\end{cor}

\begin{proof}
In the special case where $a_2$ is an isomorphism, this is \cite[Theorem 3.3.1(2)]{CH}.
The same argument applies here.  Briefly, by  \cite[Theorem 7.1]{FYZ2} we have the equality 
\[
[\cZ^r _{ \mathcal{E}_1}  ( a_1 )  ]^{\vir}  \cdot  [\cZ^r _{ \mathcal{E}_2}  ( a_2 )  ]^{\vir}  =
\sum_{  a = \left( \begin{smallmatrix} a_1 & * \\ * & a_2 \end{smallmatrix} \right) }
 [\cZ^r _{ \mathcal{E}  }  ( a )  ]^{\vir}  \in \mathrm{CH}_c^{2r}( \mathrm{Sht}^r_{\U(2)}  ),
\]
where $\cE = \cE_1 \oplus \cE_2$.  The assumption that $a_2$ is injective guarantees that every $a : \cE \to \sigma^* \cE^\vee$ appearing in the sum has corank zero or one.  Applying \cite[Theorem 1.1]{FYZ} to the corank zero terms, and Theorem \ref{thm:main corank one} to the corank one terms, we find that \eqref{final corank one} holds (with $n=2$) for every term in the sum.  The claim follows by comparing this with \eqref{eq:geometric doubling coefficients}.
\end{proof}

The remaining corollaries of  Theorem \ref{thm:main corank one} are closely related to the following special case of the Modularity Conjecture of \cite[Conjecture 4.15]{FYZ2}.

\begin{conj}\label{conj:modularity}
As $(\cE,a)$ varies over all pairs consisting of a line bundle $\cE$ on $X'$ and a Hermitian  morphism $a : \cE \to \sigma^* \cE^\vee$,  the rescaled classes 
\begin{equation}\label{eq:cycle coefficients}
 \chi(  \mathcal{E}  )   \cdot  q^{   - d(\mathcal{E} ) }   \cdot   [\cZ_{\cE}^r(a)]^{\vir} \in \mathrm{CH}^r( \mathrm{Sht}^r_{\U(2)}  )_\C
\end{equation}
 are the Fourier coefficients of an unramified automorphic form 
  \begin{equation}\label{eq:generating series}
Z^{r,\chi}  \colon  H_1(F)   \backslash H_1 (\A_F) / K_1   \to \mathrm{CH}^r( \mathrm{Sht}^r_{\U(2)}  )_\C.
\end{equation}
More precisely, there exists a function \eqref{eq:generating series} with the following property:  for any linear functional $\lambda$ on  the Chow group,  the unramified automorphic form $ \lambda(Z^{r,\chi}  ) :=   \lambda \circ Z^{r,\chi} $ has Fourier coefficients
\[
\lambda(Z^{r,\chi}  )_{(\cE,a)} = 
 \chi(  \mathcal{E}  )   \cdot  q^{   - d(\mathcal{E} ) }   \cdot    \lambda(  [\cZ^r _{ \mathcal{E}}  ( a )  ]^{\vir} ) .
\]
\end{conj}

Our second corollary of  Theorem \ref{thm:main corank one} provides evidence for the above modularity conjecture, by showing that the cycle classes \eqref{eq:cycle coefficients} become the coefficients of an automorphic form after applying particular linear functionals to them.

\begin{cor}
Fix a pair $(\cE_2,a_2)$ consisting of a line bundle $\cE_2$ on $X'$ and an injective Hermitian morphism $a_2 : \cE_2 \to \sigma^* \cE_2^\vee$.   
As $(\cE,a)$ varies over all line bundles $\cE$ on $X'$ equipped with a Hermitian morphism $a:\cE \to \sigma^* \cE^\vee$, the scalars 
\[
 \chi(  \mathcal{E}  )   \cdot  q^{   - d(\mathcal{E} ) }   \cdot  \deg\big(   [\cZ^r _{ \mathcal{E}}  ( a )  ]^{\vir}    \cdot    [\cZ^r _{ \mathcal{E}_2}  ( a_2 )  ]^{\vir}  \big) 
\]
are the Fourier coefficients of an unramified automorphic form on $H_1(\A_F)$.
\end{cor}

\begin{proof}
In the special case where $a_2$ is an isomorphism, this is \cite[Theorem 3.3.1(3)]{CH}.
The proof of the general case is identical: start with the automorphic form 
\[
\frac{1}{ (\log q)^r } \cdot \frac{ q^{   d( \cE_2)    } }{ \chi(   \cE_2 ) }
\cdot \frac{d^r}{ds^r} \Big|_{s=0} \left(   q^{2s \deg_X(\omega_X) }  \sL_2(s,\chi_0)  D ( g_1 , g_2 , s,\chi)  \right) 
\]
on $H_1(\A_F) \times H_1(\A_F)$, and  take its $(\cE_2,a_2)$-Fourier coefficient in the variable $g_2$.  
It follows from Corollary \ref{cor: doubling coefficients} that the resulting form in the variable $g_1$ has the desired Fourier coefficients.
\end{proof}

Our final corollary of Theorem \ref{thm:main corank one} is conditional on Conjecture \ref{conj:modularity}, which we now assume to hold.  
Suppose $\pi$ is an unramified irreducible cuspidal automorphic representation of $H_1(\A_F)$, and $f \in \pi$ is a $K_1$-fixed automorphic form.  Following ideas of Kudla, Conjecture \ref{conj:modularity} allows us to  define the \emph{arithmetic theta lift} 
  \begin{equation}\label{arith theta}
\vartheta^{r,\chi}(f) : = \int_{ H_1(F) \backslash H_1(\A_F)  / K_1} f(g) Z^{r,\chi}(g) \, dg \in  \Ch^r( \Sht^r_{\mathrm{U}(2)} ) _\C.
  \end{equation}
A theorem of Harder implies that $f$ is a compactly supported function on $ H_1(F) \backslash H_1(\A_F)$, so the integral is actually a finite sum.
 
 \begin{cor}\label{cor:GZformula}
Assume Conjecture \ref{conj:modularity}. Fix a pair $(\cE_2,a_2)$ consisting of a line bundle $\cE_2$ on $X'$ and an injective Hermitian morphism $a_2 : \cE_2 \to \sigma^* \cE_2^\vee$.    For any $f \in \pi$ fixed by $K_1$, we have the equality
\begin{align*}   
&  \deg \big(   \vartheta^{r,\chi} (f) \cdot [ \cZ^r_{\mathcal{E}_2}(a_2)]^{\vir} \big)  \\
& =
  f_{( \cE_2 , -a_2)  }   \cdot   \frac{    q^{ d(\cE_2) }    }{ (\log q)^r} \cdot 
    \frac{d^r}{ds^r}\Big|_{s=0}\left( q^{  2s \deg ( \omega_X)  }    L( s +1/2  ,  \mathrm{BC}(\pi) \otimes \chi )\right),
\end{align*}
where $f_{(\cE_2,-a_2)}$ is the $(\cE_2,-a_2)$-Fourier coefficient of $f$, and $ L( s   ,  \mathrm{BC}(\pi) \otimes \chi )$ is the twisted base-change $L$-function.
\end{cor}

\begin{proof}
In the special case where $a_2$ is an isomorphism, this is \cite[Theorem 3.3.1(4)]{CH}.
The proof given there shows that this result is a formal consequence of Corollary \ref{cor: doubling coefficients} and the duplication formula of \cite[Theorem 2.4.2]{CH}, and so  applies here as well.
\end{proof}

\begin{remark}
A form of the Modularity Conjecture has been announced in \cite{FYZ5}. It differs slightly from \cite[Conjecture 4.15]{FYZ2}, but implies it whenever $\cE$ has corank one, which in particular includes Conjecture \ref{conj:modularity}. Therefore, Corollary \ref{cor:GZformula} should soon become unconditional. 
\end{remark}

\subsection{Generalization with similitude factors}\label{sec: simil}

In this subsection, we will explain how the results can be generalized to shtukas for unitary groups with similitude factors. The main new feature of this case is that it supports examples where interesting identities are obtained for odd $r$. All the arguments adapt essentially verbatim, so we will omit them. 

\subsubsection{Hitchin base} Fix an integer $n \geq 1$ and let $\fL$ be a line bundle on $X$. We refer to \cite[Sections 2.3 and 2.4]{FYZ2} for a detailed discussion of the notions of Hermitian bundles and unitary shtukas with similitude factor $\fL$. In particular, the \emph{Hitchin base} $\CA_{\fL}$ is the stack whose $R$-points parametrize diagrams
    \[
    \CE \xrightarrow{a} \sigma^* \CE^{\vee} \otimes \nu^* \fL\]
where
\begin{itemize}
\item $\cE$ is rank $n-1$ vector bundle on $X_R'$, 
\item $a$ is a Hermitian morphism  that is injective fiberwise over $\Spec(R)$.
\end{itemize}

\subsubsection{Shtukas} There is a moduli space of rank $n$ unitary shtukas with similitude factor $\fL$, denoted $\Sht^r_{\U(n),\fL}$. The case $\fL = \cO_X$ recovers the $\Sht^r_{\U(n)}  = \Sht^r_{\U(n),\cO_X}$ discussed in the preceding parts of this paper. 

\subsubsection{Special cycles} For each $(\CE,a) \in \CA_{\fL}(k)$ and $r \geq 0$, there is a \emph{special cycle} 
\[\cZ^r_{\CE,\fL}(a) \to \Sht^r_{\U(n),\fL},\] equipped with a \emph{virtual fundamental class} \[[\cZ_{\cE,\fL}^r(a)]^{\vir} \in \Ch_r( \cZ_{\cE,\fL}^r(a) ), \]
    as defined in \cite[\S 4]{FYZ2}. For each $1 \leq i \leq r$, there is tautological line bundle $\ell_i$ on $\cZ^r_{\CE,\fL}(a)$ \cite[\S 4.3, Remark 4.7]{FYZ2}, whose definition is essentially identical to the case $\fL = \cO_X$ that we have discussed.

The proof of Proposition \ref{prop: corank one Z is proper} adapts to show:

\begin{prop} For any $(\CE,a) \in \CA_{\fL}(k)$, the $k$-stack $\cZ_{\cE,\fL}^r(a)$ is proper.
\end{prop}

Consequently, there is a well-defined degree morphism $\mathrm{CH}_0 ( \cZ_{\cE,\fL}^r(a)) \to \Q$.

\subsubsection{Intersection numbers} Given a vector bundle $\CE'$ on $X'$, we will use the notation 
\[
d_{\fL}(\CE') : = \rk(\CE') (\deg_X(\omega_X) + \deg_X(\fL)) - \deg_{X'}(\CE'),
\]
generalizing \eqref{eq: d}.
For any line bundle $\CE_0$ on $X'$, consider the class 

 \begin{equation}\label{simil tautological intersection}
 \left( \prod_{i=1}^r c_1(p_i^* \sigma^* \cE_0^{-1} \otimes \ell_i) \right) \cdot  [\cZ_{\cE,\fL}^r(a)]^{\vir}
 \in \mathrm{CH}_0 ( \cZ_{\cE,\fL}^r(a) )
 \end{equation}
 in the Chow group of $0$-cycles on $\cZ^r_{\CE,\fL}(a)$.
The following generalization of Theorem \ref{thm: geometric side} can be proved in an identical manner. 

 \begin{thm}\label{thm: simil key degree}
The degree of  \eqref{simil tautological intersection} is $0$ if $(-1)^r \neq \eta(\fL)^n$, and is 
 \[
 \frac{ 2 }{  (\log q)^r  }
 \frac{d^r}{ds^r}\Big|_{s=0}
  \Big[  q^{ s d_{\fL} ( \cE_0 )  + s d_{\fL}(\cE) }  
 L(2s,\eta)  \Den_\eta( q^{ 1- 2 s}  ,\cE) \Big] 
\]
if $(-1)^r = \eta(\fL)^n$.  
 \end{thm}

\subsubsection{Relation to Fourier coefficients of Eisenstein series} Finally, we can relate the expression appearing in Theorem \ref{thm: simil key degree} to Eisenstein series
on the quasi-split  unitary similitude group $\mathrm{GU}(n-1 , n-1)$ over $F$, as follows. We refer to \cite[\S 9.1 and 9.2]{FYZ2} for a detailed discussion of the relevant Eisenstein series. In particular, it depends on the choice of two Hecke characters, which (following the notation of {\it{loc. cit}}) we denote $\chi, \chi_0 :   \A_{F'}^\times \to \C^\times$.\footnote{The notation $\chi_0$ was used in the main body of the paper with a different meaning. It is unrelated to the present use of $\chi_0$, which is chosen to match the notation with \cite[\S 9]{FYZ2}.} The choice of $\chi_0$ is arbitrary, while $\chi$ is assumed to satisfy $\chi|_{\A_F^{\times}} = \eta^n$. This differs from the assumption $\chi|_{\A_F^{\times}} = \eta^{n-1}$ made in \cite[\S 9]{FYZ2}, which leads to the appearance of $\Den_{\eta}$ in equation \ref{eqn: simil Eisenstein formula} below.

Upon making the above choices, and using the normalized spherical section $\Phi$ as in \cite[Equation 9.1]{FYZ2}, the discussion in \cite[\S 9.1]{FYZ2} defines an Eisenstein series that we will denote by
\[E(g,s,(\chi,\chi_0))_{n-1}.\]

Moreover, \cite[Theorem 9.1]{FYZ2}, combined with the modification necessitated by $\chi|_{\A_F^{\times}} = \eta^n$ (as in the second part of Proposition \ref{prop:eisenstein density}), shows that for $(\CE,a) \in \CA_{\fL}(k)$, we have the following explicit formula for the Fourier coefficients of our Eisenstein series:
\begin{equation}\label{eqn: simil Eisenstein formula}
    E_{(\CE,\fL,a)}(s,(\chi,\chi_0))_{n-1} = \frac{(\chi_0\eta^{n})(\fL)\chi(\det \CE)}{\sL_{n-1}(s,\eta^n)} q^{(s-\frac{n-1}{2})d_{\fL}(\CE) - s (n-1)\deg \omega_X} \Den_{\eta}(q^{-2s},\CE).
\end{equation}

The renormalized Eisenstein series
\[
    \wt{E}_{(\CE,\fL,a)}(s,(\chi,\chi_0))_{n-1} := \frac{\sL_{n-1}(s,\eta^n)}{(\chi_0 \eta^{n})(\fL)\chi(\det \CE)} q^{(-s + \frac{n-1}{2})d_{\fL}(\CE) + s(n-1) \deg \omega_X} E_{(\CE,\fL,a)}(s,(\chi,\chi_0))_{n-1}
\]
 then satisfies
\begin{equation}\label{eq:simil Eisenstein renormalized}
    \wt{E}_{(\CE,\fL,a)}(s,(\chi,\chi_0))_{n-1} = \Den_{\eta}(q^{-2s}, \CE).
\end{equation}
We can now rewrite Theorem \ref{thm: simil key degree} in terms of $\wt{E}_{(\CE,\fL,a)}(s,(\chi,\chi_0))_{n-1}$. 

\begin{thm}\label{thm: simil main result Eisenstein}
    The degree of  \eqref{simil tautological intersection} is $0$ if $(-1)^r \neq \eta(\fL)^n$, and is
 \[
 \frac{ 2 }{  (\log q)^r  }
 \frac{d^r}{ds^r}\Big|_{s=0}
  \Big[  q^{ s d_{\fL} ( \cE_0 )  + s d_{\fL}(\cE) }  
 L(2s,\eta)  \wt{E}_{(\CE,\fL,a)}(s-1/2,(\chi,\chi_0))_{n-1} \Big] 
\]
if $(-1)^r = \eta(\fL)^n$.  
\end{thm}



\appendix


\section{Formulas for local density polynomials}

The goal of this Appendix is to explain how the sheaf-theoretic interpretation \eqref{eq: inert local density 1} of density polynomials proved in \cite{FYZ} can be used to derive more explicit formulae for them.  
We then show these explicit formulae can be used to derive further inductive formulae, amenable to computation. These formulae are not used directly in the main body of the paper, although they were useful for generating examples that aided the authors' investigation of twisted density polynomials. 

Finally, we explain how to deduce the analogous formulas for density polynomials of Hermitian lattices in the more classical $p$-adic setting.  In the $p$-adic  setting these formulas also follow from known results of Hironaka \cite{H99}.

Throughout this Appendix, $q$ always denotes an odd prime power.

\subsection{Universal density polynomials}

Let $\CO := \F_{q^2}[[\varpi]]$ and $\CO_0 := \F_{q}[[\varpi]]$ be power series rings in the variable $\varpi$, so that $\CO/\CO_0$ is an unramified quadratic extension of complete equicharacteristic discrete valuation rings. 
We recall the following combinatorial parametrization of torsion Hermitian $\cO$-modules. 

\begin{prop}
The set of isomorphism classes of torsion Hermitian $\CO$-modules of $\F_{q^2}$-length $d$ is in bijection with the set of partitions $\lambda$ of $d$, via the map associating to $Q$ its Jordan type as a torsion $\CO$-module.
\end{prop}
\begin{proof}
    This is a special case of \cite[Lemma 4.12]{FYZ}.
\end{proof}
In this context, we will denote the torsion Hermitian $\CO$-module corresponding to a partition $\lambda$ by $Q_{\lambda}$. 
Torsion $\CO_0$-modules of $\F_q$-length $d$ are also parametrized by partitions $\lambda$ of $d$, via Jordan type; for a partition $\lambda$ the corresponding torsion $\CO_0$-module will be denoted by $Q_{\lambda}^{\flat}$. 

 Given a partition $\lambda$, we will denote its length by $t(\lambda)$, and given either a torsion Hermitian $\cO$-module or torsion $\cO_0$-module $Q$ of Jordan type $\lambda$, we set 
 \[
 t(Q) := t(\lambda).
 \]
 This definition is consistent with \cite[\S 2.4]{FYZ}.
 
We  now define certain universal density polynomials. As in \ref{ss:density polynomials}, we take the formula of \cite[Theorem 2.3 (3)]{FYZ}, which is the unitary analogue of \cite[Theorem 1.1]{CY}, as our definition of density polynomials.  These are universal in the sense that, while their definitions depend on a fixed choice of odd prime power $q$, they will eventually turn out to also be polynomial in the parameter $q$.

\begin{defn}\label{def: univ-density} Let $\lambda$ be a partition of  $d \geq 0$.
\begin{enumerate}
    \item For $\epsilon \in \{\pm 1\}$ and $t\ge 0$, we define the polynomial
    \[\fm^{\epsilon}(q, t, T) := \prod_{i=0}^{t-1} (1 - (\epsilon q)^i T) \in \Z[T].\]

    \item  
    Let $Q_{\lambda}^{\flat}$ be the torsion $\CO_0$-module of Jordan type $\lambda$. 
    The \emph{split density polynomial} is 
    \[\Den^+(q,\lambda,T) := \sum_{0 \subset I_1 \subset I_2 \subset Q_\lambda^{\flat}} T^{\ell_{\CO_0}(I_1) + \ell_{\CO_0}(Q_{\lambda}^{\flat}/I_2)} \fm^+(q,t(I_2/I_1),T) \in \Z[T].\]
    Here the summation is over all chains of submodules $I_1 \subset I_2$, and $\ell_{\cO_0}(I_i)$ is the length of $I_i$ as an $\cO_0$-module. 
    
    \item   Let $Q_{\lambda}$ be the torsion Hermitian $\CO$-module of Jordan type $\lambda$. The \emph{inert density polynomial} is
    \[\Den^-(q, \lambda, T) := \sum_{0 \subset I \subset Q_\lambda} T^{2 \ell_{\CO}(I)} \fm^-(q,t(I^\perp/I),T) \in \Z[T]. \]
    Here the summation is over all isotropic submodules $I \subset Q_{\lambda}$, and $\ell_{\CO}(I)$ is the length of $I$ as an $\cO$-module. 
\end{enumerate}
\end{defn}

\begin{remark}\label{rem:universal poly comp}
As in the main body of the paper, suppose $X' \to X$ is a finite \'etale double cover  of smooth, projective, geometrically connected curves over the finite field $k=\F_q$.
Let $Q$ be a Hermitian torsion sheaf on $X'$, as in \S \ref{ss: hermitian torsion}.  For every closed point $v\in X$, the  local density polynomial  \eqref{local density def} satisfies  
\[
\Den(T,Q_v) = \Den^\pm(q_v,\lambda_v , T) ,
\]
where $q_v$ is the cardinality of the residue field of $v$.  If $v$ is inert in $X'$, the sign  is $-$, and $\lambda_v$ is the Jordan type of the completed stalk of $Q$ at the unique point of $X'$ above $v$.  If $v$ is split in $X'$, the sign  is $+$, and $\lambda_v$ is the Jordan type of the completed stalk  of $Q$ at either of two points of $X'$ above $v$.  
When $v$ is inert in $X'$ this equality is immediate from the definition \eqref{local density def}.  When $v$ is split it follows from  \cite[Remark 2.5]{FYZ}.
\end{remark}

To express $\Den^{\pm}(q,\lambda,T)$ as  polynomials in both $q$ and $T$, we introduce certain universal polynomials. 

\begin{defn}
    Given an integer $d\ge 0$, a partition $\lambda$ of $d$, and an integer $0 \leq a \leq d$, let 
    \[
    \Sub_{a,\lambda}(t) \in \Z[t]
    \]
    be the unique polynomial with the following property: for any odd prime power $q$, $\Sub_{a,\lambda}(q)$ equals the number of length $a$ submodules in the finite torsion $\F_q[[\varpi]]$-module with Jordan type $\lambda$.
\end{defn}

The existence of such polynomials is well-known, and implicit in the universality of Hall polynomials; see for example \cite[Chapter II]{macdonald1995symmetric}.
 We will express both the split and non-split density polynomials in terms of the $\Sub_{a,\lambda}(q)$. 
 
 \subsection{Split density polynomials} Let us begin with the split case.

\begin{prop}\label{prop: splitdenuniv}
    For any  partition $\lambda$ of $d \ge 0$, we have
    \[\Den^+(q,\lambda,T) = \sum_{a=0}^d \Sub_{a,\lambda}(q)  \cdot T^a.\]
\end{prop}
\begin{proof}
    This can be deduced from \eqref{eq: inert local density 1}, as in the proof of Proposition \ref{inertdenuniv} below;  however, let us give an elementary proof. 
    First we will slightly rephrase the desired equality. In the notation of Definition \ref{def: univ-density}, and abbreviating $\ell = \ell_{\cO_0}$ for length as an $\cO_0$-module, for any submodules $I_1 \subset I_2 \subset Q_\lambda^\flat$ we have 
    \begin{align*}
        T^{\ell(I_1) + \ell(Q_\lambda^\flat/I_2)} \fm^+( q , t(I_2/I_1), T) &= T^d \cdot T^{\ell(I_1) - \ell(I_2)} \cdot \prod_{i=0}^{t(I_2/I_1) - 1} (1-q^i T) \\
        &= T^d \cdot T^{\ell(I_1) - \ell(I_2) + t(I_2/I_1)} \cdot \prod_{i=0}^{t(I_2/I_1)-1} (T^{-1} - q^i) \\
        &= T^d \wt{\fm}^+(I_2/I_1, T^{-1}).
    \end{align*}
    Here, for a torsion $\CO_0$-module $I$, we have set
    \[
    \wt{\fm}^+(I,T) := T^{\ell(I) - t(I)} \prod_{i=0}^{t(I)-1}(T-q^i).
    \]
    Therefore, upon dividing both sides by $T^d$ and replacing $T$ by $T^{-1}$ (and using that the RHS is palindromic), we are reduced to showing 
    \[ \sum_{0 \subset I_1 \subset I_2 \subset Q_\lambda^\flat} \wt{\fm}^+(I_2/I_1, T) = \sum_{j=0}^d \Sub_{j,\lambda}(q) T^j.\]
    By definition, the RHS equals $\sum_{I \subset Q_\lambda^\flat} T^{\ell(I)}$. 
    For a fixed $I \subset Q_\lambda^\flat$, by grouping $T^{\ell(I)}$ on the RHS with the $I_2 = I$ terms on the LHS, we see that it suffices to prove the  identity 
    \[ \sum_{0 \subset I_1 \subset I} \wt{\fm}^+(I/I_1, T) = T^b \]
    for any finite $\cO_0$-module $I$ of length $b$.
    This last identity can be equivalently written as
    \[ \sum_{0 \subset I_1 \subset I} \wt{\fm}^+(I_1, T) = T^b,\]
because the multisets of submodules and quotients of $I$ coincide (as multisets of isomorphism classes of torsion modules). 
We will prove this equality by verifying that both sides take the same values on $T = q^a$, for any $a \ge 0$. In this case the right-hand side equals $q^{ab}$, which can be interpreted as $\# \Hom(\CO_0^a, I)$.

On the other hand, we can calculate  $\#\Hom(\CO_0^a,I)$ in the following alternative way. The data of a map $\CO_0^a \to I$ is equivalent to the data of a pair $(I_1,f)$, where $I_1 \subset I$ is its image, and $f: \CO_0^a \to I_1$ is a surjection. It suffices to verify that the number of such surjections equals
\[\wt{\fm}^+(I_1,q^a) = q^{a(\ell(I_1) - t(I_1))} \prod_{i=0}^{t(I_1)-1} (q^a - q^i).\]
By Nakayama's lemma, a $\cO_0$-module homomorphism $\CO_0^a \to I_1$ is a surjection if and only if it is a surjection upon reduction mod $\varpi$. Therefore we must prove that the number of surjections $k^a \to I_1/\varpi \cong k^{t(I_1)}$ equals $\prod_{i=0}^{t(I_1)-1}(q^a - q^i)$, as then the $q^{a(\ell(I_1)-t(I_1))}$ term simply counts all the lifts of a given surjection mod $\varpi$ to a map $\CO_0^a \to I_1$. Finally, we have
\begin{align*}
    \#\on{Surj}(k^a, k^t) &= \qbin{a}{t} \cdot |\GL_t(k)|  \\
    &= \frac{a!_q}{t!_q (a-t)!_q} q^{t(t-1)/2} \prod_{i=1}^{t} (q^i-1)  \\
    &= q^{t(t-1)/2 }\prod_{i=a-t+1}^{a} (q^i-1) \\
    &=\prod_{i=0}^{t-1} (q^a - q^i),
\end{align*}
where the first equality describes a surjection as the choice of kernel and an isomorphism between the cokernel and the target.
\end{proof}


\subsection{Inert density polynomials}


We now move towards an analogue of Proposition  \ref {prop: splitdenuniv} for the inert density polynomials.
In contrast to the proof of Proposition  \ref {prop: splitdenuniv}, the methods will be global in nature.

As in the main body of the paper, let $X' \to X$ be a finite \'etale double cover  of smooth, projective, geometrically connected curves over the finite field $k=\F_q$.
As in \cite[\S 3]{FYZ}, denote by $\Coh_d$ the moduli stack of torsion coherent sheaves on $X$ of length $d$.
Let  $\Spr_d$ be the Springer sheaf of \cite[Corollary 3.4]{FYZ},  a perverse sheaf on $\Coh_d$ with an action of the symmetric group $S_d$.

As in \S \ref{ssec: split support theorem}, for an integer  $a \leq d$, denote by $\Coh_{a \subset d}$  the moduli space parametrizing pairs $(I \subset Q)$, in which $Q \in \Coh_d$ is a length $d$ torsion sheaf on $X$, and $I \subset Q$ is a length $a$ subsheaf. The forgetful morphism $\Coh_{a\subset d} \to \Coh_d$ will be denoted by $\pi_{a \subset d}$. 
As  already noted in \eqref{eq: partial small},  we have an isomorphism of perverse sheaves on $\Coh_d$,
\begin{equation}\label{eq: SprSaSbinterpretation}
R(\pi_{a \subset d})_* (\Q_{\ell}) \cong \Spr_d^{S_a \times S_{d-a}}.
\end{equation}
Hence,  for any  $Q \in \Coh_d(k)$, the geometric cohomology of the fiber 
\[
\Sub_{a}(Q) := \pi_{a \subset d}^{-1}(Q)
\]
 is canonically identified, as a graded $\Ql$-vector space,  with the stalk of \eqref{eq: SprSaSbinterpretation} at (a geometric point above) $Q$.

\begin{prop}\label{prop: subsheafspacecohomology}
    Let $Q \in \Coh_d(k)$ be supported over a single $k$-rational point of $X$, with Jordan type a partition $\lambda$ of $d$. For any $a \leq d$, the geometric cohomology of $\Sub_a(Q)$ vanishes in odd degrees, and the Poincar\'e polynomial of $\Sub_a(Q)$ in even degrees is
    \begin{equation*}
        \sum_i \dim \rH^{2i}(\Sub_{a}(Q))  \cdot t^{i} = \Sub_{a,\lambda}(t).
    \end{equation*}
    Moreover, the Frobenius acts on $\rH^{2i}(\Sub_a(Q))$ by multiplication by $q^i$. 
\end{prop}

\begin{proof}
As in \S \ref{ss: SpringerComp}, let $\CB_Q$ be the $k$-scheme  parametrizing  complete flags of $\cO_X$-modules 
\[
0 \subset Q_1 \subset Q_2 \subset \ldots \subset Q_{d-1} \subset Q_d =  Q,
\]
so that  $\rH^*( \CB_Q) \iso (\Spr_d)_Q$, as graded $\Ql$-vector spaces.
    By \cite{Sp76}, the full Springer fiber $\CB_Q$ has a paving by affine spaces\footnote{\cite{Sp76} works over an algebraically closed field, but the construction of the paving does not use this, and works as well over $k$.}  \footnote{We can avoid invoking the affine space stratification of $\CB_Q$ by instead using the fact that the pullback map $\rH^*(\CB) \to \rH^*(\CB_Q)$ is surjective (where $\CB$ is the full flag variety). This reduces proving the odd cohomology vanishing and Frobenius action statements for the full flag variety $\CB$, where the affine space stratification can be provided easily.}
    This implies that the odd cohomology of $\CB_Q$ vanishes, and that Frobenius acts on $\rH^{2i}(\CB_Q)$ by multiplication by $q^i$. Hence the same is true of 
    \[
    \rH^*(\Sub_a(Q)) =  (\Spr_d^{S_a \times S_{d-a}} )_Q    = \rH^*(\CB_Q)^{S_a \times S_b}.
    \]
     
    It remains to prove that the Poincar\'e polynomial
    \begin{equation*}
        P(t) := \sum_i \dim \rH^{2i}(\Sub_a(Q)) \cdot t^i
    \end{equation*}
    satisfies $P(t) = \Sub_{a,\lambda}(t)$. From the above description of the Frobenius action on $\rH^*(\Sub_a(Q))$ and the Grothendieck-Lefschetz trace formula, we know that $\# \Sub_a(Q)(k) = P(q)$.
    On the other hand,  $\#\Sub_a(Q)(k) = \Sub_{a,\lambda}(q)$ holds by definition of the polynomial $\Sub_{a,\lambda}(t)$, hence we have 
    \[
    P(q) = \Sub_{a,\lambda}(q).
    \]
    
    It now suffices to show that the polynomial $P$ depends only on $a$ and $\lambda$, and not on the curve $X$ or its field of definition $k=\F_q$.  
    For then the equality $P(q) = \Sub_{a,\lambda}(q)$ must hold for any odd prime power $q$, which implies $P(t) = \Sub_{a,\lambda}(t)$. 
    
    This follows from the well-known fact that $P$ can be expressed via so-called \emph{Kostka-Foulkes polynomials} depending on $a$ and $\lambda$ only: see \cite{L81} and \cite[Chapter 8]{A2021}. 
    More precisely, as a graded $S_d$-representation, the geometric cohomology of the full Springer fiber $\cB_Q$ can be expressed as 
    \[\rH^*(\CB_Q) \cong \bigoplus_{\mu} M_{\lambda,\mu} \boxtimes V^{\mu},\]
    where the summation is over partitions $\mu$ of $d$, $V^{\mu}$ is the irreducible representation of $S_d$ corresponding to $\mu$, and the graded multiplicity space $M_{\lambda,\mu}$ has Poincar\'e polynomial equal to the modified Kostka-Foulkes polynomial $\wt{K}_{\lambda,\mu}(t)$. 
    By taking $S_a \times S_{d-a}$-invariants, one sees that  $\rH^*(\Sub_a(Q))$ similarly depends only on $a$ and $\lambda$.
\end{proof}

Recall from \eqref{HSpr} the perverse sheaf
\[
\HSpr_d :=  (\Spr_{2d}^{\Herm})^{(\Z/2\Z)^d}
\in D_c^b ( \Herm_{2d}  ),
\]
endowed with its Springer action of the subgroup $S_d$ of $W_d := (\Z/2\Z)^d \rtimes S_d$.

\begin{prop}\label{prop: HSprqminusqlemma}
Let $Q \in \Herm_{2d}(k)$ be supported above a single  $k$-rational point $v\in X$, inert in $X'$, with Jordan type given by a partition $\lambda$ of $d$. 
Then for any $0 \leq a \leq d$,
    \[
    \Tr(\Fr_Q, (\HSpr_d^{S_a \times S_{d-a}})_Q) = \Sub_{a,\lambda}(-q).
    \]
\end{prop}

\begin{proof}
Let $Q^\flat \in \Coh_d(k)$ be as in \S \ref{ss:FrobComp}.  
In particular, $Q^\flat$ is a  length $d$ torsion coherent sheaf on $X$ supported $v$, and its Jordan type is also $\lambda$.

Taking $\chi$ to be the trivial character in Proposition \ref{prop:t-twist},  the graded vector spaces
\[
(\HSpr_d^{S_a \times S_{d-a}})_Q  \iso \rH^*( \cB_Q^{\Herm} )^{ (\Z/2\Z)^d \rtimes(S_a \times S_{b-a})  } 
\]
and 
\[
(\Spr_d^{S_a \times S_{d-a}})_{Q^{\flat}} 
\iso \rH^*(\CB_{Q^\flat})^{S_a \times S_{d-a}}
\iso \rH^*(\Sub_a(Q^{\flat}))
\]
are isomorphic, but the Frobenius actions differ by an automorphism $\theta$.  Going back to the definition of $\theta$ from  \cite[Proposition 4.13]{FYZ}, our assumption that $Q$  is supported above a single inert $k$-rational point of $X$ implies that $\theta$ acts as $1$ on $\rH^{4i}(\Sub_a(Q^{\flat}))$ and as $-1$ on $\rH^{4i+2}(\Sub_a(Q^{\flat}))$. Combining this with \ref{prop: subsheafspacecohomology} yields the claim.
\end{proof}

\begin{prop}\label{inertdenuniv}
    For any odd prime power $q$ and partition $\lambda$ of $d$, there is an equality
    \[\Den^-(q,\lambda,T) = \sum_{a=0}^d (-1)^a \Sub_{a,\lambda}(-q)  \cdot T^a.\]
\end{prop}

\begin{proof}
 Choose  the \'etale double cover $X'\to X$ of smooth, projective, geometrically irreducible curves over $\F_q$  in such a way that there exists a torsion Hermitian sheaf  $Q$ on $X'$ satisfying 
\begin{itemize}
\item
$Q$   is supported  above a single  $\F_q$-rational point  $v\in X$, inert in $X'$,
\item
the  completed stalk  of $Q$ at the unique  point of $X'$ above $v$ has Jordan type $\lambda$.
\end{itemize}
Such data exist for any odd prime power $q$ and partition $\lambda$ (although for a fixed  $X'\to X$ such a $Q$ need not exist).  
For such choices,  we have 
 \begin{align*}
    \Den^-(q,\lambda,T)  & = \Den(T,Q)  \\
    & = \sum_{a+b=d} \Tr(\Fr_Q, (\HSpr_{a,b}^{S_a \times S_b})_Q)  \cdot T^a \\
    & = \sum_{a+b=d}(-1)^a  \Tr(\Fr_Q, (\HSpr_d^{S_a \times S_b})_Q)  \cdot T^a,
\end{align*}
where the first equality is by Remark \ref{rem:universal poly comp}, the second is \eqref{eq: inert local density 1}, and the third is Proposition \ref{prop:inert trace comparison}.
Combining this with Proposition \ref{prop: HSprqminusqlemma} completes the proof.
\end{proof}

Propositions \ref{prop: splitdenuniv} and \ref{inertdenuniv} immediately imply the following Corollary.

\begin{cor}\label{cor: denareunivpolyqT}
    The density polynomials $\Den^{\pm}(q,\lambda,T)$ interpolate into polynomials in both $q$ and $T$, and these polynomials are related by 
    \[\Den^-(q,\lambda,T) = \Den^+(-q,\lambda,-T).\]
\end{cor}


\subsection{Induction formulae}


In this subsection, we explain how Propositions \ref{prop: splitdenuniv} and \ref{inertdenuniv}  can be used to prove induction formulae for density polynomials. 
First we recall without proof a standard lemma on submodules in direct sums of $\CO = \F_q[[\varpi]]$-modules.

\begin{lemma}\label{lem: subs in direct sum}
    Let $A$ and $B$ be finite $\CO$-modules. The data of a submodule $I \subset A \oplus B$ of length $c$ is equivalent to the data of $(I_1, I_2, f)$, where 
    \begin{itemize}
    \item
     $I_1 \subset A$ and $I_2 \subset B$ are submodules satisfying $\ell(I_1) + \ell(I_2) = c$, 
     \item
     $f: I_1 \to B/I_2$ is a morphism of $\CO$-modules. 
     \end{itemize}
     The bijection sends $I$ to the triple $(I_1,I_2,f)$ with $I_1 = \on{im}(I \to A)$, $I_2 = I \cap (0 \oplus B)$, and the function $f(a) = p_2(\wt{a})$. Here $\wt{a}$ is any lift of $a \in I_1 \subset A$ to $I$, and $p_2$ is  projection to the second factor in $A\oplus B$.
\end{lemma}

We can now prove an inductive formula for the polynomials $\Sub_{a,\lambda}(t)$.

\begin{prop}\label{prop: induction for Salambda}
    Fix a partition $\lambda = (a_1 \geq a_2 \geq \dots \geq a_n)$. For any $m \geq a_1$ and $a \leq |\lambda|$ we have the equality
    \begin{equation}\label{eq: induction for Salambda}\Sub_{a,(m,\lambda)}(t) = \Sub_{a-1,(m-1,\lambda)}(t) + t^a \Sub_{a,\lambda}(t).
    \end{equation}
Here $(k,\lambda)$ is the partition obtained by inserting $k$ into the list $(a_1 \geq a_2 \geq \dots \geq a_n)$ at the appropriate spot to maintain the descending ordering.
\end{prop}

\begin{proof}

It suffices to prove the equality for any prime power $t = q$. Let  $Q = Q_{\lambda}$ be the length $|\lambda|$ torsion $\CO$-module of Jordan type $\lambda$, and let $Q_m$ be the cyclic $\CO$-module of length $m$.

    The left-hand side of \eqref{eq: induction for Salambda} counts the number of length $a$ submodules in $Q \oplus Q_m$. By Lemma \ref{lem: subs in direct sum}, this number equals
    \[\sum_{i=0}^a \#\{(I_i \subset Q, f:I_i \to Q_m / Q_{a-i}\}),\]
    where $I_i$ runs over all submodules $I_i \subset Q$ of length $i$, and the summation over $i$ is with the understanding that if $Q_m / Q_{a-i}$ does not make sense, i.e. $a - i > m$, the corresponding term is zero.
    
    The term $\Sub_{(a-1),(m-1,\lambda)}(q)$ on the right-hand side of \eqref{eq: induction for Salambda} counts the number of length $a-1$ submodules in $Q \oplus Q_{m-1}$. By Lemma \ref{lem: subs in direct sum} again, this number equals
    \[\sum_{i=0}^{a-1} \#\{(I_i \subset Q, f:I_i \to Q_{m-1}/Q_{a-1-i})\},\]
    with the same caveat on the summation over $i$ as above. 
    
    The terms with the same value of $i$ in both sums coincide, therefore the difference equals the $i = a$ term from the first sum, which is
    \[\Sub_{a,(m,\lambda)}(q) - \Sub_{a-1,(m-1,\lambda)}(q) = \#\{(I_a \subset Q, f:I_a \to Q_m \})
    .\]
By our  assumption that $m \geq a_1$, for  any length $a$ submodule $I_a \subset Q$ we have
    \[
    \#\Hom(I_a, Q_m) =  |I_a| = q^a.
    \]
   Therefore
    \[\Sub_{a,(m,\lambda)}(q) - \Sub_{a-1,(m-1,\lambda)}(q) = \sum_{I_a \subset Q} q^a = q^a \Sub_{a,\lambda}(q),\]
    as desired.\end{proof}
    
We can use this identity to give simple proofs of the following induction formulae for local density polynomials. Part (1) of the following corollary appears as the split case of \cite[Proposition 9.3.2]{ChenI}, with a different proof method.

\begin{prop}\label{prop: induction formula strong}
    Let $\lambda = (a_1 \geq a_2 \geq \dots \geq a_n)$. Then for any $m \geq a_1$, the following identities hold:
    \begin{enumerate}
        \item $\Den^+(q,(m,\lambda),T) = T \Den^+(q,(m-1,\lambda),T) + \Den^+(q,\lambda,qT)$.
        \item $\Den^-(q,(m,\lambda),T) = -T \Den^-(q,(m-1,\lambda),T) + \Den^-(q,\lambda,-qT).$
    \end{enumerate}
\end{prop}

\begin{proof}
    For any partition $\mu$, we have by Proposition \ref{prop: splitdenuniv} that 
    \[\Den^+(q,\mu,T) = \sum_{a=0}^{|\mu|} \Sub_{a,\mu}(q) \cdot  T^a.\]
The equality  (1) follows by applying Proposition \ref{prop: induction for Salambda}, and  (2) follows from (1) using Corollary \ref{cor: denareunivpolyqT}. 
\end{proof}

As a corollary, we obtain the following induction formula, as proved via different methods in \cite[Theorem 5.1]{T13} and the inert case of \cite[Proposition 9.3.2]{ChenI}; see also \cite[Proposition 3.7.1]{LZ1}.

\begin{cor}
\label{cor: induction formula weak}
    Given a partition $\lambda$ and $m \geq a_1 + 1$, the following identity holds:
    \[\Den^-(q,(m,\lambda),T) = T^2 \Den^-(q,(m-2,\lambda),T) + (1-T) \Den^-(q,\lambda,-qT).\]
\end{cor}
\begin{proof}
    This follows by applying Proposition \ref{prop: induction formula strong}(2) twice:
    \begin{align*}
        \Den^-(q,(m,\lambda),T) &= -T \Den^-(q,(m-1,\lambda),T) + \Den^-(q,\lambda,-qT) \\
        &= T^2 \Den^-(q,(m-2,\lambda),T) + (1-T) \Den^-(q,\lambda,-qT),
    \end{align*}
    where we used the hypothesis $m-1 \geq a_1$ to justify the second application. 
\end{proof}


\subsection{The mixed characteristic case}


Let $F/F_0$ be an unramified quadratic extension of local fields,  with rings of integers denoted $\CO/ \CO_0$, and residue fields $\F_{q^2}/\F_q$ of odd characteristic.
The goal of this subsection is to state and prove a version of Proposition \ref{inertdenuniv} for Hermitian $\CO$-lattices.  
We emphasize that we are now allowing $\CO_0$ to be either of mixed characteristic or equal characteristic.

Fix a uniformizing parameter $\varpi$ of $\CO_0$.
For any Hermitian $\CO$-lattice $L$ of rank $n$, there is a unique sequence of integers $a_1 \geq a_2 \dots \geq a_n \geq 0$ such that $L$ is the orthogonal direct sum
\[L = \oplus \langle \varpi^{a_i} \rangle,\]
where $\langle \varpi^a \rangle$ is a rank $1$ $\CO$-lattice with a generator of Hermitian norm  $\varpi^a$. 
If we discard the $a_i$ that equal $0$, the rest of the $a_i$ form a partition, which we denote by $\type(L)$. Conversely, given a partition $\lambda$ of some integer $d\ge 0$ with length of $\lambda$ at most $n$, we denote by $L_{\lambda}$ the corresponding rank $n$ lattice with $\type(L_{\lambda}) = \lambda$. This defines a bijection between the sets of isomorphism classes of rank $n$ Hermitian $\CO$-lattices and the set of partitions of length at most $n$.

    We will also use the ``type" notation in the following context. Given a finitely generated torsion $\CO$-module $M$, there is a unique partition $\lambda = a_1 \geq a_2 \geq \dots$ such that 
    \[M \cong \bigoplus_i \CO/ ( \varpi^{a_i}) .\]
    We will write $\type(M) = \lambda$. The two usages of ``type" are consistent in the following sense: given a Hermitian $\CO$-lattice $L$, we have
    \[
    \type(L) = \type(L^{\vee}/L),
    \]
where $L^\vee$ is the dual lattice of $L$.

Given a Hermitian $\cO$-lattice $L$, there is an associated local density polynomial $\Den(L,T)$ \cite[\S 2.3-2.4]{FYZ}.  If $F_0$ has nonzero characteristic, so that $\CO_0\iso \F_q[[\varpi]]$, the Cho-Yamauchi formula of \cite[Remark 2.4]{FYZ} implies 
\begin{equation}\label{CYisUniversal}
\Den(L,T) = \Den^-(q, \lambda, T),
\end{equation}
where $\lambda=\type(L)$.
    
We will give two proofs of Proposition \ref{prop:p-adic density} below, providing a coefficient-by-coefficient formula for $\Den(L,T)$.  
The first proof, assuming that $\mathrm{char}(F_0)=0$,  is based on a formula of  Hironaka \cite{H99}.
Hironaka systematically assumes that $F_0$ has characteristic $0$, and the authors have not  checked if the same  methods can be applied to the equal characteristic setting.
The second proof is entirely different.  The idea is to show that $\Den(L,T)$ is, in a suitable sense, independent of the characteristic of $F_0$.  This will allow us to deduce Proposition \ref{prop:p-adic density}, in both the equal and mixed characteristic cases,  from Proposition \ref{inertdenuniv}.

    \begin{prop}\label{prop:p-adic density}
     For any Hermitian $\CO$-lattice $L$ of type $\lambda$, we have 
    \[
    \Den(L,T) = \bigoplus_{a=0}^{|\lambda|} (-1)^a \Sub_{a,\lambda}(-q) \cdot  T^a.
    \]
    \end{prop}
    
   \begin{proof}[Proof of Proposition \ref{prop:p-adic density}, assuming $\mathrm{char}(F_0) =0$] 
   It follows from \cite[Corollary 3.3]{H99} that for any $j \geq 0$, using notation from \textit{{loc. cit.}},
    \begin{align*}
    \frac{\Den(\CO^{m+j},L)}{{\Den(\CO^{m+j},\CO^{m}})} &= \sum_{\mu \le \lambda} (-1)^{(j-1)|\mu|} q^{-j |\mu|} N^{\lambda}_\mu(-q^{-1}) \\
    &= \sum_{\mu \leq \lambda} (-1)^{|\mu|} \cdot (-q)^{-j |\mu|} \cdot N_{\mu}^{\lambda}(-q^{-1}).
    \end{align*}  
It follows from \cite[Chapters 2 and 3]{macdonald1995symmetric} that $N^{\lambda}_\mu(q^{-1})$ equals
    \[\#\{I \subset Q_{\lambda}: \on{type}(I) = \mu\},\]
where $Q_{\lambda}$ is the torsion $\F_q[[\varpi]]$-module of Jordan type $\lambda$. Therefore
    \[ \frac{\Den(\CO^{m+j},L)}{{\Den(\CO^{m+j},\CO^{m}})} =  \sum_{a=0}^{|\lambda|} (-1)^a \Sub_{a,\lambda}(-q) (-q)^{-ja}.\]
    The defining property of the polynomial $\Den(L,T)$ is that
    \[\Den(L, (-q)^{-j}) = \frac{\Den(\CO^{m+j},L)}{\Den(\CO^{m+j},\CO^{m})},\]
   completing the proof.
\end{proof}

We now begin our second proof of Proposition \ref{prop:p-adic density}.

\begin{defn}
Let $f$ be a function on the set of isomorphism classes of triples
\[
(\CO_0 \subset \CO, L)
\]
where $\CO_0 \subset \CO$ is an unramified quadratic extension of complete discrete valuation rings with finite residue fields of odd characteristic, and $L$ is a Hermitian $\CO$-lattice of rank $n$. 
We say that $f$ is \emph{characteristic-independent} if it only depends on the cardinality of the residue field of $\CO_0$ and the partition $\type(L)$.  
\end{defn}

\begin{remark}
A characteristic-independent $f$ is insensitive to the characteristic of the field of fractions of $\CO_0$ 
(hence the name)  and is determined by its values on triples of the form $(\F_q[[\varpi]]\subset \F_{q^2}[[\varpi]], L)$.
\end{remark}

Let $\CO$ be a complete discrete valuation ring with residue field $\F_{q^2}$. 
For a partition $\lambda$, let $M_{\lambda}$ be the torsion $\CO$-module with Jordan type $\lambda$. 
Given an integer $d$, a partition $\lambda$ of $d$, and a sequence of partitions $\{\mu_i\}=\{ \mu_0,\ldots, \mu_{k-1} \}$ with $\sum |\mu_i| = d$, we set
    \[ \Fl(\CO,\{\mu_i\},\lambda) = \#\{0 = I_0 \subset I_1 \subset I_2 \subset \dots \subset I_k = M_{\lambda}: \type(I_{i+1}/I_i) = \mu_i)\}.\]

    It is well-known that the numbers $\Fl(\CO,\{\mu_i\},\lambda)$ are characteristic-independent,  i.e. only depend on $q, \{\mu_i\}$ and $\lambda$: this is because they can be expressed via Hall polynomials, which are characteristic-independent, see \cite[Chapters 2]{macdonald1995symmetric}. This will be a key tool in the proof of \ref{prop: counting is independent}.

\begin{prop}\label{prop: flag counting functions lemma}
    Given integers $d \geq 2$ and $n$, consider the set $P(d, \leq n)$ of partitions of $d$ with length at most $n$. 
    Then the space of functions on $P(d, \leq n)$ is spanned by the functions $\Fl(\CO,\{\mu_i\},-)$, where $\{\mu_i\}$ runs over all sequences of partitions with total size $d$, with each $\mu_i$ of length at most $n$, and $\{\mu_i\} \neq \{(d)\}$ (i.e., there is more than one partition). 
\end{prop}

\begin{proof}
    Given a partition $\lambda \in P(d, \leq n)$, consider the dual partition $\lambda^{\vee} = n_1 \geq \dots \geq n_k$. 
    Let $\nu_i$ be the partition $1\ge 1 \ge \cdots \ge 1$ consisting of $n_i$ copies of $1$.
    The function $\Fl( \CO, \{ \nu_i\},-)$ takes the value $1$ on $\lambda$, and is supported on partitions that are equal to or below $\lambda$ in the dominance order. (Note that $n \geq n_1 \geq n_i$, therefore each $\nu_i$ has length at most $n$.) By induction on the dominance order, this allows us to express each indicator function $1_{\lambda}$ as a linear combination of functions of the form  $\Fl( \CO, \{\mu_i\},-)$. This proves the claim.
\end{proof}

\begin{prop}\label{prop: counting is independent}
    Let $\mu$ and $\nu$ be partitions of length at most $n$. 
    The function $N(\mu,\nu)$ that sends a triple $(\CO_0 \subset \CO, L)$ to the integer 
    \[\#\{L \subset L'  \subset (L')^\vee \subset L^\vee: \type(L'/L) = \mu \mbox{ and }\type((L')^{\vee}/L') = \nu\}\]
     is characteristic-independent.
\end{prop}

\begin{proof}
We will prove this by induction first on the size of $\lambda = \type(L)$, then on the size of $\mu = \type(L')$.

For the base case, consider $N(\mu,\nu)$ with $\mu = (1)$. In this case the set of all colength-1 lattices $L \subset L'$ is in bijection with $\PP(L \otimes_{\CO} \F_{q^2})$. The integral lattices $L'$ correspond to the subset $Y \subset \PP(L \otimes_{\CO} \F_{q^2})$ of lines that are isotropic under the natural Hermitian pairing on $L \otimes_{\CO} \F_{q^2}$. Moreover, one can directly check that under this bijection, one can control the type of $L'$, i.e. for any $\nu$, if we denote by $Y_{\nu} \subset Y$ the subset of lines such that the corresponding lattice $L'$ has type $\nu$, then the cardinality of $Y_{\nu}$ is characteristic-independent, i.e. only depends on $q, \nu$ and $\type(L)$. This shows the $\mu = (1)$ base case.

For the inductive step, we proceed as follows. Given a sequence of partitions $\mu_0, \dots, \mu_{k-1}$ with $\sum |\mu_i| = |\mu| = d \geq 2$, and $\{\mu_i\} \neq \{(d)\}$ (i.e. $k \geq 2$), consider the function $N(\{\mu_i\},\nu)$ that sends $(\CO_0 \subset \CO, L)$ to
\[
\#\{L = L_0 \subset L_1 \subset \dots \subset L_k = L': \type(L') = \nu \mbox{ and } \type(L_{i+1}/L_i) = \mu_i\}.
\]
Then $N(\{\mu_i\},\nu)(\CO_0 \subset \CO, L)$ can be computed as follows: we first choose $L \subset L_1$ with $\type(L_1/L) = \mu_1$, and then choose the remaining sequence of lattices containing $L_1$. Because $|\type(L_1)| < |\type(L)|$, 
the inductive hypothesis implies that the count of lattices containing $L_1$ is characteristic-independent. This shows that $N(\{\mu_i\},\nu)$ is also characteristic-independent. 

On the other hand, we can compute $N(\{\mu_i\},\nu)(\CO_0 \subset \CO,L)$
by first choosing $L \subset L_k = L'$ with $\type(L') = \nu$, and then choosing the rest of the flag in the torsion module $L'/L$. Therefore we have
\[N(\{\mu_i\},\nu)(\CO_0 \subset \CO, L) = \sum_{\mu \in P(d, \leq n)} N(\mu,\nu)(\CO_0 \subset \CO,L) \cdot \on{Fl}(\CO,\{\mu_i\},\mu).\]

By \ref{prop: flag counting functions lemma}, for any $\mu \in P(d,\leq n)$ the indicator function $1_{\mu}$ can be expressed as a linear combination of the functions $\Fl(\CO,\{\mu_i\},-)$. This implies that the numbers $N(\{\mu_i\},\nu)(\CO_0 \subset \CO,L)$ determine $N(\mu,\nu)(\CO_0 \subset \CO,L)$ uniquely. We know that the $N(\{\mu_i\},\nu)(\CO_0 \subset \CO, L)$ are characteristic-independent by induction, and $\Fl(\CO,\{\mu_i\},\mu)$ are characteristic-independent by the universality of the Hall algebra. By induction, this implies that the functions $N(\mu,\nu)$ are characteristic-independent.
\end{proof}

\begin{cor}\label{cor: Den formula in mixed char}
    The function sending a triple $(\CO_0 \subset \CO, L)$ to the polynomial $\Den(L,T)$ is characteristic-independent. 
\end{cor}
\begin{proof}
The Cho-Yamauchi formula \cite[(2.7)]{FYZ} expresses $\Den(L,T)$ as a weighted sum over lattices  $L'$ satisfying 
\[
L \subset L' \subset (L')^\vee \subset L^\vee ,
\]
where the weights depend only on the types of $L'/L$ and $(L')^\vee / L'$.
By Proposition \ref{prop: counting is independent}, this counting problem is characteristic-independent, therefore so is $\Den(L,T)$.
\end{proof}

\begin{proof}[Proof of Proposition \ref{prop:p-adic density} in general]
As in the statement of Proposition \ref{prop:p-adic density}, suppose $L$ is a Hermitian $\CO$-lattice of type $\lambda$. 
Recalling that  $\F_q \subset \F_{q^2}$ are the residue fields of $\CO_0 \subset \CO$, form the  quadratic extension $\F_q[[\varpi]] \subset \F_{q^2}[[\varpi]]$ of equicharacteristic complete discrete valuation rings, and fix a Hermitian $\F_{q^2}[[\varpi]]$-lattice $L_\lambda$ of type  $\lambda$.
Corollary \ref{cor: Den formula in mixed char} implies the first equality in 
\[
\Den(L,T) = \Den(L_\lambda,T) \stackrel{\eqref{CYisUniversal}}{=} \Den^-(q,\lambda,T) ,
\]
and so the claim follows from Proposition \ref{inertdenuniv}.
\end{proof}

\bibliographystyle{amsalpha}
\bibliography{Bibliography}

\end{document}